\documentclass[a4paper,titlepage,12pt]{amsbook}

\usepackage{amsmath, amsfonts, amsthm, amsfonts, amssymb, url,multicol,graphicx}
\usepackage[francais]{babel}
\usepackage[T1]{fontenc}
\usepackage[latin1]{inputenc}

\def\]{\textup{\mbox{]\hspace{-.15em}]}}}
\def\[{\textup{\mbox{[\hspace{-.15em}[}}}

\def\got{\mathfrak}

\newcommand{\PGL}{{\rm PGL}}

\newcommand{\WW}{\mathcal{W}}
\newcommand{\Hom}{\mathrm{Hom}}
  
\newcommand{\Gal}{\mathrm{Gal}}
\newcommand{\G}{\mathrm{G}}

\newcommand{\Qpb}{\overline{\mathbb{Q}}_p}
\newcommand{\Spec}{\mathrm{Spec}}
\newcommand{\diag}{\mathrm{diag}}
\newcommand{\GL}{\mathrm{GL}}
\newcommand{\SL}{\mathrm{SL}}

\newcommand{\ps}{\par \smallskip}

\newcommand{\Z}{\mathbb{Z}}

\newcommand{\Q}{\mathbb{Q}}

\newcommand{\Qp}{\mathbb{Q}_p}  
\newcommand{\R}{\mathbb{R}}  
\newcommand{\C}{\mathbb{C}}

\newcommand{\AAA}{\mathbb{A}}
\newcommand{\OO}{\mathcal{O}}

\newcommand{\F}{\mathbb{F}}
\newcommand{\isomo}{\overset{\sim}{\rightarrow}}

\newtheorem{definition}[subsection]{Definition}
\newtheorem{thmetoile}[subsection]{Th\'eor\`eme$^\ast$}

\newtheorem{prop}[subsection]{Proposition}

\newtheorem{cor}[subsection]{Corollaire}

\newtheorem{question}[subsection]{Question}
\newtheorem{thm}[subsection]{Th\'eoreme}
\newtheorem{conjecture}[subsection]{Conjecture}

\newenvironment{pf}
{\medskip\noindent {\it Preuve --- \ }}
{\hfill\nobreak $\Box$ \par\bigbreak} 

\makeatletter
\def\@citex[#1]#2{\leavevmode
  \let\@citea\@empty
  \@cite{\sc\small\@for\@citeb:=#2\do
    {\@citea\def\@citea{,\penalty\@m\ }%
     \edef\@citeb{\expandafter\@firstofone\@citeb\@empty}%
     \if@filesw\immediate\write\@auxout{\string\citation{\@citeb}}\fi
     \@ifundefined{b@\@citeb}{\hbox{\reset@font\bfseries ?}%
       \G@refundefinedtrue
       \@latex@warning
         {Citation `\@citeb' on page \thepage \space undefined}}%
       {\@cite@ofmt{\csname b@\@citeb\endcsname}}}}{#1}}
\makeatother

\title{Repr\'esentations galoisiennes automorphes et cons\'equences arithm\'etiques des conjectures de Langlands et Arthur}
\author{Ga\"etan CHENEVIER}

\begin{document}

\title{%
    \begin{minipage}\linewidth\centering
    \small\sc M\'emoire d'habilitation \`a diriger des recherches
    \vskip2mm
    \small Soutenu le 29 Mars 2013 \`a l'universit\'e Paris XI
    \vskip4cm
     \centering\bf\Large Repr\'esentations galoisiennes automorphes et cons\'equences arithm\'etiques des conjectures de Langlands et Arthur
    \end{minipage}
}
\maketitle

\frontmatter

\newpage

\setcounter{tocdepth}{1}
\tableofcontents

\chapter*{Introduction}

	Ce\footnote{L'auteur est financ\'e par le CNRS et soutenu par le projet ANR-10-BLAN 0114.} m\'emoire a pour but de pr\'esenter des r\'esultats
en math\'ematiques que j'ai obtenus depuis ma th\`ese, en vue de
l'obtention d'une habilitation \`a diriger des recherches.  Mes travaux
portent principalement sur l'\'etude des repr\'esentations du groupe de
Galois absolu de $\Q$ et sur divers probl\`emes arithm\'etiques qui leurs
sont attach\'es de mani\`ere plus ou moins cach\'ee. Ils sont presque tous
inspir\'es des conjectures de R. P. Langlands, ou des raffinements de ces derni\`eres propos\'es par
J. Arthur.\ps

Le texte est organis\'e en deux parties, qui bien que profond\'ement
reli\'ees peuvent \^etre lues de mani\`ere ind\'ependante.  La premi\`ere
traite des repr\'esentations galoisiennes "automorphes", c'est-\`a-dire qui
sont attach\'ees aux repr\'esentations automorphes cuspidales de $\GL(n)$,
et de leurs d\'eformations.  Ma contribution principale \`a ce th\`eme est
la construction et l'\'etude de la foug\`ere infinie, un certain
sous-ensemble aux allures fractales de la vari\'et\'e des caract\`eres
$p$-adique de ${\rm Gal}(\overline{\Q}/\Q)$.  J'en ai obtenu plusieurs
applications : \begin{itemize}\ps\ps

\item[-] Construction de repr\'esentations galoisiennes automorphes en poids peu
r\'eguliers. \ps

\item[-] D\'etermination de l'alternative symplectique/orthogonale pour les
repr\'esentations galoisiennes automorphes essentiellement autoduales sur
les corps totalement r\'eels. D\'etermination du "signe" des repr\'esentations
galoisiennes automorphes polaris\'ees sur les corps CM.  \ps

\item[-] \'Etude de l'adh\'erence Zariski des repr\'esentations galoisiennes
automorphes dans la vari\'et\'e des caract\`eres $p$-adique et analogues locaux
cristallins.  \ps \end{itemize} \ps\ps

Le second versant de mes recherches, expos\'e dans la seconde partie, 
a \'et\'e port\'e par la volont\'e de comprendre ce que disent les conjectures 
d'Arthur-Langlands sur les objets les plus familiers de la th\'eorie des nombres.\ps 
En effet, les progr\`es r\'ecents en th\'eorie des formes automorphes \`a plusieurs 
variables \'etablissent une partie de ces conjectures, et donnent du corps \`a ce r\^eve 
imagin\'e par Langlands. Je pense notamment \`a la d\'emonstration par M. Harris et R. Taylor, 
puis G. Henniart, de la conjecture de
Langlands locale, celle par G. Laumon, B.C. Ng\^o et J.-L. Waldspurger du
lemme fondamental de Langlands, ou encore la classification par Arthur du
spectre discret automorphe des groupes classiques. Malgr\'e l'incroyable 
complexit\'e des d\'emonstrations, ces \'enonc\'es sont fascinants par leur extr\^eme beaut\'e, 
et d'une certaine fa\c{c}on, par leur simplicit\'e. Ils forment aussi un outil nouveau formidable pour l'arithm\'eticien. 
La fra\^iche d\'emonstration de la conjecture de Sato-Tate en est une illustration frappante.  \ps

Dans ce texte, je vais exposer quelques autres applications plus modestes de
cet outil que j'ai d\'ecouvertes dans mes recherches, souvent guid\'e par
mon attrait pour les objets "peu ramifi\'es" en th\'eorie des nombres :
\begin{itemize}\ps\ps

\item[-] Construction de corps de nombres \`a ramification prescrite,
structure du groupe de Galois d'une extension alg\'ebrique
maximale de $\Q$ non ramifi\'ee hors d'un ensemble fini de nombres premiers fix\'es.\ps

\item[-] \'Enum\'eration des voisins de Kneser des r\'eseaux unimodulaires pairs en rang $\leq 24$. Applications aux conjectures de Harder et de Nebe-Venkov. \ps

\item[-] \'Enum\'eration conjecturale des motifs purs sur $\Q$ de conducteur $1$
(r\'eguliers, polaris\'es) en petit rang, calculs de dimensions d'espaces de formes automorphes 
pour les petits groupes r\'eductifs sur $\Z$. \ps

\item[-] Construction par voie automorphe d'extensions galoisiennes pr\'evues par les
conjectures de Bloch et Kato. \ps
\end{itemize}
\noindent

\bigskip \bigskip

\noindent {\sc Bonne lecture !}

\newpage

\begin{center}{\sc Travaux pr\'esent\'es }\end{center}\ps
\bigskip

Je vais pr\'esenter dans ce m\'emoire les r\'esultats des articles  publi\'es ou \`a para\^itre [1] \`a
[13] suivants, ainsi que ceux des pr\'epublications [A] \`a [D]. Ils sont  
disponibles sur mon site internet~\url{http://www.math.polytechnique.fr/~chenevier/}. Je renvoie \`a la bibliographie g\'en\'erale pour la liste compl\`ete de mes publications.\ps\ps

\smallskip
{\small 

[13]  {\it Construction of automorphic Galois representations II},  avec M. Harris,  \`a para\^itre \`a Cambridge Journal of Math. 1 (2013). \ps \ps

[12] {\it Sur la densit\'e des repr\'esentations cristallines du groupe de Galois absolu de $\mathbb{Q}_p$}, \`a para\^itre \`a Math. Annalen. \ps \ps 

[11]  {\it On the infinite fern of Galois representations of unitary type}, Ann. Sci. \'E.N.S. 44, 963-1019 (2011).\ps \ps 

[10] {\it The sign of Galois representations attached to automorphic forms for unitary groups}, avec J. Bella\"iche, 
Compositio Math. 147, 1337-1352 (2011).\ps \ps

[9]  {\it The $p$-adic analytic space of pseudocharacters of a profinite group, and pseudorepresentations over arbitrary rings}, \`a para\^itre aux {\it Proceedings of the LMS Durham Symposium}, {\it Automorphic forms and Galois representations} (2011).\ps\ps

[8]  {\it On the vanishing of some non semisimple orbital integrals}, avec D. Renard, 
Expo. Math. 28, 276-289 (2010). \ps \ps

[7]  {\it Families of Galois representations and Selmer groups},  avec J. Bella\"iche, 
Ast\'erisque 324, Soc. Math. France (2009), 331 pages.\ps \ps

[6] {\it Corps de nombres peu ramifi\'es et formes automorphes autoduales}, avec L. Clozel, 
Journal of the A.M.S. 22, 467-519 (2009).\ps \ps


[5]  {\it Quelques courbes de Hecke se plongent dans l'espace de Colmez},
Journal of Number Theory 128, 2430-2449 (2008). \ps \ps

[4] {\it On number fields with given ramification}, 
Compositio Math. 143, 1359-1373 (2007).\ps \ps


[3]  {\it Lissit\'e de la courbe de Hecke de GL(2) aux points Eisenstein critiques}, avec J. Bella\"iche, 
Journal Inst. Math. Jussieu 5 (2006).\ps \ps

%

[2]  {\it Formes non temp\'er\'ees pour U(3) et conjectures de Bloch-Kato}, avec J. Bella\"iche, 
Ann. Sci. \'E.N.S. 37, 611-662 (2004).\ps \ps

[1] {\it Familles $p$-adiques de formes automorphes pour GL(n)}, 
Journal f\"ur die reine und angewandte Mathematik 570, 143-217 (2004). \ps \ps

[D] {\it Level one algebraic cuspforms of classical groups of small ranks}, avec D. Renard, pr\'epublication (2012).\ps \ps

[C] {\it Kneser neighbours and orthogonal Galois representations in dimensions $16$ and $24$}, avec J. Lannes, 
Algebraic Number Theory, Oberwolfach Report 31 (2011).  \ps \ps

[B] {\it Sur la vari\'et\'e des caract\`eres $p$-adique du groupe de Galois absolu de $\mathbb{Q}_p$}, pr\'epublication (2009).\ps \ps

[A] {\it Une application des vari\'et\'es de Hecke des groupes unitaires}, article pour le tome 2 du projet de livre dirig\'e par Michael Harris, pr\'epublication (2009). \ps \ps

}

\mainmatter

\chapter{La foug\`ere infinie des repr\'esentations galoisiennes et applications}

\section{Repr\'esentations galoisiennes automorphes}\label{repgalaut}

\subsection{Notations}\label{notations} Soit $\ell$ un nombre premier et soit
$\overline{\Q}_\ell$ une cl\^oture alg\'ebrique\footnote{ On d\'esignera
g\'en\'eralement par $\overline{k}$ une cl\^oture alg\'ebrique du corps
$k$.} du corps $\Q_\ell$ des nombres $\ell$-adiques.  Nous allons dans ce
chapitre consid\'erer des repr\'esentations continues $${\rm Gal}(\overline{F}/F)
\rightarrow {\rm GL}_n(\overline{\Q}_\ell),$$ 
o\`u $F$ est soit un corps de nombres, soit un corps local, de groupe de
Galois absolu ${\rm Gal}(\overline{F}/F)$. Comme tout groupe de Galois
absolu, ce dernier n'est bien d\'efini que si l'on fixe une cl\^oture
alg\'ebrique $\overline{F}$ de $F$, sinon il ne l'est qu'\`a automorphismes
int\'erieurs pr\'es. Cette ambig\"uit\'e n'aura cependant jamais d'incidence car seules
ses classes de conjugaison d'\'el\'ements, ou encore les classes 
d'isomorphie de ses repr\'esentations, nous int\'eresserons.  \ps

On rappelle que si $v$ est une place d'un corps de nombres $F$, de
compl\'etion associ\'ee $F_v$, on dispose d'une classe de conjugaison
canonique de morphismes de groupes $${\rm Gal}(\overline{F}_v/F_v)
\rightarrow {\rm Gal}(\overline{F}/F)$$ continus et injectifs, dont les
images sont les {\it sous-groupes de d\'ecomposition en $v$} de ${\rm
Gal}(\overline{F}/F)$.  Si $\rho$ est une repr\'esentation de ${\rm
Gal}(\overline{F}/F)$, on notera ${\rho}_{|{\rm Gal}(\overline{F}_v/F_v)}$,
ou simplement $\rho_v$, la compos\'ee de $\rho$ par un tel homomorphisme.  La
repr\'esentation $\rho_v$ est bien d\'efinie \`a isomorphisme pr\`es.  \ps

Nous ferons dans ce texte un usage constant de la param\'etrisation de
Langlands des repr\'esentations admissibles irr\'eductibles de $\GL_n(F)$
quand $F$ est un corps local (i.e. pour nous $F \simeq \R, \C$ ou une extension
finie de $\Q_p$). Pour un tel corps $F$, on d\'esignera par $W_F$ le groupe de Weil de
$F$~\cite{tate} et par $W'_F$ son groupe de Weil-Deligne. Rappelons que 
$W'_F=W_F \times {\rm SU}(2)$ si $F$ est non
archim\'edien et $W'_F=W_F$ sinon. Ce sont des groupes topologiques
localement compacts. Si $\pi$ est une repr\'esentation complexe admissible
irr\'eductible de $\GL_n(F)$, il lui est associ\'e de mani\`ere naturelle un param\`etre de
Langlands, c'est-\`a-dire 
une repr\'esentation semisimple continue
$${\rm L}(\pi) : W'_{F} \rightarrow {\rm GL}_n(\C).$$ 
Cette param\'etrisation est due \`a Langlands~\cite{langlandsreel} quand
$F$ est archim\'edien, \`a Harris-Taylor~\cite{HT} et
Henniart~\cite{henniartll} sinon. Mieux, l'application $\pi \mapsto {\rm L}(\pi)$ induit une bijection entre classes
d'isomorphie de repr\'esentations irr\'eductibles admissibles de $\GL_n(F)$
et classes d'isomorphie de repr\'esentations complexes semisimples continues de
dimension $n$ de $W'_F$. \ps

Le cas \`a la fois le plus simple et le plus utile de cette correspondance est
celui o\`u $F$ est non archim\'edien et $\pi$ est
non ramifi\'ee, c'est-\`a-dire poss\`edant des invariants non-triviaux sous
$\GL_n(\OO_F)$, o\`u $\OO_F$ d\'esigne l'anneau des entiers de $F$.  Ces
repr\'esentations correspondent alors aux param\`etres, aussi dit non-ramifi\'es,  qui
se factorisent par le quotient canonique $W'_F \rightarrow \Z$ envoyant sur
$1$ tout rel\`evement du Frobenius g\'eom\'etrique du corps r\'esiduel de $F$. Ces
param\`etres sont donc uniquement d\'etermin\'es par l'image de $1$, qui est
un \'el\'ement semisimple quelconque de $\GL_n(\C)$. On retrouve ainsi la
param\'etrisation de Satake~\cite{satake} des repr\'esentations non ramifi\'ees de
$\GL_n(F)$ par les classes de conjugaison semisimples dans $\GL_n(\C)$, dans
l'esprit de Langlands~\cite{langlandsweil},~\cite{langlandsyale}.  \ps
 
Lorsque $F$ est non archim\'edien, les 
repr\'esentations complexes continues de dimension $n$ de $W'_F$ sont en bijection naturelle avec
les repr\'esentations dites de Weil-Deligne de $F$, dans lesquelles
l'action du ${\rm SU}(2)$ est remplac\'ee par un op\'erateur nilpotent
("monodromie"), et pour lesquelles nous renvoyons \`a~\cite{tate}. 
Les repr\'esentations de Weil-Deligne ont un sens \`a
coefficients dans n'importe quel corps. Lorsque ce coefficient est 
$\overline{\Q}_\ell$, et que $\ell$ n'est pas la caract\'eristique du corps
r\'esiduel de $F$, le th\'eor\`eme de mondromie $\ell$-adique de
Grothendieck met en bijection naturelle repr\'esentations continues $W_F
\rightarrow \GL_n(\overline{\Q}_\ell)$ et repr\'esentations de Weil-Deligne
de $F$ qui sont de dimension $n$ et \`a coefficients dans
$\overline{\Q}_\ell$ (voir {\it loc.cit.}). \ps

\subsection{Quelques conjectures du folklore}\label{conjfolk} Je rappelle
dans ce paragraphe quelques conjectures du folklore reliant repr\'esentations
galoisiennes et repr\'esentations automorphes. \ps

Soit $F$ un corps de nombres et soit $\pi$ une repr\'esentation automorphe cuspidale de
$\GL_n$ sur $F$. On dira que la repr\'esentation $\pi$ est {\it
alg\'ebrique} si pour tout plongement complexe $\sigma : F \rightarrow \C$, de place archim\'edienne sous-jacente $v$, le param\`etre de Langlands de $\pi_v$ restreint
\`a $W_{\overline{F_v}}$, identifi\'e via $\sigma$ \`a $W_{\rm \C}=\C^\ast$, est une somme directe de caract\`eres  $$z
\mapsto z^{p_{i,\sigma}} \overline{z}^{q_{i,\sigma}}$$ avec $p_{i,\sigma},q_{i,\sigma} \in \Z$
et $i=1,\cdots,n$.  Sous cette simple hypoth\`ese, on conjecture depuis
Weil, Eichler, Shimura, Serre et Langlands que les composantes locales aux places finies de
$\pi$ sont d\'efinies sur le sous-corps $\overline{\Q} \subset \C$ des
nombres alg\'ebriques.  Admettant cela, on conjecture de plus que pour tout
plongement $\iota : \overline{\Q} \rightarrow \overline{\Q}_\ell$ il existe
une repr\'esentation continue irr\'eductible $$\rho_{\pi,\iota}: {\rm
Gal}(\overline{F}/F) \longrightarrow {\rm GL}_n(\overline{\Q}_\ell)$$
satisfaisant la condition suivante : \begin{itemize}\medskip

\item[(C1)]  {\it Pour toute place finie $v$ de $F$ ne divisant pas $\ell$ et telle que $\pi_v$ est non ramifi\'ee, 
la repr\'esentation $(\rho_{\pi,\iota})_{|{\rm Gal}(\overline{F}_v/F_v)}$ est non ramifi\'ee et le polyn\^ome 
caract\'eristique de $\rho_{\pi,\iota}({\rm Frob}_v)$ co\"incide avec l'image par $\iota$ du polyn\^ome 
caract\'eristique du param\`etre de Satake de $\pi_v$}.\medskip
\end{itemize}

La notation ${\rm Frob}_v$ d\'esigne ici un rel\`evement quelconque du
Frobenius g\'eom\'etrique d'un groupe de d\'ecomposition en $v$ de ${\rm
Gal}(\overline{F}/F)$.  Observons que d'apr\`es le th\'eor\`eme de
Cebotarev et la th\'eorie des caract\`eres des repr\'esentations lin\'eaires, toute repr\'esentation semisimple ${\rm Gal}(\overline{F}/F)
\longrightarrow {\rm GL}_n(\overline{\Q}_\ell)$ satisfaisant (C1) est
n\'ecessairement unique \`a isomorphisme pr\`es.\ps

Il convient de signaler que certains auteurs, comme Borel dans~\cite[\S
18.2]{borelcorvallis}, et suivant Weil, utilisent le terme "de type $A_0$" pour "alg\'ebrique". 
Langlands utilise parfois aussi le terme "de type Hodge"~\cite{langlandsedin}. Mentionnons aussi l'article de Clozel~\cite{clozel}, qui traite notamment la question du corps de d\'efinition des $\pi_v$ 
avec $v$ finie sous l'hypoth\`ese que les $\pi_v$ avec $v$ archim\'ediens sont de type "cohomologique". Cela le conduit d'ailleurs \`a appeler alg\'ebrique un $\pi$ tel que $\pi \otimes |\det(\cdot)|^{\frac{n-1}{2}}$ l'est au sens ci-dessus. Nous
renvoyons \`a l'article r\'ecent de Buzzard et Gee~\cite{buzzardgee} pour des clarifications sur les diff\'erentes notions d'alg\'ebricit\'e.  \ps

Supposant l'existence (et donc unicit\'e) de $\rho_{\pi,\iota}$ satisfaisant (C1) acquise, on dispose aussi d'un
suppl\'ement conjectural pr\'ecisant cette repr\'esentation \`a toutes les
places, qu'il sera commode dans la discussion qui suit de d\'ecouper en trois
\'enonc\'es.  \ps

\begin{itemize}\medskip

\item[(C2)] {\it Pour toute place finie $v$ de $F$ ne divisant pas $\ell$,
la repr\'esentation de Weil-Deligne Frobenius-semisimplifi\'ee associ\'ee
\`a\footnote{Strictement, il aurait fallu \'ecrire "associ\'e \`a la
restriction de $(\rho_{\pi,\iota})_{|{\rm Gal}(\overline{F}_v/F_v)}$ \`a $W_{F_v} \subset {\rm
Gal}(\overline{F}_v/F_v)$...} $(\rho_{\pi,\iota})_{|{\rm Gal}(\overline{F}_v/F_v)}$ est l'image par
$\iota$ de la repr\'esentation de Weil-Deligne d\'efinie par ${\rm L}(\pi_v)$}. 
\medskip

\item[(C3)] {\it Pour toute place $v$ de $F$ divisant $\ell$, la
repr\'esentation $(\rho_{\pi,\iota})_{|{\rm Gal}(\overline{F}_v/F_v)}$ est
de De Rham au sens de Fontaine~\cite{PP}, et la repr\'esentation de Weil-Deligne
Frobenius-semisimplifi\'ee associ\'ee \`a son ${\rm D}_{\rm pst}$, est
l'image par $\iota$ de la repr\'esentation de Weil-Deligne d\'efinie par
${\rm L}(\pi_v)$ .}\medskip

\item[(C4)] {\it Pour toute place r\'eelle $v$ de $F$, la trace d'une
conjugaison complexe en $v$ est de la forme $a_v-b_v$ o\`u $a_v+b_v$ est le
nombre de caract\`eres (de dimension $1$) de ${\rm W}_\R$ apparaissant dans
${\rm L}(\pi_v)$, et o\`u $a_v$ est le nombre de ces caract\`eres qui sont de la
forme $|.|^{2m}$ ou $\varepsilon |.|^{2m+1}$avec $m \in \Z$}.  \medskip

\end{itemize}

Dans (C4), $\varepsilon$ d\'esigne le caract\`ere d'ordre $2$ de ${\rm W}_\R$. 
Pour des raisons d'exhaustivit\'e, mentionnons que l'on
conjecture d'apr\`es Tate que les repr\'esentations de Weil-Deligne
associ\'ees aux $(\rho_{\pi,\iota})_{|{\rm Gal}(\overline{F}_v/F_v)}$ sont
d\'ej\`a Frobenius-semisimples quand $\ell$ est premier \`a $v$.  De plus, la description compl\`ete de
$\rho_{\pi,\iota}$ aux places restantes, i.e.  divisant $\ell$, est l'objet
du {\it programme de Langlands $p$-adique} initi\'e par Breuil~\cite{breuil},~\cite{colmezecm}. La description des nombres de Hodge-Tate est en revanche \'el\'ementaire. 
En effet, si $v$ est une telle place, et si $\sigma : F_v \rightarrow \overline{\Q}_\ell$ est un plongement, alors $\iota$ et $\sigma$ d\'eterminent un plongement $\sigma' : F \rightarrow \C$. 
Les poids de Hodge-Tate de $(\rho_{\pi,\iota})_v$ associ\'es \`a $\sigma$ sont alors les $- p_{i,\sigma'}$, $i=1,\cdots,n$, sous la convention que le poids de Hodge-Tate du caract\`ere cyclotomique est $-1$. \ps

Dans cette g\'en\'eralit\'e, l'existence de $\rho_{\pi,\iota}$ satisfaisant
(C1) \`a (C4) n'est connue que pour $n=1$, par une construction fameuse
due \`a Weil issue de la th\'eorie du corps de classes.  Pour $n>1$, elle a
fait l'objet de tr\`es nombreux travaux, notamment d'Eichler, Shimura,
Igusa, Kuga, Deligne, Rappoport, Serre, Katz, Carayol, Mazur, Wiles, Takeshi Saito,
et Faltings dans le cas $n=2$ et $F=\Q$, auquel cas seule l'\'eventualit\'e
o\`u ${\rm L}(\pi_\infty)$ est purement scalaire \'echappe encore \`a notre
compr\'ehension ("formes de Maass alg\'ebriques").  L'action galoisienne sur la cohomologie
$\ell$-adique des courbes modulaires, et plus g\'en\'eralement des courbes
de Shimura, joue ici un r\^ole fondamental dans la construction des
$\rho_{\pi,\iota}$.  Toutes les g\'en\'eralisations qui ont suivi de ces
travaux ont consist\'e en une \'etude de plus en plus fine de l'action
galoisienne sur la cohomologie des vari\'et\'es de Shimura de dimension
sup\'erieure, qui s'est av\'er\'ee notablement plus difficile d'une part
pour des raisons g\'eom\'etriques, mais aussi pour les questions de
$L$-indistinguabilit\'e (Langlands, Kottwitz~\cite{kottwitz}).  \ps

Sans doute faut-il mentionner, avant d'entrer dans la discussion de ces
g\'en\'eralisations, que la conjecture ci-dessus admet une r\'eciproque non
moins fabuleuse (et tout aussi conjecturale), due \`a Fontaine et
Mazur~(\cite{fontainemazur}) : les repr\'esentations
$\rho_{\pi,\iota}$ ainsi obtenues sont exactement les repr\'esentations
irr\'eductibles continues ${\rm Gal}(\overline{F}/F) \longrightarrow {\rm
GL}_n(\overline{\Q}_\ell)$ qui sont non ramifi\'ees hors d'un ensemble fini
de places, et dont la restriction aux groupes de d\'ecompositions aux places
$v$ divisant $\ell$ est de De Rham au sens de Fontaine. L\`a encore,
cette conjecture n'est connue en toute g\'en\'eralit\'e qu'en dimension
$1$~\cite{fontainemazur}. En dimension $2$ et pour $F=\Q$, elle est
corrobor\'ee par la d\'emonstration de la conjecture de
Shimura-Taniyama-Weil (Wiles~\cite{wilesfermat}, Taylor-Wiles~\cite{tw},
Breuil-Conrad-Diamond-Taylor~\cite{bcdt}), et elle est maintenant
d\'emontr\'ee dans de nombreux cas (Khare-Winteberger~\cite{kw}, Kisin~\cite{kisinfm},
Emerton~\cite{emertonfm}). \ps

\subsection{Cas des repr\'esentations r\'eguli\`eres
polaris\'ees}\label{reppol} Une
repr\'esentation automorphe cuspidale alg\'ebrique $\pi$ de $\GL_n$ sur $F$
est dite {\it r\'eguli\`ere}, si pour tout place archim\'edienne $v$ de $F$,
les entiers $p_{i,\sigma}$, $i=1,\cdots,n$, d\'efinis plus haut sont
distincts.  Nous dirons aussi que $\pi$ est {\it polaris\'ee} si l'une des
deux conditions suivantes est satisfaite : \ps \begin{itemize} \item[-] Soit
$F$ est un corps totalement r\'eel et $\pi^\vee \otimes \chi$ est isomorphe
\`a $\pi$, o\`u $\chi=\prod_v \chi_v$ est un caract\`ere de Hecke de $F$ tel que le signe
$\chi_v(-1)$ est ind\'ependant\footnote{
Nous verrons plus loin que cette condition est conjecturalement automatique. C'est \'evident si $n$ est impair car $\chi^{n}$ est le carr\'e du caract\`ere central de $\pi$, 
et donc $\chi_v(-1)=1$ pour tout $v$.  La notation $\pi \otimes \chi$ d\'esigne abusivement  $\pi \otimes \chi \circ \det$.}  de la place r\'eelle $v$ de $F$.\ps \item[-]
Soit $F$ est une extension quadratique totalement imaginaire d'un corps
totalement r\'eel $F^+$ (i.e.  $F$ est un corps CM) et $\pi^\vee \otimes
|\cdot|^{1-n}$ est isomorphe \`a $\pi^c$, o\`u $c$ est le g\'en\'erateur de
${\rm Gal}(F/F^+)$ et o\`u $\pi^c$ d\'esigne le conjugu\'e ext\'erieur de $\pi$
par $c$.  \ps \end{itemize}

\begin{thm}\label{existencegalois} Soit $\pi$ une repr\'esentation automorphe cuspidale
alg\'ebrique, r\'eguli\`ere et polaris\'ee\footnote{Signalons que l'on trouve
dans la litt\'erature les acronymes RAESDC (cas totalement r\'eel) et RASDC (cas $CM$) pour les repr\'esentations de
l'\'enonc\'e. } de $\GL_n$ sur $F$.  Alors pour
tout plongement $\iota : \overline{\Q} \rightarrow \overline{\Q}_\ell$ il
existe une unique repr\'esentation semisimple $\rho_{\pi,\iota}$
satisfaisant (C1), (C2) et (C3).  \end{thm}

Ce th\'eor\`eme est la synth\`ese de travaux de nombreux math\'ematiciens
dont je vais t\^acher de citer correctement les contributions dans ce qui
suit. J'ai eu la chance, en grande partie offerte par mon directeur de th\`ese M.  Harris, de
participer \`a cet \'edifice, bien que de mani\`ere infime, par mon
travail~\cite{chgrfa} et mon article avec Harris~\cite{chharris}, que je
d\'ecrirai bri\`evement. Il faut rajouter que le th\'eor\`eme ci-dessus
repose aussi de mani\`ere fondamentale sur les travaux d'autres
auteurs non directement associ\'es \`a son \'enonc\'e, je pense notamment
\`a Deligne, Grothendieck, Kottwitz, Langlands, Laumon, Ng\^o et Walsdpurger. \ps

Quand $n \leq 3$, l'existence de $\rho_{\pi,\iota}$ satisfaisant (C1) et
partiellement (C2) est due \`a Blasius et Rogawski, dans le cadre du projet
collectif sur la cohomologie des surfaces de Picard~\cite{picard}. La classification endoscopique, par Rogawski, des repr\'esentations automorphes des groupes unitaires 
\`a trois variables joue un r\^ole tr\`es important dans ce travail~\cite{roglivre}. Quand
$n=2$, Blasius et Rogawski doivent rajouter une hypoth\`ese de r\'egularit\'e
suppl\'ementaire sur $\pi$ aux places archim\'ediennes, qui ne sera enlev\'ee que plus tard par Wiles
puis Taylor par un argument de congruences. \ps

Pour $n$ g\'en\'eral, si $\pi$ est de carr\'e int\'egrable modulo le
centre \`a au moins une place finie, l'existence de $\rho_{\pi,\iota}$
satisfaisant (C1) est due \`a Harris et Taylor~\cite{HT} (voir
aussi le compl\'ement de Taylor~\cite{taylortoulouse} dans le cas totalement
r\'eel). Elle est d\'ecoup\'ee dans la
cohomologie $\ell$-adique des "vari\'et\'es de Shimura simples" au sens de Kottwitz, qui sont
des quotients arithm\'etiques cocompacts de la boule unit\'e
hermitienne~$\sum_{i=1}^{n-1} |z_i|^2<1$ dans $\C^{n-1}$ par certains r\'eseaux de ${\rm PU}(n-1,1)$. Ces travaux g\'en\'eralisent des constructions
ant\'erieures dans des cas partiels dues \`a Clozel~\cite{clozelgal} (qui constituaient les premiers cas non triviaux connus en dimension g\'en\'erale),
bas\'ees notamment sur des travaux de Kottwitz~\cite{kottwitzinvent}. Harris et Taylor d\'emontrent
en m\^eme temps la correspondance de Langlands locale, ainsi que 
(C2) \`a semisimplification pr\`es, affaiblissement supprim\'e plus tard par Taylor et
Yoshida~\cite{TY}. 
 \ps

La construction de $\rho_{\pi,\iota}$ sans l'hypoth\`ese de "carr\'e
int\'egrabilit\'e \`a au moins une place finie", souvent g\'enante dans les
applications arithm\'etiques, par exemple dans les probl\`emes de
constructions d'extensions galoisiennes \`a la Ribet~\cite{ribet}, \'etait
notamment l'objet du projet de livre~\cite{projetlivregrfa} initi\'e par M. Harris. Il
a fait aussi l'objet de travaux ind\'ependants de Shin~\cite{shin} et de
Morel~\cite{morel}. Ces travaux ont notamment \'et\'e permis par la d\'emonstration du lemme fondamentale de Langlands par Laumon-Ng\^o~\cite{laumonngo}, Ng\^o~\cite{ngo}, et Waldspurger~\cite{waldspurger1},~\cite{waldspurger2}, qui rend inconditionnelle la stabilisation de la formule des traces pour les groupes unitaires. En effet, c'est une \'etape n\'ecessaire pour relier le calcul du nombre de points dans les corps fini des vari\'et\'es de Shimura "unitaires" et les param\`etres de Satake des formes automorphes mises en jeu (Langlands, Kottwitz). \ps

D'une part, on savait depuis les travaux de Blasius-Rogawski ("patching lemma") et Harris-Taylor que l'on peut ramener
le cas g\'en\'eral du Th\'eor\`eme~\ref{existencegalois} au cas CM par changements de bases quadratiques
(Arthur-Clozel). Harris et Labesse ont m\^eme observ\'e que l'on peut
se ramener, dans le cas CM, \`a la situation o\`u $F/F^+$ est non ramifi\'ee
\`a toutes les places finies et d\'ecompos\'ee aux places finies o\`u $\pi$
est ramifi\'ee, hypoth\`ese simplificatrice notamment pratique pour les questions de descente et de changement de base~\cite{labesse} entre groupes unitaires et groupes lin\'eaires.  Dans ce cas, le th\'eor\`eme est d\^u ind\'ependamment \`a
Shin~\cite{shin} et \`a Clozel-Harris-Labesse~\cite{chla},~\cite{chlb}, en toute
g\'en\'eralit\'e quand $n$ est impair, et sous une hypoth\`ese parasite
quand $n$ est pair, de mani\`ere analogue au cas $n\leq 3$. \ps
Clozel, Harris et Labesse demandent que pour tout $\sigma : F \rightarrow \C$
on ait $|p_{i,\sigma} -p_{j,\sigma}|\geq 2$ si $i \neq j$,
et pour Shin il est m\^eme suffisant d'imposer cette condition \`a une seule
place archim\'edienne. Il serait trop long d'expliquer ici en d\'etail
la raison de l'apparition de cette condition dans le cas $n$ pair.  En deux
mots, cela vient de l'obstruction au principe de Hasse pour les groupes
unitaires \`a un nombre pair de variables et des contraintes donn\'ees par la formule de
multiplicit\'e \`a la Labesse-Langlands pour le transfert temp\'er\'e de ${\rm U}(n)
\times {\rm U}(1)$ vers ${\rm U}(n+1)$. On verra par ailleurs de multiples
analogues de ces contraintes dans ce m\'emoire quand il sera question de mes
travaux avec Lannes et avec Renard.\ps

Ma contribution principale est d'avoir donn\'e avec Harris
dans~\cite{chharris} un argument permettant de d\'eduire des cas
pr\'ec\'edents le cas g\'en\'eral, i.e.  de construire $\rho_{\pi,\iota}$
satisfaisant (C1), ainsi que (C2) \`a semisimplification pr\`es quand $n$
est pair\footnote{Dans ce cas on a tout de m\^eme une majoration de la
taille de la monodromie.  Ainsi que nous l'avions observ\'e avec J. 
Bella\"iche dans notre livre\cite{bchlivre}, c'est d'ailleurs tout ce dont
on a besoin pour les applications aux constructions d'extensions
galoisiennes \`a la Ribet~\cite{ribet}, ou pour les d\'emonstrations \`a la
Wiles de la conjecture principale d'Iwasawa~\cite{wilesmc}.}, sans autres
hypoth\`eses que celles de l'\'enonc\'e.  La propri\'et\'e (C2) a \'et\'e
par la suite v\'erifi\'ee en toute g\'en\'eralit\'e par
Caraiani~\cite{caraiani1}.  Ce cas $n$ pair "sans hypoth\`ese de
r\'egularit\'e suppl\'ementaire" pr\'esentait un certain int\'er\^et \`a
l'\'epoque, car c'\'etait un ingr\'edient manquant important pour faire
marcher la d\'emonstration de la conjecture de Sato-Tate
(\cite{hsbt},\cite{cht},\cite{taylorsatotate}) dans le cas des courbes
elliptiques sur $\Q$ de $j$-invariant dans $\Z$.  Une raison simple est que
les ${\rm Sym}^{2k+1}$ des courbes elliptiques ne satisfont aucune des
conditions de r\'egularit\'e ci-dessus.  \ps

Dans l'article avec Harris, nous commen\c{c}ons par nous ramener, en faisant
des changements de base cycliques ad\'equats~\cite{arthurclozel} et via une
g\'en\'eralisation du "patching Lemma" de Blasius-Rogawski d\^u \`a Harris
et Sorensen~\cite{sorensen}, au cas o\`u $\pi$ admet des vecteurs invariants
par un sous-groupe d'Iwahori en toutes les places de $F$ divisant $\ell$, chacune de ces places \'etant de plus d\'ecompos\'ee sur $F^+$.  Sous cette hypoth\`ese, nous utilisons alors un
argument de d\'eformation $\ell$-adique pour trouver des congruences avec
des repr\'esentations $\pi'$ satisfaisant les m\^emes hypoth\`eses que $\pi$
mais suffisament r\'eguli\`eres aux places archim\'ediennes.  J'avais
construit de telles congruences dans ma th\`ese~\cite{chcrelle} (\`a la
demande de Harris!) \`a ceci pr\`es que mon corps de base $F$ \'etait choisi
quadratique imaginaire pour simplifier.  \ps

Dans l'article~\cite{chgrfa}, j'ai repris ce travail en m'affranchissant de
cette hypoth\`ese sur $F$, une t\^ache d'ordre essentiellement technique. 
Les d\'emonstrations passent par une descente au groupe unitaire \`a $n$
variables relativement \`a $F/F^+$ qui est compact \`a toutes les places
r\'eelles et quasi-d\'eploy\'e aux places finies (cela existe si $[F^+:\Q]$
est pair, ce que l'on peut toujours supposer).  Les congruences cherch\'ees
sont alors cons\'equences triviales de ma construction des vari\'et\'es de
Hecke dans ce contexte, i.e.  de la foug\`ere infinie, sur laquelle je
reviendrai en d\'etail au~\S\ref{fougereunitaire}.  La propri\'et\'e (C2) \`a semisimplification pr\`es
s'obtient par un argument simple, d\'ecouvert dans~\cite{chgrfa}, bas\'e sur
la construction d'un pseudo-caract\`ere universel de ${\rm W}_E$ dans le
centre de Bernstein de $\GL_n(E)$ quand $E$ est un corps $p$-adique. 
L'in\'egalit\'e susmentionn\'ee concernant la monodromie d\'ecoule
d'observations simples sur les familles de repr\'esentations au sens des
traces faites dans mon livre avec Bella\"iche~\cite{bchlivre}.

\ps

Dans tous les cas o\`u la construction de $\rho_{\pi,\iota}$ est faite par
voie g\'eom\'etrique, on obtient "gratuitement" une version affaiblie de
(C3) sur $\rho_{\pi,\iota}$ comme cons\'equence des th\'eor\`emes de
comparaisons en g\'eom\'etrie arithm\'etique, notamment la propri\'et\'e
d'\^etre de De Rham ainsi que le polyn\^ome caract\'eristique du Frobenius
cristallin dans le cas non ramifi\'e.  Une des innovations de mon
article~\cite{chgrfa} est d'avoir compris comment pr\'eserver la
propri\'et\'e (C3) dans l'argument de d\'eformation.  En effet, il est bien
connu que les repr\'esentations galoisiennes construites par limites dans
des familles $p$-adiques comme celles apparaissant dans la foug\`ere infinie
ne sont pas en g\'en\'eral de De Rham, loin s'en faut, m\^eme quand elles
sont de Hodge-Tate.  Mon astuce est de d\'edoubler tout d'abord par
changement de base quadratique r\'eel de $F^+$ toutes les places de $F^+$
divisant $\ell$, puis de proc\'eder par un argument de d\'eformations en
fixant les poids attach\'es \`a toutes les places archim\'ediennes sauf
une\footnote{C.  Skinner m'a inform\'e qu'il avait eu ind\'ependamment une
id\'ee similaire}.  Je conclus par les travaux de
Berger-Colmez~\cite{bergercolmez} sur les familles analytiques de
repr\'esentations "\`a poids de Hodge-Tate fix\'es", dont le lieu de De Rham
est un ferm\'e analytique comme le d\'emontrent ces auteurs.  Je note
d'ailleurs qu'il est important dans cet argument d'utiliser les travaux de
Shin~\cite{shin} (qui ne n\'ecessite que de la r\'egularit\'e \`a une place archim\'edienne)
plut\^ot que ceux de~\cite{chlb}.  \ps

La propri\'et\'e (C3) a ensuite \'et\'e d\'emontr\'ee en toute
g\'en\'eralit\'e par Barnet-Lamb, Gee, Geraghty et Taylor~\cite[I \&
II]{blgt} et Caraiani~\cite{caraiani2}.  Mentionnons que
l'irr\'eductibilit\'e de $\rho_{\pi,\iota}$ n'est pas connue en g\'en\'eral,
sauf si $n \leq 3$ dans le cas CM (Ribet, Blasius-Rogaswki), et si $n\leq 5$
dans le cas totalement r\'eel par Calegari et Gee~\cite{calegarigee}.  Des
travaux r\'ecents de Taylor et Patrikis montrent qu'en g\'en\'eral elle vaut
pour un ensemble de densit\'e $1$ de "choix de $\iota$" (cela n'a de sens
que si l'on fait intervenir un corps de coefficients de $\pi$).\ps

Peut-\^etre paradoxallement, la propri\'et\'e (C4) concernant la
repr\'esentation $\rho_{\pi,\iota}$ a \'et\'e \'egalement r\'esistante. 
Ce n'est pas \'etonnant si l'on se rappelle que quand $F$ est totalement
r\'eel, la repr\'esentation $\rho_{\pi,\iota}$ est construite par
r\'eduction au cas CM, c'est-\`a-dire par recollement \`a partir de ses
restrictions aux ${\rm Gal}(\overline{F}/E)$ o\`u $E/F$ est totalement
imaginaire quadratique : la seule perte d'information concerne les places
archim\'ediennes.  Avant de donner un \'enonc\'e sur (C4) je voudrais
discuter d'une autre propri\'et\'e cach\'ee mais tr\`es utile de la
repr\'esentation $\rho_{\pi,\iota}$.\ps

\subsection{L'alternative symplectique-orthogonale}\label{signegalois} 

Soit $F$ un corps de nombres totalement r\'eel. Supposons que $\pi$ est une
repr\'esentation automorphe cuspidale de $\GL_n$ sur $F$ qui est
alg\'ebrique, r\'eguli\`ere et polaris\'ee.  Pour chaque $\iota :
\overline{\Q} \rightarrow \overline{\Q}_\ell$ on dispose donc d'une
repr\'esentation $\rho_{\pi,\iota}$ comme dans le
Th\'eor\`eme~\ref{existencegalois}.  Soit $\chi$ un caract\`ere de Hecke de
$F$ tel que $\pi^\vee \otimes \chi \simeq \pi $ et tel que le signe
$$\chi_\infty(-1):=\chi_v(-1)$$ ne d\'epende pas de la place archim\'edienne
$v$.  L'alg\'ebricit\'e de $\pi$ entra\^ine alors celle de $\chi$, et l'on
dispose donc d'un caract\`ere galoisien associ\'e
$\chi_{\iota}:=\rho_{\chi,\iota}$.  La condition $\pi \simeq \pi^\vee
\otimes \chi$ se traduit alors, via la propri\'et\'e (C1), en l'existence
d'une accouplement non d\'eg\'en\'er\'e et ${\rm
Gal}(\overline{F}/F)$-\'equivariant $$ \rho_{\pi,\iota} \otimes
\rho_{\pi,\iota} \longrightarrow \chi_{\iota}$$ Si $\rho_{\pi,\iota}$ est
irr\'eductible, ce qui est conjectur\'e, un tel accouplement est unique \`a
un scalair pr\`es.  La question qui nous int\'eresse ici est de d\'eterminer
\`a quelle condition cet accouplement est sym\'etrique ou altern\'e.  \ps

Le th\'eor\`eme suivant est d\'emontr\'e dans mon article~\cite{bchsigne},
qui est un travail en commun avec Jo\"el Bella\"iche.  Soit $c_v \subset
{\rm Gal}(\overline{F}/F)$ la classe de conjugaison form\'ee des
conjugaisons complexes de $F$ en la place archim\'edienne $v$. 
L'ind\'ependance de $\chi_v(-1)$ en une telle $v$ entra\^ine que le signe
$\chi_\iota(c):=\chi_\iota(c_v)$ ne d\'epend pas non plus de $v$.  On
prendra garde que $\chi_\iota(c)$ ne co\"incide pas n\'ecessairement avec
$\chi_\infty(-1)$ (voir plus bas).  \ps

\begin{thm}\label{signetotreel} Si $n$ est pair, et si $\chi_\iota(c)=-1$, alors il existe un accouplement
non-d\'eg\'en\'er\'e, ${\rm Gal}(\overline{F}/F)$-\'equivariant, 
$\rho_{\pi,\iota} \otimes \rho_{\pi,\iota} \longrightarrow \chi_\iota$,
qui est altern\'e. Dans tous les autres cas, il existe un tel accouplement qui est sym\'etrique.
\end{thm}

Autrement dit, la repr\'esentation $\rho_{\pi,\iota}$ peut \^etre conjugu\'e
\`a une repr\'esentation \`a valeurs dans ${\rm GSp}(n,\overline{\Q}_\ell)$
dans le premier cas, dans ${\rm GO}(n,\overline{\Q}_\ell)$ dans le second,
pour le caract\`ere de similitude $\chi_\iota$.  Rappelons qu'il est important de pr\'eciser le caract\`ere de similitudes dans ce genre
de situation, car une repr\'esentation donn\'ee peut \^etre symplectique
pour un certain caract\`ere de similitudes, orthogonale pour un autre. 
L'exemple le plus \'evident est celui de la repr\'esentation tautologique du
groupe ${\rm O}(2,\overline{\Q}_\ell)$, qui est symplectique pour le
caract\`ere d\'eterminant et orthogonale pour le caract\`ere trivial.  \ps

Je voudrais pr\'eciser que ce r\'esultat n'est pas vide en dimension $n$
impaire, car l'irr\'eductibilit\'e de $\rho_{\pi,\iota}$ n'est pas connue. 
En retour, c'est m\^eme un ingr\'edient important pour d\'emontrer cette
irr\'eductibilit\'e quand $n\leq 5$ dans les travaux de Calegari et
Gee~\cite{calegarigee}.  \ps

Regardons de plus pr\`es la condition $\chi_\iota(c)=-1$. Le corps $F$
\'etant totalement r\'eel, rappelons que $\chi$ s'\'ecrit de mani\`ere
unique sous la forme $\chi=|.|^w\psi$ o\`u $w \in \Z$ et $\psi$ est un
caract\`ere de Hecke d'image finie de $F$ (raffinement classique d\^u \`a
Chevalley du th\'eor\`eme des unit\'es de Dirichlet).  Par la construction
de Weil, on a l'identit\'e $$\chi_\iota=\omega_\ell^w \cdot \iota(\psi)
\circ {\rm rec} $$ o\`u $\omega_\ell$ d\'esigne le caract\`ere cyclotomique
$\ell$-adique de $F$, $\iota(\psi)$ est simplement le compos\'e de $\iota$
par $\psi$ (ce qui a un sens car ${\rm Im}(\psi) \subset \overline{\Q}$), et
${\rm rec}$ est l'isomorphisme de r\'eciprocit\'e de la th\'eorie du corps
de classes globale.  De $\psi(-1)=\chi_\infty(-1)$ on d\'eduit l'identit\'e
\begin{equation}\label{relation-1}\chi_{\iota}(c)=(-1)^w
\chi_\infty(-1).\end{equation} \ps

Je propose maintenant d'expliquer en quoi le Th\'eor\`eme~\ref{signetotreel}
est attendu du point de vue des conjectures de Langlands.  En effet,
Langlands associe \`a $\pi$ une repr\'esentation irr\'eductible $\rho :
\mathcal{L}_F \rightarrow \GL_n(\C)$ de son groupe conjectural
$\mathcal{L}_F$, telle que pour toute place $v$ de $F$ la restriction de
$\rho$ \`a l'homomorphisme structural ${\rm W}'_{F_v} \rightarrow
\mathcal{L}_F$ est isomorphe \`a ${\rm L}(\pi_v)$.  Le caract\`ere $\chi$
peut \^etre vu comme un caract\`ere de $\mathcal{L}_F^{\rm
ab}=F^\ast\backslash \AAA^\ast_F$.  L'hypoth\`ese $\pi \simeq \pi^\vee
\otimes \chi$ se traduit en une hypoth\`ese de type $\rho \simeq \rho^\vee
\otimes \chi$ par la propri\'et\'e de Cebotarev de $\mathcal{L}_F$.  Il
existe donc un unique accouplement "global" $\mathcal{L}_F$-\'equivariant
non d\'eg\'en\'er\'e $\rho \otimes \rho \rightarrow \chi$.  Il est donc
n\'ecessairement sym\'etrique si $n$ est impair, on suppose donc $n$ pair. 
\ps

Cet accouplement global induit en une place archim\'edienne donn\'ee $v$ un
accouplement non-d\'eg\'en\'er\'e ${\rm L}(\pi_v) \otimes {\rm L}(\pi_v)
\rightarrow \chi_v$.  Comme $\pi$ est alg\'ebrique r\'eguli\`ere, ${\rm
L}(\pi_v)$ est sans multiplicit\'e et somme directe de repr\'esentations de
la forme $$I_{p,q}={\rm Ind}_{{\rm W}_\C}^{{\rm W}_\R} z^p\overline{z}^q.$$
Le lemme de puret\'e de Clozel~\cite{clozel} entra\^ine d'ailleurs que
l'entier $p+q$ est le m\^eme pour tous les constituants $I_{p,q}$.  Il est
m\^eme ind\'ependant de $v$ : c'est l'entier $w$ d\'efini plus haut par la
relation $\pi^\vee \otimes \chi \simeq \pi$.  Par multiplicit\'e $1$, on en
d\'eduit que l'accouplement global induit sur chaque facteur $I_{p,q}$ de
${\rm L}(\pi_v)$ un accouplement non d\'eg\'en\'er\'e de facteur de
similitude $\chi_v$.  Observons que le type sym\'etrique/altern\'e de cet
accouplement est fix\'e ind\'ependamment de $v$ par l'accouplement global. 
Mais un tel accouplement sur $I_{p,q}$ est n\'ecessairement sym\'etrique si
$\chi_v(-1)=1$, altern\'e sinon.  En effet, si $\varepsilon_{\C/\R}$ et
$x=|.|\varepsilon_{\C/\R}$ d\'esignent respectivement les caract\`eres
d'ordre $2$ et identit\'e de ${\rm W}_\R^{\rm ab}=\R^\ast$, alors
$$\det(I_{p,q}) \otimes \chi_v^{-1} = \varepsilon_{\C/\R}x^w \chi_v^{-1},$$
par la formule du d\'eterminant d'une induite, qui vaut donc $1$ (cas
"altern\'e") si et seulement si $\chi_v(-1)=(-1)^{w+1}$.  \ps

On en d\'eduit tout d'abord que tous les $\chi_v(-1)$ doivent
n\'ecessairement avoir m\^eme signe, autrement dit : l'hypoth\`ese sur le
signe des $\chi_v(-1)$ apparaissant dans la d\'efinition de $\pi$ d\'ecoule
conjecturalement des autres hypoth\`eses.  Cette observation r\'epond
d'ailleurs \`a mon sens aux questions sur ce sujet soulev\'ees dans~\cite[\S
4.3]{cht}.  On d\'eduit enfin de la formule~\eqref{relation-1} que
l'alternative symplectique-orthogonale de la repr\'esentation $\rho$ doit
satisfaire exactement la m\^eme dichotomie que dans le
Th\'eor\`eme~\ref{signetotreel}.  Je laisse au lecteur le plaisir de se
convaincre que $\rho$ et $\rho_{\pi,\iota}$ doivent satisfaire la m\^eme
alternative.  Il pourra par exemple consid\'erer le cas particulier o\`u
$\pi$ est discr\`ete \`a une place finie.  On notera tout de m\^eme que
l'argument ci-dessus n'admet bien s\^ur pas d'analogue direct si l'on
remplace $\mathcal{L}_F$ par ${\rm Gal}(\overline{F}/F)$.  \ps

Nous d\'eduisons en fait le Th\'eor\`eme~\ref{signetotreel} d'un
th\'eor\`eme analogue dans le cas $CM$.  Supposons donc maintenant que $F$
est un corps $CM$ et que $\pi$ est une repr\'esentation automorphe cuspidale
alg\'ebrique r\'eguli\`ere polaris\'ee de $\GL_n$ sur $F$.  En particulier,
on dispose d'un accouplement non-d\'eg\'en\'er\'e et
$\Gal(\overline{F}/F)$-\'equivariant $$\rho_{\pi,\iota} \otimes
\rho_{\pi,\iota}^c \longrightarrow \omega_\ell^{1-n}.$$ Concr\`etement,
c'est une forme bilin\'eaire non-d\'eg\'en\'er\'ee
$\langle\cdot,\cdot\rangle$ sur l'espace $V$ de $\rho_{\pi,\iota}$ telle que
pour tout $x,y \in V$ et tout $g \in \Gal(\overline{F}/F)$ on ait la
relation $$\langle
\rho_{\pi,\iota}(g)x,\rho_{\pi,\iota}(\widetilde{c}g\widetilde{c}^{-1})y
\rangle = \omega_\ell(g)^{1-n} \langle x, y\rangle.$$  Dans cette formule,
$\widetilde{c} \in {\rm Gal}(\overline{F}/F^+)$ d\'esigne une conjugaison
complexe fix\'ee, qui est donc un \'el\'ement de carr\'e $1$ et d'image $c
\in {\rm Gal}(F/F^+)$.  L\`a encore, si $\rho_{\pi,\iota}$ est
irr\'eductible un tel accouplement est unique \`a un scalair pr\`es, et donc
soit sym\'etrique, soit altern\'e.  Le r\'esultat principal
de~\cite{bchsigne} est le r\'esultat suivant.

\begin{thm}\label{signecm} Soit $\pi$ une repr\'esentation automorphe
cuspidale alg\'ebrique r\'eguli\`ere et polaris\'ee de $\GL_n$ sur le corps de nombres $F$ suppos\'e CM.  Il existe un accouplement sym\'etrique, non-d\'eg\'en\'er\'e, et
$\Gal(\overline{F}/F)$-\'equivariant, $\rho_{\pi,\iota} \otimes
\rho_{\pi,\iota}^c \longrightarrow \omega_\ell^{1-n}$.  \end{thm}

Des cas tr\`es particuliers de ce th\'eor\`eme avaient \'et\'e d\'emontr\'e
par Clozel-Harris-Taylor dans~\cite{cht}, comme cons\'equences de leurs
th\'eor\`emes de type $R=T$; ces auteurs ont notamment vu l'importance de
cette propri\'et\'e de symm\'etrie dans les questions de d\'eformations. Nous la retrouverons au~\S\ref{fougeredimsup}. \ps

Observons que $\rho_{\pi,\iota}$ est sans multiplicit\'e comme on le voit
sur ses poids de Hodge-Tate (r\'egularit\'e).  Il suit que sur chaque
facteur irr\'eductible $r$ de $\rho_{\pi,\iota}$ tel que $r^\vee \otimes
\omega_\ell^{1-n} \simeq r^c$, l'unique accouplement du type ci-dessus qui
est non d\'eg\'en\'er\'e est sym\'etrique.  (Dans la terminologie de notre
article nous disons aussi que $\rho_{\pi,\iota}$ est de signe $+1$.  )\ps

Le fait que le Th\'eor\`eme~\ref{signecm} entra\^ine le
Th\'eor\`eme~\ref{signetotreel}, bien que peut-\^etre contre-intuitif, est
sans difficult\'e.  Le seul point non-trivial est d'observer,
suivant~\cite[\S 4]{cht}, que si $\pi$ est comme dans le
Th\'eor\`eme~\ref{signetotreel}, alors le changement de base de $\pi$ \`a
une extension quadratique CM convenable de $F$ admet une torsion par un
caract\`ere de Hecke qui est polaris\'ee au sens CM : c'est de toutes
fa\c{c}ons ainsi que $\rho_{\pi,\iota}$ est construite par recollement dans
le cas totalement r\'eel.  \ps

On pourrait v\'erifier comme plus haut que le formalisme de Langlands
sugg\`ere \'egalement le th\'eor\`eme ci-dessus.  Une telle \'etude
montrerait qu'\`a cause de l'hypoth\`ese sur $\pi_\infty$, le param\`etre de
Langlands global d'un $\pi$ comme dans l'\'enonc\'e du th\'eor\`eme est
toujours le changement de base d'un param\`etre global, discr\^et aux places
archim\'ediennes, du groupe unitaire quasi-d\'eploy\'e \`a $n$ variables
associ\'e \`a $F/F^+$~(voir par exemple~\cite[A.11.7]{bchlivre}).  C'est d'ailleurs un fait
important dans notre d\'emonstration que la fonctorialit\'e associ\'ee \`a
ce cas est en fait connue, du moins sous des hypoth\`eses suppl\'ementaires
auxquelles on peut toujours se ramener, par les travaux de
Labesse~\cite{labesse}.  Ceci \'etant dit, le th\'eor\`eme ci-dessus peut
alors \^etre vu comme un cas particulier de la conjecture de "parit\'e" de
Gross~\cite{grossodd}.  \ps

Terminons, enfin!,  par une br\`eve discussion de notre d\'emonstration du
th\'eor\`eme.  L'id\'ee est de se ramener par des congruences, plus
pr\'ecis\'ement des familles $\ell$-adiques, \`a d\'emontrer le r\'esultat
dans le cas particulier o\`u $n$ est impair et $\rho_{\pi,\iota}$ est
irr\'eductible...  auquel cas il est \'evident !  On v\'erifie en effet
d'abord que l'alternative sym\'etrique/altern\'ee est constante dans une
famille $\ell$-adique de repr\'esentations $\rho : {\rm Gal}(\overline{F}/F)
\rightarrow \GL_n(\overline{\Q}_\ell)$ satisfaisant toutes $\rho^\vee
\otimes \omega_\ell^{1-n} \simeq \rho$, d\`es lors que la base de cette
famille est un affino\"ide connexe.  \ps

Pour des raisons techniques, on se ram\`ene par changement de base
r\'esoluble bien choisi au cas o\`u $F/F^+$ est non ramifi\'ee \`a toutes
les places finie, $[F^+:\Q]$ est pair, $\pi$ est non ramifi\'ee aux places
non d\'ecompos\'ees au dessus de $F$, et o\`u $\pi$ admet des invariants par
un sous-groupe d'Iwahori \`a toutes les places divisant $\ell$.  On
d\'esigne par $U(n)$ le groupe unitaire \`a $n$ variables sur $F$ associ\'e
\`a $F/F^+$ qui est quasi-d\'eploy\'e \`a toutes les places finies et
compact aux places archim\'ediennes.  \ps

Si $n$ est impair, on commence par descendre $\pi$ en une repr\'esentation
automorphe $\pi_0$ de $U(n)$.  C'est possible par les travaux de
Labesse~\cite{labesse}.  Ensuite, nous utilisons la foug\`ere infinie pour
ce groupe unitaire, et g\'en\'eralisons une observation d\'ej\`a vue dans ma
th\`ese en dimension $3$~\cite[\S 9.1]{bchens} : en se d\'epla\c{c}ant dans
cette foug\`ere, on peut connecter $\rho_{\pi,\iota}=\rho_{\pi_0,\iota}$ \`a
un $\rho_{\pi',\iota}$ qui est de plus irr\'eductible restreinte \`a toutes
les places divisant $\ell$.  Il suffit pour cela de se d\'eplacer
intelligemment de sorte \`a prescrire suffisament les pentes des Frobenius
cristallins des ${\rm D}_{\rm cris}((\rho_{\pi',\iota})_v)$ pour $v$
divisant $\ell$ : cet argument sera expliqu\'e au~\S\ref{appraffinee}. Le r\'esultat est alors \'evident pour $\rho_{\pi',\iota}$
et on conclut.\ps

Dans le cas $n$ pair, on utilise un transfert endoscopique temp\'er\'e du
type ${}^L{\rm U}(n) \times {}^L{\rm U}(1) \rightarrow {}^L{\rm U}(n+1)$
pour se ramener \`a la dimension $n+1$ impaire et proc\'eder comme
pr\'ec\'edemment : on d\'eforme cette fois-ci irr\'eductiblement
$\rho_{\pi,\iota} \oplus \nu_\iota$ pour un caract\`ere de Hecke
alg\'ebrique $\nu$ de $F$ bien choisi.  Il y a une obstruction \`a
l'existence de tels transferts en g\'en\'eral provenant de la formule de
multiplicit\'e \`a la Labesse-Langlands, d\'emontr\'ee dans ce cadre par
Clozel, Harris et Labesse dans~\cite{chlb}.  Il faut qu'aux places
archim\'ediennes $v$ de $F$ les poids de $\nu_v$ et de $\pi_v$ soient
entrelac\'es dans certains ordres.  Cela n'a pas d'importance pour notre
d\'emonstration car nous pouvons choisir $\nu$ et modifier au pr\'ealable
\`a loisir les poids de $\pi$ quitte \`a se d\'eplacer d'abord un peu dans
la foug\`ere infinie de $U(n)$.  \ps

Je voudrais mentionner qu'une strat\'egie plus g\'eom\'etrique pour
d\'emontrer le Th\'eor\`eme~\ref{signecm} consisterait \`a utiliser
l'accouplement de Poincar\'e sur l'espace de cohomologie dans lequel
$\rho_{\pi,\iota}$ est construite.  On se heurte cependant notamment \`a la
difficult\'e que la conjugaison complexe n'agit pas naturellement sur cet
espace, car la vari\'et\'e sous-jacente n'est que d\'efinie sur $F$.  Une
mani\`ere de contourner ce probl\`eme a \'et\'e propos\'ee par Taylor
dans~\cite{taylorcomplex}.  Il serait int\'eressant de savoir si cette
m\'ethode permet de re-d\'emontrer le th\'eor\`eme ci-dessus.  Pour
terminer, je voudrais discuter comme promis la condition~(C4).

\begin{thm}~\cite{taylorcomplex},~\cite{taibi} Supposons $F$ totalement
r\'eel et soient $\pi$ et $\iota$ comme dans le
Th\'eor\`eme~\ref{existencegalois}.  Si $n$ est pair et si
$\chi_\iota(c)=1$, on suppose de plus que $\chi_\infty(-1)=1$.  Alors
$\rho_{\pi,\iota}$ satisfait (C4).  \end{thm}

Ce th\'eor\`eme est d\^u \`a Taylor~\cite{taylorcomplex} dans le cas
particulier o\`u $n$ est impair et o\`u $\rho_{\pi,\iota}$ est
irr\'eductible.  Le r\'esultat de Taylor a ensuite \'et\'e \'etendu au cadre
g\'en\'eral ci-dessus par mon \'etudiant Olivier  Ta\"ibi dans~\cite{taibi}.  La
strat\'egie de Ta\"ibi est de se ramener par des congruences au r\'esultat
de Taylor, un peu \`a la mani\`ere de la d\'emonstration ci-dessus.  De
nouvelles difficult\'es apparaissent cependant, qui sont inh\'erentes aux
propri\'et\'es de la foug\`ere infinie pour les groupes classiques. 
Notamment, ainsi que l'a observ\'e Ta\"ibi, il n'est plus tout-\`a-fait
possible de d\'eformer irr\'eductiblement tout $\rho_{\pi,\iota}$ dans la
foug\`ere infinie de ${\rm Sp}(n-1)$ quand $n$ est impair.  Ses
d\'emonstrations utilisent au final de mani\`ere astucieuse des groupes
classiques annexes et plusieurs cas de fonctorialit\'e d\'emontr\'es dans
les travaux r\'ecents d'Arthur~\cite{arthur}.

\newpage
\section{Interlude : vari\'et\'es de caract\`eres et analogues $p$-adiques} 

\subsection{R\'esum\'e et perspectives}	Avant d'aborder mes travaux sur la
foug\`ere infinie, il me semble opportun de faire quelques rappels sur les
familles de repr\'esentations d'un groupe g\'en\'eral et sur la notion de
vari\'et\'e de caract\`eres.  J'exposerai notamment des r\'esultats de mon
article~\cite{chdet} dans lequel je d\'efinis une notion de vari\'et\'e de
caract\`eres param\'etrant les repr\'esentations continues $p$-adiques d'un
groupe profini donn\'e.  Cela fournira un cadre agr\'eable, \`a d\'efaut
d'\^etre indispensable, pour l'\'etude des familles de repr\'esentations
galoisiennes et de la foug\`ere infinie, dans les parties suivantes.  Il
ouvre aussi des questions int\'eressantes qui me semblent inexplor\'ees. 
Afin de justifier que l'objet que je d\'efinis est bien un analogue de la
vari\'et\'e des caract\`eres usuelle pour les groupes discrets, je
commencerai par des rappels sur ces derni\`eres.  \ps

	Je profiterai \'egalement de l'occasion pour exposer des
r\'esultats, obtenus en collaboration avec Jo\"el
Bella\"iche dans notre livre~\cite[Ch.  1]{bchlivre}, concernant l'\'etude
g\'en\'erale des vari\'et\'es de caract\`eres au voisinage des points
param\'etrant des repr\'esentations r\'eductibles sans multiplicit\'es. 
C'est un sujet certainement classique, au moins depuis~\cite{lubmag},
surlequel la litt\'erature semble assez vaste.  Je ne suis malheureusement
pas comp\'etent pour faire l'historique et la synth\`ese de ces travaux
(voir par exemple~\cite{sikora} pour de nombreuses r\'ef\'erences sur des
points de vue diff\'erents du notre). Notre motivation initiale, issue de la th\'eorie des nombres,
\'etait de g\'en\'eraliser en dimension sup\'erieure la m\'ethode initi\'ee
par Ribet~\cite{ribet} pour construire des extensions entre
repr\'esentations galoisiennes, et aussi de comprendre les limites de cette
m\'ethode; je renvoie au rapport r\'ecent de Mazur~\cite{mazurribet} \`a ce
sujet.  Nos r\'esultats jouent par exemple un r\^ole technique important
dans le th\'eor\`eme principal de notre livre (voir~\S\ref{thmlivre} Ch. 2); ils devraient
s'appliquer aussi \`a des questions analogues concernant les groupes discrets, par exemple dans l'analogie entre noeuds et nombres premiers, que l'on
rediscutera au~\S\ref{dehnar} Ch. 2.  \ps

\subsection{Pr\'eliminaire : un point de vue modulaire sur les vari\'et\'es de caract\`eres}\label{prelimvcar} Soit $G$ un groupe.  Si $A$ est un
anneau commutatif unitaire, une famille de repr\'esentations de $G$ de
dimension $d$ et {\it \`a valeurs dans $A$} (ou encore {\it param\'etr\'ee
par ${\rm Spec}(A)$}), est traditionnellement un homomorphisme de groupes
${\rm G} \rightarrow \GL_d(A)$. Notons $${\rm Hom}(G,\GL_d)$$ le
foncteur $A \mapsto
\Hom(G,\GL_d(A))$, des anneaux commutatifs unitaires vers les ensembles, param\'etrant
les familles de repr\'esentations de dimension $d$ de $G$.  Ce foncteur est \'evidemment repr\'esentable ("matrices
g\'en\'eriques"), par un anneau $R_d(G)$ qui est m\^eme de type fini sur $\Z$ si le groupe $G$ est
de type fini, ce que l'on supposera d\'esormais pour simplifier.  \ps

	Pour de nombreuses questions, il est souvent plus naturel, et parfois
m\^eme obligatoire, de consid\'erer des repr\'esentations \`a isomorphisme pr\`es,
i.e.  modulo l'action naturelle de ${\rm GL}_d$ par conjugaison sur ${\rm
Hom}(G,\GL_d)$.  Cependant, il est bien connu que le foncteur $$A
\rightarrow {\rm Hom}({\rm G},\GL_d(A))/{\rm GL}_d(A)$$ n'est pas en
g\'en\'eral repr\'esentable, en consid\'erant par exemple des repr\'esentations
non d\'efinies sur leur corps des traces (mais il y a pire!). On se restreint d\'esormais aux anneaux $A$ qui sont des $k$-alg\`ebres (commutatives unitaires) o\`u
$k$ est un corps de caract\'eristique $0$. L'approche traditionnelle est
alors de d\'efinir la vari\'et\'e des $k$-caract\`eres de dimension $d$ de
$G$ \`a la Mumford, comme \'etant le quotient GIT  de
$\Hom(G,\GL_d)$ par $\GL_d$~(\cite{lubmag}). On la notera ${\rm Car}_d(G)$.
C'est par d\'efinition le
$k$-sch\'ema affine de $k$-alg\`ebre $$k[{\rm Car}_d(G)]:=R_d(G)^{\GL_d(k)} \subset R_d(G).$$  \ps

Il est alors non trivial mais connu depuis M. Artin~\cite{artin}, voir aussi~\cite{Proc0}, que les points de ${\rm Car}_d(G)$ dans une extension
alg\'ebriquement close $L$ de $k$ param\`etrent exactement les classes d'isomorphie de repr\'esentations semisimples $G \rightarrow \GL_d(L)$. La terminologie "vari\'et\'e des caract\`eres" vient alors du fait classique que ces classes d'isomorphies
sont uniquement d\'etermin\'ees par leurs caract\`eres (i.e. leurs
fonctions traces), car $k$ est de caract\'eristique nulle. \ps

Cependant, cette description ne s'\'etend pas en g\'en\'eral aux points \`a valeurs dans une $k$-alg\`ebre quelconque, ce qui complique sensiblement l'\'etude de ${\rm Car}_d(G)$. D'un point de vue na\"if, c'est simplement le passage au quotient (ou aux invariants sur les fonctions) qui fait perdre cette interpr\'etation modulaire. C'est particuli\`erement embarrassant pour
les questions infinit\'esimales, comme les calculs d'espaces tangents (de
Zariski).~\footnote{Traditionnellement, on s'en sort toutefois sur l'ouvert Zariski de ${\rm
Car}_d(G)$ param\'etrant les repr\'esentations absolument irr\'eductibles, sur
lequel on dispose d'une interpr\'etation modulaire en terme de
repr\'esentations dans des alg\`ebres d'Azumaya, voir par exemple~\cite{Proc3}.} \ps
	
	Il se trouve que l'on peut contourner ces probl\`emes en utilisant
les travaux de Procesi sur les invariants de similitudes simultan\'es d'une
famille de matrices dans $M_d$~\cite{Proc}.  Procesi d\'emontre en
effet que ces invariants sont d'une part engendr\'es par les traces en les
mon\^omes en ces matrices, et d'autre part que les relations entre ces
traces se d\'eduisent toutes de l'identit\'e de Cayley-Hamilton, faits bien
connus dans le cas d'une seule matrice.  Cette derni\`ere identit\'e joue
alors un r\^ole combinatoire central dans ces questions, car elle conduit
comme nous allons le rappeler \`a une th\'eorie des invariants explicite.  \ps

	Soit $\rho : G \rightarrow \GL_d(A)$ une repr\'esentation et soit
$T={\rm Trace} \circ \rho : G \rightarrow A$.  On sait depuis Frobenius~\cite[p. 50, formule (21)]{frobenius} que
$T$ satisfait l'identit\'e suivante\footnote{Si $M=A^d$, cela vient d'une
part de ce que $\Lambda_A^{d+1}M=0$, et d'autre part de ce que si
$g_1,\cdots,g_n \in \GL_A(M)$ et $\sigma \in \got{S}_n$, la trace de
$(g_1,\cdots,g_n)\sigma$ sur $M^{\otimes_A n}$ co\"incide avec
$T^{\sigma}(g_{1},\dots,g_{n})$ (Kostant).} :
\begin{equation}\label{npseudo} \forall g_1,\, g_2,\,\dots ,g_d, g_{{d+1}}
\in G, \,\,\, \, \sum_{\sigma \in \got S_{d+1}}
\varepsilon(\sigma)T^\sigma(g_1,\,g_2,\,\dots,\,g_{d+1})=0.  \end{equation}
Explications : si $n\geq 1$ est un entier, et si $\sigma \in \got S_n$, on a
pos\'e $T^{\sigma}(g_1,g_2,\dots ,g_n)=T(g_{i_1}g_{i_2}\dots g_{i_r})$ si
$\sigma$ est le cycle $(i_1 i_2 \dots i_r)$, et en g\'en\'eral
$T^\sigma=\prod T^{c_i}$ si $\sigma=c_1\dots c_s$ est la d\'ecomposition en
cycles de $\sigma$.  Ainsi que Procesi l'a remarqu\'e~\cite{Proc}, cette
identit\'e peut \'egalement \^etre vue comme la polarisation totale de
l'identit\'e de Cayley-Hamilton, vue sous la forme ${\rm
Trace}(\chi_M(M)N)=0$ o\`u $M,N \in M_d(A)$ et $\chi_M$ est le polyn\^ome
caract\'eristique de $M$. Cela a conduit Taylor~\cite{tay} \`a consid\'erer
des fonctions g\'en\'erales $T : G \rightarrow A$ telles que $T(1)=d$,
$T(gh)=T(hg)$ pour tout $g,h \in G$, et satisfaisant
l'identit\'e~(\ref{npseudo}).  Ces fonctions ont \'et\'e depuis largement
utilis\'ees en th\'eorie des nombres sous le nom de {\it
pseudo-repr\'esentations},\footnote{Le premier \`a avoir introduit cette
terminologie et vu son utilit\'e, en dimension $2$ et sous une forme
l\'eg\`erement diff\'erente, est Wiles dans~\cite{wiles}.} ou {\it
pseudo-caract\`eres}, de dimension $d$, et nous adopterons cette derni\`ere
terminologie, h\'erit\'ee de Rouquier~\cite{rouquier}.  \ps

\newcommand{\PS}{\mathrm{PCar}} D\'esignons par $\PS_d(G)(A)$ l'ensemble des
pseudo-caract\`eres de dimension $d$ de $G$ \`a valeurs dans $A$.  Il est
\'evident que $A \mapsto \PS_d(G)(A)$ d\'efinit un foncteur $\PS_d(G)$ des
$k$-alg\`ebres commutatives vers les ensembles, et que ce foncteur est
repr\'esentable, par une $k$-alg\`ebre que nous noterons $k[\PS_d(G)]$. 
L'observation de Frobenius montre que la trace fournit un morphisme de
foncteurs $\Hom(G,\GL_d)/\GL_d \rightarrow \PS_d(G)$.  Il se trouve que ce
morphisme induit un isomorphisme ${\rm Car}_d(G) \isomo \PS_d(G)$.  C'est en
effet une cons\'equence imm\'ediate des r\'esultats de
Procesi~\cite{Proc},~\cite{Proc2}, je redonne l'argument ci-dessous par manque de r\'ef\'erence
ad\'equate.  Je rappelle que $k$ est un corps de caract\'eristique nulle. 
\ps

\begin{prop} La trace de la repr\'esentation universelle $G \rightarrow
\GL_d(R_d(G))$ induit un isomorphisme de $k$-alg\`ebres $$k[\PS_d(G)] \isomo
k[{\rm Car}_d(G)].$$ Autrement dit, $\PS_d(G)$ est canoniquement le foncteur
des points de ${\rm Car}_d(G)$.  \end{prop}

\begin{pf} Soit $\rho^{\rm univ} : G \rightarrow {\rm GL}_d(R_d(G))$ la
 repr\'esentation universelle et $T^{\rm univ} : G \rightarrow k[\PS_d(G)]$ le pseudo-caract\`ere universel.  Il est \'evident que ${\rm
 Trace} \circ \rho^{\rm univ}(G) \subset k[{\rm Car}_d(G)]$.  Le r\'esultat
 de Frobenius suscit\'e fournit donc un morphisme canonique $f :
 k[\PS_d(G)] \rightarrow k[{\rm Car}_d(G)]$ envoyant $T^{\rm univ}(g)$, pour $g
 \in G$, sur ${\rm Tr} \circ \rho^{\rm univ}(g)$.  Ce morphisme est
 surjectif si $G$ est un groupe libre de type fini par~\cite{Proc}, il le
 reste donc pour $G$ de type fini quelconque par r\'eductivit\'e lin\'eaire
 de ${\rm GL}_d$ sur $k$.  Autrement dit la $k$-alg\`ebre $k[{\rm
 Car}_d(G)]$ est engendr\'ee par les ${\rm Trace}(\rho^{\rm univ}(g))$, $g
 \in G$ (l'assertion analogue pour $k[\PS_d(G)]$ et les $T^{\rm
 univ}(g)$ est \'evidente).  \ps R\'eciproquement, d'apr\`es~\cite{Proc2}
 appliqu\'e \`a l'une quelconque des alg\`ebres de Cayley-Hamilton
 associ\'ees \`a $T^{\rm univ} $ (voir la note~\ref{notepcar} ci-apr\`es ou~\cite[\S 1.2]{bchlivre}), il existe
 un morphisme injectif de $k$-alg\`ebres $f' : k[\PS_d(G)] \rightarrow B$ et
 une repr\'esentation $\rho : G \rightarrow \GL_d(B)$ tels que ${\rm Tr}
 \circ \rho = f' \circ T^{\rm univ}$.  On en d\'eduit un $k$-morphisme $f''
 : R_d(G) \rightarrow B$ tel que $\rho = f'' \circ \rho^{\rm univ}$.  La
 derni\`ere assertion du paragraphe pr\'ec\'edent assure que $f''(k[{\rm
 Car}_d(G)]) \subset f'(k[\PS_d(G)])$.  L'injectivit\'e de $f'$ montre que
 $(f')^{-1} \circ f'' : k[{\rm Car}_d(G)] \rightarrow k[\PS_d(G)]$ est bien
 d\'efini.  Il envoit par construction ${\rm Tr} \circ \rho^{\rm univ}(g)$
 sur ${\rm T}^{\rm univ}(g)$ : c'est donc un inverse de $f$.  \end{pf}

Je voudrais mentionner que si l'on ne s'int\'eresse qu`\`a des
repr\'esentations de dimension $2$ et de d\'eterminant $1$, la notion
adapt\'ee de pseudo-caract\`ere est simplement une fonction $T : G
\rightarrow A$ telle que $T(1)=2$ et telle que pour tout $g,h \in G$,
$T(gh)=T(hg)$ et $$T(gh)+T(gh^{-1})=T(g)T(h).$$ Cette identit\'e de traces
est omnipr\'esente dans la litt\'erature sur les vari\'et\'es de
caract\`eres \`a valeurs dans $\SL_2$, comme par exemple
dans~\cite{cullershalen} (voir aussi l'exposition par Goldman~\cite{goldman}
sur la d\'etermination par Fricke et Vogt de ${\rm Car}_\Q({\rm F}_2)$ (fin
${\rm XIX}$\`eme), ${\rm F}_g$ d\'esignant le groupe libre \`a $g$
g\'en\'erateurs).  Toutefois, sauf erreur de ma part, ces auteurs ne me
semblent pas aller jusqu'\`a affirmer que dans ce cas la vari\'et\'e des
caract\`eres est simplement l'espace de modules fin de ces fonctions (ce qui
est vrai).  Cette propri\'et\'e me semble n\'eanmoins int\'eressante.  Par
exemple, elle permet de voir imm\'ediatement que l'espace tangent de Zariski
en la repr\'esentation triviale est l'espace des fonctions $t : G
\rightarrow k$ telles que $\forall g,h \in G$, $t(gh)=t(hg)$ et
$t(gh)+t(gh^{-1})=2(t(h)+t(g))$ ("identit\'e du parall\'elogramme"), ce qui
a des cons\'equences amusantes pour lesquelles je renvoie \`a ma
note~\cite{chcarpad}.  \ps

Pour revenir au cas g\'en\'eral, il suit par ailleurs de ce que l'on a dit
que si $k$ est un corps alg\'ebriquement clos de caract\'eristique nulle,
les pseudo-caract\`eres $G \rightarrow k$ de dimension $d$ sont exactement
les traces des (classes d'isomorphie de) repr\'esentations semisimples $G
\rightarrow k$.  Cela avait originalement \'et\'e observ\'e par Taylor
dans~\cite{tay} par une approche similaire \`a celle ci-dessus.  Rouquier a
par la suite obtenu dans~\cite{rouquier} une d\'emonstration assez directe
de ce r\'esultat (dans l'esprit de~\cite{Proc0}), qu'il a m\^eme \'etendue au
cas des corps $k$ alg\'ebriquement clos tels que $d!$ est inversible dans
$k$.  La morale de cette affaire est que nous aurions directement pu
d\'efinir la vari\'et\'e des $k$-caract\`eres de dimension $d$ de $G$ comme
\'etant le $k$-sch\'ema affine $\PS_d(G)$.  C'est ce point de vue que nous
allons g\'en\'eraliser dans le cas profini.\ps
\newcommand{\Car}{{\rm Car}}
\newcommand{\CH}{\mathcal{R}}
Pour conclure ce paragraphe de g\'en\'eralit\'es sur les vari\'et\'es des
caract\`eres, je voudrais mentionner quelques r\'esultats g\'en\'eraux que
j'ai obtenus avec Bella\"iche dans~\cite[Ch. 1]{bchlivre} (voir
aussi~\cite{bckl},~\cite{bchjimj}). Notre objectif \'etait d'\'etudier ${\rm Car}_d(G)$ au
voisinage d'un point $x \in {\rm Car}_d(G)$ param\'etrant une repr\'esentation semisimple
$\rho_x : G \rightarrow \GL_d(\overline{k})$ \'eventuellement r\'eductible mais sans multiplicit\'es, i.e. 
$\rho_x \simeq \oplus_{i=1}^s r_i$ avec les $r_i$ irr\'eductibles et $r_i \not\simeq r_j$ si $i \neq j$. Notre travail
comprend : \begin{itemize}\ps

\item[(i)] Une \'etude sch\'ematique des diff\'erents lieux de r\'eductibilit\'e du pseudo-caract\`ere universel 
au voisinage du point $x$ (\'equations, propri\'et\'es locales des morphismes
naturels $\prod_i  \Car_{d_i}(G) \rightarrow \Car_{d}(G)$ o\`u $d=\sum_i d_i$).\ps

\item[(ii)] Une \'etude des liens entre les ${\rm
Ext}_{k[G]}(r_i,r_j)$ et la g\'eom\'etrie en $x$ de $\Car_d(G)$ et des lieux de r\'eductibilit\'e
ci-dessus.\ps

\item[(iii)] Des crit\`eres pour qu'un morphisme ${\rm Spec}(A) \rightarrow
\Car_d(G)$ d'image contenant $x$ soit la trace d'une repr\'esentation $G
\rightarrow \GL_d(A)$, du moins apr\`es restriction \`a un voisinage \'etale
assez petit de $x$. Par exemple, nous montrons qu'il suffit que
$\Car_d(G)$ soit lisse, ou m\^eme simplement de henselis\'e strict factoriel, en $x$. \ps
\end{itemize}

Je renvoie \`a~\cite{bchlivre} pour plus de d\'etails sur ces sujets, qu'il
serait trop long de d\'evelopper ici. Ils reposent de mani\`ere importante sur un
th\'eor\`eme de structure concernant l'alg\`ebre de Cayley-Hamilton
universelle $\CH^{\rm univ}$ sur $\Car_d(G)$.\footnote{\label{notepcar}Soit $T : G \rightarrow A$ un pseudo-caract\`ere de dimension
$d$, que l'on prolonge $A$-lin\'eairement en $T: A[G] \rightarrow A$. Les relations de Newton permettent
de former, pour tout $x \in A[G]$, un polyn\^ome "caract\'eristique" abstrait 
$P_{x,T}=X^d-T(x)X^{d-1}+\frac{T(x)^2-T(x^2)}{2}X^{d-2}+\cdots  \in A[X]$. 
On consid\`ere alors la $A$-alg\`ebre $\CH(T)$ quotient de
l'alg\`ebre du groupe $A[G]$ par l'id\'eal bilat\`ere engendr\'e par les
$P_{g,T}(g) \in A[G]$ pour tout $g \in G$. L'identit\'e~\eqref{npseudo}
assure que $T$ se factorise par $\CH(T)$ et que tout $x \in \CH(T)$ satisfait l'identit\'e de Cayley-Hamilton $P_{x,T}(x)=0$. Le couple
$(\CH(T),T)$ est donc une alg\`ebre de Cayley-Hamilton au sens de
Procesi~\cite{Proc2}. On dispose alors tautologiquement d'une repr\'esentation $G \rightarrow
\CH(T)^\ast$ "de trace" $T$. La formation de $(\CH(T),T)$ commute \`a tout changement de base $A
\rightarrow B$. Si $A=k[\Car_d(G)]$ et $T=T^{\rm
univ}$, on pose $\CH^{\rm univ}=\CH(T)$. Je renvoie \`a~\cite{Proc2},~\cite{dcprr} et~\cite[Ch. 
1]{bchlivre} pour la th\'eorie des alg\`ebres de Cayley-Hamilton.} Nous d\'emontrons en effet que
si $\OO_x$ d\'esigne le hens\'elis\'e strict de l'anneau local de
$\Car_d(G)$ en $x$, alors $$\CH^{\rm univ} \otimes_{k[\Car_d(G)]} \OO_x$$ est
une {\it $\OO_x$-alg\`ebre de matrices g\'en\'eralis\'ee de type $(\dim
r_1,\cdots,\dim r_s)$}. C'est une classe d'alg\`ebres \`a traces que nous avons
introduite et \'etudi\'ee dans~\cite[Ch. 1]{bchlivre}. Dans le cas particulier o\`u $\rho_x$ est
irr\'eductible, cette alg\`ebre est simplement isomorphe \`a
$M_d(\OO_x)$, un r\'esultat d\'ej\`a d\'emontr\'e par
Nyssen~\cite{nyssen}, Rouquier~\cite{rouquier} et Procesi~\cite{Proc3}. Si 
$${\rm Car}_d(G)^{\rm irr} \subset {\rm Car}_d(G)$$ d\'esigne l'ouvert
Zariski constitu\'e des points $x$ tels que $\rho_x$ est irr\'eductible, 
la restriction de $\CH^{\rm univ}$ \`a ${\rm Car}_d(G)^{\rm irr}$ est donc une
alg\`ebre d'Azumaya de rang $d^2$. Mentionnons pour finir que nos r\'esultats permettent
par exemple de d\'evisser de mani\`ere simple
 l'espace tangent de $\Car_d(G)$ en un point $x$ tel que $\rho_x$ est sans
multiplicit\'e~\cite{joeldef} (quand $\rho_x$ est irr\'eductible, il est
bien connu que cet espace est isomorphe \`a ${\rm H}^1(G,{\rm
Ad}(\rho_x))$~\cite{lubmag}).  \ps

\subsection{La vari\'et\'e des caract\`eres $p$-adique}\label{carpadique} On suppose
maintenant que $G$ est un groupe profini et l'on fixe un nombre premier $p$
ainsi qu'un entier $d\geq 1$. On s'int\'eresse dans ce cas aux homomorphismes continus
$G \rightarrow \GL_d(A)$ o\`u $A$ est une $\Q_p$-alg\`ebre de Banach
affino\"ide au sens de Tate~\cite{taterig}, auxquels on pensera comme \`a
des familles analytiques de repr\'esentations continues $G \rightarrow
\GL_d(\overline{\Q}_p)$ param\'etr\'ees par l'espace analytique
affino\"ide ${\rm Max}(A)$.  On consid\`erera plus g\'en\'eralement des
espaces analytiques sur $\Q_p$ au sens de Tate~\cite{BGR}.  Si $X$ est un
tel espace, on d\'esignera par $\OO(X)$ la $\Q_p$-alg\`ebre des fonctions
analytiques globales de $X$.  C'est un $\Q_p$-espace de Fr\'echet pour la
topologie de la convergence uniforme sur les ouverts affino\"ides.  
\ps
\newcommand{\PSC}{\mathrm{Car}^p}

On suppose que $G$ v\'erifie la propri\'et\'e suivante : pour tout sous-groupe
ouvert $H \subset G$, il n'y a qu'un nombre fini d'homomorphismes continus
$H \rightarrow \Z/p\Z$.  Cette condition est par exemple satisfaite si $G$
poss\`ede un sous-groupe dense de type fini. On consid\`ere alors le foncteur $\PSC_d(G)$,
des espaces analytiques sur $\Q_p$ vers les ensembles, associant \`a $X$
l'ensemble des pseudo-caract\`eres $G \rightarrow \OO(X)$ qui sont de
dimension $d$ et continus (en tant que fonction).  Je d\'emontre
dans~\cite{chdet} le r\'esultat suivant.  \ps

\begin{prop} Le foncteur $\PSC_d(G)$ est repr\'esentable. \end{prop}

C'est l'espace analytique repr\'esentant $\PSC_d(G)$ que j'appelle la vari\'et\'e des caract\`eres
$p$-adique de $G$ en dimension $d$. Je vais bri\`evement expliquer ses propri\'et\'es
principales ainsi que sa construction. \ps

Je dois commencer par un rappel classique sur la notion de {\it repr\'esentation
r\'esiduelle}. Soit $L$ une extension finie de $\Q_p$ et soit $\rho : G \rightarrow
\GL_d(\overline{L})$ une repr\'esentation semisimple telle que ${\rm
Trace}(\rho(G)) \subset L$. Alors $\rho$ est continue si et seulement si
${\rm Trace} \circ \rho$ l'est, auquel cas les polyn\^omes
caract\'eristiques $\det(t-\rho(g))$, $g \in G$, sont \`a coefficients
dans l'anneau des entiers $\OO_L$ de $L$. Si $k_L=\OO_L/\pi_L$ d\'esigne le
corps fini r\'esiduel de $\OO_L$, on peut alors associer \`a $\rho$ une unique repr\'esentation semisimple continue
$$\overline{\rho} : G \rightarrow \GL_d(\overline{k_L}),$$ dite
repr\'esentation r\'esiduelle de $\rho$, telle que $\det(t-\overline{\rho}(g))
\equiv \det(t-\rho(g))$ modulo l'id\'eal maximal de $\OO_L$. 
L'unicit\'e est une observation classique de Brauer-Nesbitt, et pour
l'existence on raisonne traditionnellement en se pla\c{c}ant d'abord dans une extension finie
$L'/L$ sur laquelle $\rho$ est d\'efinie, et on d\'efinit ensuite $\overline{\rho}$ en
semisimplifiant la r\'eduction modulo $\pi_{L'}$ d'un $\OO_{L'}$-r\'eseau stable par $G$.
Soit $\mathcal{I}$ l'ensemble des orbites de ${\rm
Gal}(\overline{\F}_p/\F_p)$ agissant sur l'ensemble des classes d'isomorphie de
repr\'esentations semisimples et continues $G \rightarrow \GL_d(\overline{\F}_p)$ (action sur les coefficients).
La repr\'esentation $\overline{\rho}$ d\'efinit en particulier un \'el\'ement de
$\mathcal{I}$.
\ps

Si $X$ est un espace analytique, si $T \in \PSC_d(G)(X)$, et si
$x \in X$ est de corps r\'esiduel $k(x)$ (une extension finie de
$\Q_p$), on d\'esignera par $$\rho_x : G \rightarrow
\GL_d(\overline{k(x)})$$ la repr\'esentation semisimple continue, unique \`a
isomorphisme pr\`es, de trace la compos\'ee $T : G \rightarrow \OO(X) \rightarrow
k(x)$~(\S\ref{prelimvcar}). On commence par observer que l'application $x
\mapsto \overline{\rho_x}$, $X \mapsto \mathcal{I}$, est constante sur les
composantes connexes de $X$ (en particulier, ne prend qu'un nombre
fini de valeurs si $X$ est affino\"ide). Il suit que $\PSC_d(G)$ est "r\'eunion disjointe" 
des sous-foncteurs $\PSC_d(G)[r]$, $r \in
\mathcal{I}$, param\'etrant les $T$ dont toutes les repr\'esentations
r\'esiduelles associ\'ees d\'efinissent le m\^eme \'el\'ement de
$\mathcal{I}$. On d\'emontre ensuite la repr\'esentabilit\'e de
$\PSC_d(G)[r]$ en l'identifiant \`a la fibre g\'en\'erique au sens de
Berthelot~\cite{berthelot} d'une certaine $\Z_p$-alg\`ebre locale noeth\'erienne compl\`ete
$R(r)$. \ps

Si $p>d$, on d\'efinit $R(r)$ comme l'anneau de d\'eformation
universel du pseudo-caract\`ere continu ${\rm Trace}(r) : G \rightarrow \overline{\F}_p$
(bien d\'efini modulo action de ${\rm Gal}(\overline{\F}_p/\F_p)$ \`a l'arriv\'ee). Il est bien
noeth\'erien sous notre hypoth\`ese sur $G$ par les r\'esultats de
Taylor~\cite{tay}. Il ne reste alors qu'\`a v\'erifier que la fibre
g\'en\'erique de $R(r)$ repr\'esente bien $\PSC_d(G)[r]$, ce qui est un
exercice sans difficult\'e, notamment compte tenu de l'interpr\'etation g\'en\'erale de
ce type de fibres  donn\'ee dans~\cite{dj}. Mentionnons que quand $r$ est irr\'eductible, on s'en sort m\^eme quelque soit $p$ en consid\'erant plut\^ot pour $R(r)$ l'anneau de d\'eformation universel de la repr\'esentation continue $r$ au sens de Mazur~\cite{mazurdef}, anneau par ailleurs isomorphe \`a l'anneau consid\'er\'e pr\'ec\'edemment si $p>d$~\cite{rouquier}.
\ps

Le cas $p\leq d$ et $r$ quelconque est plus vache, et occupe en fait la majeure partie de mon
article~\cite{chdet}. Dans ce cas, et g\'en\'eralement quand $d!$ n'est pas
inversible dans $A$, la th\'eorie des pseudo-caract\`eres $G \rightarrow A$
est plus pathologique, car l'identit\'e~\eqref{npseudo} n'entra\^ine pas l'identit\'e de
Cayley-Hamilton.  Cela m'a conduit \`a consid\'erer une variante
multiplicative de la th\'eorie des pseudo-caract\`eres qui a des bonnes
propri\'et\'es sur $\Z$, que j'appelle les {\it lois-d\'eterminants}.  Ces
lois sont d\'ej\`a essentiellement pr\'esentes sous forme de "normes" dans
les travaux pionniers de Procesi~\cite{Proc0}.  La d\'efinition que je
propose, bas\'ee sur la notion de loi polynomiale multiplicative au sens
de Roby (\cite{roby},~\cite{roby2}), est la suivante.

\begin{definition} Soient $A$ un anneau commutatif unitaire, $H$ un groupe,
et $d\geq 1$ un entier. Une loi-d\'eterminant sur $H$, de dimension $d$ et
\`a valeurs dans $A$, est une loi polynomiale $D: \Z[H] \rightarrow A$ qui est
multiplicative, homog\`ene de degr\'e $d$, et telle que
$D(1)=1$.  
\end{definition}

Le cas typique de telle loi est la compos\'ee d'un homomorphisme d'anneaux
$\Z[H] \rightarrow M_d(A)$ par $\det : M_d(A) \rightarrow A$.  En g\'en\'eral, une loi-d\'eterminant peut
\^etre vue comme la donn\'ee d'une collection de $d$ fonctions $H
\rightarrow A$, "les coefficients du polyn\^ome caract\'eristique", satisfaisant un ensemble d'identit\'es polynomiales en g\'en\'eral assez
complexes, dont l'identit\'e d'Amitsur sur le d\'eterminant d'une somme~\cite{amitsur},~\cite{rs}. En\footnote{Voir ma
note~\cite{chserre} pour une discussion \'el\'ementaire de cette th\'eorie en dimension $2$.} dimension $2$, il s'agit simplement de
couples $(t,d)$ de fonctions $H \rightarrow A$ telles que $d : H \rightarrow
A^\ast$ est un homomorphisme de groupes et $t : H \rightarrow A$ satisfait
$t(1)=2$ et pour tout $g,h \in H$, $t(hg)=t(hg)$ et
$d(g)t(g^{-1}h)=t(g)t(h)-t(gh)$. \ps

Le foncteur associant \`a un anneau commutatif $A$ l'ensemble des lois-d\'eterminants sur $H$, de dimension $d$ et \`a valeurs dans $A$, est trivialement
repr\'esentable.\footnote{Les travaux de Roby d\'emontrent m\^eme qu'il est
repr\'esent\'e par l'ab\'elianis\'e de l'anneau des invariants  $(\Z[H]^{\otimes
d})^{\got{S}_d}$, la loi-d\'eterminant universelle \'etant simplement $x \mapsto
x^{\otimes d}$.} Des r\'esultats importants de Donkin~\cite{Do}, Haboush, Seshadri~\cite{seshadri},
Vaccarino (\cite{vac2}, \cite{vac0}, \cite{vac1}), et Zubkov~\cite{zubkov}, assurent que cet anneau universel est de type fini sur $\Z$ si $H$ est un groupe de type fini. 
Je renvoie \`a mon article~\cite{chdet} pour plusieurs crit\`eres assurant qu'une loi-d\'eterminant $\Z[H] \rightarrow A$ est le d\'eterminant d'une repr\'esentation $\Z[H] \rightarrow {\rm M}_d(A)$; c'est par exemple vrai si $A$ est un corps alg\'ebriquement clos (sans hypoth\`ese de caract\'eristique). \ps

Pour revenir \`a la question de la repr\'esentabilit\'e de $\PSC_d(G)[r]$ quand $p \leq d$, on
d\'efinit cette fois-ci $R(r)$ comme l'anneau de d\'eformation universel de la
loi-d\'eterminant continue $\det \circ r : \Z[G] \rightarrow \overline{\F}_p$. C'est un anneau noeth\'erien par les remarques pr\'ec\'edentes, et on
v\'erifie que sa fibre g\'en\'erique au
sens de Berthelot est bien $\PSC_d(G)[r]$.  \ps

Il suit de cette longue discussion que si $X_d=\PSC_d(G)$ d\'esigne la vari\'et\'e des caract\`eres
$p$-adique de $G$ en dimension $d$, alors $X_d$ est r\'eunion disjointe admissible de sous-espaces
ouverts ferm\'es $X_d(r)$ ind\'ex\'es par les repr\'esentations
r\'esiduelles $r \in \mathcal{I}$. Chaque $X_d(r)$ est de plus isomorphe \`a
un ferm\'e dans la boule unit\'e ouverte d'une certaine dimension, dont
l'id\'eal de d\'efinition est engendr\'e par un nombre fini d'\'el\'ements de $R(r)$. Si $x \in X_d$, on notera en g\'en\'eral $$\rho_x : G \rightarrow
\GL_d(\overline{k(x)})$$ 
l'unique repr\'esentation semisimple continue dont la trace est d\'efinie
par $x$. Si $X_d^{\rm irr} \subset X_d$ d\'esigne l'ensemble des
$x$ tels que $\rho_x$ est irr\'eductible, il n'est pas difficile de
v\'erifier que $X_d^{\rm irr}$ est un ouvert Zariski de $X_d$. Une propri\'et\'e int\'eressante de $X_d^{\rm irr}$ est qu'il repr\'esente le foncteur associant \`a un affino\"ide
$Y={\rm Max}(A)$ l'ensemble des classes d'isomorphie de couples $(S,\rho)$ tels que $S$ est une $A$-alg\`ebre d'Azumaya de rang $d^2$ et $\rho : G \rightarrow S^\ast$ est un homomorphisme continu dont la trace r\'eduite induit un morphisme $Y \rightarrow X_d^{\rm irr}$. \ps

En guise d'exemple, je propose de clore cette section par le th\'eor\`eme
suivant, qui est le r\'esultat principal de ma note~\cite{chcarpad}.
Sa d\'emonstration, par r\'ecurrence sur la dimension $d$, utilise notamment les r\'esultats de Tate sur la cohomologie des repr\'esentations $p$-adiques de
${\rm Gal}(\overline{\Q}_p/\Q_p)$, des techniques de
d\'eformations \`a la Mazur~\cite{mazurdef}, ainsi que la plupart des id\'ees
introduites dans cette partie. \ps

\begin{thm} Soit $X_d$ la vari\'et\'e des caract\`eres $p$-adique, en
dimension $d$, du groupe de Galois absolu ${\rm Gal}(\overline{F}/F)$ d'une extension finie $F$ de $\Q_p$.
\begin{itemize}\ps 
\item[(i)] $X_d$ est \'equidimensionnelle de dimension $[F:\Q_p]d^2+1$.  \ps
\item[(ii)] L'ouvert Zariski $X_d^{\rm irr}$ est dense dans $X_d$.  \ps
\item[(iii)] Si $d \neq 2$ ou si $F \neq \Q_p$, $X_d^{\rm irr}$ co\"incide exactement avec
le lieu r\'egulier de $X_d$.  \ps
\item[(iv)] Quand $F=\Q_p$, le lieu singulier de $X_2$ est celui
param\'etrant les torsions par un caract\`ere de $1\oplus \omega$ o\`u
$\omega$ est le caract\`ere cyclotomique.\ps
\end{itemize}
\end{thm}

On aimerait avoir plus d'outils pour \'etudier les vari\'et\'es de
caract\`eres $p$-adique en g\'en\'eral. On voudrait aussi en savoir un peu
plus sur les $X_d(r)$ dans le contexte de l'\'enonc\'e ci-dessus. Hormis le cas $d=1$, le seul cas
vraiment facile est celui o\`u $r$ est irr\'eductible et satisfait $r
\not\simeq r\otimes \omega$ o\`u $\omega$ est le caract\`ere cyclotomique
modulo $p$ (condition automatique si $p-1$ ne divise pas $d$ quand $F=\Q_p$).  Dans ce cas,
les r\'esultats de Tate et Mazur montrent que $X_d(r)$ est simplement une boule ouverte.\footnote{Plus exactement, c'est
la restriction \`a $\Q_p$ d'une boule ouverte d\'efinie sur l'extension non
ramifi\'ee de $\Q_p$ de corps r\'esiduel le corps de d\'efinition de $r$.} 
M\^eme dans le cas particulier $K=\Q_p$ et $d=2$, je ne connais pas $X_2(r)$ dans tous les cas, notamment quand 
$p=2$ et $r$ est la repr\'esentation triviale. Pour des r\'esultats, un peu \'eparpill\'es, sur ces th\`emes, voir aussi~\cite{kisinfm},~\cite{kisinfern},~\cite{paskunas},~\cite{chcarpad} et~\cite{boeckle}.

\newpage

\section{La foug\`ere infinie et la densit\'e des points modulaires}\label{sectionfougereglobale}

\subsection{Le lieu g\'eom\'etrique des vari\'et\'es de caract\`eres
$p$-adique du groupe de Galois absolu de $\Q$}\label{introlieugeom}

Soient $p$ un nombre premier, $S$ un ensemble fini de nombres premiers
contenant $p$, et soit $G_{\Q,S}$ le groupe Galois d'une extension
alg\'ebrique maximale de $\Q$ non ramifi\'ee hors de $S$ (et de l'infini). 
On fixe un entier $d\geq 1$ et on consid\`ere la vari\'et\'e des
caract\`eres $p$-adique $X_d$ de $G_{\Q,S}$ en dimension
$d$~(\S\ref{carpadique}).  La structure m\^eme de $X_d$ est tr\`es loin
d'\^etre comprise en g\'en\'eral, et nous n'aurons rien \`a apporter sur
cette question.  Cela n'aura pas tellement d'incidence sur les questions qui
vont nous int\'eresser ici.  En guise d'exemple, disons simplement que si
$$r : G_{\Q,S} \rightarrow \GL_d(\overline{\F}_p)$$ est une repr\'esentation
continue irr\'eductible, et disons si $p>2$, alors $X_d(r)$ est de dimension
$\geq 1+2i(d-i)$, o\`u $i=\dim {\rm Ker}(r(c)-{\rm id})$ et $c \in G_{\Q,S}$
d\'esigne une conjugaison complexe.  Cela suit en effet des propri\'et\'es de la cohomologie
galoisienne de $G_{\Q,S}$ (Poitou, Tate~\cite{milneadt}) et de la th\'eorie de
Mazur~\cite{mazurdef}. On dira que $r$ est {\it r\'eguli\`ere} si
$H^2(G_{\Q,S},{\rm ad}(r))=0$.  Dans ce cas, assez fr\'equent en pratique,
$X_d(r)$ est simplement la $K$-boule unit\'e ouverte\footnote{Nous entendons
par l\`a la boule unit\'e ouverte de l'espace affine analytique sur $K$ que
l'on voit par restriction comme un espace analytique sur $\Q_p$.  C'est
aussi la fibre g\'en\'erique de l'anneau local complet
$\OO_K[[t_1,\cdots,t_n]]$.} de dimension $1+2i(d-i)$, o\`u $K$ est
l'extension non ramifi\'ee de $\Q_p$ dont le corps r\'esiduel est le corps
de d\'efinition de $r$.  \ps

Soit $x \in X_d$ et soit $\rho_x : G_{F,S} \rightarrow
 \GL_d(\overline{k(x}))$ la repr\'esentation semisimple continue associ\'ee. 
 On dira que le point $x$ est {\it g\'eom\'etrique} s'il existe une
 vari\'et\'e projective lisse $Y$ d\'efinie sur $\Q$, et des entiers $m \in
 \Z$ et $i \geq0$, tels que $\rho_x$ est un sous-quotient de $$H^i_{\rm
 et}(Y_{\overline{\Q}},\Q_p(m)) \otimes_{\Q_p} \overline{k(x)}.$$ Le
 sous-ensemble $X_d^{\rm geom} \subset X_d$ des points g\'eom\'etriques est
 assez myst\'erieux, et sujet de nombreux probl\`emes ouverts. Par exemple, d'apr\`es les conjectures du
folklore rappel\'ees au~\S\ref{conjfolk}, c'est aussi l'ensemble des $x$ tels que $\rho_x $ est de la forme $\rho_{\pi,\iota}$ pour une repr\'esentation automorphe
cuspidale alg\'ebrique $\pi$ de $\GL_d$ sur $F$, ou encore tel que $\rho_x$
soit g\'eom\'etrique au sens de Fontaine-Mazur, i.e.  simplement de De Rham
aux places divisant $p$. \ps

Il est naturel de se demander si $X_d^{\rm geom}$ poss\`ede une quelconque structure particuli\`ere en tant que partie de $X_d$.  Observons d\'ej\`a que
$X_d^{\rm geom}$ est un ensemble d\'enombrable (donc petit), car il en va de m\^eme des classes d'isomorphie de vari\'et\'es projectives sur $\Q$, il ne
contient donc pas de sous-vari\'et\'e de dimension $>0$.  La question qui va nous int\'eresser principalement dans cette partie est la suivante. 
 
 \begin{question} Que peut-on dire de l'adh\'erence Zariski de $X_d^{\rm geom}$ dans $X_d$ ? \end{question}

Je rappelle que l'adh\'erence Zariski d'une partie $A$ d'un espace analytique s\'epar\'e $X$ est le plus petit ferm\'e analytique r\'eduit de $X$ contenant $A$. Il existe par noeth\'erianit\'e des affino\"ides. On dit notamment que $A$ est Zariski-dense dans $X$ si son adh\'erence Zariski est la nilreduction de $X$ (voir~\cite{conradirr}).\ps

C'est un exercice de d\'eterminer $X_1$ et $X_1^{\rm geom}$ lorsque $d=1$.
Le th\'eor\`eme de Kronecker-Weber montre en effet que $X_1$ est une r\'eunion
disjointe d'espaces analytiques de la forme $Z \times \mathcal{W}$ o\`u $Z$
est fini et o\`u $$\mathcal{W}={\rm Hom}(\Z_p^\ast,\mathbb{G}_m)$$ est la
vari\'et\'e des caract\`eres $p$-adique du groupe $\Z_p^\ast$ en dimension
$1$.  Cet espace $\mathcal{W}$, parfois appel\'e {\it espace des poids} en
th\'eorie des formes modulaires $p$-adiques, joue un r\^ole fondamental dans
ce sujet.\footnote{Par d\'efinition, ses points dans une $\Q_p$-alg\`ebre
affino\"ide $A$ sont simplement les homomorphismes de groupes continus
$\Z_p^\ast \rightarrow A^\ast$.} Il est isomorphe \`a $\widehat{\mu} \times
\mathcal{B}$ o\`u $\widehat{\mu}$ d\'esigne la vari\'et\'e (finie) des
caract\`eres de dimension $1$ du sous-groupe de torsion $\mu \subset
\Z_p^\ast$, et $\mathcal{B}=\{t \in \AAA^1, |t|<1\}$ est la boule unit\'e
ouverte de dimension $1$ sur $\Q_p$, via l'application $\chi \mapsto
(\chi_{|\mu(\Z_p^\ast)},\chi(1+2p)-1)$.  \ps

Les points g\'eom\'etriques de $X_1$ sont alors ceux dont la composante dans
$\WW$ est un caract\`ere d'ordre fini fois un caract\`ere de $\Z_p^\ast$ de
la forme $x \mapsto x^k$ avec $k \in \Z$.  Ces points sont notamment
Zariski-denses dans $X_1$, et ce
dans un sens tr\`es fort : c'est aussi une partie d'accumulation. On dira qu'une partie $A$ d'un espace analytique $X$ s'accumule
en $x \in X$ si $x$ admet un syst\`eme de voisinages ouverts affino\"ides $U$
tels que $A \cap U$ est Zariski-dense dans $U$; on dira que $A$ est
d'accumulation si elle s'accumule en chacun de ses points. \ps
 
Le premier r\'esultat int\'eressant concerne le cas $d=2$ de la question ci-dessus, il est d\^u \`a Gouv\^ea et Mazur~\cite{gm}. En effet, 
ces auteurs ont mis en \'evidence une partie de  $X_2$ aux allures de type fractal, qu'ils l'ont appell\'ee la {\it foug\`ere infinie}.  Il est
 remarquable qu'elle leur ait \'et\'e sugg\'er\'ee par leur recherche
 num\'erique dans~\cite{gmnum} de familles analytiques de formes modulaires
 g\'en\'eralisant la construction pionni\`ere de Hida~\cite{hida}.  Les observations de Gouv\^ea-Mazur ont \'et\'e valid\'ees
 ensuite par la construction de ces familles par Coleman
 dans~\cite{coleman1} et~\cite{coleman2}.  Une cons\'equence notable de la
 foug\`ere infinie est qu'elle permet de montrer que les points modulaires
 sont Zariski-denses dans certaines composantes connexes de
 $X_2$.  \ps

\subsection{La foug\`ere infinie de Gouv\^ea et Mazur}\label{arcgm} Nous allons commencer par faire des rappels sur les observations de Gouv\^ea et Mazur. On suppose d\'esormais
$S=\{p\}$ pour simplifier et on fixe un plongement $\iota : \overline{\Q} \rightarrow
\overline{\Q}_p$.  Soit $M_k({\rm SL}(2,\Z))$ l'espace des formes modulaires
de poids $k$ pour le groupe ${\rm SL}(2,\Z)$.  On rappelle que si
$$f=\sum_{n\geq 0} a_n q^n \in M_k(\SL(2,\Z))$$ est propre pour tous les
op\'erateurs de Hecke, et normalis\'ee par $a_1=1$, on peut lui associer
d'apr\`es Deligne une repr\'esentation continue $\rho_{f} :
G_{\Q,S} \rightarrow \GL_2(\overline{\Q}_p)$ telle que le polyn\^ome
caract\'eristique d'un Frobenius g\'eom\'etrique en $\ell \neq p$ soit
$X^2-\iota(a_p)X+p^{k-1}$ : c'est un cas tr\`es particulier des
constructions discut\'ees au~\S\ref{conjfolk}. En particulier une telle
forme d\'efinit un point $x_f \in X_2$ : il est m\^eme g\'eom\'etrique
d'apr\`es Deligne et Scholl~\cite{scholl}.\ps

Rappelons comment Coleman associe en g\'en\'eral deux  germes d'arcs canoniques issus de $x_f$ dans
$X_2$.  On part du polyn\^ome $P_f(X)=X^2-\iota(a_p)X+p^{k-1}$. Soit
$\alpha$ une racine de $P_f$ : c'est un entier alg\'ebrique\footnote{... en
fait, un $p$-nombre de Weil de poids $k-1$ (Deligne),} dont la valuation $p$-adique est
comprise entre $0$ et $k-1$. Coleman associe au couple $(f,\alpha)$ une sous-courbe (sans point isol\'e) d'un voisinage ouvert de $x_f$
dans $X_2$, dont le germe en $x_f$ est canonique, que l'on notera
$$C_{(f,\alpha)} \subset X_2.$$ Cette courbe a plusieurs propri\'et\'es
remarquables. D'une part, ses points de la forme $x_{f'}$ y sont
Zariski-denses et d'accumulation (ce qui force bien entendu le poids des
$f'$ apparaissant \`a \^etre non born\'e). D'autre part, si $x_{f'} \in
C_{(f,\alpha)}$ alors\footnote{Au sens strict, cette propri\'et\'e ne d\'ecoule pas imm\'ediatement de la construction de Coleman, il faut y rajouter le r\'esultat principal de~\cite{kisinoc}. Je l'utiliserai cependant ici pour simplifier l'exposition.} $P_{f'}$ admet une racine ayant m\^eme valuation que
$\alpha$. Quand la valuation de $\alpha$ est nulle, la construction de
$C_{(f,\alpha)}$ avec ces propri\'et\'es remonte \`a Hida~\cite{hida}, le premier exemple remontant \`a Serre~\cite{serrepadic} ("famille d'Eisenstein ordinaire").\ps

Nous donnerons quelques indications sur cette construction de Coleman au~\S\ref{constcoleman}.  Je voudrais d'abord expliquer comment la foug\`ere infinie de Gouv\^ea-Mazur en jaillit.
Sa d\'efinition donn\'ee dans~\cite{gm} est simplement la r\'eunion ensembliste de tous les arcs de Coleman
$$\mathcal{F}_{GM} := \bigcup_{(f,\alpha)} C_{(f,\alpha)} \subset X_2.$$

Si l'on part de $f$, et si $P_f$ admet deux racines $\alpha,\beta$ de valuations $v(\alpha)\neq v(\beta)$ distinctes, Gouv\^ea et Mazur observent que l'intersection $C_{(f,\alpha)} \cap C_{(f,\beta)}$ est finie. En effet, il suit de la construction de Coleman que l'on peut supposer que $C_{(f,\alpha)}$ et $C_{(f,\beta)}$ sont deux ferm\'es d'un m\^eme voisinage affino\"ide de $x_f$, de sorte que si l'intersection est infinie, elle contient une composante irr\'eductible de cette courbe, et donc une infinit\'e de points de la forme $x_{f'}$. Pour un tel $f'$, $P_{f'}$ admet alors une racine de valuation $v(\alpha)$ et une autre de valuation $v(\beta) \neq v(\alpha)$... donc $f'$ est de m\^eme poids $k=v(\alpha)+v(\beta)+1$ que la forme $f$, ce qui contredit l'infinit\'e des $f'$. Il n'est en fait pas trop difficile de d\'emontrer que $C_{(f,\alpha)}$ et $C_{(f,\beta)}$ sont m\^eme transverses en $x_f$ d\`es que $0<v(\alpha),v(\beta) < k-1$~\cite{chens}. \ps

\begin{figure}[htp]
\centering
\includegraphics[scale=0.75]{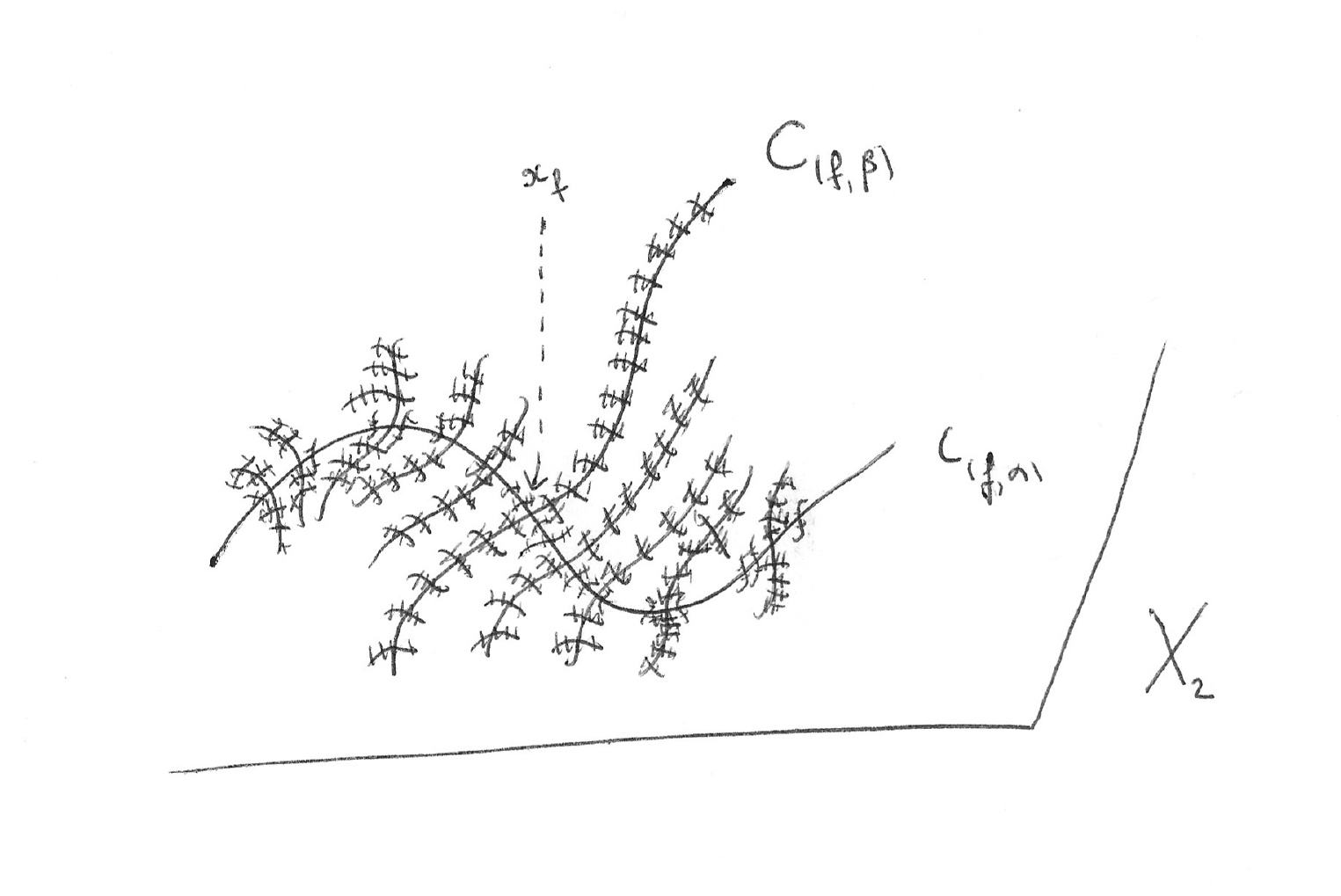}
\caption{La foug\`ere infinie de Gouv\^ea-Mazur}
\label{figfougere1}
\end{figure}

Observons enfin que si l'hypoth\`ese que $P_f$ admet deux racines de valuations distinctes n'est pas satisfaite, elle le sera de toutes fa\c{c}ons pour tous les point de la forme $x_{f'}$ avec $f'$ de poids $k' \neq k$ qui sont dans un arc de Coleman passant par $x_f$, \`a cause de la condition $v(P_{f'}(0))=k'-1$. La structure fractale de $\mathcal{F}_{GM}$ est alors bien claire, ainsi que le terme de {\it foug\`ere infinie} ! (voir la figure~\ref{figfougere1} ci-dessous). \ps

Au final, il est alors ais\'e de voir que l'adh\'erence Zariski dans $X_2$ des points de la forme $x_f$, disons $W$, a toutes ses composantes irr\'eductibles de dimension au moins $2$. Notons $x_{f,m} \in X_2$ le point de repr\'esentation galoisienne $\rho_x=\rho_f \otimes \omega^m$, $\omega$ \'etant le caract\`ere cyclotomique $p$-adique et appelons {\it modulaire} un tel point. La th\'eorie de Sen~\cite{sen} montre que $W$ reste dans un ferm\'e partout strict de $X_2$, un nombre de Hodge-Tate g\'en\'eralis\'e y \'etant constant et \'egal \`a $0$, le th\'eor\`eme suivant d'en d\'eduit.
\bigskip

\begin{thm} (Gouv\^ea-Mazur, Coleman) L'adh\'erence Zariski des points modulaires dans $X_2$ a toutes ses composantes irr\'eductibles de dimension $\geq 3$. \par
Si $r$ est irr\'eductible, r\'eguli\`ere\footnote{La philosophie des valeurs sp\'eciales de fonctions $L$ sugg\`ere que cette r\'egularit\'e est tr\`es fr\'equente. Par exemple, si $\Delta \in M_{12}({\rm SL}(2,\Z))$ est la fonction discriminant, alors $\overline{\rho_{\Delta}}$ est irr\'eductible et r\'eguli\`ere d\`es que $p>13$ et $p \neq 691$ (Mazur, Weston~\cite{weston}).}, et de la forme $\overline{\rho_f}$ pour une certaine forme $f$, auquel cas $X_2(r)$ est une boule de dimension $3$, alors les points modulaires sont Zariski-denses et d'accumulation dans $X_2(r)$. 
\end{thm}

Rappelons que d'apr\`es la conjecture de modularit\'e de Serre~\cite{serremod1} (voir aussi \cite{serremod2}), d\'emontr\'ee par Khare et Wintenberger~\cite{kw}, l'hypoth\`ese $r \simeq  \overline{\rho_f}$ pour un certain $f$ est automatique d\`es que $r$ est impaire, i.e. ${\rm Tr}(r(c))=0$. L'argument ci-dessus montre plus g\'en\'eralement que si $F \subset X_2$ est une composante irr\'eductible contenant un point lisse qui est modulaire, alors les points modulaires sont Zariski-denses et d'accumulation dans $F$. On conjecture que l'ouvert ferm\'e $X_2^{\rm odd} \subset X_2$ param\'etrant les $x$ tels que ${\rm Trace}(\rho_x(c))=0$ est \'equi-dimensionnel de dimension $3$, et les r\'esultats de Wiles~\cite{wilesfermat} et Taylor-Wiles~\cite{tw} vont dans cette direction. En combinant notamment ces id\'ees, Boeckle a pu montrer que les points modulaires sont Zariski-denses et d'accumulation dans la plupart des composantes connexes de $X_2^{\rm odd}$. Je renvoie \`a l'expos\'e Bourbaki d'Emerton~\cite{emerton} pour plus de r\'ef\'erences sur ces questions. \ps

\subsection{Retour sur les travaux de Coleman~\cite{coleman1},\cite{coleman2}}\label{constcoleman} Avant d'aller plus loin, revenons sur l'approche de Coleman pour construire les arcs $C_{(f,\alpha)}$. La forme $f$ engendre un sous-espace de dimension $2$ dans l'espace $M_k(\Gamma_0(p))$ des formes modulaires pour le sous-groupe de congruence
$\Gamma_0(p)$ usuel, \`a savoir $\langle f(\tau), f(p\tau) \rangle$, sur lequel l'op\'erateur de Hecke $U_p$ d'Atkin-Lehner admet le polyn\^ome caract\'eristique $P_f$. Le choix de $\alpha$ correspond donc au choix d'une forme propre $f_\alpha$ pour $U_p$ dans cet espace. Elle est n\'ecessairement propre pour les $T_\ell$ avec $\ell \neq p$, de m\^eme valeurs propres que $f$. Quand $P_f$ admet deux racines $\alpha \neq \beta$, les formes $f_\alpha$ et $f_\beta$ sont parfois appel\'ees {\it formes jumelles} associ\'ees \`a $f$.  \ps

On consid\`ere le mod\`ele rationnel usuel de la courbe modulaire ${\rm X}_0(p)$, et plus exactement son extension des scalaires \`a $\Q_p$, que l'on verra comme un courbe analytique $p$-adique. 
Cette courbe vient avec un fibr\'e en droites canonique\footnote{On mettra sous le tapis
dans cette discution les probl\`emes li\'es \`a la torsion dans
$\Gamma_0(p)$, qui ne sont pas difficiles \`a contourner en introduisant des niveaux auxilliaires.}$\omega$. Le choix de $\iota$ et le th\'eor\`eme de type GAGA de Kiehl permet de voir $M_k(\Gamma_0(p))$ comme des sections globales de $\omega^k$ sur $X_0(p)$. Si $k \in \Z$, on consid\`ere alors l'espace
$M_k^\dagger$ des sections du faisceau $\omega^k$ convergentes sur un voisinage strict d'un certain ouvert affino\"ide fixe $Z \subset X_0(p)$. Cet ouvert $Z$ est la composante connexe de la pointe infini dans $X_0(p) \backslash ]S[$, o\`u $]S[$ d\'esigne le tube supersingulier, c'est-\`a-dire l'image inverse par la r\'eduction canonique $X_0(p) \rightarrow X_0(p)_{\overline{\F}_p}$ 
de l'ensemble des points param\'etrant les courbes elliptiques supersinguli\`eres.  Bien s\^ur, $M_k(\Gamma_0(p))\subset_\iota M_k^\dagger$, mais ce dernier est un espace gigantesque, une r\'eunion croissante de Banach $p$-adiques; les \'el\'ements de $M_k^\dagger$ sont appel\'es formes modulaires $p$-adiques surconvergentes de poids $k$. \ps

Coleman construit alors un faisceau $\mathcal{M}$ de modules
topologiques sur l'espace des poids $\WW$, qui est localement libre dans un
certain sens, et dont la sp\'ecialisation en tout $k \in \Z$, vu comme un
\'el\'ement de $\WW$ par le caract\`ere $x \mapsto x^k$, co\"incide
avec $M_k^\dagger$. Un point clef est que les $M_k^\dagger$, en tant qu'espaces vectoriels topologiques, ne d\'ependent pas vraiment de $k$ car le fibr\'e $\omega$ est trivial sur $Z$.  Coleman v\'erifie que l'action naturelle sur les $M_k^\dagger$ des op\'erateurs de Hecke $T_\ell$ avec $\ell \neq p$,
ainsi que l'op\'erateur d'Atkin-Lehner $U_p$, s'\'etend en une action $\WW$-lin\'eaire
continue sur $\mathcal{M}$. \ps

Un point important est que $U_p$ agit comme un op\'erateur compact sur les $M_k^\dagger$ (Dwork, Katz~\cite{katzpmod}). Coleman introduit dans~\cite{coleman2} une notion d'op\'erateurs compacts en famille et \'etablit une th\'eorie spectrale \`a la Fredholm pour ces derniers, g\'en\'eralisant~\cite{serrecompact}.
Elle s'applique \`a $(U_p,\mathcal{M})$ et lui permet de d\'efinir notamment une s\'erie de Fredholm $$\det(1-T {U_p}|{\mathcal{M}}) \in \OO( \WW \times \AAA^1) \subset \OO(\WW)[[T]].$$ En \'etudiant la factorisation de cette s\'erie localement au voisinage d'un point de $\WW$, il d\'emontre que toute forme modulaire $f \in M_k^\dagger$ qui est propre pour les $T_\ell$, et propre pour $U_p$ de valeur propre non nulle, fait partie d'une famille analytique de formes modulaires $p$-adiques surconvergentes de poids variant dans un petit voisinage de $k \in \WW$, \'egalement propres pour les $T_\ell$ et $U_p$. Le fait que ces familles contiennent beaucoup de formes modulaires
classiques suit alors de son autre travail~\cite{coleman1} ("crit\`ere de
classicit\'e de Coleman"). Je renvoie \`a~\cite{coleman2} pour la d\'efinition exacte de familles de formes modulaires $p$-adiques. Disons simplement ici que les arcs introduits ci-dessus ne sont que des
traductions galoisiennes de ces d\'efinitions et constructions. \ps

Mentionnons en particulier que si $f \in M_k^\dagger$ est une forme modulaire $p$-adique surconvergente, disons $f=\sum_{n\geq 0} a_n q^n$ avec $a_1=1$, et si $f$ est propre pour $U_p$ et les $T_\ell$ avec $\ell \neq p$, les valeurs propres respectives \'etant alors n\'ecessairement $a_p$ et les $a_\ell$, la th\'eorie de Coleman permet d'associer \`a $f$ une unique repr\'esentation semisimple continue $\rho_f : G_{\Q,S} \rightarrow \GL_2(\overline{\Q}_p)$ telle que $\det(X-\rho_f({\rm Frob}_\ell))=X^2-a_\ell X + \ell^{k-1}$ pour tout $\ell \neq p$. Cette repr\'esentation est obtenue par un proc\'ed\'e formel d'interpolation $p$-adique \`a partir de la construction de Deligne. \ps

\newcommand{\cX}{\mathcal{X}}
\newcommand{\cE}{\mathcal{E}}
\newcommand{\Dc}{{\mathrm D}_{\rm cris}}
\newcommand{\Ref}{\Phi}

\subsection{The eigencurve}\label{eigencurvefern}L'histoire de la foug\`ere infinie ne s'arr\`ete pas tout \`a fait ici. En effet, quelques ann\'ees plus tard, Coleman et Mazur ont revisit\'e dans~\cite{colemanmazur} les constructions de Coleman et les ont "recoll\'ees" sous la forme d'un object unique tout aussi merveilleux que la foug\`ere infinie qu'ils appellent~{\it the eigencurve}. Sa d\'efinition la plus simple est de consid\'erer l'adh\'erence Zariski $$\mathcal{C} \subset X_2 \times \mathbb{G}_m$$ des couples de la forme $(x_f,\alpha^{-1})$ o\`u $f$ parcourt les formes modulaires propres, normalis\'ees, pour ${\rm SL}(2,\Z)$, et o\`u $P_f(\alpha)=0$. Ils d\'emontrent que $\mathcal{C}$ est une courbe \'equidimensionnelle. C'est un fait absolument remarquable car nous avons vu pr\'ecis\'ement que l'adh\'erence Zariski des $x_f$ dans $X_2$ est de dimension $\geq 2$. \ps

Cette courbe est d'abord construite ind\'ependamment de consid\'erations galoisiennes dans~\cite{colemanmazur}, o\`u elle appara\^it comme \'etant l'espace des param\`etres $\mathcal{D}$ de la famille $p$-adique universelle de formes surconvergentes propres, de valeur propre inversible pour $U_p$. Par construction, elle est munie d'un morphisme fini et surjectif 
$$\nu : \mathcal{D} \rightarrow \mathcal{Z}(U_p)$$
o\`u $\mathcal{Z}(U_p) \subset \WW \times \mathbb{G}_m$ est le ferm\'e d\'efini par $\det(1-TU_p|\mathcal{M})=0$ ("hypersurface de Fredholm de $U_p$"), et donc en particulier d'un
 morphisme $$\kappa : \mathcal{D} \rightarrow \WW$$ obtenu en composant $\nu$ par la premi\`ere projection. Le morphisme $\kappa$ est plat et localement fini : tout point de $\mathcal{D}$ admet un voisinage ouvert affino\"ide $V$ tel que $\kappa(V) \subset \WW$ soit un ouvert affino\"ide et tel que $\kappa_{|V} : V \rightarrow \kappa(V)$ soit fini et plat. En revanche, $\kappa$ est en g\'en\'eral de fibres infinies. \ps
 
 Des arguments de Coleman et Mazur d\'emontrent que les points de la courbe $\mathcal{D}$  param\'etrant les formes modulaires classiques sont Zariski-denses et d'accumulation. Un argument~\`a la Wiles permet alors de d\'efinir un unique pseudo-caract\`ere $G_{\Q,S} \rightarrow \OO(\mathcal{D})$ dont les \'evaluations en ces points "classiques" sont ceux de Deligne ($\mathcal{D}$ est r\'eduite). Il en r\'esulte un morphisme canonique $\mathcal{D} \rightarrow X_2$. Il n'est pas difficile de voir sur la construction de $\mathcal{D}$ que le morphisme induit $\mathcal{D} \rightarrow X_2 \times \mathbb{G}_m$, le second facteur \'etant la compos\'ee de $\nu$ par la projection vers $\mathbb{G}_m$, est une immersion ferm\'ee. On en d\'eduit alors que $\mathcal{D} \simeq \mathcal{C}$  (d'ailleurs cet argument me semble plus simple que celui donn\'e dans~\cite{colemanmazur}). \ps

Soit $\mu : \mathcal{C} \rightarrow X_2$ le morphisme d\'eduit de la premi\`ere projection. L'arc de Coleman $C_{(f,\alpha)}$ appara\^it  simplement dans ce langage comme l'image par $\mu$ d'un voisinage ouvert affino\"ide assez petit de $(x_f,\alpha)$ dans $\mathcal{C}$. J'appellerai foug\`ere infinie {\it compl\`ete}  l'ensemble 
			$$\mathcal{F}=\mu(\mathcal{C}) \subset X_2.$$
C'est une sorte de prolongement analytique canonique de tous les arcs de Coleman constituant $\mathcal{F}_{GM}$, et dont la courbe $\mathcal{C}$ est une sorte de d\'epliage en tous les points doubles de la forme $x_f$. Cet ensemble admet une d\'efinition plus directe : c'est l'ensemble des points $x \in X_2$ tels que $\rho_x$ est la repr\'esentation galoisienne associ\'ee \`a une forme modulaire $p$-adique surconvergente, propre pour les op\'erateurs de Hecke, de valeur propre non nulle pour $U_p$ (voir~\cite{colemanmazur} pour une d\'efinition, ici n\'ecessaire, de telles formes de poids non seulement dans $\Z$ mais dans $\WW$).\ps 

	Malheureusement, tr\`es peu de choses sont connues sur la g\'eom\'etrie globale de $\mathcal{C}$. Coleman et Mazur conjecturent  qu'elle n'admet qu'un nombre fini de composantes irr\'eductibles, ce qui serait assez formidable : cela signifierait qu'en prolongeant analytiquement un seul arc de Coleman bien choisi dans $X_2$ on pourrait parcourir de mani\`ere Zariski-dense une composante irr\'eductible donn\'ee de $X_2$ ! \`A ma connaissance, on ne conna\^it aucun exemple  de ce ph\'enom\`ene. \ps
	
\subsection{Quelques contributions sur la structure de la courbe de Coleman-Mazur} Je voudrais terminer cette discussion sur la courbe $\mathcal{C}$ de Coleman-Mazur en mentionnant quelques r\'esultats (assez modestes !) que j'ai obtenus sur sa structure. Le premier, un travail en collaboration avec J. Bella\"iche~\cite{bchjimj}, concerne les points {\it Eisenstein critiques}. Ce sont les points $e_k=(x_f,\alpha^{-1})$ o\`u $f$ est la s\'erie d'Eisenstein de poids $k\geq 4$ pour $\SL(2,\Z)$ et $\alpha=p^{k-1}$. 
Nous nous sommes notamment int\'eress\'es \`a ces points car nous avions cr\^u un temps que l'\'etude de la repr\'esentation galoisienne port\'ee par leurs voisinages dans $\mathcal{C}$ permettrait de d\'emontrer la non-nullit\'e de la fonction z\^eta $p$-adique de Kubota-Leopold aux entiers impairs positifs, un fameux probl\`eme ouvert dans l'esprit de la conjecture de Leopoldt. Nous avons finalement montr\'e plut\^ot le r\'esultat suivant.
		
\begin{thm} La courbe $\mathcal{C}$ est lisse en $e_k$ pour tout $k\geq 4$. De plus, le morphisme $\kappa$ est \'etale en $e_k$ si et seulement si $\zeta_p(k-1) \neq 0$ o\`u $\zeta_p$ d\'esigne la fonction z\^eta $p$-adique de Kubota-Leopold. 
\end{thm}

	Le second r\'esultat que je voudrais mentionner concerne l'\'etude de $\mathcal{C}$ pour des petites valeurs de $p$ (je rappelle que je travaille encore sous l'hypoth\`ese $S=\{p\}$, i.e. que le "niveau mod\'er\'e est $1$"). Au moins quand $p=2$, certains ouverts explicites de $\mathcal{C}$ ont \'et\'e mis en \'evidences par plusieurs auteurs, dont Emerton, Smithline, Buzzard, Kilford, Jacobs, Calegari et  Loeffler. Je renvoie \`a~\cite{BCa} et \cite{BK} pour des r\'ef\'erences sur ce sujet. Buzzard et Kilford d\'eterminent notamment compl\`etement $\mathcal{C}$ pour $p=2$ au dessus du compl\'ementaire d'un certain disque ferm\'e de $\WW$. Ces auteurs proc\`edent en g\'en\'eral par calculs explicites de l'action de la correspondance de Hecke $U_2$ sur certaines r\'egions de $X_0(2)$. Buzzard et Calegari d\'emontrent aussi que le morphisme $\kappa$ est propre quand $p=2$~\cite{BCa2} (voir aussi~\cite{calegaripropre}). \ps
	
	Ma contribution est un peu diff\'erente. Tout d'abord, j'ai d\'emontr\'e dans~\cite{chdet} que pour $p=2$ alors $X_2^{\rm odd}$ est la boule unit\'e ouverte de dimension $3$ sur $\Q_2$. Cela suit du fait, d\^u \`a Tate~\cite{tatemod2}, que la seule repr\'esentation r\'esiduelle dans ce cas est la triviale, et des m\'ethodes que j'ai d\'ej\`a introduites au~\S\ref{carpadique}. C'\'etait d'ailleurs initiallement dans ce but que j'avais commenc\'e \`a m'int\'eresser aux lois-d\'eterminants. \ps
	
	J'ai aussi observ\'e, toujours pour $p=2$, que la seule composante irr\'eductible
$F$ de $\mathcal{C}$ telle que $\kappa_{|F} : F \rightarrow \WW$ soit finie est la composante Eisenstein ordinaire. C'est l'occasion ici de donner l'argument, qui est tr\`es simple (mais non publi\'e). Soit $F \subset  \mathcal{C}$ une telle composante, et soit $f \in \OO(F)$ la fonction analytique associ\'ee \`a $U_2$, i.e. l'inverse du param\`etre de $\mathbb{G}_m$ restreint \`a $\mathcal{C}$. Soit $g={\rm Norme}_{F/\WW}(f) \in \OO(\WW)$. C'est une fonction qui ne s'annule pas sur $\WW$ et qui est partout born\'ee par $1$, car $f$ l'est sur $F$, c'est donc un \'el\'ement  inversible de $\Z_2[[\Z_2^\ast]][1/2]$. En particulier, sa valuation est constante sur $\WW$. Mais il suit de~\cite{BK} que $v(g(w))$ tends vers $0$ quand $w$ s'approche du bord de $\WW$ : l'\'el\'ement $g$ est donc de valuation constante \'egale \`a $0$. Il suit que $v(f(x))=0$ pour tout $x \in F$, i.e. $F$ est inclus dans le lieu ordinaire, Q.E.D. Si l'on combine cet argument \`a la propret\'e de $\kappa$ d\'emontr\'ee par Buzzard et Calegari, il semblerait\footnote{Je parle ici au conditionnel car la notion de propret\'e employ\'ee par ces auteurs semble un peu plus restrictive que celle dont on a besoin ici pour conclure, et je n'ai pas v\'erifi\'e que leurs arguments s'\'etendent.} que le seul facteur irr\'eductible polynomial de $\det(1-T{U_2}|\mathcal{M})$ soit $1-T$, les autres facteurs irr\'eductibles  \'etant alors de degr\'e infini. \ps
	
Enfin, ma derni\`ere contribution est inspir\'ee par les travaux de
Colmez~\cite{colmeztri} et~\cite{colmezinvL}.  Je consid\`ere l'\'eclat\'e
de $\WW \times \mathbb{G}_m$ en tous les points {\it sp\'eciaux},
i.e.  de la forme $(k,\lambda^{-1})$ o\`u $k\geq 2$ est un entier pair et
$\lambda^2=p^{k-2}$.  Les points $x$ de $\mathcal{C}$ tels que
$\nu(x)=(k,\lambda^{-1})$ est sp\'ecial sont exactement ceux param\'etrant
les formes modulaires propres $M_k(\Gamma_0(p))$ qui sont nouvelles en $p$
au sens d'Artkin-Lehner, i.e.  non de la forme $f_\alpha$ avec $f \in
M_k({\rm SL}_2(\Z))$ et $P_f(\alpha)=0$.  D\'esignons par
$\widetilde{\mathcal{Z}(U_p)}$ le transform\'e strict de l'hypersurface de
Fredholm $\mathcal{Z}(U_p)$ dans cet \'eclat\'e.  Le r\'esultat suivant est
le th\'eor\`eme principal de mon article~\cite{chhecke}.
	
\begin{thm}\label{eclatementeigen} Pour $p \in \{2,3,5,7\}$, le morphisme structural $\mathcal{C} \rightarrow \mathcal{Z}(U_p)$ se rel\`eve en un isomorphisme $\mathcal{C} \isomo \widetilde{\mathcal{Z}(U_p)}$. \end{thm}

Ces valeurs de $p$ sont celles telles que ${\rm X}_1(p)$ est de genre nul. Ces $p$ \'etant \'egalement r\'eguliers au sens de Kummer, cela entra\^ine que les $\rho_f$ sont uniquement d\'etermin\'ees par leurs restrictions \`a un groupe de d\'ecomposition en $p$, une propri\'et\'e que m'avait d\'ej\`a fait observer Emerton pour $p=2$. La d\'emonstration consiste ensuite \`a comprendre pourquoi  $\widetilde{Z(U_p)}$ contient toute l'information n\'ecessaire pour param\'etrer ces repr\'esentations locales. Je le d\'emontre en utilisant de mani\`ere essentielle les propri\'et\'es "triangulines"  (y compris au sens infinit\'esimal) des repr\'esentations galoisiennes port\'ees par $\mathcal{C}$, sur lesquelles je reviendrai un peu plus loin.

\subsection{G\'en\'eralisations en dimension sup\'erieure}\label{fougeredimsup} Il me semble peu probable qu'en dimension $d\geq 3$, $X_d^{\rm geom}$ soit Zariski-dense dans des composantes connexes de $X_d$ sur lesquelles la trace des conjugaisons complexes vaut $\pm 1$ (qui sont de dimension $\geq 5$).  Il serait int\'eressant de le d\'emontrer. Un indicateur de cela est la difficult\'e \`a construire des repr\'esentations automorphes alg\'ebriques $\pi$ de $\GL_d$ sur $\Q$ quand $d \geq 3$ car leurs composantes archim\'ediennes $\pi_\infty$ ne sont pas discr\`etes restreintes \`a ${\rm SL}_d(\R)$.\footnote{Ash et Pollack ont m\^eme avanc\'e dans~\cite{ashpollack} l'hypoth\`ese qu'en conducteur $1$ et pour $d=3$, toutes ces repr\'esentations seraient \`a torsion pr\`es des carr\'e sym\'etriques de formes modulaires pour ${\rm SL}(2,\Z)$.} En revanche, il est bien connu qu'il existe pl\'ethore de repr\'esentations $\pi$ comme plus haut qui sont polaris\'ees au sens du~\S\ref{reppol} : je renvoie par exemple au~\S\ref{dehnar} Ch. 2 o\`u ce probl\`eme sera d'ailleurs \'etudi\'e en d\'etail. Cela vient au fond de ce que les composantes archim\'ediennes des repr\'esentations autoduales alg\'ebriques r\'eguli\`eres sont discr\`etes dans le
 dual unitaire autodual de ${\rm SL}_d(\R)$. \ps

Ceci sugg\`ere de se restreindre en dimension $d\geq 3$ au ferm\'e de $X_d$ constitu\'e des $x$ param\'etrant les repr\'esentations polaris\'ees, i.e. tels que $\rho_x^\ast \simeq \rho_x \otimes \chi_x$ pour un certain caract\`ere $\chi_x$. Pour des raisons techniques ce n'est pas l'approche que j'ai suivie. En effet, le cadre automorphe naturel pour travailler avec ce type de repr\'esentations est celui des groupes classiques ${\rm SO}$ et ${\rm Sp}$, mais on ne disposait pas des travaux r\'ecents d'Arthur au moment o\`u j'ai commenc\'e \`a \'etudier ces questions (ces travaux sont n\'ecessaires pour relier repr\'esentations automorphes de ces groupes classiques et des $\GL_m$, $m\geq 1$, et donc par exemple pour leur associer des repr\'esentations galoisiennes). 
J'ai plut\^ot \'etudi\'e la contribution \`a $X_d$ des repr\'esentations galoisiennes associ\'ees aux repr\'esentations automorphes cuspidales alg\'ebriques r\'eguli\`eres et polaris\'ees dans le cas CM (cf. \S\ref{reppol}). Une autre bonne raison de se focaliser sur ce cas est que les d\'emonstrations dont j'ai parl\'ees au~\S\ref{repgalaut} utilisent aussi de mani\`ere essentielle les propri\'et\'es de la foug\`ere infinie dans ce contexte. \ps

Fixons un corps de nombres $F$ qui est CM, extension quadratique totalement imaginaire du corps totalement r\'eel $F^+$, et fixons $c \in {\rm Gal}(\overline{F}/F^+)$ une conjugaison complexe. La notation $S$ d\'esignera d\'esormais un ensemble fini de places finies de $F^+$ contenant les places au dessus de $p$ et les places ramifi\'ees dans $F$. Soit $G_{F,S}$ le groupe de Galois d'une extension alg\'ebrique maximale de $F$ non ramifi\'ee hors de $S$. Consid\'erons  
la vari\'et\'e des caract\`eres $p$-adique $Y_d$ de $G_{F,S}$ en dimension $d$. Nous allons nous int\'eresser principalement \`a des repr\'esentations $\rho : G_{F,S} \rightarrow \GL_d(A)$, $A$ \'etant disons une $\Z_p$-alg\`ebre, telles que $\rho^\ast \simeq \rho^c \otimes \omega^{d-1}$, o\`u $\omega$ d\'esigne le caract\`ere cyclotomique $p$-adique. Une telle repr\'esentation sera dite {\it de type unitaire}.\footnote{On peut justifier cette terminologie en reliant les repr\'esentations de type unitaire de $G_{F,S}$ avec les morphismes de $G_{F^+,S}$ dans un groupe alg\'ebrique non connexe qui appara\^it comme le dual de Langlands alg\'ebrique  du groupe unitaire quasi-d\'eploy\'e, \`a $d$ variables, associ\'e \`a $F/F^+$ : voir~\cite{cht} et~\cite{buzzardgee}.} De m\^eme, un pseudo-caract\`ere $G_{F,S} \rightarrow A$ sera dit de type unitaire si $T(g^{-1})=T(cgc^{-1})\omega(g)^{d-1}$ pour tout $g \in G_{F,S}$. On consid\`ere $$Y_d^\bot \subset Y_d$$
le ferm\'e de $Y_d$ param\'etrant les pseudo-caract\`eres de type unitaire.  En particulier, un point $x \in Y_d$ est dans $Y_d^\bot$ si et seulement si $\rho_x$ est de type unitaire.  Bien entendu, l'induction fournissant un morphisme naturel de $Y_d \rightarrow X_{d[F:\Q]}$ (pour des ensembles $S$ convenables des deux c\^ot\'es),  \'etudier $Y_d$ c'est aussi un peu \'etudier $X_d$. Nous allons nous focaliser sur $Y_d^\bot$.  \ps

\begin{definition}\label{defptaut} On dira qu'un point $x \in Y_d^\bot$ est automorphe s'il existe une repr\'esentation automorphe cuspidale alg\'ebrique r\'eguli\`ere polaris\'ee $\pi$ de $\GL_d$ sur $F$, et un plongement $\iota : \overline{\Q} \rightarrow \overline{k(x)}$, tels que : \begin{itemize}\ps
\item[(i)] $\pi_v$ est non ramifi\'ee si $v|p$ ou si $v$ n'est pas au-dessus de $S$, \ps
\item[(ii)] $\rho_x  \simeq \rho_{\pi,\iota}$. \ps
\end{itemize}
\end{definition}

 La conjecture optimiste suivante me semble n\'eanmoins raisonnable. Une
repr\'esentation semisimple $r : G_{F,S} \rightarrow
\GL_d(\overline{\F}_p)$ sera dite automorphe si elle est de la forme
$\overline{\rho_x}$ avec $x \in Y_d^\bot$ automorphe.  On pose
$Y_d^\bot(r):=Y_d(r) \cap Y_d^\bot$.  \ps

\begin{conjecture} Les points automorphes sont Zariski-denses et d'accumulation dans $Y_d^\bot(r)$ pour toute repr\'esentation r\'esiduelle automorphe $r$. \end{conjecture}

On pourrait sans doute supprimer la donn\'ee du $r$ dans l'\'enonc\'e ci-dessus, via une g\'en\'eralisation convenable de la conjecture de Serre, \`a condition de se restreindre aux composantes irr\'eductibles de $Y_d^\bot$ param\'etrant des repr\'esentations satisfaisant une condition de "signe" convenable dans l'esprit du Th\'eor\`eme~\ref{signecm}. \ps

Le th\'eor\`eme suivant, r\'esultat principal de mon article~\cite{chens}, est un premier pas en direction de la conjecture ci-dessus. Cet article contient aussi un \'enonc\'e similaire concernant les formes modulaires de Hilbert, que je ne d\'etaillerai pas ici. 

\begin{thm}\label{densite3} On suppose $d\leq 3$ et que $F$ est totalement d\'ecompos\'e en $p$. Alors l'adh\'erence Zariski des points automorphes dans $Y_d^\bot$ a toutes ses composantes 
irr\'eductibles de dimension $\geq [F^+:\Q]\frac{d(d+1)}{2}$.
\end{thm}

Je rappelle que pour $d=3$, les $\rho_x$ avec $x \in Y_3^\bot$ automorphe appara\^issent dans la cohomologie \'etale $p$-adique de puissances de la vari\'et\'e ab\'elienne universelle sur les surfaces de Picard associ\'ees \`a $F$ (ils sont m\^eme g\'eom\'etriques par construction si $[F^+:\Q]>1$), voir le~\S\ref{repgalaut} et~\cite{picard}.\ps

Mon but dans les paragraphes qui suivent est de donner des indications sur la d\'emonstration de ce th\'eor\`eme. Je veux mentionner auparavant que $[F^+:\Q]\frac{d(d+1)}{2}$ est \'egalement la dimension probable des $Y_d^\bot(r)$ avec $r$ automorphe. Les premiers indicateurs de ce fait remontent \`a~\cite{cht}. L'ingr\'edient crucial pour le comprendre est le Th\'eor\`eme~\ref{signecm}. Partons pour simplifier d'une repr\'esentation irr\'eductible continue $r : G_{F,S} \rightarrow \GL_d(\overline{\F}_p)$ de type unitaire. De m\^eme que la th\'eorie des d\'eformations usuelle de $r$ est gouvern\'ee par la repr\'esentation ${\rm ad}(r)$~\cite{mazurdef}, celle des d\'eformations de $r$ qui sont de type unitaire  est gouvern\'ee par une certaine repr\'esentation que je noterai ici $${\rm ad}^+ (r)$$
du groupe de Galois $G_{F^+,S}$ (le groupe de Galois d'une extension alg\'ebrique maximale de $F^+$ non ramifi\'ee hors de $S$).  \ps

Elle a les propri\'et\'es suivantes. D'une part, si $V$ d\'esigne l'espace de $r$ alors l'espace sous-jacent \`a ${\rm ad}^+(r)$ est ${\rm End}_{\overline{\F}_p}(V)$, et la restriction \`a $G_{F,S}$ 
de ${\rm ad}^+(r)$ est simplement ${\rm ad}(r)$. D'autre part, si $\langle,\rangle$ est un $G_{F,S}$-accouplement non-d\'eg\'en\'er\'e $V \otimes V^c \rightarrow \omega^{1-d}$, alors l'\'el\'ement $c$ agit dans ${\rm ad}^+(r)$ comme l'oppos\'e de l'adjonction par rapport \`a cet accouplement. Si $p\neq 2$, et si cet accouplement est sym\'etrique, ce qui est notamment le cas pour les $r$ automorphes d'apr\`es le Th\'eor\`eme~\ref{signecm}, on constate que 
$$\dim ({\rm End}_{\overline{\F}_p}V)^{{\rm ad}^+(r)(c)=-{\rm id}}=\frac{d(d+1)}{2}$$ (au lieu de $\frac{d(d-1)}{2}$ dans le cas altern\'e).  Si ${\rm H}^2(G_{F^+,S},{\rm ad}^+(r))=0$, le foncteur des d\'eformations continues de type unitaire de $r$ est formellement lisse~\footnote{...sur l'anneau des vecteurs de Witt du corps des coefficients de $r$,} et d'espace tangent de dimension\footnote{Cela r\'esulte de la formule de la caract\'eristique d'Euler de Tate~\cite{milneadt}. On observera que les $G_{F,S}$-invariants de ${\rm ad}^+(r)$ sont nuls, car ils sont inclus dans les homoth\'eties par irr\'eductibilit\'e de $r$, sur lesquels $c$ agit par $-{\rm id}$.}  $[F^+:\Q]\dim ({\rm End}_{\overline{\F}_p}V)^{{\rm ad}^+(r)(c)=-{\rm id}}$, et en particulier $Y_d(r)$ est une boule ouverte de cette dimension. 

\begin{cor} Supposons $p>2$ et que $r : G_{F,S} \rightarrow \GL_3(\overline{\F}_p)$ est irr\'eductible, automorphe, telle que $H^2(G_{F^+,S},{\rm ad}^+(r))=0$. Alors les points automorphes sont Zariski-denses et d'accumulation dans $Y_3^\bot(r)$, qui est une boule ouverte de dimension $6[F^+:\Q]$.
\end{cor}

Voici un exemple concret d'application pour lequel je renvoie \`a l'appendice de~\cite{chens}. Soit $E$ une courbe elliptique sur $\Q$ de conducteur $N$, soit $F$ un corps quadratique imaginaire, et soit $S$ l'ensemble des nombres premiers divisant $2pN$. La repr\'esentation $r=({\rm Sym}^2 E[p])^\ast_{|G_{F,S}}$ est bien de type unitaire, elle est m\^eme automorphe d'apr\`es Wiles et ses continuateurs~\cite{bcdt}, Gelbart-Jacquet~\cite{GJ}, et Arthur-Clozel~\cite{arthurclozel}. C'est un exercice de v\'erifier que $${\rm ad}^+(r)=\varepsilon_{F/\Q} \oplus r \otimes \omega \oplus ({\rm Sym}^4 E[p])^\ast \otimes \varepsilon_{F/\Q}\omega^2$$ o\`u $\varepsilon_{F/\Q}$ est le caract\`ere d'ordre $2$ de $G_{\Q,S}$ d\'efini par $F/\Q$. J'ai pu v\'erifi\'e {\it loc. cit.} que si $p=5$ et $F=\Q(i)$, et si la classe d'isog\'enie de $E$ est du type
$$17A, \, \, 21A, \, \, 37B, \, \, 39A, \, \, 51A, \, \, 53A, \, \, 69A, \, \, 73A, \, \, 83A, \, \, 91B$$ 
dans les tables de Cremona,  alors $r$ est irr\'eductible et satisfait $H^2(G_{F^+,S},{\rm ad}^+(r))=0$. C'est par exemple le cas de la courbe $y^2+xy+y = x^3-x^2-x$ de conducteur $17$.  Ces calculs d\'ependent de calculs de nombre de classes faits par Pari~\cite{GP} qui sont conditionnels \`a GRH.

\subsection{La foug\`ere infinie de type unitaire}\label{fougereunitaire} Le
premier ingr\'edient crucial dans la d\'emonstration du
Th\'eor\`eme~\ref{densite3} est la pr\'esence d'une g\'en\'eralisation de la
th\'eorie de Coleman, ainsi que de la foug\`ere infinie de Gouv\^ea-Mazur,
dans le contexte des repr\'esentations galoisiennes de type unitaire.  Pour
l'essentiel, et quand $F^+=\Q$, j'avais construit de telles
g\'en\'eralisations dans mon article~\cite{chcrelle} issu de ma th\`ese. Le
lecteur ne perdrait d'ailleurs par grand chose en premi\`ere approche \`a supposer $F=\Q^+$ et
$d=3$ dans ce qui suit.  Au prix d'une perte d'un peu de clart\'e, je me sens
n\'eanmoins oblig\'e de me placer dans le cadre d'un corps totalement r\'eel
$F^+$ g\'en\'eral, qui est celui dont on a besoin pour d\'emontrer les
th\'eor\`emes \'enonc\'es dans le~\S\ref{repgalaut}, et pour lequel je
renvoie \`a mon article~\cite{chgrfa}.\ps

On conserve les notations du~\S\ref{fougeredimsup}. On d\'esignera de plus
par $S_p$ l'ensemble des places de $F^+$ divisant $p$. Mes hypoth\`eses de travail seront d\'esormais les suivantes
: \begin{itemize}\ps
\item[(i)] Si $d$ est pair, alors $d[F^+:\Q] \equiv 0 \bmod 4$.\ps
\item[(ii)] Si $d>3$, alors $F$ est non ramifi\'e au dessus de toutes les
places finies de $F^+$, et $S$ ne contient que des places d\'ecompos\'ees dans $F$. \ps
\item[(iii)] Chaque $v \in S_p$ est d\'ecompos\'ee dans $F$.\ps
\end{itemize}

Expliquons ces hypoth\`eses. Sous la condition (i), le principe de Hasse
pour les formes hermitiennes assure qu'il existe une forme hermitienne sur
$F^d$ relativement \`a $F/F^+$, telle que l'espace hermitien $F^d
\otimes_{F^+} F_v^+$ soit d'indice maximal $[d/2]$ pour toute place finie $v$ de $F^+$
et d\'efini positif pour toute place r\'eelle.  On d\'esignera alors par
$$U(d)$$ le groupe unitaire associ\'e, qui est un groupe alg\'ebrique sur
$F^+$.  \ps

La condition (ii) est la condition sous laquelle il est
actuellement connu que les repr\'esentations automorphes de $U(d)$ se
transf\`erent par "changement de base" \`a $U(d) \times_{F^+} F \simeq
\GL(d)$ selon les recettes de Langlands et Arthur.  Si $d\leq 3$, c'est le th\'eor\`eme
principal de Rogawski~\cite{roglivre} (voir aussi~\cite{picard}), et sous
(ii) c'est d\^u \`a Labesse~\cite{labesse}.  Il est raisonnable d'esp\'erer
que dans un futur proche l'hypoth\`ese (ii) ne soit plus n\'ecessaire,
\'etant donn\'es notamment les r\'esultats r\'ecents de
Moeglin~\cite{moeglin}, et ceux de Mok~\cite{mok}, sur la classification
d'Arthur pour les groupe unitaires. Il faut rajouter que du point de vue des astuces "\`a la Blasius-Rogawski" permettant de raisonner apr\`es changement de base, l'hypoth\`ese (ii) n'est pas aussi farfelue qu'elle n'y para\^it. Enfin, l'hypoth\`ese (iii) est une hypoth\`ese
simplificatrice que nous avons faite dans notre travail mais qui pourrait
sans doute \^etre lev\'ee. \ps

Notre objectif dans les paragraphes qui suivent est d'abord de d\'efinir un analogue
de la courbe de Coleman-Mazur.  Pour que nos \'enonc\'es soient corrects il sera n\'ecessaire d'\'etendre l\'eg\`erement
la notion de point automorphe de $Y_d^\bot$ d\'efinie pr\'ec\'edemment, en une notion relative au groupe unitaire $U(d)$. \ps

 Si $\pi$ est une repr\'esentation automorphe de $U(d)$, d\'esignons par
$\widetilde{\pi}$ la torsion par $|\cdot|^{\frac{1-d}{2}}$ du changement de base
de $\pi$ \`a $\GL(d)$ sur $F$.  Il existe alors une \'ecriture $d=\sum_{i=1}^k d_i$,
o\`u $k$ et les $d_i$ sont des entiers $\geq 1$, ainsi que des
repr\'esentations automorphes discr\`etes alg\'ebriques r\'eguli\`eres
$\Pi_i$ des $\GL(d_i)$ sur $F$ telles que $\Pi_i^\vee \simeq \Pi_i^c
|\cdot|^{\frac{d-1}{2}}$, tels que $\widetilde{\pi}$ soit l'induite
parabolique normalis\'ee de $\Pi_1 \times \cdots \times \Pi_k$. Les r\'esultats\footnote{Lorsque $d$ et $d_i$ sont
de m\^eme parit\'e, alors $\Pi_i |\cdot|^{\frac{d-d_i}{2}}$ est polaris\'ee. 
Sinon, c'est $\Pi_i |\cdot|^{\frac{d-d_i-1}{2}}\mu$ qui l'est, $\mu$ \'etant
n'importe quel caract\`ere de Hecke alg\'ebrique de $F$ tel que
$\mu^{-1}=\mu^c|\cdot|$.}
du~\S\ref{repgalaut}, et la recette de Moeglin et
Walspurger~\cite{moeglinwaldspurger}, permettent d'associer \`a ces $\Pi_i$ des repr\'esentations galoisiennes
$\rho_{\Pi_i,\iota}$ satisfaisant (C1), (C2) et (C3) de la mani\`ere \'evidente, et donc de poser
$$\rho_{\pi,\iota}=\oplus_{i=1}^k \rho_{\Pi_i,\iota}.$$ 
Cette repr\'esentation de $G_{F,S}$ d\'etermine d'ailleurs en retour l'entier $k$ et les  $(\Pi_i,d_i)$ \`a permutation pr\`es~\cite{jasha}.

\begin{definition}\label{defptUaut} Un point $x \in Y_d^\bot$ sera dit $U(d)$-automorphe s'il existe une repr\'esentation automorphe $\pi$ de $U(d)$, non ramifi\'ee
en tout $v \notin S$ et tout $v|p$, ainsi qu'un plongement $\iota :
\overline{\Q} \rightarrow \overline{k(x)}$, tels que $\rho_x \simeq
\rho_{\pi,\iota}$.
\end{definition}

C'est un fait essentiel que tous les points automorphes de $Y_d^\bot$ au sens de la D\'efinition~\ref{defptaut} sont
$U(d)$-automorphes (\cite{roglivre},\cite{labesse}).  La caract\'erisation
des points $U(d)$-automorphes g\'en\'eraux en terme des $\Pi_i$ associ\'ees
est en revanche assez subtile : elle est d\'ecrite par la formule de
multiplicit\'e d'Arthur, qui \`a l'heure actuelle est encore non
d\'emontr\'ee dans ce contexte d\`es que $d>3$.  Nous n'en aurons cependant
pas besoin.  Il nous reste \`a d\'efinir la notion de point
$U(d)$-automorphe {\it raffin\'e}.  Cette notion de raffinement sera
l'analogue du choix d'une racine $\alpha$ de $P_f$ dans la construction des
arcs de Coleman passant par $x_f$ rappel\'ee au~\S\ref{constcoleman}.  \ps

Soit $x$ un point $U(d)$-automorphe de $Y_d^\bot$. La repr\'esentation
galoisienne $\rho_x$ est cristalline en toutes les places $v$ divisant $p$
par la propri\'et\'e~(C3).  Si $v$ est une telle place, regardons le
$\varphi$-module filtr\'e ${\rm D}_{\rm cris}((\rho_x)_v)$ d\'efini par
Fontaine.  Il admet un op\'erateur Frobenius $\varphi$, et nous
d\'esignerons par $P_{x,v} \in \overline{k(x)}[T]$ le polyn\^ome
caract\'eristique de $\varphi_{x,v}=\varphi^{f_v}$, o\`u $p^{f_v}$ est le
cardinal du corps r\'esiduel de $F_v$.  Inspir\'es par Mazur~\cite{Maz},
nous appellerons {\it raffinement de $x$ en $v$} la donn\'ee d'une
\'enum\'eration\footnote{Lorsque le point $x$ param\`etre une
repr\'esentation automorphe de $U(d)$ dont les $\Pi_i$ associ\'es par la
recette de changement de base rappel\'ee plus haut ne sont pas cuspidaux,
nous faisons en fait une hypoth\`ese suppl\'ementaire sur $\Ref_v$ que nous
mettons ici sous le tapis : je renvoie \`a~\cite[\S 6.4]{bchlivre} pour la
d\'efinition correcte dans ce cas (cela correspond \`a la notion de
raffinement {\it accessible} au sens de {\it loc.  cit.}).  Cette mise en
garde ne concerne ni les points automorphes, ni les poids $U(d)$-automorphes
param\'etrant les $\pi$ dont le caract\`ere infinit\'esimal est assez
r\'egulier. Elle jouera par contre un r\^ole au Ch. 2~\S\ref{thmlivre}.} $$\Ref_v=(\phi_{1,v},\phi_{2,v},\cdots,\phi_{d,v}) \in
\overline{k(x)}^d$$ des racines de $P_{x,v}$ dans $\overline{k(x)}$
(compt\'ees avec multiplicit\'es).  \ps

D'apr\`es la propri\'et\'e (C3) de $(\rho_x)_v$, il est aussi \'equivalent
de choisir un ordre sur les valeurs propres du param\`etre de Satake de
$\pi_u$ si $\pi$ est la repr\'esentation automorphe de $U(d)$ associ\'ee \`a
$x$ et si $u$ est la place de $F_+$ au dessous de $v$.  Dans le cas
"g\'en\'erique" o\`u les racines de $P_{x,v}$ sont distinctes, alors $x$
admet exactement $d!$ raffinements en $v$.  \ps

Il sera commode de choisir\footnote{Ce choix n'aura qu'une influence b\'enigne sur les
consid\'erations qui suivent, mais nous \'evitera des p\'eriphrases
p\'enibles.} une fois pour toutes, pour chaque $v \in S_p$, l'une
des deux places de $F$ au dessus de $v$, et de d\'esigner par
$\widetilde{S}_p$ l'ensemble de ces places choisies de $F$.

\begin{definition} Un point $U(d)$-automorphe raffin\'e est un
couple $(x,\Ref)$ o\`u $x \in Y_d^\bot$ est un point $U(d)$-automorphe, et 
$\Ref=\{\Ref_v, v \in \widetilde{S}_p\}$ est une collection de raffinements
$\Ref_v$ de $x$ en chaque $v \in \widetilde{S}_p$.
\end{definition}
\newcommand{\uk}{\underline{k}}

Nous allons pour finir associer \`a chaque point $U(d)$-automorphe raffin\'e
$(x,\Ref)$ un point canonique \begin{equation} \label{defxref}
x_\Ref=(x,\delta) \in Y_d^\bot \times \mathcal{T}^d,\end{equation} o\`u
$\mathcal{T}$ d\'esigne l'espace analytique $p$-adique param\'etrant les
caract\`eres continus du groupe multiplicatif\footnote{Cela signifie que se
donner un morphisme $X \rightarrow \mathcal{T}$, $X$ \'etant un
$\Q_p$-affino\"ide, \'equivaut \`a se donner un homomorphisme continu de
groupes $F_p^\ast \rightarrow \OO(X)^\ast$.  } de $(F^+_p)^\ast$, o\`u
$F_p^+=F^+ \otimes_\Q \Q_p$.  L'espace $\mathcal{T}$ est tr\`es concret :
le choix d'une uniformisante de $F_v^+$ pour tout $v \in S_p$ l'identifie au
produit direct $$\mathcal{T} \simeq \mathbb{G}_m^{S_p} \times
\mathcal{T}_0,$$ o\`u $\mathcal{T}_0$ d\'esigne la vari\'et\'e des
caract\`eres $p$-adique en dimension $1$ du groupe compact $\prod_{v \in
S_p} \OO_{F^+_v}^\ast$ (une r\'eunion disjointe finie de boules ouvertes de
dimension $[F^+:\Q]$).  Quand $F^+=\Q$ on a simplement $\mathcal{T}_0=\WW$
et $\mathcal{T}\simeq \mathbb{G}_m \times \WW$.  \ps

Tout \'el\'ement $\uk=(k_{i,\sigma}) \in (\Z^d)^{{\rm
Hom}(F^+,\overline{\Q}_p)}$ pourra \^etre vu comme un \'el\'ement de
$\mathcal{T}(\overline{\Q}_p)$ en consid\'erant le caract\`ere
"alg\'ebrique" $$(x_{i,v}) \mapsto \prod_{i,\sigma}
\sigma(x_{i,v_\sigma})^{k_{i,\sigma}},$$ o\`u $\sigma$ parcourt ${\rm
Hom}(F^+,\overline{\Q}_p)$ et $v_\sigma \in S_p$ d\'esigne la place induite
par $\sigma$.\ps \ps

Il reste \`a d\'efinir  le caract\`ere $\delta$ dans la
formule~\eqref{defxref}.  C'est le produit $\delta=\delta_{nr}\delta_w$ de
deux caract\`eres, o\`u $\delta_{nr}$ ne d\'epend que de $\Ref$, et o\`u
$\delta_w$ ne d\'epend que des poids de Hodge-Tate des $(\rho_x)_v$ avec $v
\in \widetilde{S}_p$.  Ces derniers, naturellement ind\'ex\'es par les
plongements $\sigma : F^+ \rightarrow \overline{k(x)}$, sont distincts par
d\'efinition \`a $\sigma$ fix\'e, et seront num\'erot\'es dans l'ordre
croissant : $$k_{1,\sigma} < k_{2,\sigma} < \cdots < k_{d,\sigma}.$$ On pose
$\uk=(k_{i,\sigma})$ et $\delta_w=\uk^{-1}$.  Pour d\'ecrire enfin
$\delta_{nr}$, on fixe une place $v \in S_p$, ainsi qu'un entier $i \in
\{1,\cdots,d\}$, et il suffit de donner sa restriction au facteur
$(F_v^+)^\ast$ mis \`a la place $i$ (vu comme sous-groupe avec des $1$ aux
autres coordonn\'ees).  Soit $u \in \widetilde{S}_p$ l'unique place divisant
$v$.  Alors $\delta_{nr}$ est trivial sur $\OO_{F_v^+}^\ast$ et envoit une
uniformisante de $F_v^+$ sur le $i$-\`eme \'el\'ement $\phi_{i,u}$ de
$\Ref_u$.  \ps

\begin{definition} La vari\'et\'e de Hecke de $U(d)$ est l'adh\'erence
Zariski $$\mathcal{E}_d \subset Y_d^\bot \times \mathcal{T}^d$$ de tous les
points $x_\Ref$, o\`u $(x,\Ref)$ parcourt les points $U(d)$-automorphes
raffin\'es de $Y_d^\bot$.  \end{definition}

C'est la g\'en\'eralisation naturelle de la courbe de Coleman-Mazur. La
terminologie anglaise d'usage pour $\mathcal{E}_d$ est ${\it eigenvariety}$. 
La traduction fran\c{c}aise litt\'erale \'etant ambig\"ue nous lui avions
pr\'ef\'er\'e ce terme dans~\cite{chcrelle} car elle appara\^it aussi de
mani\`ere naturelle comme un spectre d'une alg\`ebre de Hecke.  Observons
que la seconde projection, ainsi que la projection naturelle $\mathcal{T}
\rightarrow \mathcal{T}_0$, d\'efinissent des morphismes naturels $\nu :
\mathcal{E}_d \rightarrow \mathcal{T}^d$ et $\kappa : \mathcal{E}_d
\rightarrow \mathcal{T}_0^d$ ("morphisme poids").  \ps

\begin{thm}\label{vhecke} \begin{itemize}\item[(i)] $\mathcal{E}_d$ est d'\'equidimension $\dim(\mathcal{T}_0^d)=d[F^+:\Q]$. \ps
\item[(ii)] Le morphisme $\nu$ est fini et le morphisme $\kappa$ est localement fini et ouvert. \ps
\item[(iii)] Les points automorphes raffin\'es $x_\Ref$ sont Zariski-denses et d'accumulation dans~$\cE_d$.\end{itemize}
\end{thm}

Ce th\'eor\`eme est le r\'esultat principal de mes articles~\cite{chcrelle}
et \cite{chgrfa}, \'etendant notamment des travaux de Hida~\cite{hidasln}
(voir aussi~\cite{ashstevens},~\cite{ashstevens2},
\cite{hidabook},\cite{emerton}, \cite{buzzardeigen}, \cite{loeffler},
\cite{urban}).  Je renvoie \`a~\cite[\S 7]{bchlivre}, ainsi qu'\`a mon cours
Peccot~\cite{cpeccot}, pour une discussion
d\'etaill\'ee du cas particulier $F^+=\Q$.  La propri\'et\'e pr\'ecise de $\kappa$ sous-entendue dans le (ii) ci-dessus est la suivante : $\mathcal{E}_d$ est admissiblement recouvert par des ouverts affino\"ides $\Omega$ tels que $\kappa(\Omega)$ est ouvert et  $\kappa_{|\Omega} : \Omega \rightarrow \kappa(\Omega)$ est fini. \ps

Il ne me semble pas n\'ecessaire de revenir sur la d\'emonstration de ce
th\'eor\`eme dans ce m\'emoire, car il est plut\^ot issu de mon travail de
th\`ese et rallongerait un peu trop cette discussion d\'ej\`a longuette. 
Disons simplement que nous d\'eveloppons une th\'eorie ad hoc de formes
automorphes $p$-adiques pour le groupe unitaire $U(d)$ ayant des
propri\'et\'es formelles similaires aux constructions de Coleman
rappell\'ees au~\S\ref{constcoleman}.\footnote{D\'efinition d'une notion de
forme automorphe $p$-adique englobant la notion classique, construction d'un
certain module sur l'espace des poids $\mathcal{T}_0$ interpolant les
espaces de telles formes, avec action des op\'erateurs de Hecke, application
de la th\'eorie spectrale de Coleman en famille \`a un certain op\'erateur
compact, etc...} La compacit\'e du groupe $U(d)(F^+_v)$ pour toute place $v$
archim\'edienne joue un r\^ole simplificateur crucial dans mes
constructions.  En effet, toutes les repr\'esentations automorphes de $U(d)$
ayant alors de la cohomologie en degr\'e $0$, on est simplement ramen\'e \`a
interpoler $p$-adiquement les repr\'esentations alg\'ebriques
irr\'eductibles du $\Q$-groupe $G={\rm Res}^{F^+}_\Q U(d)$, ce que l'on fait
\`a l'aide de la s\'erie principale localement analytique des
repr\'esentations d'un sous-groupe d'Iwahori de $G(\Q_p)$.  \ps La m\'ethode
de Coleman-Mazur, \'etendue par Buzzard~\cite{buzzardeigen}, appliqu\'ee aux
constructions \'evoqu\'ees ci-dessus, fournit d'abord un analoque purement
automorphe de leur construction $\mathcal{D}$~(\S\ref{eigencurvefern}) dans
ce contexte qui a les propri\'et\'es (i) et (iii) du Th\'eor\`eme, et qui
vient par construction avec des morphismes $\nu$ et $\kappa$ satisfaisant
(ii) \'egalement.  Le (iii) suit d'un crit\`ere beaucoup plus pr\'ecis
assurant qu'un point $x \in \mathcal{D}$ est $U(d)$-automorphe d\`es que $\nu(x)$
satisfait un crit\`ere explicite \`a la
Coleman~\cite{coleman1}~\cite{chgrfa}.  L'existence de repr\'esentations
galoisiennes attach\'ees aux repr\'esentations automorphes de $U(d)$
par la recette de changement de base d\'ecrite plus haut, qui repose sur
le~\S\ref{repgalaut}, et un argument de pseudo-caract\`eres \`a la
Wiles~\cite{wiles},
permet alors de d\'efinir un morphisme canonique $\eta : \mathcal{D}
\rightarrow Y_d^\bot$.  Le morphisme qui s'en d\'eduit $$\eta \times \nu :
\mathcal{D} \longrightarrow Y_d^\bot \times \mathcal{T}^d$$ est par
construction une immersion ferm\'ee dont l'image contient tous les points
$U(d)$-automorphes raffin\'es, cette image est donc simplement $\mathcal{E}_d$, d'o\`u le th\'eor\`eme.  Je renvoie
\`a~\cite[\S 2.1]{chens} pour le d\'etail de cette partie de l'argument. 
\ps

Je veux rajouter ici que pour pouvoir d\'efinir le morphisme $\eta$, il
suffit de savoir associer des repr\'esentations galoisiennes \`a un
sous-ensemble Zariski-dense dans $\mathcal{D}$ de points $U(d)$-automorphes, par
exemple ceux param\'etrant les repr\'esentations automorphes ayant un
caract\`ere infinit\'esimal assez r\'egulier.  On obtient alors gratuitement
des repr\'esentations galoisiennes associ\'ees aux autres points
automorphes\footnote{En fait, l'espace $\mathcal{D}$ contient aussi des
points param\'etrant des repr\'esentations automorphes de $U(d)$ qui sont
ramifi\'ees aux places divisant $p$, notamment par construction toutes celles qui
ont des invariants par un sous-groupe d'Iwahori en toutes ces places.  C'est
un fait important dans la d\'emonstration du
Th\'eor\`eme~\ref{existencegalois} car \'etant donn\'ee une repr\'esentation
automorphe cuspidale $\pi$ de $\GL_d$ sur $F$, on peut toujours trouver une
extension $CM$ r\'esoluble $F'/F$ telle que le changement de base de $\pi$
\`a $F'$ ait des invariants par un sous-groupe d'Iwahori en chaque place
("potentielle semi-stabilit\'e").  } de $\mathcal{D}$, et m\^eme plus
g\'en\'eralement \`a tous ses points, qui param\`etrent par d\'efinition les
formes automorphes $p$-adiques pour $U(d)$ qui sont propres et de "pentes
finies" : c'est un fait simple mais essentiel dans la d\'emonstration de ma
contribution au Th\'eor\`eme~\ref{existencegalois}.  \ps

\`A ce stade la d\'efinition de la foug\`ere infinie unitaire compl\`ete est
similaire \`a celle du~\S\ref{eigencurvefern}.  La projection $Y_d^\bot
\times \mathcal{T}^d \rightarrow Y_d^\bot$ induit un morphisme naturel $\mu
: \mathcal{E}_d \rightarrow Y_d^\bot$.

\begin{definition} La foug\`ere infinie unitaire est l'ensemble
$\mathcal{F}_d=\mu(\mathcal{E}_d) \subset Y_d^\bot$.  \end{definition}

	Cette foug\`ere pr\'esente une structure fractale dans l'esprit de
celle de Gouv\^ea et Mazur.  En effet, soit $x \in Y_d^\bot$ un point
automorphe, ou m\^eme plus g\'en\'eralement $U(d)$-automorphe.  Le choix
d'un raffinement $\Phi$ de $x$ d\'efinit un point $x_{\Ref} \in
\mathcal{E}_d$.  Soit $\Omega$ un voisinage ouvert affino\"ide de $x_\Ref$
dans $\mathcal{E}_d$.  Quitte \`a r\'etr\'ecir $\Omega$, la popri\'et\'e
(ii) assure que $$B_{x,\Ref}:=\mu(\Omega) \subset Y_d^\bot$$ est un ferm\'e
d'un voisinage ouvert affino\"ide de $x$ dans $Y_d^\bot$.  C'est l'analogue
des arcs de Coleman.  Le germe en $x$ de $B_{x,\Ref}$ est canonique.  Nous
dirons en g\'en\'eral que c'est la branche de la foug\`ere infinie
$\mathcal{F}_d$ en $x$ associ\'ee au raffinement $\Ref$.  \ps

\begin{figure}[htp]
\centering
\includegraphics[scale=0.75]{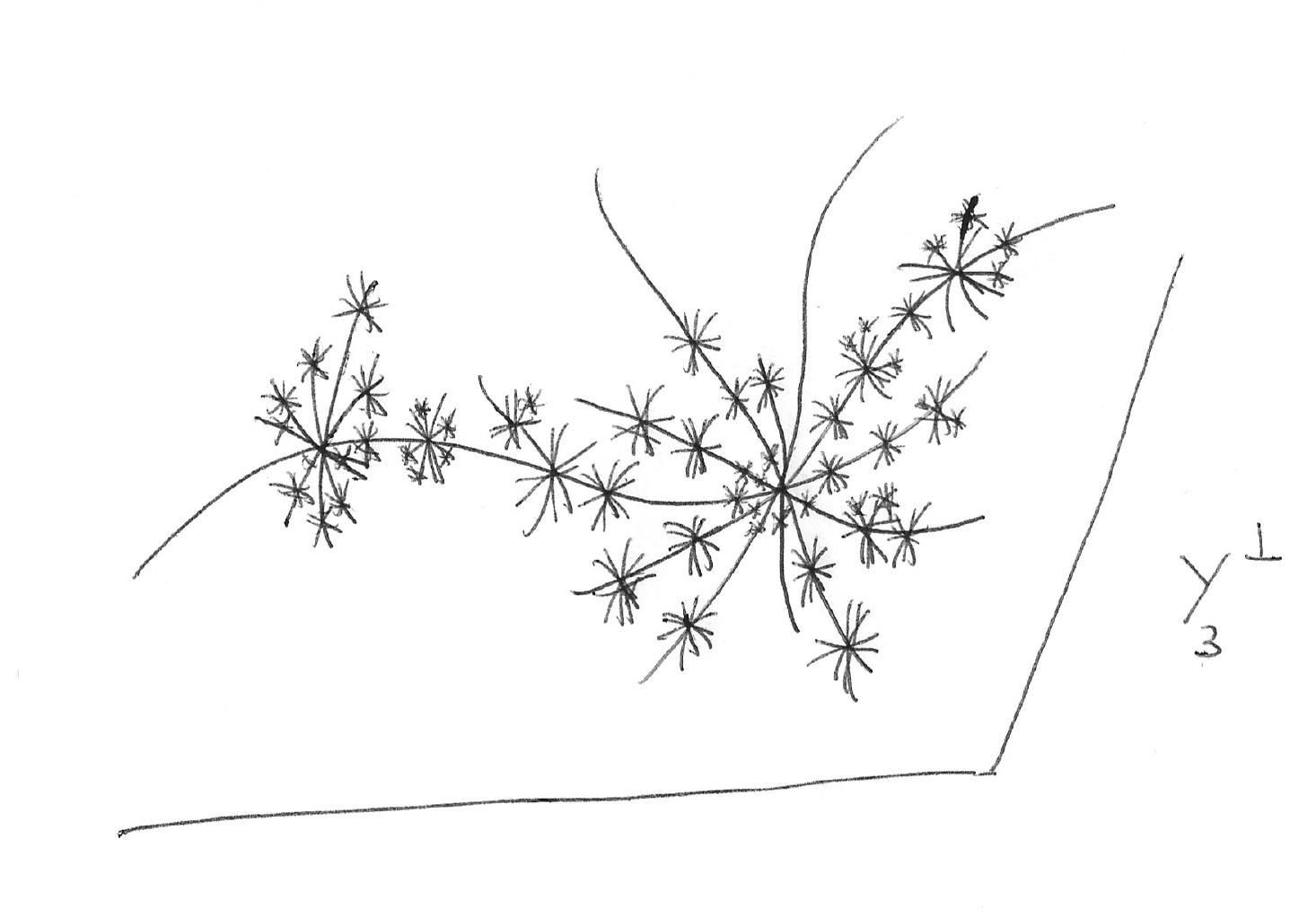}
\caption{La foug\`ere infinie unitaire dans $Y_3^\bot$}
\label{figfougere2}
\end{figure}

Comme $x$ poss\`ede jusqu`\`a $(d!)^{|S_d|}$ raffinements, il peut y avoir
tout autant de branches de la foug\`ere passant par $x$.  Mais les points
$U(d)$-automorphes s'accumulent en $x$ dans $B_{x,\Ref}$ par la
propri\'et\'e (iii).  La structure fractale de $\mathcal{F}_d$ devient alors
tout \`a fait transparente !  La figure~\ref{figfougere2} ci-dessous est une
illustration quand $d=3$ et $F^+=\Q$.  C'est ici le point de d\'epart de mon
article~\cite{chens}.\ps

Ainsi, la foug\`ere infinie $\mathcal{F}_d$ appara\^it comme un "espace de dimension $d[F^+:\Q]$" qui \'evolue dans l'espace $Y_d^\bot$ de dimension conjecturale beaucoup plus grande $\frac{d(d+1)}{2}[F^+:\Q]$. Il poss\`ede un ensemble Zariski-dense et d'accumulation de points hautement multiples, \`a savoir les points $U(d)$-automorphes. La num\'erologie rend assez d\'elicate la g\'en\'eralisation directe de l'argument de Gouv\^ea et Mazur. Mon approche pour minorer la dimension de son adh\'erence Zariski sera de comparer les positions relatives des diff\'erentes branches $B_{x,\Phi}$ en un point automorphe $x$. J'utiliserai pour cela les propri\'et\'es en $v \in \widetilde{S}_p$ des $(\rho_z)_v$ quand $z$ varie dans ces branches.


\subsection{Propri\'et\'es aux places divisant $p$ des repr\'esentations galoisiennes dans la foug\`ere infinie} \label{appraffinee} Commen\c{c}ons par quelques observations assez simples mais essentielles qui d\'ecoulent des propri\'et\'es directes de $\mathcal{F}_d$.  En guise de pr\'etexte, je me propose de les illustrer par le r\'esultat suivant~(\cite[\S 9]{bchens},\cite[\S 7.7]{bchlivre},\cite[\S 3]{bchsigne}), ingr\'edient clef dans notre d\'emonstration du Th\'eor\`eme~\ref{signecm} esquiss\'ee au~\S\ref{signegalois}.\ps

\begin{prop}\label{propirr} Pour tout point $U(d)$-automorphe $x \in Y_d^\bot$, et tout voisinage ouvert affino\"ide $V$ de $x$ dans $Y_d^\bot$, il existe un point automorphe $y \in V$  tel que 
$(\rho_y)_v$ est irr\'eductible pour tout $v \in \widetilde{S}_p$.
\end{prop}

L'id\'ee de la d\'emonstration est de partir du point $x$ de cette
proposition et de se d\'eplacer dans la foug\`ere infinie $\mathcal{F}_d$ en
choisissant bien les branches par lesquelles on passe \`a chaque fois qu'un
choix de raffinement est possible !  \ps

Comme je vais l'expliquer ci-dessous, on dispose d'un certain contr\^ole sur
les polygones de Hodge et de Newton des Frobenius cristallins des
$(\rho_{y_{\Ref'}})_v$ quand $y_{\Ref'}$ est dans un petit voisinage d'un
$x_\Ref \in \mathcal{E}_d$ donn\'e.  Ce contr\^ole va nous permettre
d'appliquer des crit\`eres du type suivant : pour qu'une
$\overline{\Q}_p$-repr\'esentation cristalline $W$ de dimension $d$ de ${\rm
Gal}(\overline{\Q}_p/\Q_p)$ soit irr\'eductible, il suffit que les
polyg\^ones de Newton et de Hodge de ${\rm D}_{\rm cris}(W)$ n'ait pas de
"sous-polygones" respectifs dont les extr\'emit\'es soient communes.  En
effet, si $W' \subset W$ est une sous-repr\'esentation, n\'ecessairement
cristalline, les polyg\^ones de Hodge et Newton de ${\rm D}_{\rm cris}(W')
\subset {\rm D}_{\rm cris}(W)$ fournissent de tels
sous-polygones.\footnote{Cette assertion devrait indiquer au lecteur le sens
que je donne au terme "sous-polygone".}\ps
	
Pour simplifier, je supposerai d\'esormais que $F^+=\Q$, de sorte qu'il n'y
	ait qu'une seule place $v$ dans $\widetilde{S}_p$, qui satisfait de
	plus $F_v=\Q_p$.  L'uniformisante $p \in \Q_p$ fournit une
	identification $\mathcal{T}=\mathbb{G}_m \times \WW$.  D\'esignons
	par $$F_1,\cdots,F_d \in \OO(\mathcal{E}_d)^\ast$$ les restrictions
	\`a $\mathcal{E}_d$ par l'application $\nu$ des coordonn\'ees des
	facteurs $\mathbb{G}_m$ de $\mathcal{T}^d$.  Soit
	$\Ref_v=(\phi_1,\phi_2,\cdots,\phi_d)$ un raffinement de $x$ en $v$,
	et $\Ref=\{\Ref_v\}$.  Si $k_1 < k_2 < \cdots < k_d$ sont les poids
	de Hodge-Tate de $(\rho_{x_\Ref})_v$, on a
	$$(F_1(x_\Ref),F_2(x_\Ref),\cdots,F_d(x_\Ref)) = ( \phi_1\,
	p^{-k_1}, \phi_2 \, p^{-k_2},\cdots,\phi_d\, p^{-k_d})$$
par d\'efinition. Les fonctions $F_i$ \'etant analytiques sur $\mathcal{E}_d$ par construction, on constate que la collection des repr\'esentations cristallines $(\rho_{x_\Ref})_v$, quand $x_\Ref$ varie dans $\mathcal{E}_d$, a la propri\'et\'e que les {\it valeurs propres du Frobenius cristallin $(\rho_{x_\Ref})_v$, prises dans l'ordre de $\Ref_v$, et renormalis\'ees par les poids de Hodge-Tate pris dans l'ordre croissant, varient analytiquement}. C'est un fait absolument remarquable, \'etant donn\'e que typiquement la fonction $k \mapsto p^k$ n'est bien s\^ur pas continue pour la topologie $p$-adique sur $\Z$, et que ces valeurs propres de Frobenius sont g\'en\'eralement divisibles par des grandes puissances de $p$. \ps

\begin{figure}[htp]
\centering
\includegraphics[scale=0.5]{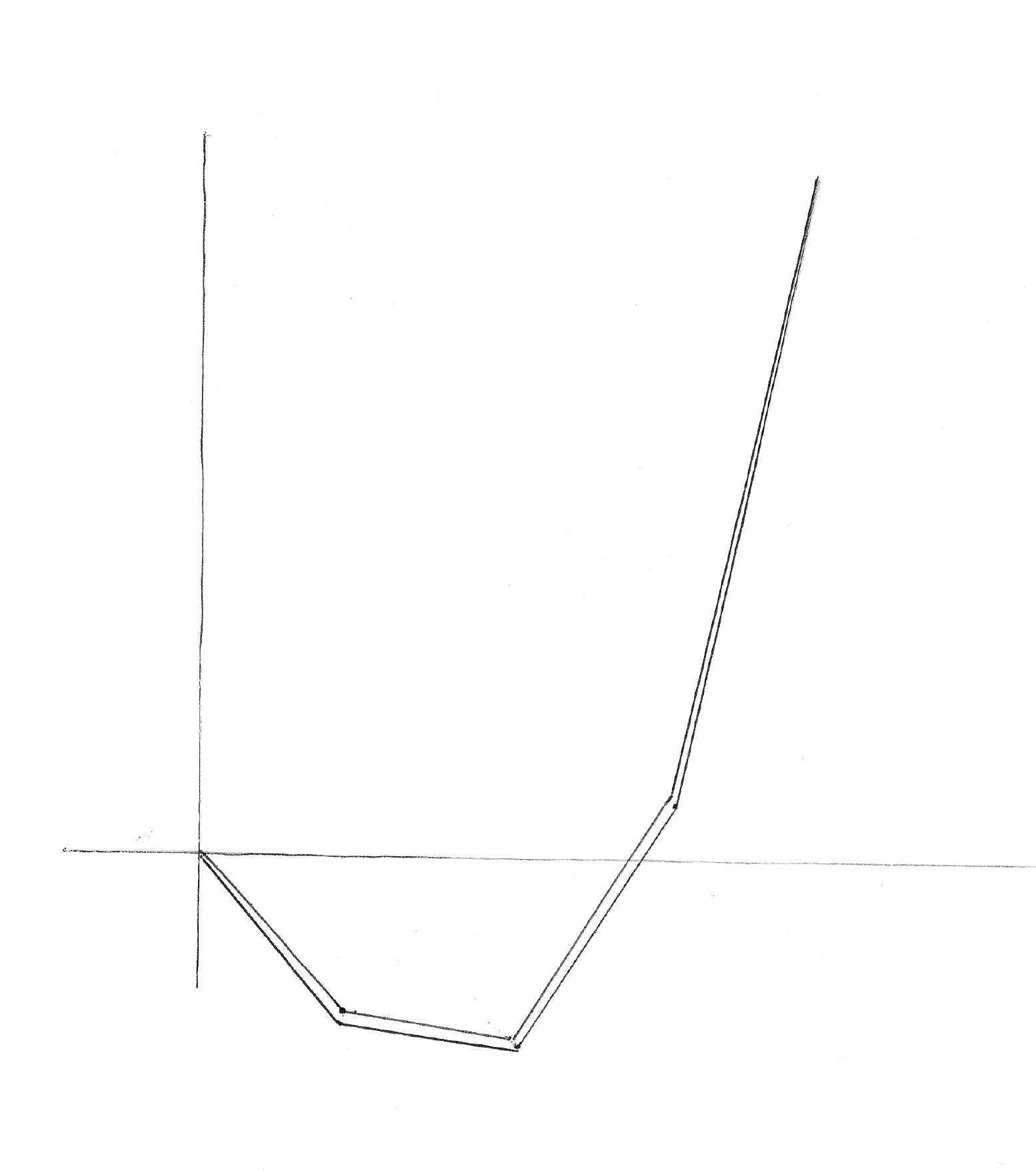}
\caption{{\small Le ${\rm D}_{\rm cris}((\rho_{x_\Phi})_v)$ d'un point $U(d)$-automorphe raffin\'e typique $x_\Phi \in \mathcal{E}_d$}}
\label{polygonenp}
\end{figure}

La situation typique est donc qu'au voisinage $\Omega$ assez petit d'un
point $x_\Ref \in \mathcal{E}_d$, les fonctions $x \mapsto v(F_i(x))$ sont
des constantes $v_i \in \Q$, alors que la variation des poids donn\'ees par
$\kappa$ est aussi g\'en\'erique que possible : si $y_{\Ref'} \in \Omega$ a
pour poids de Hodge-Tate en $v$ les entiers $k'_1 <\!\!< k'_2 <\!\!< \cdots
<\!\!< k'_d$ (les tels points s'accumulent en $z$) les pentes du
polyg\^one de Newton de ${\rm D}_{\rm cris}((\rho_{y_{\Ref'}})_v)$ sont les
$$k'_1 + v_1 <\!\!< k'_2 + v_2 <\!\!< \cdots <\!\!< k'_d+v_d.$$ Il en
r\'esulte que "vu de tr\`es loin", les polygones de Hodge et de Newton des
${\rm D}_{\rm cris}((\rho_{y_{\Ref'}})_v)$ avec $y_{\Ref'} \in \Omega$ vont
"presque co\"incider", avec un \'ecart mesur\'e par les $v_i$ !  La figure~\ref{polygonenp} ci-dessous illustre cette propri\'et\'e. Elle est d'ailleurs mieux comprise lorsque l'on introduit la notion
de repr\'esentation trianguline~\cite{colmeztri}, qui forme une
g\'en\'eralisation simultan\'ee des repr\'esentations totalement
r\'eductibles et des repr\'esentations cristallines de ${\rm
Gal}(\overline{\Q}_p/\Q_p)$ : nous reviendrons sur ce point dans le
paragraphe qui suit.  \ps

Observons enfin que si l'on \'etait parti d'un m\^eme $x$ mais d'un
raffinement $\Ref' \neq \Ref$ diff\'erent, les constantes $v_i$ ci-dessus
seraient en g\'en\'eral tout-\`a-fait diff\'erentes : ce serait les
$v_{\sigma(i)}+k_{\sigma(i)}-k_i$ si
$\Ref'_v=(\varphi_{\sigma(1)},\cdots,\varphi_{\sigma(d)})$ et $\sigma \in
\got{S}_d$ .  Il en suit une combinatoire plut\^ot amusante qui nous permet,
en partant du $x$ de l'\'enonc\'e ci-dessus et en bougeant un peu dans la
foug\`ere, de produire un $y \in V$ avec des polyg\^ones de Hodge et de
Newton profond\'ement modifi\'es, satisfaisant par exemple le crit\`ere
d'irr\'eductibilit\'e plus haut.  On v\'erifie dans~\cite{bchsigne}, auquel
je renvoie pour les d\'etails, qu'il suffit pour cela de se promener dans
$3$ branches bien choisies de la foug\`ere \`a partir de $x$. 
$\square$\ps\ps

Les observations  ci-dessus d\'evoilent des propri\'et\'es
int\'eressantes des repr\'esentations galoisiennes $(\rho_{x_\Ref})_v$ pour
$x_\Ref \in \mathcal{E}_d$ et $v \in \widetilde{S}_p$.  Pour aller plus loin
dans la compr\'ehension de la foug\`ere infinie, il est
n\'ecessaire de comprendre avec plus de pr\'ecision les propri\'et\'es de
ces $(\rho_x)_v$ pour un point g\'en\'eral $x \in \mathcal{E}_d$, ou du
moins dans un voisinage d'un $x_\Phi$.  Rappelons que ces $\rho_x$ sont
associ\'ees aux formes automorphes $p$-adiques propres de pente finie de
$U(d)$.  Elles n'ont pas \'et\'e construites par voie g\'eom\'etrique, mais
plut\^ot par interpolation $p$-adique \`a partir des points
$U(d)$-automorphes, de surcro\^it par un argument global utilisant
Cebotarev, ce qui rend peu \'evident le probl\`eme de comprendre
$(\rho_x)_v$.  \ps \newcommand{\Ro}{\mathcal{R}}

	Un progr\`es important dans cette direction a \'et\'e obtenu par
Kisin. Il a observ\'e dans~\cite{kisinoc} que dans les familles $p$-adiques
de repr\'esentations de ${\rm Gal}(\overline{\Q}_p/\Q_p)$ satisfaisant, dans
un cadre abstrait, des propri\'et\'es similaires \`a celles mises en
\'evidence ci-dessus, tous les $\rho_x$ h\'eritent automatiquement de
p\'eriodes cristallines par interpolation. Kisin consid\`ere pr\'ecis\'ement
le cadre abstrait suivant. \ps

Soient $Y$ un espace analytique $p$-adique r\'eduit et s\'epar\'e, $F \in
\OO(Y)^\ast$ une fonction analytique inversible, $M$ un $\OO_Y$-module
libre de rang $d$ muni d'une repr\'esentation $\OO_Y$-lin\'eaire
continue de ${\rm Gal}(\overline{\Q}_p/\Q_p)$ et $Z \subset Y$ un
sous-ensemble Zariski-dense, ayant les deux propri\'et\'es suivantes.  D'une
part, on suppose que pour tout $z \in Z$, la $k(z)$-repr\'esentation $M_z$
obtenue par \'evaluation de $M$ en $z$ est cristalline et satisfait ${\rm
D}_{\rm cris}(M_z)^{\varphi=F(z)} \neq 0$.  D'autre part, on demande que
pour tout r\'eel $C>0$ l'ensemble des $z \in Z$ tels que $M_z$ admet $0$ pour
poids de Hodge-Tate simple, et a tous ses autres poids $>C$, est Zariski-dense
dans $Y$.  Kisin d\'emontre alors dans~\cite{kisinoc} que 
\begin{equation}\label{ckisin1} \forall y
\in Y, \, \, {\rm D}_{\rm cris}(M_z)^{\varphi=F(y)} \neq 0.\end{equation} Il d\'emontre
\'egalement que si ${\rm D}_{\rm cris}(M_y^{\rm ss})^{\varphi=F(y)}$ est de
$k(y)$-dimension $1$, $M_y^{\rm ss}$ d\'esignant la semisimplifi\'ee de la
$k(y)$-repr\'esentation $M_y$, alors pour tout id\'eal $I \subset \OO_y$ de
codimension finie alors 
\begin{equation}\label{ckisin2}{\rm D}_{\rm cris}(M/IM)^{\varphi=F}\end{equation} 
est libre de rang $1$ sur $\OO_y/I$.  Ces r\'esultats de Kisin sont \'el\'ementaires quand la
fonction $y \mapsto v(F(y))$ est identiquement nulle, car la famille $M$
poss\`ede alors une famille de sous-repr\'esentations de rang $1$, mais le cas
g\'en\'eral repose sur des propri\'et\'es plus fines des anneaux de
Fontaine.  \ps 

Je voudrais mentionner deux extensions de ces r\'esultats de Kisin.  D'une part
Nakamura les a \'etendu dans~\cite{nakamurakisin} au cas o\`u le
corps de base $\Q_p$ est remplac\'e par une extension finie de $\Q_p$. 
D'autre part, dans~\cite[Ch. 3]{bchlivre}, nous les avons \'etendu 
aux $\OO_Y$-modules coh\'erents $M$ simplement suppos\'es sans
torsion, ce qui est le cadre naturel dont on a besoin pour les applications
aux vari\'et\'es de Hecke $Y$. Quand $Y$ est une courbe, un argument simple
de normalisation permet de se ramener au cas de Kisin, alors que dans notre
contexte avec Bella\"iche, nous devons raisonner par platification du module
$M$ sur un \'eclatement bien choisi, ce qui introduit un certain nombre de
difficult\'es techniques. Comme je viens de le dire, Kisin applique ses r\'esultats
\`a la normalisation de la courbe de $\mathcal{C}$ de Coleman-Mazur, et
obtient notamment que si $f=\sum_{n\geq 0} a_n q^n$ est une forme modulaire
$p$-adique surconvergente de poids quelconque (\'eventuellement non entier),
normalis\'ee, propre pour tous les op\'erateurs de Hecke, et telle que $a_p
\neq 0$, alors \begin{equation}\label{condkisin}{\rm D}_{\rm
cris}((\rho_f)_p)^{\varphi=a_p} \neq 0.\end{equation} (Quand $f \in M_k({\rm
SL}(2,\Z))_\iota$ est une forme modulaire usuelle, on sait m\^eme que
$(\rho_f)_p$ est cristalline en $p$ et que le polyn\^ome caract\'eristique
du Frobenius cristallin de ${\rm D}_{\rm cris}((\rho_f)_p)$ est le
polyn\^ome $P_f=X^2-a_p X + p^{k-1}$ du~\S\ref{arcgm}).  \ps

 Les repr\'esentations continues $\rho : {\rm Gal}(\overline{\Q}_p/\Q_p)
\rightarrow \GL_2(\overline{\Q}_p)$ telles que ${\rm D}_{\rm cris}(\rho)
\neq 0$ ont depuis suscit\'e une attention toute
particuli\`ere.\footnote{Une raison suppl\'ementaire \`a cela est l'espoir
formul\'e par Kisin que la foug\`ere infinie dans $X_2$ puisse \^etre
simplement caract\'eris\'ee comme \'etant le sous-ensemble des $x \in X_2$
tels que $(\rho_x)_{|{\rm Gal}(\overline{\Q}_p/\Q_p)}$ ait cette
propri\'et\'e (voir~\cite{emertonbki} pour les progr\`es r\'ecents \`a ce
sujet).  } \renewcommand{\fg}{(\varphi,\Gamma)} Colmez a notamment observ\'e que la
condition de Kisin s'interpr\^ete de mani\`ere tr\`es simple sur le
$\fg$-module ${\rm D}_{\rm rig}(\rho)$ \`a coefficients dans l'anneau de
Robba, ce sur quoi nous reviendrons dans la
partie suivante. Il montre qu'il existe un caract\`ere continu $\chi :
{\rm Gal}(\overline{\Q}_p/\Q_p) \rightarrow \overline{\Q}_p^\ast$ tel que
${\rm D}_{\rm cris}(\rho \otimes \chi) \neq 0$ si et seulement si ${\rm
D}_{\rm rig}(\rho)$ est r\'eductible, et il obtient une classification
compl\`ete de ces repr\'esentations. Colmez appelle plus g\'en\'eralement
{\it triangulines} les $\overline{\Q}_p$-repr\'esentations $V$ de ${\rm
Gal}(\overline{\Q}_p/\Q_p)$ telles que ${\rm D}_{\rm rig}(V)$ soit extension
successive de $\fg$-modules de rang $1$ sur l'anneau de Robba.
Dans~\cite[Ch. 2]{bchlivre}, nous \'etendons l'approche de Colmez en donnant des interpr\'etations "triangulines" de la
condition de nature infinit\'esimale~\eqref{ckisin2} d\'emontr\'ee par Kisin,
ce qui fournit une variante utile de l'approche de Kisin
dans~\cite{kisinoc}. C'est par exemple l'un des ingr\'edients dans ma d\'emonstration
du Th\'eor\`eme~\ref{eclatementeigen} (\cite{chhecke}). Je renvoie
\`a l'article d'exposition de Berger~\cite{bergertri} pour un survol de la
th\'eorie des repr\'esentations triangulines.
\footnote{Mentionnons tout de m\^eme l'extension par Nakamura~\cite{nakamuratri} des
r\'esultats de Colmez~\cite{colmeztri} au cas des extensions finies de
$\Q_p$.} \ps 

Nous avons poursuivi cette voie dans notre livre~\cite{bchlivre} dans
le but de comprendre les propri\'et\'es des $(\rho_x)_v$ pour $x \in
\mathcal{E}_d$ et $d$ g\'en\'eral. Comme nous l'avions observ\'e depuis longtemps~\cite{bchens}, les
r\'esultats de Kisin susmentionn\'es s'appliquent aux puissances ext\'erieures de la famille
de (pseudo)-repr\'esentations galoisiennes port\'ee par $\mathcal{E}_d$, et montrent   
dans les notations plus haut que pour tout $x \in \mathcal{E}_d$, $${\rm  
D}_{\rm cris}(\Lambda^i (\rho_x)_v \otimes \chi_1 \cdots
\chi_i)^{\varphi=F_1(x)\cdots F_i(x)} \neq 0, \, \, \, \, \, \forall
i=1,\cdots,d,$$ pour certain caract\`eres $\chi_i$ que je ne pr\'ecise pas
ici mais qui sont d\'etermin\'es par $\kappa(x)$. En dimension $d>2$, cette
propri\'et\'e n'est pas suffisante pour entra\^iner que $(\rho_x)_v$ est
trianguline au sens de Colmez,\footnote{Consid\'erer la somme directe d'un
caract\`ere et d'une repr\'esentation non trianguline de dimension $2$.}  
m\^eme si cela l'indique fortement. Nous nous attendions ainsi \`a ce que
ce soit le cas, et \'egalement que les sous-$\fg$-modules associ\'es des ${\rm
D}_{\rm rig}((\rho_x)_v)$ varient essentiellement analytiquement sur
$\mathcal{E}_d$. \ps

Modulo quelques acrobaties, cette m\'ethode des puissances ext\'erieures
nous a tout de m\^eme permis de
d\'emontrer dans~\cite{bchlivre} les propri\'et\'es de
triangulinit\'e attendues dans un voisinage formel des points
$U(d)$-automorphes raffin\'es "non-critiques" et "r\'eguliers", qui est le
cas crucial pour les applications arithm\'etiques de notre livre, et
\'egalement un ingr\'edient pour d\'emontrer le Th\'eor\`eme~\ref{densite3}. Je renvoie au Th\'eor\`eme~\ref{critlocal} ci-dessous pour un \'enonc\'e pr\'ecis. 
\ps

Les difficult\'es d'apparence technique pour \'etudier les points
g\'en\'eraux se sont av\'er\'ees assez s\'erieuses, et n'ont d'ailleurs
\'et\'e finalement r\'esolues dans la g\'en\'eralit\'e qu'elles m\'eritent
que tout r\'ecemment, ind\'ependamment par Hellman~\cite{hellmann},
Kedlaya-Pottharst-Xiao~\cite{kpx}, Liu~\cite{liu2}.  Je renvoie d'ailleurs
\'egalement \`a~\cite{chhecke} pour d'autres mises en gardes d\'ej\`a
visibles dans le cas de la courbe de Coleman-Mazur lorsque l'on
s'int\'eresse aux propri\'et\'es de triangulinit\'e en famille, qui ne
peuvent \^etre vrai au sens strict au voisinage de tous les points de
$\mathcal{E}_d$.  \ps

\subsection{D\'eformations triangulines des repr\'esentations cristallines raffin\'ees}\label{deftricris} Dans cette partie, nous oublions temporairement la foug\`ere infinie pour \'etudier un probl\`eme purement local, ce qui donnera \'egalement sans doute un peu d'air au lecteur. \ps

Soit $L$ une extension finie de $\Q_p$, et soit $V$ une $L$-repr\'esentation cristalline de dimension $d$  de $G_p:={\rm Gal}(\overline{\Q_p}/\Q_p)$. Il sera commode de faire les hypoth\`eses suivantes : \ps \begin{itemize}\ps
\item[(i)] ${\rm End}_{G_p}(V)$ est r\'eduit aux homoth\'eties $L$.\ps
\item[(ii)] Les poids de Hodge-Tate de $V$ sont distincts. \ps
\item[(iii)] Le polyn\^ome caract\'eristique du Frobenius cristallin de ${\rm D}_{\rm cris}(V)$ est scind\'e dans $L$ \`a racines distinctes. De plus, si $\phi$ et $\phi'$ sont deux telles racines, alors $\phi^{-1}\phi' \neq p$. \ps
\end{itemize}

Soit $\mathcal{C}$ la cat\'egorie des $L$-alg\`ebres locales artiniennes de corps r\'esiduel isomorphe \`a $L$.  On s'int\'eresse au foncteur
$$\mathcal{X}_V : \mathcal{C} \rightarrow {\rm Ens}$$ des d\'eformations de
$V$ \`a $\mathcal{C}$~(\cite{mazurdef},\cite{kisinoc}).  Par d\'efinition,
pour tout objet $A$ de $\mathcal{C}$, on d\'esigne par $\mathcal{X}_V(A)$
l'ensemble des classes d'isomorphie de repr\'esentations continues de
$G_p$ sur un $A$-module libre de rang fini $V_A$, telles que $V_A \otimes_A L
\simeq V$ en tant que $A[G_p]$-module.  Les conditions (i) et (iii), ainsi que des r\'esultats de Tate,
assurent que $\mathcal{X}_V$ est formellement lisse de dimension
$d^2+1=\dim(H^1(G_p,{\rm End}(V)))$, i.e.  $$\mathcal{X}_V \simeq {\rm
Spf}(L[[X_0,\cdots,X_{d^2}]].$$

Nous allons maintenant introduire et comparer toute une collection de sous-foncteurs naturels de $\mathcal{X}_V$. Comme plus haut, nous appellerons {\it raffinement de $V$} la donn\'ee d'un ordre $\Phi=(\phi_1,\cdots,\phi_d)$ sur les valeurs propres de $\varphi$ agissant sur ${\rm D}_{\rm cris}(V)$. Nous allons notamment associer \`a $(V,\Phi)$, suivant~\cite[Ch. 2]{bchlivre}, un sous-foncteur $$\mathcal{X}_{V,\Phi} \subset \mathcal{X}_V$$
appel\'e foncteur des d\'eformations $\Phi$-triangulines de $V$. La d\'efinition de ce foncteur est un peu d\'etourn\'ee : elle passe par les $\fg$-modules sur l'anneau de Robba $\Ro_L$. \ps

On rappelle que ce dernier d\'esigne l'anneau des s\'eries de Laurent
$$f(z)= \sum_{n\in \Z} a_n (z-1)^n$$ telles que $a_n \in L$ pour tout $n \in
\Z$, et qui convergent pour tout $z \in \mathbb{C}_p$ dans une couronne de
la forme $r_f \leq |z-1|<1$, le rayon int\'erieur $r_f$ d\'ependant de $f$. 
Il est muni d'actions commutantes de la lettre $\varphi$ et du groupe
$\Gamma=\Z_p^\ast$ par les formules $$\varphi(f)(z)=f(z^p), \, \, \, \,
\gamma(f)(z)=f(z^\gamma) \, \, \, \, \, \forall \gamma \in \Gamma.$$ \ps Un
$\fg$-module sur $\Ro_L$ est un $\Ro_L$-module libre de rang fini muni
d'actions semi-lin\'eaires commutantes de $\varphi$ et du groupe $\Gamma$
satisfaisant les deux conditions suivantes : $\varphi(D)$ engendre $D$ comme
$\Ro_L$-module et l'action de $\Gamma$ est continue dans un certain sens.  Les
$\fg$-modules sur $\Ro_L$ forment une cat\'egorie $L$-lin\'eaire tensorielle
$\fg/L$ de mani\`ere naturelle. Nous noterons aussi ${\rm Rep}_L$ la
cat\'egorie des repr\'esentations $L$-lin\'eaires continues de ${\rm
Gal}(\overline{\Q_p}/\Q_p)$ sur un $L$-module libre de rang fini.  Ainsi que
l'ont observ\'e Berger et Colmez (voir~\cite{colmeztri}), les travaux de
Fontaine~\cite{fontainegr}, combin\'es \`a ceux de
Cherbonnier-Colmez~\cite{chercol} et de Kedlaya~\cite{kedlaya}, assurent
l'existence d'un $\otimes$-foncteur exact, $L$-lin\'eaire, et pleinement
fid\`ele, $${\rm D}_{\rm rig} : {\rm Rep}_L \rightarrow \fg/L.$$ Un
$\fg$-module sur $\Ro_L$ est dans l'image essentielle de ${\rm D}_{\rm rig}$
si et seulement si il est {\it \'etale}, une condition portant uniquement
sur l'action de $\varphi$ (voir~{\rm loc.  cit.}). Ces d\'efinitions
et \'enonc\'es s'\'etendent alors verbatim lorsque $L$ est partout remplac\'e par un
objet $A$ quelconque de la cat\'egorie $\mathcal{C}$~ (\cite[Ch. 
2]{bchlivre}).\ps \newcommand{\Fil}{\mathrm{Fil}} 

\ps Retournons \`a la repr\'esentation $V$ plus haut et posons $D={\rm D}_{\rm rig}(V)$.  Soit
$t={\rm log}(z) \in \Ro_L$, les travaux de
Berger~\cite{berger1},\cite{berger2} d\'emontrent l'existence d'un
isomorphisme canonique $L[\varphi]$-\'equivariant $$\beta: {\rm D}_{\rm
cris}(V) \isomo (D[1/t])^\Gamma.$$ Dans cette bijection, on observe que les
sous-espaces $L[\varphi]$-stables $W \subset {\rm D}_{\rm cris}(V)$
correspondent bijectivement aux sous-$\fg$-modules $D' \subset D$ qui sont
facteurs directs comme $\Ro_L$-modules, via $$W \mapsto
(\Ro[1/t]\beta(W))\cap D,$$ 
le rang de $(\Ro[1/t]\beta(W))\cap D$ sur $\Ro_L$ \'etant $\dim_L W$. Mais les raffinements de $V$ correspondent
canoniquement aux $L$-drapeaux complets $\varphi$-stables de ${\rm D}_{\rm cris}(V)$
d'apr\`es l'hypoth\`ese (iii) sur $V$. L'observation pr\'ec\'edente les met
donc \'egalement en bijection avec les filtrations croissantes
$(\Fil_i(D))_{0\leq i \leq d}$ de $D$ par des sous-$\fg$-modules $\Fil_i(D)$
de rang $i$ qui sont facteurs directs comme $\Ro_L$-modules.  En
particulier, $V$ est trianguline au sens de Colmez, et ce d'autant de
fa\c{c}ons diff\'erentes qu'il y a de raffinements de $V$, \`a savoir $d!$. 
\ps

Je peux enfin d\'efinir le foncteur $\mathcal{X}_{V,\Phi} : \mathcal{C} \rightarrow {\rm Ens}$. Soit $(\Fil_i(D))$ la filtration de $D$ associ\'ee \`a $\Phi$ par la recette ci-dessus. Si $A$ est un objet de $\mathcal{C}$,  on d\'esigne par $\mathcal{X}_{V,\Phi}(A)$
l'ensemble des classes d'isomorphie de triplets $(V_A,(\Fil_i(D_A)),\pi)$ o\`u : \ps \begin{itemize} 
\item[-] $V_A$ est un $A$-module libre de rang $d$ muni d'une action $A$-lin\'eaire continue de $G_p$,\ps
\item[-] $(\Fil_i(D_A))_{0\leq i \leq d}$ est une filtration croissante de $D_A={\rm D}_{\rm rig}(V_A)$  par des sous-$\fg/A$-modules de rang $i$ facteurs directs comme $\Ro_A$-modules, \ps
\item[-] $\pi$ est un isomorphisme $A[G_p]$-\'equivariant $V_A \otimes_A L \isomo V$ tel que l'isomorphisme ${\rm D}_{\rm rig}(\pi) :  D_A \otimes_A L \isomo D$ envoie $\Fil_i(D_A)$ sur $\Fil_i(D)$ pour tout $0 \leq i \leq d$.\ps
\end{itemize}
\medskip
\noindent On d\'emontre dans~\cite{bchlivre} que le morphisme naturel $\mathcal{X}_{V,\Phi} \rightarrow \mathcal{X}_V$ fait de $\mathcal{X}_{V,\Phi}$ un sous-foncteur de $\mathcal{X}_V$, et que ce sous-foncteur est de plus pro-repr\'esentable et formellement lisse de dimension $$\dim(\mathcal{X}_{V,\Phi})=\frac{d(d+1)}{2}+1.$$ On utilise pour cela des r\'esultats de Colmez~\cite{colmeztri}, compl\'et\'es par Liu~\cite{liu}, sur la cohomologie \`a la Fontaine-Herr des $\fg$-modules de rang $1$ sur $\Ro_L$. \ps

Outre cette collection de $d!$ sous-foncteurs canoniques de
$\mathcal{X}_{V}$ ainsi construits on dispose \'egalement des sous-foncteurs
$$\mathcal{X}_{V,{\rm cris}} \subset \mathcal{X}_{V,{\rm HT}} \subset
\mathcal{X}_V$$ param\'etrant les d\'eformations $V_A$ de $V$ qui sont
respectivement cristallines et de Hodge-Tate, vues comme
$L$-repr\'esentations de $G_p$ par oubli de la $A$-structure.  Le polyn\^ome
de Sen de la d\'eformation universelle d\'efinit \'egalement un morphisme de
foncteurs $\Delta : \mathcal{X}_V \rightarrow \mathbb{G}_a^d$ tel que
l'\'equation $\Delta = 0$ d\'ecoupe exactement le ferm\'e
$\mathcal{X}_{V,{\rm HT}} \subset \mathcal{X}_V$.  \ps

Le r\'esultat suivant est un de nos r\'esultats clefs dans~\cite[Ch.
2]{bchlivre}.  On dit que le raffinement $\Phi$ de $V$ est non-critique si
le drapeau complet associ\'e de ${\rm D}_{\rm cris}(V)$ est en position
g\'en\'erale avec la filtration de Hodge sur ${\rm D}_{\rm cris}(V)$.
Observons que cette filtration admet exactement un sous-espace de chaque
dimension $\leq d$ par la condition (ii) sur $V$.

\begin{thm}\label{deftricrit} Supposons que $\Phi$ est un raffinement non critique de $V$. Alors $\mathcal{X}_{V,{\rm cris}} \subset \mathcal{X}_{V,\Phi}$, et pour tout objet $A$ de $\mathcal{C}$ on a 
$\mathcal{X}_{V,{\rm cris}}(A)=\mathcal{X}_{V,{\rm HT}}(A) \cap \mathcal{X}_{V,\Phi}(A)$.  
En particulier, l'application de Sen induit une suite exacte sur les espaces tangents
$$ 0 \longrightarrow \mathcal{X}_{V,{\rm cris}}(L[\varepsilon]) \longrightarrow \mathcal{X}_{V,\Phi}(L[\varepsilon]) \overset{\Delta} \longrightarrow L^d \longrightarrow 0.$$
\end{thm}

L'hypoth\`ese de non-criticit\'e de $\Phi$ est n\'ecessaire. Sous cette
hypoth\`ese, le foncteur $\mathcal{X}_{V,\Phi}$ appara\^it comme \'etant une
l\'eg\`ere g\'en\'eralisation de $\mathcal{X}_{V,{\rm cris}}$.  Comme nous
allons le voir il est \'egalement particuli\`erement pertinent dans
l'\'etude des vari\'et\'es de Hecke. \ps

Repla\c{c}ons-nous en effet dans les hypoth\`eses du~\S\ref{fougereunitaire}. 
Fixons $v \in \widetilde{S}_p$ telle que $F_v=\Q_p$.  Soit $(x,\Phi)$ un point
automorphe raffin\'e, ce qui d\'efinit un point $x_\Phi \in \mathcal{E}_d$. 
Supposons que $V=(\rho_x)_v$ est irr\'eductible.  La restriction \`a un
groupe de d\'ecomposition en $v$ fournit alors un morphisme canonique
\begin{equation}\label{locglobnat}\widehat{\mathcal{E}_d}^x
\longrightarrow \mathcal{X}_V\end{equation} o\`u $\widehat{\mathcal{E}_d}^x$ d\'esigne
le compl\'et\'e formel de $\mathcal{E}_d$ en $x$.  La question soulev\'ee
\`a la fin du~\S\ref{appraffinee} est de
qualifier son image. \ps

On dira que $\Ref_v=(\phi_{1,v},\cdots,\phi_{d,v})$ est {\it r\'egulier} si les $\varphi_{i,v}$ sont deux \`a deux distincts,
tels que $\varphi_{i,v}\varphi_{j,v}^{-1} \neq p$ pour tout $i,j$, et si de
plus pour tout entier $i=1,\cdots,d$, l'\'el\'ement $\prod_{j\leq
i}\varphi_{j,v}$ est valeur propre simple de $\varphi$ sur ${\rm D}_{\rm
cris}(\Lambda^i (\rho_x)_v)$.  Cette derni\`ere hypoth\`ese est purement technique
dans ce qui suit, et li\'ee \`a la m\'ethode des puissances ext\'erieures
que nous employons : elle pourrait d'ailleurs \^etre supprim\'ee dans
l'\'enonc\'e qui suit \`a l'aide des travaux de~\cite{kpx}.  En revanche, la
premi\`ere est importante.  Elle assure notamment que $V$ satisfait les
conditions g\'en\'erales de ce chapitre et permet \'egalement de voir
$\Phi_v$ comme un raffinement de $(\rho_x)_v$ au sens ci-dessus.  On dira
alors que $\Phi_v$ est non-critique si ce dernier l'est.  Le r\'esultat
ci-dessous est d\'emontr\'e dans~\cite[Ch.  4]{bchlivre}.

\begin{thm}\label{critlocal} Soit $x_\Phi \in \mathcal{E}_d$ un point automorphe raffin\'e.
On suppose que $V=(\rho_{x_\Phi})_v$ est irr\'eductible et que $\Ref_v$ en est un
raffinement r\'egulier et non-critique. Alors l'application
naturelle~\eqref{locglobnat} se factorise par $\mathcal{X}_{V,\Ref_v}
\subset \mathcal{X}_V$.
\end{thm}

Cet \'enonc\'e met bien en valeur l'int\'er\^et des foncteurs de
d\'eformations triangulines. Il rend de plus le
Th\'eor\`eme~\ref{deftricrit} particuli\`erement int\'eressant, notamment lorsqu'on le confronte au
th\'eor\`eme suivant, conjectur\'e dans~\cite{bchlivre} et d\'emontr\'e dans
mon article~\cite{chens}.

\begin{thm}\label{poidsetale} Soit $x_\Phi \in \mathcal{E}_d$ un point
automorphe raffin\'e. On suppose que pour tout $v \in \widetilde{S}_p$,
$F_v=\Q_p$, $V=(\rho_{x_\Phi})_v$ est irr\'eductible, et que $\Ref_v$ en est un
raffinement r\'egulier et non-critique. Alors l'application $\kappa : \mathcal{E}_d
\longrightarrow \mathcal{T}_0^d$ est \'etale en $x_\Phi$.
\end{thm}

La d\'emonstration de ce th\'eor\`eme utilise les propri\'et\'es les plus fines des
vari\'et\'es de Hecke. Elle repose \'egalement sur divers \'enonc\'es de
multiplicit\'e $1$ en th\'eorie des formes automorphes pour $U(d)$ (\cite{roglivre},
\cite{labesse}), ainsi que sur une g\'en\'eralisation au cas de $U(d)$ de la
th\'eorie des formes compagnons et de l'op\'erateur "$\Theta^{k-1}$" de la th\'eorie des formes
modulaires $p$-adiques (voir notamment~\cite{jones} pour le complexe BGG localement analytique). \ps

Retournons au contexte purement local pr\'ec\'edent. Ma derni\`ere contribution sur le sujet des d\'eformations triangulines est
un th\'eor\`eme comparant les divers $\mathcal{X}_{V,\Ref}$ \`a
l'int\'erieur de $\mathcal{X}_V$. Le r\'esultat suivant est peut-\^etre
l'innovation principale de mon article~\cite{chens}.

\begin{thm}\label{cltri} Supposons que les $d!$ raffinements de $V$ soient non-critiques, alors dans l'espace tangent $\mathcal{X}_V(L[\varepsilon])$ on a l'\'egalit\'e
$$\sum_{\Phi} \mathcal{X}_{V,\Phi}(L[\varepsilon]) = \mathcal{X}_V(L[\varepsilon]).$$
\end{thm}

Autrement dit, toute d\'eformation de $V$ \`a l'ordre $1$ est une
combinaison lin\'eaire de d\'eformations triangulines. Je d\'emontre plus
pr\'ecis\'ement {\it loc. cit.} qu'il suffit que $d$ raffinements bien choisis
$\Phi_1,\cdots,\Phi_d$ de $V$ soient non-critiques pour que $\sum_{i=1}^d
\mathcal{X}_{V,\Phi_i}(L[\varepsilon]) = \mathcal{X}_V(L[\varepsilon])$. Si
$d\leq 3$ et $V$ est irr\'eductible, il existe en fait toujours trois tels
raffinements.
\ps

Ma d\'emonstration de ce th\'eor\`eme est une r\'ecurrence sur la dimension
$d$ de $V$, ce qui est peut-\^etre un peu surprenant \'etant donn\'e que $V$
peut tr\`es bien \^etre irr\'eductible !  Le point est en fait de
d\'emontrer un \'enonc\'e plus g\'en\'eral valable pour tous les
$\fg$-modules cristallins non n\'ecessairement \'etales, une notion que nous
avions introduite avec Bella\"iche dans~\cite[Ch. 2]{bchlivre}. L'argument
fonctionne alors comme suit. On part d'un raffinement non-critique
$\Phi$ de $V$. On v\'erifie qu'un foncteur de d\'eformation
"paraboline de type $(1,d-1)$" annexe bien choisi de $V$, que je
dois \'etudier au m\^eme titre que les foncteurs pr\'ec\'edents, a la
propri\'et\'e que son espace tangent engendre $\mathcal{X}_V(L[\varepsilon])$
avec $\mathcal{X}_{V,\Phi}(L[\varepsilon])$. On raisonne ensuite simplement par
r\'ecurrence dans le bloc de taille $d-1$.

\subsection{D\'emonstration du Th\'eor\`eme~\ref{densite3}}\label{preuvedensite3} Indiquons enfin comment terminer la d\'emonstration du Th\'eor\`eme~\ref{densite3}. On suppose que $F_v=\Q_p$ pour tout $v \in \widetilde{S}_p$. \ps
	L'id\'ee de la d\'emonstration est la suivante. Soit $W$ l'adh\'erence Zariski des points automorphes dans $Y_d^\bot$, et soit $W_0$ une composante irr\'eductible de $W$.
Le lieu singulier de $W$ est un ferm\'e strict par les propri\'et\'es d'excellence des affino\"ides~\cite{conradirr}, de sorte que $W_0$ contient au moins un point automorphe $x$ r\'egulier dans $W$.  Quitte \`a se d\'eplacer un peu dans la foug\`ere infinie en partant de $x$, on peut supposer que $(\rho_{x})_v$ est irr\'eductible pour tout $v \in \widetilde{S}_p$ d'apr\`es la Proposition~\ref{propirr}. On peut \'egalement supposer que ces $(\rho_{x})_v$ satisfont les hypoth\`eses (ii) et (iii) du~\S\ref{deftricris}. \ps

Nous allons regarder l'espace tangent $T_x(W_0)$ de $W_0$ en $x$. Soit $L$ une extension finie de $\Q_p$ assez grande de sorte que pour tout $v \in \widetilde{S}_p$, la repr\'esentation $V_v=(\rho_{x})_v$ soit d\'efinie sur $L$. On dispose d'une application $k(x)$-lin\'eaire naturelle 
$$T_x(W_0)\otimes_{k(x)}L  \longrightarrow \prod_{v \in \widetilde{S}_p} \mathcal{X}_{V_v}(L[\varepsilon]).$$
Il sera commode de noter $T_{V_v,?}$ le $L$-espace vectoriel tangent $\mathcal{X}_{V_v,?}(L[\varepsilon])$. Regardons l'application naturelle suivante d\'eduite du morphisme ci-dessus :
$$f : T_x(W_0)\otimes_{k(x)}L  \longrightarrow \prod_{v \in \widetilde{S}_p} T_{V_v}/T_{V_v,{\rm cris}}.$$

	Soit $\Phi=\{\Phi_v\}$ un raffinement de $x$. La branche $B_{x,\Phi}$ est incluse dans $W_0$ au voisinage de $x$. En particulier, ${\rm Im}(f)$ contient l'image naturelle de $T_{x_\Ref}(\mathcal{E}_d)$ dans le terme de droite. Mais si les $\Phi_v$ sont des raffinements r\'eguliers et non-critiques des $V_v$, la combinaison des Th\'eor\`emes~\ref{deftricrit}.~\ref{critlocal}, et~\ref{poidsetale} assure que cette image est simplement $$\prod_v T_{V_v,\Phi_v}/T_{V_v,{\rm cris}}.$$ En particulier, si les $V_v$ ont tous leurs raffinements non-critiques et r\'eguliers, l'application $f$ est surjective d'apr\`es le Th\'eor\`eme~\ref{cltri} ! Cela vaut plus g\'en\'eralement si ces $V_v$ ont suffisamment de raffinements non critiques, et en particulier si $d\leq 3$ par les remarques suivant ce m\^eme th\'eor\`eme. \ps
Observons enfin que l'espace d'arriv\'ee de $f$ est de dimension 
$[F^+:\Q]\frac{d(d+1)}{2}$ : c'est un calcul purement local, par exemple cons\'equence du Th\'eor\`eme~\ref{deftricrit}. Comme $W_0$ est lisse en $x$ cela entraine que $\dim(W_0) \geq [F^+:\Q]\frac{d(d+1)}{2}$, ce que l'on voulait d\'emontrer. $\square$ 
\ps \medskip
Il n'aura pas \'echapp\'e au lecteur la co\"incidence num\'erique r\'ealisant \`a la fois le nombre $[F^+:\Q]\frac{d(d+1)}{2}$ comme la codimension du lieu cristallin local et \'egalement comme la dimension globale attendue de $Y_d^\bot$. Cette co\"incidence n'est pas sans rappeler celles observ\'ees dans~\cite{cht} ! Ce qui m'emp\^eche de conclure en dimension $d>3$ est qu'actuellement je n'arrive pas \`a d\'emontrer l'existence de points automorphes dans $W_0$ ayant suffisamment de raffinements non critiques. 
Nous verrons en revanche dans la partie suivante que ces id\'ees s'appliquent parfaitement \`a un analogue purement local des questions consid\'er\'ees dans cette partie. \ps

\newpage
\section{Analogues locaux cristallins}

Dans cette partie, je me propose de pr\'esenter un analogue purement local des questions \'etudi\'ees aux~\S\ref{sectionfougereglobale}, et notamment de la foug\`ere infinie. J'expose les r\'esultats de mon article~"{\it Sur la densit\'e des repr\'esentations cristallines du groupe de Galois absolu de $\Q_p$}".\ps

Fixons $p$ un nombre premier et $d\geq 1$ un entier. Soit $X_d$ la vari\'et\'e des caract\`eres $p$-adique en dimension $d$ du groupe $G={\rm Gal}(\overline{\Q}_p/\Q_p)$ (\S\ref{carpadique}). Nous dirons qu'un point $x \in X_d$ est cristallin, si la repr\'esentation associ\'ee $\rho_x : G \rightarrow \GL_d(\overline{k(x)})$ l'est au sens de Fontaine. Une question certainement naturelle \`a ce stade dans ce m\'emoire est de s'interroger sur les propri\'et\'es du lieu des points cristallins de $X_d$. Un point $x \in X_d$ sera dit irr\'eductible si $\rho_x$ l'est.  Le r\'esultat principal de~\cite{chmannalen} est le suivant. Il est d\^u \`a ind\'ependamment \`a Kisin~\cite{kisinfern} et Colmez~\cite{colmeztri} en dimension $d=2$, et \'evident si $d=1$. Quand $d=2$, il joue notamment un r\^ole important dans la d\'emonstration par Colmez~\cite{colmezecm} de la correspondance de Breuil-Langlands $p$-adique pour $\GL_2(\Q_p)$. 

\begin{thm}\label{denscrisloc} Les points cristallins de $X_d$ sont Zariski-denses et d'accumulation dans toutes les composantes irr\'eductibles de $X_d$ contenant un point cristallin irr\'eductible. \end{thm}

Faisons quelques remarques avant de donner les id\'ees de la d\'emonstration. Si $r : G \rightarrow \GL_d(\overline{\F}_p)$ est une repr\'esentation continue et irr\'eductible telle que $r \not\simeq r \otimes \omega$, $\omega$ d\'esignant le caract\`ere cyclotomique modulo $p$, alors $X_d(r)$ est une boule ouverte de dimension $d^2+1$. Il n'est pas difficile de voir que cette boule contient des points cristallins en consid\'erant  par exemple les induites d'un produit de caract\`eres de Lubin-Tate de l'extension non ramifi\'ee de degr\'e $d$ de $\Q_p$. Le Th\'eor\`eme~\ref{denscrisloc} entra\^ine donc que les points cristallins sont Zariski-denses et d'accumulation dans un tel $X_d(r)$. Bien entendu, les repr\'esentations induites mentionn\'ees pr\'ecedemment ne suffisent pas \`a entra\^iner ce r\'esultat d\`es que $d>1$, leur adh\'erence Zariski \'etant plut\^ot de dimension $d+1$. \ps

Il me semble naturel d'esp\'erer que toutes les composantes irr\'eductibles de $X_d$ contiennent des points cristallins irr\'eductibles, et donc que les points cristallins sont Zariski-denses dans $X_d$ d'apr\`es le th\'eor\`eme. D'une certaine mani\`ere, cela constituerait un analogue local de la conjecture de modularit\'e de Serre.

\begin{conjecture} Toute composante irr\'eductible de $X_d$ contient un point cristallin irr\'eductible. \end{conjecture}

Les m\'ethodes d\'evelopp\'ees dans~\cite{chmannalen} m'ont permis de v\'erifier {\it loc. cit.} que pour toute repr\'esentation semisimple continue $r : G \rightarrow \GL_d(\overline{\F}_p)$, alors $X_d(r)$ contient en effet un point irr\'eductible cristallin. La conjecture ci-dessus semble cependant sensiblement plus fine. Quand $r$ est scalaire, je crois que l'on ne sait toujours pas si elle est v\'erifi\'ee pour les composantes de $X_2(r)$, c'est sans doute pour beaucoup d\^u au fait que l'on ne sait pas encore d\'ecrire ces derni\`eres comme je l'ai rappel\'e au~\S\ref{carpadique}, bien que cela ne semble pas insurmontable. \ps

Mentionnons que l'on peut se demander si les \'enonc\'es ci-dessus restent raisonnables si l'on ne s'int\'eresse qu'aux repr\'esentations cristallines apparaissant dans la cohomologie \'etale $p$-adique des vari\'et\'es projectives lisses sur $\Q_p$ ayant bonne r\'eduction sur $\Z_p$. Cela semble probable, mais il serait int\'eressant de le d\'emontrer.  En revanche, ma d\'emonstration du Th\'eor\`eme~\ref{denscrisloc} a r\'ecemment \'et\'e \'etendue par Nakamura~\cite{nakamuradensity} dans le cas o\`u $G$ est remplac\'e par le groupe de Galois absolu d'une extension finie quelconque de $\Q_p$.\ps

La m\'ethode que j'ai suivie pour d\'emontrer le Th\'eor\`eme~\ref{denscrisloc} est tr\`es proche de celle du~\S\ref{preuvedensite3} utilis\'ee pour d\'emontrer le Th\'eor\`eme~\ref{densite3}. On dira que le point $x \in X_d$ est {\it triangulin}, si la repr\'esentation $\rho_x$ est trianguline au sens de Colmez~\cite{colmeztri}, c'est-\`a-dire si le $\fg$-module ${\rm D}_{\rm rig}(\rho_x)$ sur l'anneau de Robba est extension successive de $\fg$-modules de rang $1$ (\S\ref{deftricris}). La d\'efinition de la foug\`ere infinie dans ce contexte est simplissime.

\begin{definition} La foug\`ere infinie $\mathcal{F}_d \subset X_d$ est l'ensemble des points triangulins. \end{definition}

L'inconv\'enient avec cette d\'efinition ensembliste est qu'elle ne dit pas grand chose sur la structure de $\mathcal{F}_d$. Dans le contexte global, cette structure provenait de la structure analytique naturelle de la vari\'et\'e de Hecke $\mathcal{E}_d$, dont les propri\'et\'es principales d\'ecoulaient de constructions globales en th\'eorie des formes automorphes. Mon but principal dans ce qui suit est d'expliquer comment construire de mani\`ere directe un analogue purement local des vari\'et\'es de Hecke, poss\'edant des propri\'et\'es similaires. Il me para\^it raisonnable de noter \'egalement $\mathcal{E}_d$ cet analogue local, ce qui ne devrait g\^ener qu'un lecteur \'etourdi. \ps

L'ensemble sous-jacent \`a $\mathcal{E}_d$ sera simplement l'ensemble des paires $(x,\mathcal{T})$ o\`u $x \in X_d$ est un un point triangulin et o\`u $\mathcal{T}$ d\'esigne un drapeau complet de ${\rm D}_{\rm rig}(\rho_x)$ constitu\'e de sous-$\fg$-modules qui sont facteurs directs comme modules sur l'anneau de Robba. On dispose en particulier tautologiquement d'une application canonique
$$\mu : \mathcal{E}_d \rightarrow X_d$$
d'image $\mathcal{F}_d$. Il suit \'egalement de la discussion du~\S\ref{deftricris} que les points cristallins sont dans $\mathcal{F}_d$. Mieux, les ant\'ec\'edents par $\mu$ d'un point cristallin $x \in X_d$ sont en bijection canonique avec l'ensemble des drapeaux $\varphi$-stables de ${\rm D}_{\rm cris}(\rho_x)$, i.e. des raffinements de $\rho_x$ au sens de Mazur~\cite{Maz}. \ps

Une grosse partie du travail consiste \`a comprendre comment munir l'ensemble $\mathcal{E}_d$ d\'efini ci-dessus d'une structure int\'eressante d'espace analytique $p$-adique. Quand $d=2$, c'est le point de vue adopt\'e par Colmez dans~\cite{colmeztri}. La structure sur $\mathcal{E}_2$ d\'efinie par Colmez est ad hoc en terme de sa param\'etrisation des repr\'esentations triangulines. L'espace $\mathcal{E}_2$ est alors lisse de dimension $4$. Il est naturellement construit au dessus de l'espace $\mathcal{T}^2$, o\`u $$\mathcal{T} \simeq \mathbb{G}_m \times \WW$$ 
d\'esigne l'espace analytique param\'etrant les caract\`eres $p$-adiques continus du groupe $\Q_p^\ast$ (\S\ref{fougereunitaire}). Cet espace $\mathcal{T}$ joue un r\^ole pr\'edominant car il param\`etre \'egalement de mani\`ere naturelle les $\fg$-modules de rang $1$ sur l'anneau de Robba. \ps

Pour cette structure, Colmez parvient \`a d\'emontrer que l'application $\mu : \mathcal{E}_2 \rightarrow X_2$ est localement analytique\footnote{Si $X$ et $Y$ sont des espaces analytiques $p$-adiques, et si $f : X \rightarrow Y$ est une application ensembliste, je dirai que $f$ est localement analytique en $x \in X$ s'il existe un ouvert affino\"ide $U \subset X$ contenant $x$ tel que $f_{|U} : U \rightarrow X$ soit induit par une fonction analytique de $U$ dans $X$. Cette d\'efinition conserve d'ailleurs son sens si $X$ est simplement suppos\'e localement analytique, en un sens \'evident.} en tous les points "suffisamment g\'en\'eriques" de $\mathcal{E}_2$, c'est d'ailleurs l'un des points les plus d\'elicats de ses constructions : voir~\cite[\S 5.1]{colmeztri}. Ainsi, $\mathcal{F}_2$ est un espace de dimension $4$ \'evoluant dans l'espace $X_2$ de dimension $5$. Les points cristallins sont en g\'en\'eral des points doubles de $\mathcal{F}_2$ et il est facile de voir sur la classification des repr\'esentations triangulines qu'ils sont Zariski-denses et d'accumulation dans $X_2$, il est alors ais\'e d'en d\'eduire le Th\'eor\`eme~\ref{denscrisloc} dans ce cas.\ps

Avant d'expliquer mon travail, je voudrais dire un mot de l'approche de Kisin~\cite{kisinfern}, qui est techniquement diff\'erente mais au fond similaire. L'espace $\mathcal{E}_2$ 
est remplac\'e chez Kisin par son {\it finite slope subspace} $$X_{\rm fs} \subset X_2 \times \mathbb{G}_m$$ d\'ej\`a d\'efini dans~\cite{kisinoc}. Par d\'efinition, c'est essentiellement l'adh\'erence Zariski des points cristallins raffin\'es, \`a la mani\`ere de Coleman-Mazur~\S\ref{eigencurvefern}. Kisin explique notamment pourquoi $X_{\rm fs}$ est de dimension $\leq 3$, ce qui est assez d\'elicat,
alors que Colmez construit explicitement $\mathcal{E}_2$, ce qui est peut-\^etre plus naturel par rapport \`a notre probl\`eme. Bien entendu, ne perdons pas de vue que la construction par Kisin de $X_{\rm fs}$ est bien ant\'erieure \`a celle de Colmez ! \ps
\newcommand{\cT}{\mathcal{T}}
Revenons au cas de la dimension $d$ quelconque. Je vais maintenant \'enoncer le r\'esultat clef de mon article~\cite{chmannalen}. Si $X$ est un affino\"ide sur $\Q_p$, je noterai  $\Ro_X$ l'anneau de Robba \`a coefficients dans $\OO(X)$. Je renvoie \`a~\cite{chmannalen} pour la d\'efinition pr\'ecise ainsi que pour la notion plus ou moins \'evidente de $\fg$-module sur $\Ro_X$. Si $\delta \in \cT(X)$, c'est-\`a-dire un homomorphisme continu $\Q_p^\ast \rightarrow \OO(X)^\ast$, on dispose d'un $\fg$-module $$\Ro_X(\delta)$$ de rang $1$ sur $\Ro_X$. 
Il poss\`ede une $\Ro_X$-base $e$ telle que $\varphi(e)=\delta(p)e$ et $\gamma(e)=\delta(\gamma) e$ pour tout $\gamma \in \Gamma=\Z_p^\ast$. Nous appellerons "$\fg$-module triangulaire rigidifi\'e de rang $d$ sur $\Ro_X$" la donn\'ee d'un quadruplet $$(D,(\Fil_i(D)),(\delta_i),({\rm gr}_i))$$ o\`u $D$ est un $\fg$-module de rang $d$ sur $\Ro_X$, $(\Fil_i(D))_{i=0,\cdots,d}$ est une filtration croissante de $D$ par des sous-$\fg$-modules de rang $i$ et facteurs directs comme $\Ro_X$-modules, $(\delta_i) \in \mathcal{T}(X)^d$, et o\`u pour tout entier $i=1,\cdots,d$, $${\rm gr}_i: \Fil_i(D)/\Fil_{i-1}(D) \isomo \Ro_X(\delta_i)$$ est un isomorphisme de $\fg$-modules sur $\Ro_X$. Un tel $\fg$-module sera dit {\it r\'egulier} si pour tout $i<j$, on a $\delta_i\delta_j^{-1} \in \mathcal{T}^{\rm reg}$, o\`u $\mathcal{T}^{\rm reg} \subset \mathcal{T}$ d\'esigne l'ouvert compl\'ementaire de l'ensemble des points param\'etrant les caract\`eres de $\Q_p^\ast$ de la forme $x \mapsto x^{-i}$ ou $x \mapsto |x|x^{i+1}$ pour en entier $i\geq 0$.\ps

Soit ${\rm Aff}$ la cat\'egorie des espaces affino\"ides et morphismes analytiques. On consid\`ere le foncteur $\triangle_d^\square : {\rm Aff} \longrightarrow {\rm Ens}$
associant \`a $X$ l'ensemble des classes d'isomorphie de $\fg$-modules triangulaires, rigidifi\'es, et r\'eguliers, de rang $d$ sur $\Ro_X$. \ps

\begin{thm}\label{reptri} Le foncteur $\bigtriangleup_d^{\square}$ est repr\'esentable par un espace analytique $p$-adique\footnote{Pr\'ecisons que cet espace n'est pas un affino\"ide, bien qu'uniquement d\'etermin\'e par ses points \`a valeurs affino\"ides.} qui est 
irr\'eductible et lisse sur $\Q_p$, de dimension $\frac{d(d+3)}{2}$. Les $\fg$-modules cristallins\footnote{Suivant~\cite{bchlivre}[Ch. 2], on dit qu'un $\fg$-module $D$ sur $\Ro_L$ est dit cristallin si $(D[1/t])^\Gamma$ est de $L$-dimension ${\rm rang}_{\Ro_L}(D)$. Lorsque $D={\rm D}_{\rm rig}(V)$, les th\'eor\`emes de Berger rappel\'es au~\S\ref{deftricris} montrent que $D$ est cristallin en ce sens si et seulement si $V$ l'est au sens de Fontaine.} sont Zariski-denses et d'accumulation dans $\bigtriangleup_d^\square$.
\end{thm}

Ce r\'esultat est en fait nouveau m\^eme pour $d=2$, o\`u il compl\`ete les r\'esultats de Colmez et de Kisin susmentionn\'es. La notion de rigidification utilis\'ee fait que l'espace
$\bigtriangleup_d^\square$ admet moralement $d-1$ dimensions de trop, cela n'aura pas vraiment d'incidence.\ps

Le coeur technique de sa d\'emonstration 
est un ensemble de r\'esultats sur les $\fg$-modules sur
$\Ro_X$, notamment sur leur cohomologie \`a la Fontaine-Herr (\cite{herr1},\cite{herr2}). Un point d\'elicat  consiste \`a v\'erifier  dans ce 
contexte que les complexes
$C_{\varphi,\gamma}$ et $C_{\psi,\gamma}$ sont quasi-isomorphes comme dans
la th\'eorie classique de Herr, ce qui avait notamment \'et\'e conjectur\'e
par Kedlaya dans~\cite[\S 2.6]{kedlayaseul}. Nous y parvenons pour les $\fg$-modules triangulaires sur $\Ro_X$. Nous calculons enfin la cohomologie des $\Ro_X(\delta)$, \'etendant des r\'esultats de Colmez (en degr\'es $0$ et $1$,
\cite{colmeztri})) et de Liu (en degr\'e $2$, \cite{liu}) dans le cas
particulier o\`u $X$ est un point, ainsi que des r\'esultats partiels de Bella\"iche~\cite{joeltrans}. Notre preuve, bien qu'inspir\'ee de celle de
Colmez, est en fait un peu plus simple, m\^eme dans le cas d'un point. Elle repose sur un d\'evissage cher \`a Colmez de l'anneau de Robba~\cite[Thm. 0.1]{colmezgros}. Voici un \'echantillon des r\'esultats obtenus. \ps

\begin{thm}\label{calccohomo} Soit $D$ un $\fg$-module triangulaire de rang $d$ sur $\Ro_X$. Alors $H^i(D)$ est de type fini sur $\OO(X)$ pour tout $i$. De plus, dans le groupe de
Grothendieck des $\OO(X)$-modules de type fini on a la relation
$$[H^0(D)]-[H^1(D)]+[H^2(D)]=-[\OO(X)^d].$$ Enfin, la formation des $H^i(D)$ commute \`a
tout changement de base plat affino\"ide. \ps
	Si de plus $D$ est r\'egulier, alors $H^0(D)=H^2(D)=0$ et $H^1(D)$ est libre de
rang $d$ sur $\OO(X)$. La formation des $H^i(D)$ commute dans ce cas \`a tout changement
de base affino\"ide. 
\end{thm}

Mentionnons \`a ce propos les r\'esultats r\'ecents de Kedlaya-Pottharst-Xiao~\cite{kpx}, qui constituent une vaste g\'en\'eralisation de ce dernier th\'eor\`eme. \ps

Pour revenir \`a notre probl\`eme initial, il reste \`a passer des familles de $\fg$-modules aux familles de repr\'esentations galoisiennes. Cette extension en famille des propri\'et\'es du foncteur ${\rm D}_{\rm rig}$ a \'et\'e \'etudi\'ee en d\'etail par Berger-Colmez dans un sens~\cite{bergercolmez}, et par Kedlaya-Liu~\cite{kedliu} dans l'autre. Elle n'est pas sans subtilit\'e. Les r\'esultats de Kedlaya et Liu montrent que pour tout $x \in  \bigtriangleup_d^\square$ param\'etrant un $\fg$-module $D_x$ \'etale, il existe un ouvert affino\"ide $U \subset  \bigtriangleup_d^\square$ contenant $x$ tel que : \begin{itemize}\ps
\item[(i)]  $D_y$ est \'etale pour tout $y \in U$, \ps
\item[(ii)] il existe une repr\'esentation continue $G \rightarrow \GL_d(\OO(U))$, unique \`a isomorphisme pr\`es, dont le $\fg$-module sur $\Ro_U$ associ\'e par Berger-Colmez soit isomorphe au $\fg$-module sur $\Ro_U$ d\'efini par l'inclusion $U \subset \bigtriangleup_d^\square$.
\end{itemize}
\ps
Notons $\mathcal{E}_d^\square \subset \bigtriangleup_d^\square$  le sous-ensemble constitu\'e des points $x \in\bigtriangleup_d^\square $ tels que le $\fg$-module $D_x$ param\'etr\'e par $x$ soit \'etale. La propri\'et\'e (i) ci-dessus munit $\mathcal{E}_d^\square$ d'une structure d'espace localement analytique $p$-adique.\footnote{On pourrait m\^eme en faire un vrai espace analytique de mani\`ere ad hoc en consid\'erant une r\'eunion disjointe d'ouverts affino\"ides $U \subset \bigtriangleup_d^\square$ satisfaisant les deux propri\'et\'es ci-dessus et formant une partition ensembliste de $\mathcal{E}_d^\square$. Cela n'est cependant pas tr\`es canonique et n'apporterait rien. Le lecteur choqu\'e par ce d\'enoument pourra consulter consulter \`a profit le travail d'Hellmann~\cite{hellmann}. } La propri\'et\'e (ii) assure enfin, apr\`es passage aux traces, qu'il existe une unique application localement analytique
$$\tau : \mathcal{E}_d^\square \longrightarrow X_d,$$
envoyant tout $x \in \mathcal{E}$ sur la trace $\tau(x)$ de la repr\'esentation $V_x$ telle que ${\rm D}_{\rm rig}(V_x) \simeq D_x$. C'est bien ce que nous voulions d\'emontrer, en version rigidifi\'ee, ce qui n'affecte en rien la strat\'egie. On a bien s\^ur $\mathcal{F}_d=\tau(\mathcal{E}_d^\square)$.
L'argument expliqu\'e au~\S\ref{preuvedensite3} conduit alors de m\^eme au Th\'eor\`eme~\ref{denscrisloc}. Il n'y a pas de difficult\'e particuli\`ere dans ce contexte local pour trouver des repr\'esentations cristallines ayant tous leurs raffinements non-critiques, cela d\'ecoule notamment facilement des travaux de Kisin~\cite{kisinjams}. $\square$ \ps

Pour finir, je voudrais mentionner que le Th\'eor\`eme~\ref{reptri} est \'egalement le point de d\'epart des travaux r\'ecents de Hellmann~\cite{hellmann} visant \`a g\'en\'eraliser la construction $X_{\rm fs}$ de Kisin \`a la dimension $d$. \ps

\chapter[Des cons\'equences arithm\'etiques des conjectures d'Arthur-Langlands]{Quelques cons\'equences des conjectures de Langlands et Arthur}

\section[Le lemme de Dehn arithm\'etique]{Corps de nombres \`a ramification prescrite : le lemme de Dehn arithm\'etique.}\label{dehnar} Dans cette partie, 
j'expose les r\'esultats de mes articles~\cite{chnf} {\it "On number 
fields with given ramification"} et~\cite{chclo} {\it "Corps de nombres peu ramifi\'es 
et repr\'esentations automorphes autoduales"}. Ce dernier est un 
travail en commun avec Laurent Clozel. \ps

\subsection{R\'esum\'e et perspectives}

Soient $S$ un ensemble fini de nombres premiers, $\Q_S$
une extension alg\'ebrique maximale de $\Q$ non ramifi\'ee hors de $S$ (et
de l'infini) et $$\G_S=\Gal(\Q_S/\Q).$$ \ps

Un r\'esultat bien connu de Minkowski affirme que $\G_\emptyset=\{1\}$. En revanche, si $S$ est non
vide, ce que l'on supposera d\'esormais, la structure de ces groupes $\G_S$ est tr\`es mal connue, et ce
malgr\'e leur omnipr\'esence en g\'eom\'etrie arithm\'etique. 
Par exemple, un r\'esultat d'Hermite assure que $\G_S$ n'a qu'un nombre fini de sous-groupes ferm\'es 
d'indice donn\'e, mais on ne sait pour aucun $S\neq \emptyset$ si $\G_S$ est topologiquement engendr\'e 
par un nombre fini d'\'el\'ements, une question pos\'ee par Shaffarevich au
congr\`es international de Stockholm en 1962. Nous renvoyons au livre de
Neukirch, Schmidt et Wingberg~\cite[Ch. X, \S 11]{nsw} pour un \'etat de
l'art sur la structure de ${\rm G}_S$. \ps

Un autre probl\`eme du folklore, pos\'e \`a ma connaissance par Ralph Greenberg, consiste \`a d\'eterminer les sous-groupes de d\'ecomposition de $\G_S$. 
Malgr\'e le peu d'indices dont nous disposons pour appr\'ehender cette question, 
il semble commun\'ement esp\'er\'e que ces groupes soient aussi gros que la restriction 
impos\'ee sur la ramification le permette : les r\'esultats que j'ai obtenus vont dans cette
direction. Ainsi que me l'a racont\'e Greenberg, sa question fait echo \`a une question
ant\'erieure que lui avait pos\'ee James Milne au sujet de la
d\'etermination du pro-cardinal de ${\rm G}_S$. Comme l'a remarqu\'e Milne,
cette question intervient notamment dans la d\'etermination de la cohomologie galoisienne des corps de
nombres (voir~\cite[Chap. I \S 4]{milneadt}). \ps
\ps

Pr\'ecis\'ement, je me place dans le contexte suivant. Soit $p \in S$ un nombre premier. La donn\'ee d'un plongement de corps $\iota : \Q_S \longrightarrow \Qpb$ definit
un homomorphisme continu
\begin{equation}\label{decomap} \Gal(\Qpb/\Qp) \longrightarrow \G_S
\end{equation}
dont la classe de conjugaison est ind\'ependante du plongement choisi. La
question de Greenberg est de savoir si ce morphisme est injectif. C'est un exercice que de v\'erifier l'\'equivalence entre 
les propri\'et\'es suivantes : \ps \begin{itemize}

\item[(a)] l'application~\eqref{decomap} est injective, \ps
\item[(b)] $\iota(\Q_S)$ est dense dans $\Qpb$ pour la topologie $p$-adique, \ps
\item[(c)] pour toute extension finie $F$ de $\Q_p$, il existe un corps de nombres $K$ non ramifi\'e hors de $S$ et poss\'edant une place $v$ au dessus de $p$ telle que $F$ se plonge dans $K_v$. \ps
\end{itemize}

Autrement dit, la question de l'injectivit\'e de \eqref{decomap} revient \`a
savoir si il existe suffisamment de corps de nombres non ramifi\'es en
dehors de $S$ pour qu'ils engendrent tout $\overline{\Q}_p$ une fois
plong\'es dans ce dernier. \`A ma connaissance, avant les travaux que j'ai
effectu\'es cela n'\'etait connu pour aucun couple $(p,S)$. Bien
entendu, plus $S$ est petit, plus le probl\`eme est difficile, et
ultimement nous aimerions savoir si (\ref{decomap}) est injective pour
$S=\{p\}$. Dans tous les cas, une premi\`ere difficult\'e r\'eside en ce qu'il est
d\'ej\`a assez d\'elicat de construire de tels corps de nombres non
triviaux, $S$ \'etant fix\'e.\footnote{Le seul \'enonc\'e r\'eellement \'el\'ementaire venant \`a
l'esprit \`a ce sujet est un corollaire standard du lemme de Krasner : si
$F$ est une extension finie de $\Q_p$, il existe un corps de nombres $K$ de
m\^eme degr\'e que $F$ et poss\'edant une place $v$ au dessus de $p$ telle
que $K_v$ est isomorphe \`a $F$.  Bien entendu on ne contr\^ole en rien le
discriminant de $K$.  } Mon th\'eor\`eme principal est le suivant.\ps

\begin{thm}\label{chcl}  Si $|S|\geq 2$, alors (\ref{decomap}) est injective. \end{thm}

\begin{cor}\label{corchcl} Supposons $|S|\geq 2$. Alors le pro-cardinal de
${\rm G}_S$ est divisible par tous les entiers. 
\end{cor} 

Autrement dit, pour tout entier $m\geq 1$ il existe un corps de nombres non ramifi\'e hors de $S$ dont le degr\'e est un multiple de
$m$, ce qui \'etait la question de Milne d\'ej\`a cit\'ee. Le corollaire d\'ecoule simplement du th\'eor\`eme, par exemple en 
appliquant la condition (c) aux extensions non ramifi\'ees de $\Q_p$. Le cas
$S=\{p\}$ est par contre toujours ouvert : j'y reviendrai plus loin. \ps

Bien que les \'enonc\'es ci-dessus soient de formulation \'el\'ementaire,
notre d\'emonstration reposera sur une s\'erie de r\'esultats appartenant
\`a la th\'eorie arithm\'etique des formes modulaires et \`a ses
g\'en\'eralisations en dimensions sup\'erieures, qui eux ne le
sont pas du tout. Un survol de la litt\'erature concernant ${\rm G}_S$ montre que cette disproportion apparente est
en r\'ealit\'e devenue assez commune. Il n'est pas le lieu ici de faire un
historique de ce sujet, mais parmi les \'enonc\'es ayant jou\'e un
r\^ole impotrant dans ces questions il convient de mentionner l'invention
par Grothendieck de l'homologie $\ell$-adique ainsi que le faisceau de conjectures issu de la
conjecture de Shimura-Tanyama-Weil : conjecture de r\'eciprocit\'e de Langlands, conjecture de
modularit\'e de Serre, conjecture de Fontaine-Mazur. \ps

Un exemple simple
concernant notre probl\'ematique est la d\'emonstration du fait que si $S
\neq \emptyset$ alors ${\rm G}_S$ n'est pas r\'esoluble, un probl\`eme
pos\'e par Benedict Gross dans~\cite{grossconj}. Pour $p\geq 11$ cela se
d\'eduit de la construction par Deligne~\cite{deligneram} des repr\'esentations $\ell$-adiques
attach\'ees aux formes modulaires propres pour le groupe ${\rm SL}(2,\Z)$,
comme l'a d\'emontr\'e Swinnerton-Dyer~\cite{swdcong} par une m\'ethode initi\'ee par Serre
dans~\cite{serrecongram}. Le cas $S=\{p\}$ pour $p\leq 7$ n'a \'et\'e d\'emontr\'e que tr\`es recemment 
par Demb\'el\'e, Dieulefait, M. Greenberg et Voight, voir \cite{dembele2},~\cite{dgv35}
et~\cite{dieulefait7}, aussi par des
constructions automorphes. \ps

Avant d'expliquer la d\'emonstration, je voudrais mentionner que la question
de l'injectivit\'e de~\eqref{decomap} a aussi un sens g\'eom\'etrique
int\'eressant, que l'on d\'egage en r\'ealisant ${\rm G}_S$ suivant
Grothendieck comme groupe fondamental \'etale de ${\rm Spec}(\Z)\backslash
S$.  Il s'agit en effet de savoir si un "lacet" non homotope \`a $1$ dans un
voisinage infinit\'esimal du point manquant $\{p\}$ dans
$\Spec(\Z)\backslash S$ reste non homotope \`a $1$ dans $\Spec(\Z)\backslash
S$. Ce point de vue a particuli\`erement de sens dans l'analogie, presque
un guide, entre noeuds et nombres premiers, initi\'ee par Mazur dans~\cite{mazuralexander} et poursuivie
ensuite notamment par Kapranov, Resznikov et Morishita~\cite{morishita}. Cette analogie a \'et\'e enrichie encore tout r\'ecemment par Mazur
dans~\cite{mazurpo}, auquel nous renvoyons par ailleurs pour un expos\'e synth\'etique s\'eduisant de
ce th\`eme. \ps

Dans cette analogie, rappelons que $\Spec(\Z)$ correspond \`a la sph\`ere
$\mathbb{S}^3$, que les nombres premiers correspondent aux noeuds dans cette
sph\`ere, l'ensemble $S$ devenant alors un entrelac.  Si $K$ est un tel noeud, on
peut en consid\'erer un voisinage tubulaire assez petit, \`a savoir un tore
plein, et consid\'erer le bord de ce tore plein, qui est un vrai tore $T
\simeq {\mathbb S}^1 \times {\mathbb S}^1$ plong\'e dans $\mathbb{S}^3
\backslash K$.  On dispose alors d'un homomorphisme (ignorant toujours les
points bases) $$\pi_1(T) \longrightarrow \pi_1(\mathbb{S}^3\backslash K),$$
qui constitue l'analogue formel de l'homomorphisme~\eqref{decomap} dans le
cas $S=\{p\}$. Il se trouve que cet homomorphisme est injectif d\`es que
$K$ est non trivial : c'est en effet une cons\'equence bien connue du {\it
lemme de Dehn}, et plus pr\'ecis\'ement de sa g\'en\'eralisation connue sous
le nom de "th\'eor\`eme du lacet" de Papakyriakopoulos~\cite{papak}.  Comme tout indique
que les nombres premiers correspondent \`a des noeuds non triviaux, j'aime
bien penser \`a la question de Greenberg, du moins pour $S=\{p\}$, comme \`a
un "lemme de Dehn arithm\'etique".  \ps

\subsection{Id\'ee de la d\'emonstration}

L'id\'ee de la d\'emonstration est la suivante. Soit $g \in {\rm
Gal}(\overline{\Qp}/\Q_p)$ un \'el\'ement non trivial.  Observons que l'on
peut trouver un entier $n\geq 1$ et une repr\'esentation continue (d'image
finie) et irr\'eductible $$r : {\rm Gal}(\overline{\Q_p}/\Q_p) \rightarrow
{\rm GL}(n,\C)$$ telle que $r(g) \neq 1$.  En effet, $g$ est d'image
non-triviale dans un quotient fini $G$ de  ${\rm Gal}(\overline{\Q_p}/\Q_p)$.  Il agit donc non trivialement dans la
repr\'esentation r\'eguli\`ere de $G$, ainsi donc que dans l'une des
repr\'esentations irr\'eductibles de ce dernier.  \ps

On peut bien entendu voir $r$ comme une repr\'esentation continue du
groupe de Weil de $\Q_p$, et donc consid\'erer la repr\'esentation
supercuspidale $\varpi$ de $\GL(n,\Q_p)$ associ\'ee \`a $r$ par la
correspondance de Langlands locale (Harris-Taylor, Henniart). Supposons que l'on puisse
construire une repr\'esentation automorphe cuspidale $\Pi$ de $\GL_n$ sur $\Q$ telle
que : \ps

\begin{itemize}
\item[(i)] $\Pi$ est essentiellement autoduale, \ps 
\item[(ii)] $\Pi_p$ est isomorphe \`a $\varpi$, \ps
\item[(iii)] $\Pi_v$ est non ramifi\'ee si $v$ est un nombre premier non
dans $S$, \ps
\item[(iv)] $\Pi_\infty$ est alg\'ebrique r\'eguli\`ere.\footnote{La combinaison de (i) et (iv) signifie que $\Pi$ est alg\'ebrique, r\'eguli\`ere et polaris\'ee au sens du~\S\ref{repgalaut} Ch. 1.} \ps
\end{itemize}

\noindent 

Encore d'apr\`es Harris et Taylor~\cite{HT} (voir
aussi~\cite{taylortoulouse}), les hypoth\`eses
(i), (iii) et (v) assurent qu'une telle repr\'esentation admet pour tout
nombre premier $\ell$ une
r\'ealisation $\ell$-adique $$\rho_{\Pi,\ell} : {\rm Gal}(\overline{\Q}/\Q)
\rightarrow \GL(n,\overline{\Q}_\ell)$$ compatible en les premiers distincts
de $\ell$ \`a la correspondance de Langlands locale (voir~\S\ref{repgalaut} Ch. I pour plus de d\'etails). On choisit alors $\ell
\in S\backslash \{p\}$, ce qui est loisible car $|S|\geq 2$. Cela entra\^ine que
$\rho_{\Pi,\ell}$ se factorise par ${\rm G}_S$ par (iv), et que sa restriction
\`a ${\rm Gal}(\overline{\Q_p}/\Q_p)$ est isomorphe\footnote{Cela n'a de
sens bien s\^ur qu'une fois que l'on a fix\'e des plongements $\overline{\Q}
\rightarrow \C$ et $\overline{\Q} \rightarrow \overline{\Q}_\ell$, car $r$
est alors canoniquement d\'efinie sur $\overline{\Q}$ et peut \^etre vue sur
$\overline{\Q}_\ell$.} \`a $r$.  En particulier, $\rho_{\Pi,\ell}(g) \neq 1$
: ce que l'on voulait.  En passant, observons que nous avons m\^eme
d\'emontr\'e quelque chose de beaucoup plus fort que l'\'enonc\'e initial ! 
\ps

Cette d\'emonstration marche ... presque. En effet, la construction de $\Pi$
satisfaisant (ii) est en g\'en\'eral impossible, car $r^\ast$ n'est pas
n\'ecessairement isomorphe \`a une torsion de $r$ par un
caract\`ere.\footnote{Soient $p$ un nombre premier impair, $V$ un
$\F_p$-espace vectoriel de dimension finie $>1$, $u \in \GL(V)$ un
\'el\'ement d'ordre $m$ impair et ne poss\'edant pas la valeur propre $1$,
et soit $G$ une extension de $\Z/m\Z$ par $V$ donn\'ee par $u$.  Il n'est
pas difficile de construire une extension galoisienne de $\Q_p$ de groupe de
Galois isomorphe \`a un tel $G$, par exemple avec $m$ premier \`a $p$ et
$\dim(V)=m-1$.  Soit $\chi : V \rightarrow \C^\ast$ un caract\`ere d'ordre
$p$ et soit $W$ une repr\'esentation irr\'eductible de $G$ dont la
restriction \`a $V$ contient $\chi$ (Frobenius).  Alors $W^\ast$ n'est pas
une torsion de $W$ par un caract\`ere de $G$.  En effet, d'une part
$(W^\ast)_{|V}$ contient $\chi^{-1}$, qui n'est pas dans l'orbite de $\chi$
sous $\langle u \rangle$ car $m$ et $p$ sont impairs, et d'autre part
$G^{\rm ab}=\Z/m\Z$ car $1 \notin {\rm Spec}(u)$.}  J'avais contourn\'e
ce probl\`eme dans mon article~\cite{chnf}, dans lequel j'avais introduit la
strat\'egie ci-dessus, en me pla\c{c}ant au dessus d'un corps quadratique
imaginaire $E$ non ramifi\'e hors de $S$ et dans lequel $p$ est
d\'ecompos\'e.  \ps

Il s'agissait pr\'ecis\'ement de construire une repr\'esentation automorphe
cuspidale $\Pi$ de $\GL(n)$ sur $E$ telle que $\Pi^\vee \simeq \Pi^c$
(conjugu\'e ext\'erieur par $\langle c \rangle = {\rm Gal}(E/\Q)$) qui soit
cohomologique, non ramifi\'ee hors des places divisant $S$, et isomorphe \`a
$\varpi$ en l'une des deux places au dessus de $p$. Dans ce travail, je d\'emontre
l'existence d'une telle $\Pi$ en construisant tout d'abord par une simple
m\'ethode de s\'eries de Poincar\'e\footnote{Ici, on n'utilise
essentiellement que le th\'eor\`eme de Peter-Weyl car $G(\R)$ est compact!}
une repr\'esentation auxiliaire $\Pi'$ sur un groupe unitaire $G$ \`a $n$
variables relatif \`a $E$, choisi non ramifi\'e hors des premiers au dessus
de $S$, co\"incidant avec $\GL(n,\Q_p)$ au dessus de $p$, et tel que $G(\R)$
est le groupe unitaire compact (principe de Hasse). La repr\'esentation
$\Pi$ est alors obtenue comme changement de base de $\Pi'$ \`a $\GL(n)$ sur
$E$, gr\^ace \`a des r\'esultats de Clozel, Harris et
Labesse~\cite{clozellabesse},\cite{harrislabesse}.  \ps

L'argument pr\'ec\'edent conduit au r\'esultat principal de mon
article~\cite{chnf}, \`a savoir l'injectivit\'e de~\eqref{decomap} dans de
nombreux cas, ce qui en faisait les premiers cas connus. Par exemple, cela
suffit pour traiter le cas $S=\{p,\ell\}$ quand $l \equiv 3 \bmod 4$, $p>2$ et
$\left(\frac{p}{\ell}\right)=1$, en consid\'erant $E=\Q(\sqrt{-\ell})$. Cela
ne fournit cependant qu'un "quart" des cas du Th\'eor\`eme~\ref{chcl}. \ps

On peut se demander \`a ce stade d'o\`u sortent les corps de nombres peu
ramifi\'es que nous avions promis. Ils sortent bien entendu de la
construction de $\rho_{\Pi,\ell}$, c'est-\`a-dire de l'action galoisienne,
d\'efinie par Grothendieck, sur l'homologie $\ell$-adique des vari\'et\'es
de Shimura simples au sens de Kottwitz. Ce sont les quotients $\Gamma \backslash
\mathcal{B}$ de la boule unit\'e hermitienne
        $$\mathcal{B}=\{(z_i) \in \C^r, \sum_{i=1}^r |z_i|^2<1\}$$
par certains r\'eseaux cocompacts $\Gamma$ du groupe ${\rm PU}(r,1)$ des
bijections biholomorphes de $\mathcal{B}$. Il est remarquable, si l'on suit
la preuve du th\'eor\`eme de Harris et Taylor, qu'ils consid\`erent des
exemples de telles vari\'et\'es qui ont en g\'en\'eral mauvaise r\'eduction
en des places auxilliaires. Cette mauvaise r\'eduction est en fait ensuite
\'elimin\'ee dans les repr\'esentations galoisiennes par l'argument de recollement
de Blasius-Rogawski d\'ej\`a mentionn\'e au \S~\ref{reppol}. Ceci est d'autant plus remarquable que cet argument n'a pas
d'avatar praticable au sein des corps de nombres. \ps

\subsection{Constructions de repr\'esentations automorphes \`a ramification
prescrites} Comme nous l'avons d\'ej\`a dit, la m\'ethode ci-dessus ne d\'emontre pas le
Th\'eor\`eme~\ref{chcl} dans tous les cas. Je m'\'etais rendu compte que
pour y parvenir, il suffirait par le m\^eme argument de disposer du th\'eor\`eme suivant.

\begin{thm}\label{chcl2} Soit $\varpi$ une repr\'esentation
supercuspidale\footnote{Dans~\cite{chclo} nous supposons de plus pour
simplifier que $\varpi^\vee$ n'est pas isomorphe \`a une torsion non
ramifi\'ee de $\varpi$ dans cet \'enonc\'e. Pour l'application au
Th\'eor\`eme~\ref{chcl} on peut toujours se ramener \`a ce cas par torsion par un
caract\`ere ramifi\'e, et j'ignorerai ce d\'etail ici.} de $\GL(n,\Q_p)$,
de param\`etre de Langlands\footnote{Voir~\S\ref{repgalaut} Ch. I.} ${\rm L}(\varpi)$. Il existe une repr\'esentation automorphe $\Pi$ de $\GL(2n,\Q_p)$ autoduale
et cohomologique telle que : \begin{itemize}\ps
\item[(i)] le param\`etre de Langlands de $\Pi_p$ est ${\rm L}(\varpi)\otimes |\cdot |^s \oplus
{\rm L}(\varpi)^\vee \otimes |\cdot|^{-s}$ pour un
$s \in \C$, \ps
\item[(ii)] $\Pi_q$ est non ramifi\'ee si $q \notin \{p,\ell\}$.\ps
\end{itemize}
\end{thm}

Le param\`etre du (i) est trivialement auto-dual cette fois-ci. Notons que l'on ne peut pas se passer
de l'apparition du $s$ dans cet \'enonc\'e car on ne peut pas prescrire les
nombres de Weil ! Cependant, un argument \'el\'ementaire bas\'e sur la structure du quotient mod\'er\'e de ${\rm
Gal}(\overline{\Q}_p/\Q_p)$ permet de voir que cela ne trouble pas la
strat\'egie pr\'ec\'edente (voir~\cite[Lemma 4. (i)]{chnf}). La d\'emonstration du Th\'eor\`eme~\ref{chcl2} est l'objectif
principal de mon travail en commun avec Laurent Clozel~\cite{chclo}, dont je
vais donner les grandes lignes dans ce qui suit, ce qui ach\`evera donc la
d\'emonstration du Th\'eor\`eme~\ref{chcl}. \ps

Il sera commode dans cette discussion de consid\'erer temporairement le
probl\`eme plus g\'en\'eral suivant.  Fixons, $p$ et $\ell$ deux nombres
premiers distincts, $m\geq 1$ un entier et $\mathfrak{c}_p$ une composante
de Bernstein de $\GL(2m,\Q_p)$ contenant des repr\'esentations autoduales. 
Par exemple, $\mathfrak{c}_p$ peut \^etre la composante $\mathfrak{c}_p(\varpi)$ de
classe inertielle $$(\GL(m,\Q_p)\times \GL(m,\Q_p),\varpi \times
\varpi^\vee)$$ o\`u $\varpi$
est une supercuspidale de $\GL(m,\Q_p)$.  On s'int\'eresse au probl\`eme de
l'existence d'une repr\'esentation automorphe cuspidale de $\GL(2m)$ sur
$\Q$ telle que : \ps \begin{itemize} 
\item[(a)] $\Pi$ est autoduale,\ps
\item[(b)] $\Pi_p$ est dans la composante $\mathfrak{c}_p$,\ps 
\item[(c)] $\Pi_\ell$ est de carr\'e int\'egrable,\ps 
\item[(d)] $\Pi_v$ est non ramifi\'ee si $v \notin \{\ell,p,\infty\}$,\ps 
\item[(e)] $\Pi_\infty$ est cohomologique.\ps \end{itemize} \ps

Ce probl\`eme g\'en\'eral, de par sa rigidit\'e, s'est av\'er\'e
extr\^emement fin.  Une mani\`ere de le deviner est de consid\'erer comme au~\S\ref{signegalois} Ch. 1 la
repr\'esentation conjecturale du groupe de Langlands $\mathcal{L}_\Q$ de
$\Q$, disons $$\rho_\Pi : \mathcal{L}_\Q \rightarrow \GL(2m,\C),$$ associ\'e
\`a $\Pi$.  En effet, cette repr\'esentation doit \^etre autoduale par (a),
elle est donc exclusivement symplectique ou orthogonale.  En particulier,
les param\`etres des composantes locales de $\Pi$ sont soit tous
orthogonaux, soit tous symplectiques (exclusivement \`a cause de (c)), ce
qui est une contrainte notable. Cet \'enonc\'e, non conjectural, est
d'ailleurs cons\'equence des r\'esultats r\'ecents d'Arthur~\cite{arthur}. En fait, sous l'hypoth\`ese (e) le
param\`etre de Langlands de $\Pi_\infty$ ne pr\'eserve en fait aucun
accouplement orthogonal non d\'eg\'en\'er\'e, i.e.  $\rho_\Pi$ est
n\'ecessairement symplectique.  Notons que cela ne va pas \`a l'encontre de
l'\'enonc\'e du Th\'eor\`eme~\ref{chcl2}, car la composante
$\mathfrak{c}_p(\varpi)$ contient manifestement des repr\'esentations \`a la
fois symplectiques et orthogonales (de param\`etre comme dans le (i) du
Th\'eor\`eme~\ref{chcl2}).  \ps

Il convient de mentionner que le probl\`eme de construire des
repr\'esentations automorphes cuspidales autoduales pour $\GL(n)$ avec des
composantes locales prescrites a \'et\'e \'egalement consid\'er\'e
r\'ecemment par d'autres auteurs.  Par exemple, il est discut\'e par Prasad
et Ramakrishnan dans~\cite[\S 3]{prasadramakri}, auquel mes r\'esultats avec Clozel
apportent d'ailleurs certaines r\'eponses.  Il est aussi consid\'er\'e par
Khare, Larsen et Savin dans~\cite{kls} dans leur \'etude du probl\`eme de
Galois inverse.  Ces derniers expliquent comment construire un $\Pi$ comme
plus haut dans le cas o\`u $\mathfrak{c}_p$ est supercuspidale autoduale de
param\`etre symplectique (sans toutefois pouvoir assurer la condition (c)). 
Ils construisent un tel $\Pi$ par transfert \`a $\GL(2m)$ d'une
repr\'esentation automorphe cuspidale g\'en\'erique $\Pi'$ sur le
$\Q$-groupe d\'eploy\'e ${\rm SO}(2m+1)$. L'existence de ce type de
transfert repose sur une s\'erie de travaux difficiles de Cogdell, Kim,
Piatetski-Shapiro, Shahidi, Jiang et Soudry. La repr\'esentation $\Pi'$ est
construite par un argument de s\'erie de Poincar\'e~(voir
aussi~\cite{henniartgl3},~\cite{semcartan1},~\cite{semcartan2} \`a ce
sujet).\ps

Retournons maintenant \`a notre d\'emonstration, dans le cas o\`u $\mathfrak{c}_p=\mathfrak{c}_p(\varpi)$.
Notre m\'ethode pour construire la repr\'esentation $\Pi$ repose sur la
formule des traces d'Arthur-Selberg pour le groupe $\GL(2m)$ tordu par le
$\Z$-automorphisme $$\theta(g)={}^t\!g^{-1}.$$ Rappelons que cette formule
est une identit\'e de distributions ${\rm I}_{\rm spec}(\cdot)={\rm I}_{\rm
geom}(\cdot)$ dont l'utilisation pour ce type de probl\`emes est bien
connue: appliqu\'ee \`a des fonctions tests $$f = \otimes_v' f_v$$ ne
tra\c{c}ant que dans les $\Pi$ qui nous int\'eressent (ou presque), il
s'agit de montrer la non nullit\'e de son c\^ot\'e g\'eom\'etrique ${\rm
I}_{\rm geom}(f)$.  Tout le piquant \'etait alors de comprendre comment les
obstructions plus haut se manifesteraient dans l'\'etude de la formule
des traces. \ps

Si $v \notin \{\infty,\, \ell,\, p\}$, la fonction $f_v$ est naturellement la
fonction caract\'eristique de $\GL(2m,\Z_v)$.  Pour simplifier, nous
cherchons $\Pi$ telle que $\Pi_\ell$ est la repr\'esentation de Steinberg
(elle aussi de param\`etre autodual symplectique!), et on prend pour
$f_\ell$ un pseudo-coefficient tordu de cette derni\`ere.  Nous prenons pour
l'instant pour $f_p$ une fonction continue \`a support compact dont la trace
tordue s'annule hors de la composante de Bernstein $\mathfrak{c}_p$.  En la
place $\infty$, nous prenons pour $f_\infty$ un pseudo-coefficient de
s\'eries $\theta$-discr\`etes cohomologiques.  Ces repr\'esentations sont en
fait naturellement param\'etr\'ees par le poids extr\'emal $\lambda$ d'une
repr\'esentation irr\'eductible $V_\lambda$ du groupe compact
${\rm SO}(2m+1,\R)$ (car leur param\`etre est symplectique, comme on l'a d\'ej\`a
dit). Ce poids $\lambda$ est en fait essentiellement l'unique lattitude dont nous
disposons dans notre probl\`eme pour choisir la fonction $f$, que nous
noterons plut\^ot $f^\lambda$ au lieu de $f$ pour mettre en \'evidence cette
d\'ependance en $\lambda$ dans le choix de $f_\infty$.  \ps

Pour ces fonctions tests $f^\lambda$, une version simplifi\'ee de la formule
des traces due \`a Arthur \cite{ITF} s'applique, dont le c\^ot\'e
g\'eom\'etrique se r\'eduit alors aux termes port\'es par les \'elements
$\theta$-semisimples et $\Q$-elliptiques.  Nous d\'emontrerons alors
ultimement que pour un $f_p$ convenablement choisi, il existe une constante
$C \neq 0$ telle que $${\rm I}_{\rm geom}(f^\lambda) \sim C \dim(V_\lambda)$$ lorsque $\lambda$ tend vers l'infini en s'\'eloignant des
murs de Weyl pour ${\rm SO}(2m+1,\R)$. La constante $C$ est \`a un scalaire
non nul pr\`es l'int\'egrale orbitale tordue ${\rm TO}_{\gamma_0}(f_p)$ de
$f_p$ en un certain \'el\'ement $\theta$-semisimple $\Q$-elliptique $$\gamma_0 \in
\GL(2m,\Q)$$ dont le centralisateur tordu est le groupe symplectique
${\rm Sp}(2m)$ sur $\Q$.  \ps

Un premier ingr\'edient dans la d\'emonstration de l'asymptotique ci-dessus
est le fait que si $\gamma \in {\rm SO}(2m+1,\R)$ est non
trivial, alors $\frac{{\rm Trace}(\gamma,V_\lambda)}{\dim(V_\lambda)}$ tend
vers $0$ quand $\lambda$ tend s'\'eloigne des murs. Nous d\'eduisons ce fait
d'une version de la formule du caract\`ere de Weyl s'appliquant aussi aux
\'el\'ements non r\'eguliers (que nous \'etablissons pour tous les groupes de Lie compacts
connexes). Cette formule jouera d'ailleurs un r\^ole important dans mon travail ult\'erieur avec Renard
sur les formes automorphes de niveau $1$.\ps 

Le second ingr\'edient, en fait l'essentiel de la difficult\'e, consiste \`a montrer que si $f_p$ est bien choisie
alors l'int\'egrale orbitale ${\rm TO}_{\gamma_0}(f_p)$ est non nulle. Si l'on remplace
$\got{c}_p$ par la composante de Bernstein d'une repr\'esentation
supercuspidale autoduale de param\`etre symplectique, l'\'enonc\'e analogue
est un th\'eor\`eme d\^u \`a Shahidi~\cite{shahidi}, modulo une identit\'e
conjecturale entre fonctions $L$ d\'emontr\'ee par la suite par
Henniart~\cite{henniartcarresym}. L'hypoth\`ese que nous
mettons sur $f_p$ est que sa trace tordue ne prend que des valeurs $>0$ sur
les repr\'esentations temp\'er\'ees autoduales dans $\mathfrak{c}_p(\varpi)$
(th\'eor\`eme de Paley-Wiener-Rogawski). \ps

On commence par v\'erifier directement que les int\'egrales orbitales
tordues \underline{stables} de $f_p$ ne sont pas identiquement nulles.  Un
point crucial est que ceci nous permet, apr\`es stabilisation des termes
elliptiques de la formule des traces, de construire une repr\'esentation
$\Pi'$ qui a toutes les propri\'et\'es requises, sauf qu'elle est ramifi\'ee
en une place $q$ auxiliaire.  La cl\'e de l'histoire se passe ici : il n'est
pas possible de r\'eduire \`a un seul terme le c\^ot\'e g\'eom\'etrique de
la formule des traces tordue par un argument de support, car il reste
toujours au moins la classe de conjugaison stable d'un \'el\'ement.  \ps

On utilise ensuite cette repr\'esentation $\Pi'$ pour montrer la
non-annulation de ${\rm TO}_{\gamma_0}(f_p)$. Notre m\'ethode est bas\'ee sur un argument de positivit\'e
qui semble nouveau consistant \`a faire varier $\lambda$ dans la formule des traces
tordue stablilis\'ee et d'utiliser le lemme suivant. Soit $G$ un groupe de Lie compact connexe,
soient $\gamma_1,\cdots,\gamma_r$ des classes de conjugaisons non triviales
et distinctes de $G$, et soient $c_0,c_1,\cdots,c_r$ des nombres complexes. On suppose que pour tout
charact\`ere irr\'eductible $\chi$ de $G$, le nombre $c_0 \chi(1) +
\sum_{i=1}^r c_i \chi(\gamma_i)$ est un r\'eel $\geq 0$, et qu'il soit non nul
pour au moins un $\chi$. Alors $c_0$ est un r\'eel $>0$.  \ps

Nos arguments ont d'autres applications. Par exemple, Clozel en a 
trouv\'e une \`a la d\'etermination de certains signes
apparaissant dans la th\'eorie du changement base des repr\'esentations des
groupes de Lie r\'eels~\cite{projetlivregrfa}. Dans notre article, nous en donnons aussi
quelques-unes concernant l'alternative symplectique/orthogonal, dont
le th\'eor\`eme suivant. \ps

\begin{thm}\label{chclalt} Soient $F$ un corps totalement
r\'eel et $\pi$ une repr\'esentation automorphe cuspidale de $\GL(2n,\AAA_F)$. On suppose que  
$\pi$ est autoduale, de carr\'e int\'egrable modulo le centre en au moins une
place finie $v$, et cohomologique \`a toutes les places archim\'ediennes.
Alors pour toute place $v$ de $F$, le $L$-param\`etre de $\pi_v$ pr\'eserve
une forme bilin\'eaire symplectique non d\'eg\'en\'er\'ee.\ps\ps
En particulier, les repr\'esentations galoisiennes $\ell$-adiques
associ\'ees, avec $\ell \neq v$, sont aussi de similitudes symplectiques.
\end{thm}

La d\'emonstration de ce th\'eor\`eme donn\'ee dans~\cite{chclo} n'\'etait en
fait pas tout \`a fait compl\`ete car elle n\'ecessitait une propri\'et\'e
alors inconnue des int\'egrales orbitales tordues des pseudocoefficients $f_\infty$, \`a
savoir leur annulation aux \'el\'ements non semisimples elliptiques, et ceci
afin d'appliquer la formule des traces simples d'Arthur. J'ai par la
suite \'etabli cette propri\'et\'e manquante dans mon travail~\cite{chrenard1} en
commun avec David Renard. Je dois mentionner que ce th\'eor\`eme est de toutes fa\c{c}ons obsol\`ete. En effet, c'est un cas tr\`es particulier du th\'eor\`eme avec J. Bella\"iche que j'ai racont\'e au~\S\ref{signegalois} Ch. 1, m\^eme si la d\'emonstration est ici tr\`es diff\'erente. De plus, le cas particulier ci-dessus se d\'eduit aussi imm\'ediatement des r\'esultats r\'ecents d'Arthur~\cite{arthur}. \ps

\subsection{Quelques questions ouvertes} La question ouverte la plus
naturelle est de d\'emontrer l'injectivit\'e de~\eqref{decomap} dans le cas
o\`u $S=\{p\}$. Je renvoie \`a~\cite[\S 4.2,\S 4.3]{chnf} pour une
discussion des difficult\'es auxquelles on se heurte par la m\'ethode
employ\'ee ci-dessus, qui consisterait \`a utiliser des r\'ealisations
galoisiennes $p$-adiques, plut\^ot que $\ell \neq p$ adiques, des
repr\'esentations automorphes $\Pi$ construites. Disons simplement qu'il
n'est sans doute pas trop difficile de construire des $\Pi$ comme aux
paragraphes pr\'ec\'edents qui soient non ramifi\'es hors $p$ et l'infini. 
Le probl\`eme principal est plut\^ot de nature purement locale : si $\rho$ est une repr\'esentation $p$-adique
de ${\rm Gal}(\overline{\Q}_p/\Q_p)$ qui est potentiellement semi-stable, il
ne me semble pas clair en g\'en\'eral comment comparer l'extension (gigantesque) de $\Q_p$
d\'ecoup\'ee par $\rho$ et celle (pourtant petite) d\'ecoup\'ee par la
repr\'esentation de Weil-Deligne associ\'ee \`a $\rho$. \ps

Terminons par quelques questions ouvertes qui me semblent int\'eressantes. Dans ces questions, $S$ est un
ensemble fini non vide de nombres premiers, le cas le plus int\'eressant
\'etant le cas $|S|=1$. La premi\`ere, due \`a Greenberg, est toujours ouverte.\ps

\begin{question} (Greenberg) Soit $\ell$ un nombre premier qui n'est pas dans $S$, quel est le 
pro-ordre d'un \'el\'ement de Frobenius en $\ell$ (comme \'el\'ement de ${\rm G}_S$).
\end{question}

Il est tentant de penser que ce pro-ordre est divisible par tout entier, mais je ne vois pas de raison particuli\`ere \`a cela (seule la divisibilit\'e par $p^\infty$ pour $p \in S$ est imm\'ediate). Du point de vue de l'analogie entre noeuds et nombres premiers, cela dirait que deux nombres premiers dans ${\rm Spec}(\Z)$ sont autant entrelac\'es que possible. \ps

\begin{question} Soit ${\rm G}_S^{\rm res}$ le plus grand quotient
pro-r\'esoluble de ${\rm G}_S$ et soit $p \in S$. Est-ce que le morphisme ${\rm
Gal}(\overline{\Q}_p/\Q_p) \longrightarrow {\rm G}_S^{\rm res}$
induit par~\eqref{decomap} est injectif ?
\end{question}

Je pense que non, mais je ne sais pas le d\'emontrer. La question se pose car on aimerait bien savoir si il est raisonnable 
d'esp\'erer red\'emontrer le Th\'eor\`eme~\ref{chcl} sans utiliser tout l'arsenal automorphe, par exemple uniquement 
\`a l'aide de la th\'eorie du corps de classes. 
\ps

\begin{question} Est-il vrai que pour tout nombre premier $p$ il existe  un 
corps de nombres totalement r\'eel, non ramifi\'e hors de $p$, et poss\'edant au moins deux places au
dessus de $p$ ?
\end{question}

Je n'ai pas d'avis tranch\'e sur ce probl\`eme. Il est ais\'e de fabriquer
des $p$ pour lesquels cela fonctionne.  Par exemple, si $p$ est un nombre
premier de la forme $4a^3-27b^2$ avec $a,b \in \Z$ et $ p \pm 1 \neq 3a$, alors le corps de d\'ecomposition du polyn\^ome
$X^3-aX+b$ a les propri\'et\'es cherch\'ees\footnote{L'hypoth\`ese $p\pm1 \neq 3a$ permet de s'assurer de l'irr\'eductibilit\'e de ce polyn\^ome dans $\Q[X]$, i.e. d'\'eviter les polyn\^omes de la forme $(X+k)(X^2-kX+ (k^2-p)/4)$, qui sont de discriminant $p$ si  $p \pm 4$ est de la forme $9k^2$.} 
(exemples :  $p=229, 257,733, 761,1129,1229,\dots$). Il contient $\Q(\sqrt{p})$ et poss\`ede $3$ places au dessus de $p$.
 \ps

\newpage
\section[Voisins de Kneser des r\'eseaux de Niemeier]{Voisins de Kneser des r\'eseaux de Niemeier et formes automorphes pour ${\rm O}(24)$}\label{niemeierchlannes}

\subsection{R\'esum\'e et perspectives} Dans cette partie, j'expose un travail 
en commun avec Jean Lannes dont la version finale est toujours en cours de r\'edaction. Pour cette raison, 
je donnerai plus de d\'etails dans cette partie que dans les autres. Je renvoie \`a ma pr\'epublication~\cite{cl} pour des 
r\'esultats tr\`es partiels, ainsi qu'\`a mon site~\url{http://www.math.polytechnique.fr/~chenevier/niemeier/niemeier.html} 
pour un ensemble de tables. \ps

Ce travail porte sur les r\'eseaux unimodulaires pairs de rang $\leq 24$ et
particuli\`erement sur le comptage de leurs $p$-voisins au sens de Kneser. 
Nous d\'emontrons, en utilisant notamment les travaux r\'ecents d'Arthur~\cite{arthur} sur 
le spectre discret automorphe des groupes classiques, une formule explicite 
pour ces nombres en termes de certaines formes modulaires de genre $1$ et $2$. Cette formule est actuellement encore 
conditionnelle \`a certains r\'esultats de la th\'eorie d'Arthur que nous pr\'eciserons plus loin. \ps

Ces formules ouvrent un nouveau point de vue sur des probl\`emes classiques de ce sujet, 
notamment lorsqu'on les confronte \`a la th\'eorie des syst\`emes de racines dans 
ces dimensions. Elles permettent par exemple de comprendre et d'ordonner tout un ensemble de constructions 
du r\'eseau de Leech, et de d\'eterminer le graphe des $p$-voisins de Kneser des r\'eseaux de Niemeier 
pour tout nombre premier $p$ (le cas $p=2$ \'etait connu de
Borcherds~\cite{borcherdsthese}). Elles 
ont d'autres applications int\'eressantes : d\'emonstration de la conjecture de Nebe-Venkov~\cite{nebevenkov} sur les combinaisons lin\'eaires des s\'eries th\'eta des r\'eseaux de Niemeier, 
d\'etermination des valeurs propres de Hecke de certaines formes de genre $2$, 
d\'emonstration de la conjecture de Harder~\cite{harder}. \ps

Le probl\`eme du comptage des voisins des r\'eseaux de Niemeier se ram\`ene \`a
la question de d\'eterminer les repr\'esentations automorphes $\pi$ d'une certaine
$\Z$-forme $G$ semisimple du groupe sp\'ecial orthogonal euclidien ${\rm SO}(24,\R)$
qui sont telles que $\pi_\infty$ est triviale et qui sont partout non
ramifi\'ees. Les travaux de Niemeier~\cite{niemeier} montrent qu'il y a $25$ telles
repr\'esentations. Les r\'esultats r\'ecents d'Arthur proposent un angle d'attaque
pour cette question. {\it Ma motivation premi\`ere pour ce travail \'etait de comprendre
ce qu'ils disent dans ce contexte}. \ps

Je d\'ecrirai les travaux d'Arthur plus en d\'etail dans les paragraphes qui
suivent, mais je voudrais tout d'abord d\'egager de mani\`ere informelle les
int\'errogations qui me paraissaient les plus piquantes avant d'entamer ce
travail. Grosso modo, la th\'eorie d\'evelopp\'ee par Arthur associe \`a chaque repr\'esentation $\pi$ comme plus haut un "param\`etre global"
qui est une \'ecriture symbolique de la forme $$\pi_1[d_1]\oplus \pi_2[d_2]
\oplus \cdots \oplus \pi_k[d_k]$$
o\`u $1\leq k \leq 24$ est un entier, $\pi_i$ est une repr\'esentation
automorphe cuspidale "demi-alg\'ebrique" et partout non ramifi\'ee de ${\rm PGL}(r_i)$, et les $r_i,d_i$ sont
des entiers tels que $$\sum_{i=1}^k r_i d_i =24.$$ Le caract\`ere
infinit\'esimal des $\pi_i$ satisfait de plus une condition tr\`es forte que
je n'expliciterai pas ici. Ce param\`etre d\'etermine les $\pi_v$ essentiellement enti\`erement, et de mani\`ere
transparente, en termes des $(\pi_i)_v$. R\'eciproquement, la c\'el\`ebre
formule de multiplicit\'e d'Arthur donne une condition n\'ecessaire et
suffisante pour qu'un $\pi$ de param\`etre donn\'e existe. \ps

Si l'on essaie de "deviner" les param\`etres possibles, on se heurte
imm\'ediatement \`a l'examen suivant, constat d'un survol de la
litt\'erature en th\'eorie des formes automorphes : les seules
repr\'esentations automorphes demi-alg\'ebriques partout non-ramifi\'ees de
$\GL(n)$ que l'on "connaisse bien" sont essentiellement celles attach\'ees aux formes modulaires
pour ${\rm SL}(2,\Z)$ (i.e. $n=2$), et dans une moindre mesure celles attach\'ees aux formes
de Siegel de genre $g \leq 3$ (pour lesquelles il est par ailleurs extr\^emement
difficile de calculer les param\`etres de Satake). Je reviendrai plus en
d\'etail sur ces questions quand je d\'ecrirai mon travail avec Renard~\cite{chrenard2}, dont
c'est le point de d\'epart. Ainsi, si par malheur l'une des $25$ repr\'esentations $\pi$ contient 
un $\pi_i$ qui n'est pas dans la petite liste connue, alors vraisemblablement il restera inconnu
pour nous, ainsi donc que le comptage des $p$-voisins des r\'eseaux de Niemeier.
\footnote{Pour \^etre honn\^ete, je pensais initiallement pour cette raison que cette m\'ethode \'etait
vou\'ee \`a l'echec, dans le sens o\`u elle ne permettrait de d\'eterminer
qu'une petite partie des $25$ repr\'esentations cherch\'ees.} \ps

La petite liste en question, ici, est tr\`es r\'eduite : seules
les $5$ formes modulaires pour ${\rm SL}(2,\Z)$ de poids $12, 16, 18, 20$ et
$22$, les $4$ formes de Siegel de genre $2$ de poids $(6,8)$, $(4,10)$, $(8,8)$
et $(12,6)$, ainsi que le carr\'e sym\'etrique de la forme modulaire discriminant $\Delta$, ont les
propri\'et\'es de caract\`ere infinit\'esimal requises. Un premier miracle
est qu'il existe exactement $24$ param\`etres possibles, i.e. avec la bonne
condition sur le caract\`ere infinit\'esimal, que l'on peut former avec cet
ensemble de repr\'esentations. Le second miracle est que dans tous les cas
la formule de mulitplicit\'e d'Arthur donne une multiplicit\'e non nulle (et
deux dans un cas).\ps

Autant le premier miracle pouvait \^etre une co\"incidence num\'erique, autant
pour le second c'est parfaitement invraisemblable. En effet, \`a tout
param\`etre $\oplus_{i=1}^k \pi_i[d_i]$ comme plus haut Arthur associe un
$2$-groupe ab\'elien \'el\'ementaire de cardinal $2^k$ ou $2^{k-1}$. Il
associe aussi deux caract\`eres de ce groupe, l'un \'etant de nature en
g\'en\'eral globale ("caract\`ere $\varepsilon$), l'autre plut\^ot local (d\'etermin\'e par un accouplement archim\'edien). La formule de
multiplicit\'e d'Arthur dit qu'un param\`etre donnn\'e est celui d'une
repr\'esentation automorphe si et seulement si ces deux caract\`eres
co\"incident. Dans un monde al\'eatoire il y a donc \`a priori une chance
sur $2^k$ ou $2^{k-1}$ que cela soit possible, et dans les exemples en
question on a $1\leq k \leq 5$. Cependant dans tous les cas les caract\`eres
co\"incident, comme on le v\'erifie "au cas par cas" : la probabilit\'e
na\"ive de cet \'ev\`enement etait d'une chance sur $2^{59}$...  Je
reviendrai plus tard sur cette co\"incidence. \ps

Avant d'entamer une description plus pr\'ecise de ce travail, je voudrais
mentionner un difficult\'e importante si l'on veut appliquer la formule de
multiplicit\'e d'Arthur.  En effet, cette formule, bien qu'annonc\'ee par
Arthur, n'est pas encore d\'emontr\'ee.  Pire, \`a ma connaissance elle
n'est pas encore \'enonc\'ee !  Ce travail, ainsi que le travail que j'ai
men\'e de front avec Renard d\'ej\`a mentionn\'e, m'ont permis d'en deviner
une formulation sans ambigu\"it\'e dans ce contexte, et de la comprendre
plus g\'en\'eralement pour les $\Z$-groupes semisimples $G$ tels que $G(\R)$ est un
groupe compact (je la donnerai au~\S\ref{formmultarth}).  Je peux dire maintenant que
c'est une formule tr\`es simple ! En revanche, son \'etablissement reste
encore conditionnel notamment \`a certains r\'esultats annonc\'es par Arthur
dans son travail sur les formes int\'erieures des groupes
classiques~\cite[Ch.  9]{arthur}.  Il s'en suit que certains des r\'esultats
qui vont suivre seront conditionnels, nous l'indiquerons alors par une
\'etoile dans leur \'enonc\'e.  \ps

 \ps

\subsection{R\'eseaux unimodulaires pairs de petit rang} Fixons $n\geq 1$ un entier et consid\'erons l'espace euclidien 
$\R^n$, muni de son produit scalaire standard $(x_i)\cdot (y_i)=\sum_i x_i y_i$. 
On dit qu'un r\'eseau $L \subset \R^n$ est unimodulaire si il est de covolume est $1$, 
et qu'il est pair si $x \cdot x \in 2\mathbb{Z}$ pour tout $x \in L$. Un tel r\'eseau est entier, i.e. $x \cdot y \in \Z$ pour tout $x,y \in L$. \ps

Les r\'eseaux unimodulaires pairs, et la question de leur classification, apparaissent dans plusieurs 
probl\`emes math\'ematiques :  on pense aux tores iso-spectraux de Milnor, \`a la dualit\'e de Poincar\'e sur la cohomologie 
singuli\`ere en degr\'e m\'edian, aux probl\`emes d'empilements de sph\`eres 
r\'eguliers \cite{conwaysloane}, ou encore \`a la th\'eorie des formes modulaires pour ${\rm SL}(2,\Z)$, 
et plus g\'en\'eralement ${\rm Sp}(2g,\Z)$, via leurs s\'eries th\'eta. Leur classification est aussi 
\'equivalente \`a celles des formes quadratiques non d\'eg\'en\'er\'ees, et d\'efinies positives, 
sur l'anneau des entiers $\Z$, une surprenante question ouverte ! \`A chaque tel r\'eseau $L$ on 
peut en effet associer la forme quadratique $q_L : L \rightarrow \Z$ d\'efinie par 
$$q_L(x) =\frac{x \cdot x}{2}$$
qui a la propri\'et\'e que la forme bilin\'eaire sous-jacente $x \cdot y = q_L(x+y)-q_L(x)-q_L(y)$ est de d\'eterminant $1$. 
L'exemple peut-\^etre le plus 
simple est le r\'eseau de rang $8k$ $${\rm E}_{8k} = {\rm D}_{8k} + \Z e$$ o\`u
${\rm D}_n \subset \Z^n$ d\'esigne le sous-groupe des \'el\'ements $(x_i)$ tels
que $\sum_i x_i \equiv 0 \bmod 2$, et o\`u $e= \frac{1}{2}(1,1,\cdots,1)$. 
En fait, ${\rm E}_8$ est isom\'etrique au "r\'eseau des racines" de l'alg\`ebre de Lie complexe du m\^eme nom, 
muni de son unique produit scalaire invariant pour lequel les racines sont de carr\'e scalaire $2$. \ps

L'action naturelle de ${\rm O}(n,\R)$ sur $\R^n$ induit une action sur l'ensemble de ses r\'eseaux 
unimodulaires pairs. Nous d\'esignerons par $${\mathfrak X}_n$$
l'ensemble de ses orbites, i.e. des classes d'isom\'etrie de r\'eseaux unimodulaires pairs de rang $n$. La th\'eorie de la r\'eduction 
assure que $\mathfrak{X}_n$ est fini pour tout entier $n\geq 1$, et par exemple la th\'eorie de Hasse-Minkowski 
entra\^ine que $$\mathfrak{X}_n \neq \emptyset \Leftrightarrow n \equiv 0 \bmod 8.$$ \ps

Il est connu que $\mathfrak{X}_8=\{{\rm E}_8\}$ (Mordell) et $\mathfrak{X}_{16}=\{{\rm E}_{16},{\rm E}_8 \oplus {\rm E}_8\}$ 
(Witt). De plus, d'apr\`es Niemeier~\cite{niemeier} l'ensemble $\mathfrak{X}_{24}$ admet $24$ \'el\'ements le plus fameux \'etant le r\'eseau de Leech~\cite{leech}. J'en dirai plus sur cette classification ci-apr\`es.
Un r\'eseau unimodulaire pair de rang $24$ est appel\'e r\'eseau de Niemeier. La formule de masse de 
Minkowski-Siegel-Smith,  montre en revanche que $\mathfrak{X}_{32}$ a plus de $80$ millions 
d'\'el\'ements \cite{serre}, il y en a m\^eme plus d'un milliard  d'apr\`es~\cite{king}, 
et leur nombre exact semble actuellement inconnu. \ps

Le r\'eseau de Leech est, \`a isom\'etries pr\`es, le seul r\'eseau de Niemeier sans {\it racine}, 
c'est-\`a-dire sans \'el\'ement $x$ tel que $x \cdot x =2$. \`A l'oppos\'e, 
si le r\'eseau $L$, unimodulaire pair de rang $24$, poss\`ede au moins une racine, 
l'ensemble $R(L)$ de ses racines est tr\`es gros : c'est un syst\`eme de racines\footnote{Au sens usuel, par exemple au sens de Bourbaki, en particulier 
ici il engendre $\R^{24}$.} 
dans $\R^{24}$. Une d\'emonstration assez simple de cette propri\'et\'e 
a \'et\'e donn\'ee par Venkov~\cite{venkov}. Pour des raisons qui tiennent du miracle, 
l'application $L \mapsto R(L)$ induit alors une bijection entre classes d'isom\'etrie de r\'eseaux de Niemeier 
poss\'edant des racines et classes d'isomorphie de syst\`emes de racines de rang $24$ dont toutes les composantes 
irr\'eductibles sont de type ADE et ont m\^eme nombre de Coxeter.\footnote{Je rappelle que si $R$ 
est un syst\`eme de racines r\'eduit irr\'eductible de rang $\ell$, son nombre de Coxeter est 
le nombre entier $h(R)=|R|/\ell$. concr\`etement, on a $h({\rm A}_n)=n+1$, $h({\rm D}_n)=2(n-1)$ et les nombres 
de Coxeter de ${\rm E}_6$, ${\rm E}_7$ et ${\rm E}_8$ sont respectivement $12, 18$ et $30$.} C'est un exercice sans difficult\'e que de v\'erifier qu'il y a exactement $23$ tels syst\`emes de racines, list\'es dans la table~\ref{listeniemeier}.
\ps

\begin{table}[htp]
\caption{Les $23$ syst\`emes de racines \'equi-Coxeter de type ADE et de rang $24$}
\begin{tabular}{|c|c|c|c|c|c|c|c|c|c|c|c|c|}
\hline R & ${\rm D}_{24}$ & ${\rm D}_{16}{\rm E}_{8}$ & ${\rm E}_8^3$ & ${\rm A}_{24}$ & ${\rm D}_{12}^2$ & ${\rm D}_{10}{\rm E}_7^2$ & ${\rm A}_{17}{\rm E}_7$ & ${\rm A}_{15}{\rm D}_9$ & ${\rm D}_8^3$ & ${\rm A}_{12}^2$ & ${\rm A}_{11}{\rm D}_7{\rm E}_6$ & ${\rm E}_6^4$ \\
\hline h(R) & 46 & 30 & 30 & 25 & 22 & 18 & 18 & 16 & 14 & 13 & 12 & 12 \\
\hline R & ${\rm A}_9^2{\rm D}_{6}$ & ${\rm D}_{6}^4$ & ${\rm A}_8^3$ & ${\rm A}_{7}^2{\rm D}_5^2$ & ${\rm A}_{6}^4$ & ${\rm A}_{5}^4{\rm D}_4$ & ${\rm D}_4^6$ & ${\rm A}_{4}^6$ & ${\rm A}_3^8$ & ${\rm A}_2^{12}$ & ${\rm A}_1^{24}$  &  \\
\hline h(R) & 10 & 10 & 9 & 8 &  7 & 6 & 6 & 5 & 4 & 3 & 2 &  \\
\hline
\end{tabular}
\label{listeniemeier}
\end{table}
\ps

Suivant certains auteurs nous utiliserons la notation $R^+$ pour d\'esigner la classe d'isom\'etrie de 
l'unique r\'eseau de Niemeier de syst\`eme de racines $R$. Par exemple, on a ${\rm D}_{24}^+={\rm E}_{24}$, 
$({\rm E}_8^3)^+={\rm E}_8^3$ et $({\rm E_8}{\rm D}_{16})^+={\rm E}_8\oplus {\rm E}_{16}$. Pour certains $R$, la construction de $R^+$ est nettement plus subtile. \ps

concr\`etement, un syst\`eme de racines $R \subset \R^{24}$ de la liste ci-dessus \'etant donn\'e, il s'agit de trouver un r\'eseau $L$ coinc\'e entre le r\'eseau $Q(R)$ engendr\'e par $R$ et son r\'eseau dual  $Q(R)^\sharp$, 
aussi appel\'e "r\'eseau des poids" de $R$, qui soit d'une part unimodulaire pair {\it et} qui d'autre part ne contienne pas plus de racines que $R$. La premi\`ere condition \'equivaut \`a dire que l'image $I$ de $L$ dans le groupe ab\'elien fini quadratique $$q_R : Q(R)^\sharp/Q(R) \rightarrow \Q/\Z, \, \,\, \, \, x \mapsto \frac{x \cdot x}{2},$$
soit "un lagrangien", i.e. satisfasse $q_R(I) \equiv 0$ et $I=I^\bot$.  Pour qu'un tel $I$ existe, il est par exemple n\'ecessaire que $|Q(R)^\sharp/Q(R)|$, qui n'est autre que l'indice de connexion du syst\`eme
de racines $R$ au sens de Bourbaki~\cite{bourbaki}, soit un carr\'e pour les $23$ r\'eseaux $R$ de la table ci-dessus, ce qui est facile \`a v\'erifier mais constitue un premier miracle ! Venkov v\'erifie {\it au cas par cas} que pour tout $R$, il existe 
de tels lagrangiens $I$, qui de surcro\^it ne forment qu'une seule orbite sous l'action du groupe orthogonal de $Q(R)$ (le groupe $A(R)$ de Bourbaki). Le cas le plus difficile est celui de ${\rm A}_1^{24}$,  
le lagrangien recherch\'e  \'etant alors donn\'e par le code de Golay dans $(\Z/2\Z)^{24}$ ! \ps

\subsection{Voisins de Kneser et graphes $\mathfrak{X}_n(p)$} Soient $L$ et $M$ des r\'eseaux unimodulaires pairs de m\^eme rang $n$ et 
soit $p$ un nombre premier. Suivant M. Kneser, nous dirons que $L$ et $M$ sont $p$-voisins si $L\cap M$ est 
d'indice $p$ dans $L$ (et donc dans $M$). Il est ais\'e de fabriquer tous les $p$-voisins d'un $L$ donn\'e
: ils sont en bijection naturelle avec les $\F_p$-points de la quadrique
$C_L$ projective lisse sur $\Z$ d\'efinie par la forme quadratique $q_L$. \ps

concr\`etement, si $P  \in
C_L(\F_p)$, on peut trouver un \'el\'ement $v \in L \backslash pL$ de la
droite isotrope $P$ satisfaisant $\frac{v \cdot v}{2} \equiv 0 \bmod p^2$
(Hensel). Si $H$ d\'esigne le sous-r\'eseau des $x \in L$ tels que $x \cdot v \equiv 0 \bmod p$,
alors le r\'eseau $$L(P)=H + \Z \frac{v}{p}$$ est visiblement unimodulaire et pair
: c'est un $p$-voisin de $L$. Il n'est pas difficile de voir que $L(P)$ ne
d\'epend pas du choix de $v$, et que tout $p$-voisin de $L$ est de la forme $L(P)$ pour un et un seul $P \in
{\rm C}_L(\F_p)$. Observons que la quadrique $C_L$ \'etant hyperbolique sur $\Z_p$ pour
tout premier $p$, le r\'eseau $L$ admet donc exactement $$|C_L(\F_p)| = 1+p+p^2+\cdots +
p^{n-2} + p^{n/2-1}=: c_n(p)$$
$p$-voisins. 
\ps

Soient $L$ et $M$ deux r\'eseaux unimodulaires pair de rang $n$. Nous d\'esignerons par $$N_p(L,M)$$
le nombre des $p$-voisins de $L$ qui sont isom\'etriques \`a $M$. {\it C'est ce nombre myst\'erieux qui nous 
int\'eressera principalement dans la suite. Notre r\'esultat principal est d'en avoir trouv\'e 
une formule explicite quand $n\leq 24$}. Autrement dit, $L$ et $p$ \'etant donn\'e, on s'int\'eresse \`a la partition de la quadrique 
${\rm C}_L(\F_p)$ d\'efinie par la classe d'isom\'etrie du $p$-voisin associ\'e. Tentons tout d'abord de motiver cette \'etude et la notion de $p$-voisins. \ps

Tout d'abord, la relation de $p$-voisinage munit $\mathfrak{X}_n$ d'une structure de graphe d\'ependant du choix de $p$, 
que l'on notera $\mathfrak{X}_n(p)$. Kneser avait d\'emontr\'e la connexit\'e de ce graphe comme cons\'equence 
de son c\'el\`ebre th\'eor\`eme d'approximation forte. Cela fournit un proc\'ed\'e th\'eorique pour construire 
$\mathfrak{X}_n$ \`a partir du r\'eseau ${\rm E}_n$ et du choix d'un nombre premier $p$ : on part de ${\rm E}_n$, 
on d\'etermine les classes d'isom\'etrie de ses $p$-voisins, et on recommence tant que le nombre total de classes 
d'isom\'etrie construites augmente strictement. Cela a par exemple permis \`a Kneser~\cite{kneser16} de recalculer tr\`es 
simplement $\mathfrak{X}_8$ et $\mathfrak{X}_{16}$ et \`a Niemeier de d\'eterminer $\mathfrak{X}_{24}$, tous 
deux en utilisant le nombre premier $2$ et les nombreuses sym\'etries en jeu.\footnote{Par exemple, il est ais\'e de voir que si $x \in C_L(\F_p)$ et $\gamma \in {\rm O}(L)$, alors les $p$-voisins 
de $L$ associ\'es \`a $x$ et $\gamma(x)$ sont isom\'etriques. Comme l'application naturelle 
${\rm O}({\rm E}_8) \rightarrow {\rm O}({\rm E}_8 \otimes \F_2)$ est surjective, il vient que tous les $2$-voisins 
de ${\rm E}_8$ sont isom\'etriques, et en fait isom\'etriques \`a ${\rm E}_8$ lui-m\^eme comme on le v\'erifie 
facilement en calculant un $2$-voisin quelconque par exemple. Cela montre que $\mathfrak{X}_8=\{{\rm E}_8\}$. }
Il est ais\'e de voir que le graphe $\mathfrak{X}_{16}(2)$ est le graphe connexe trivial, i.e. de 
diam\`etre $1$, \`a deux sommets. Par contre le graphe $\mathfrak{X}_{24}(2)$, d\'etermin\'e par
Borcherds~\cite{conwaysloane},  
ne l'est pas du tout ! Il est de diam\`etre $5$ et cette page \url{http://en.wikipedia.org/wiki/Niemeier_lattice} 
de la Wikipedia en donne une jolie repr\'esentation (aussi due \`a Borcherds).\ps

Une seconde motivation est le probl\`eme de mieux comprendre certaines des diff\'erentes constructions du r\'eseau de Leech~\cite{conwaysloane}, en g\'en\'eral assez subtiles, qui apparaissent souvent comme des cas particuliers de constructions de $p$-voisins. La difficult\'e de construire le r\'eseau de Leech est par exemple manifeste sur le graphe ${\rm X}_{24}(2)$. En effet, 
il n'est  reli\'e dans ce dernier qu'au r\'eseau $({\rm A}_1^{24})^+$, le r\'eseau de Niemeier de la forme 
$R^+$ le plus d\'elicat \`a construire ! Cette propri\'et\'e est en fait assez simple \`a comprendre : si 
le r\'eseau de Leech est un $2$-voisin du r\'eseau $R^+$, alors un sous-groupe d'indice $2$ de ce dernier 
ne contient pas de racine. En particulier, $R$ a la propri\'et\'e que la somme de deux racines n'est jamais 
une racine : tous ses constituants irr\'eductibles sont donc de rang $1$, i.e. $R=A_1^{24}$. C'est un cas particulier du fait g\'en\'eral suivant que nous avons d\'emontr\'e, et qui est le point de d\'epart d'une s\'erie de nos
calculs bas\'es sur les propri\'et\'es du r\'eseau de Leech.

\begin{thmetoile}\label{voisinleech} Le r\'eseau de Leech est un $p$-voisin du r\'eseau $R^+$ si et seulement si $p \geq h(R)$. 
\end{thmetoile}

La condition n\'ecessaire, plus \'el\'ementaire (et inconditionnelle), est
bas\'ee sur des propri\'et\'es des groupes de Weyl affines. C'est un
analogue formel du r\'esultat, d\^u \`a Kostant~\cite{kostantsl2}, affirmant
que l'ordre minimal d'un \'el\'ement r\'egulier et d'ordre fini dans un
groupe de Lie connexe adjoint co\"incide avec le nombre de Coxeter de son
syst\`eme de racines. La condition suffisante, "construisant Leech",
n\'ecessite en revanche tout l'arsenal de notre travail. Un nombre
particuli\`erement int\'eressant est la quantit\'e $$N_p({\rm Leech},{\rm
E}_{24})$$ qui traduit la difficult\'e de construire Leech \`a partir du
plus simple des r\'eseaux de Niemeier, \`a savoir ${\rm E}_{24}$. On
constate sur ${\rm X}_{24}(2)$ qu'il faut le nombre maximal de $5$ \'etapes
pour construire le r\'eseau de Leech par voisinages \`a partir du r\'eseau
${\rm E}_{24}$ ! En revanche, on trouve dans~\cite{conwaysloane} une r\'ef\'erence
\`a l'observation suivante due \`a Thompson : il est assez facile de
construire le r\'eseau de Leech comme $47$-voisin de ${\rm E}_{24}$. En
effet, c'est le $47$-voisin associ\'e au point
$(0,1,2,\cdots,23) \in {\rm E}_{24}$. Observons simplement ici, v\'erification du pauvre, que ce point se r\'eduit 
bien en un point de la quadrique $C_{{\rm E}_{24}}(\F_{47})$ car 
$$\sum_{i=1}^{23} i^2 = \frac{1}{6} 23\cdot 24\cdot 47 \equiv 0 \bmod  47.$$ 
Ceci est bien en accord avec le th\'eor\`eme ci-dessus, car $h({\rm E}_{24})=46$ est le plus grand nombre de Coxeter intervenant, et $47$ est m\^eme le plus petit premier pour lequel une telle construction est possible ! 
\ps  
 
Un des \'enonc\'es les plus simples \`a formuler que nous obtenons avec Lannes est le suivant :

\begin{thmetoile}\label{graphep} Le graphe $\mathfrak{X}_{24}(p)$ est de diam\`etre $1$ si, et seulement si, $p\geq 47$. La liste des graphes $\mathfrak{X}_{24}(p)$ pour $p<47$ est celle donn\'ee dans~\cite{cl}.\end{thmetoile}

Par exemple, les Figures~\ref{graphex243} et~\ref{graphex247} d\'ecrivent les graphes $\mathfrak{X}_{24}(p)$ pour $p=3$ et $7$, de diam\`etres respectifs
$4$ et $2$. Les r\'eseaux de Niemeier sont num\'erot\'es par l'ordre choisi dans la table~\ref{listeniemeier}, un 
raffinement de l'ordre d\'ecroissant des nombres de Coxeter, le premier \'etant donc ${\rm E}_{\rm 24}$, 
l'avant dernier $({\rm A}_1^{24})^+$, et le r\'eseau de Leech est mis en dernier.
L'examen de la derni\`ere colonne montre par exemple bien que les deux seuls
$3$-voisins du r\'eseau de Leech sont $({\rm A}_1^{24})^+$ et $({\rm
A}_2^{12})^+$. On voit aussi que $({\rm A}_{12}^2)^+$ et $({\rm E}_6^4)^+$ ne
sont pas $3$-voisins. 

\begin{figure}[phtn]

{\tiny 

$\left[\begin{array}{cccccccccccccccccccccccc}
 1&1&0&1&1&1&0&1&0&0&0&0&0&0&0&0&0&0&0&0&0&0&0&0\cr
 1&1&1&1&1&1&1&1&1&1&1&0&1&0&0&0&0&0&0&0&0&0&0&0\cr
 0&1&1&0&0&1&1&0&1&0&1&1&0&0&1&0&0&0&0&0&0&0&0&0\cr
 1&1&0&1&1&1&0&1&0&1&1&0&1&0&1&0&0&0&0&0&0&0&0&0\cr
 1&1&0&1&1&1&1&1&1&1&1&0&1&1&1&1&0&0&0&0&0&0&0&0\cr
 1&1&1&1&1&1&1&1&1&1&1&1&1&1&1&1&1&1&0&0&0&0&0&0\cr
 0&1&1&0&1&1&1&1&1&1&1&1&1&1&1&1&1&1&0&0&0&0&0&0\cr
 1&1&0&1&1&1&1&1&1&1&1&0&1&1&1&1&1&1&0&0&0&0&0&0\cr
 0&1&1&0&1&1&1&1&1&1&1&1&1&1&1&1&1&1&1&1&0&0&0&0\cr
 0&1&0&1&1&1&1&1&1&1&1&0&1&1&1&1&1&1&0&1&0&0&0&0\cr
 0&1&1&1&1&1&1&1&1&1&1&1&1&1&1&1&1&1&1&1&1&0&0&0\cr
 0&0&1&0&0&1&1&0&1&0&1&1&1&1&1&1&1&1&1&1&0&1&0&0\cr
 0&1&0&1&1&1&1&1&1&1&1&1&1&1&1&1&1&1&1&1&1&0&0&0\cr
 0&0&0&0&1&1&1&1&1&1&1&1&1&1&1&1&1&1&1&1&1&0&0&0\cr
 0&0&1&1&1&1&1&1&1&1&1&1&1&1&1&1&1&1&1&1&1&1&0&0\cr
 0&0&0&0&1&1&1&1&1&1&1&1&1&1&1&1&1&1&1&1&1&1&0&0\cr
 0&0&0&0&0&1&1&1&1&1&1&1&1&1&1&1&1&1&1&1&1&1&0&0\cr
 0&0&0&0&0&1&1&1&1&1&1&1&1&1&1&1&1&1&1&1&1&1&1&0\cr
 0&0&0&0&0&0&0&0&1&0&1&1&1&1&1&1&1&1&1&1&1&1&1&0\cr
 0&0&0&0&0&0&0&0&1&1&1&1&1&1&1&1&1&1&1&1&1&1&1&0\cr
 0&0&0&0&0&0&0&0&0&0&1&0&1&1&1&1&1&1&1&1&1&1&1&0\cr
 0&0&0&0&0&0&0&0&0&0&0&1&0&0&1&1&1&1&1&1&1&1&1&1\cr
 0&0&0&0&0&0&0&0&0&0&0&0&0&0&0&0&0&1&1&1&1&1&1&1\cr
 0&0&0&0&0&0&0&0&0&0&0&0&0&0&0&0&0&0&0&0&0&1&1&1
\end{array}\right]$   }

\caption{La matrice d'adjacence du graphe $\mathfrak{X}_{24}(3)$}
\label{graphex243}
\bigskip\bigskip \bigskip

{\tiny 

$\left[\begin{array}{cccccccccccccccccccccccc}
 1&1&0&1&1&1&1&1&1&1&1&0&1&1&1&1&1&0&0&0&0&0&0&0\cr
 1&1&1&1&1&1&1&1&1&1&1&1&1&1&1&1&1&1&1&1&0&0&0&0\cr
 0&1&1&1&1&1&1&1&1&1&1&1&1&1&1&1&1&1&1&1&0&0&0&0\cr
 1&1&1&1&1&1&1&1&1&1&1&1&1&1&1&1&1&1&1&1&1&0&0&0\cr
 1&1&1&1&1&1&1&1&1&1&1&1&1&1&1&1&1&1&1&1&1&0&0&0\cr
 1&1&1&1&1&1&1&1&1&1&1&1&1&1&1&1&1&1&1&1&1&1&0&0\cr
 1&1&1&1&1&1&1&1&1&1&1&1&1&1&1&1&1&1&1&1&1&1&0&0\cr
 1&1&1&1&1&1&1&1&1&1&1&1&1&1&1&1&1&1&1&1&1&1&0&0\cr
 1&1&1&1&1&1&1&1&1&1&1&1&1&1&1&1&1&1&1&1&1&1&1&0\cr
 1&1&1&1&1&1&1&1&1&1&1&1&1&1&1&1&1&1&1&1&1&1&1&0\cr
 1&1&1&1&1&1&1&1&1&1&1&1&1&1&1&1&1&1&1&1&1&1&1&0\cr
 0&1&1&1&1&1&1&1&1&1&1&1&1&1&1&1&1&1&1&1&1&1&1&0\cr
 1&1&1&1&1&1&1&1&1&1&1&1&1&1&1&1&1&1&1&1&1&1&1&0\cr
 1&1&1&1&1&1&1&1&1&1&1&1&1&1&1&1&1&1&1&1&1&1&1&0\cr
 1&1&1&1&1&1&1&1&1&1&1&1&1&1&1&1&1&1&1&1&1&1&1&0\cr
 1&1&1&1&1&1&1&1&1&1&1&1&1&1&1&1&1&1&1&1&1&1&1&0\cr
 1&1&1&1&1&1&1&1&1&1&1&1&1&1&1&1&1&1&1&1&1&1&1&1\cr
 0&1&1&1&1&1&1&1&1&1&1&1&1&1&1&1&1&1&1&1&1&1&1&1\cr
 0&1&1&1&1&1&1&1&1&1&1&1&1&1&1&1&1&1&1&1&1&1&1&1\cr
 0&1&1&1&1&1&1&1&1&1&1&1&1&1&1&1&1&1&1&1&1&1&1&1\cr
 0&0&0&1&1&1&1&1&1&1&1&1&1&1&1&1&1&1&1&1&1&1&1&1\cr
 0&0&0&0&0&1&1&1&1&1&1&1&1&1&1&1&1&1&1&1&1&1&1&1\cr
 0&0&0&0&0&0&0&0&1&1&1&1&1&1&1&1&1&1&1&1&1&1&1&1\cr
 0&0&0&0&0&0&0&0&0&0&0&0&0&0&0&0&1&1&1&1&1&1&1&1
\end{array}\right]$   }

\caption{La matrice d'adjacence du graphe $\mathfrak{X}_{24}(7)$}
\label{graphex247}
\end{figure}
\ps

Je vais expliquer d'ici peu comment nous avons obtenu ces graphes ainsi que
l'id\'ee de la d\'emonstration de ces th\'eor\`emes.  Je voudrais rajouter
simplement ici qu'il semble difficile, pour calculer les $N_p(L,M)$,  d'op\'erer par calcul "brutal" \`a
l'aide de l'ordinateur (ce que ne faisait pas non plus Borcherds pour $p=2$ par
ailleurs!).  En effet, il faut 6 ms \`a PARI et mon ordinateur pour
d\'eterminer la classe d'isom\'etrie d'un $3$-voisin d'un r\'eseau de
Niemeier $L$ d\'efini par un point donn\'e de $C_L(\F_3)$ : on calcule le voisin ainsi que le nombre de ces racines, ce qui d\'etermine d\'ej\`a son nombre de Coxeter, et s'il le faut on d\'etermine l'indice de son sous-r\'eseau engendr\'e par les racines.  Il
faudrait donc au moins $$\frac{c_{24}(3)\cdot (0.005)}{3600\cdot 24 \cdot
(365.25)} \simeq 8.955$$ ann\'ees pour d\'eterminer avec leurs multiplicit\'es tous les $3$-voisins
d'un seul r\'eseau de Niemeier ! \footnote{Bien entendu, on peut r\'eduire
sensiblement ces calculs en d\'eterminant en pr\'eliminaire les orbites de
${\rm O}(L)$ sur $C_L(\F_p)$, par exemple en utilisant des propri\'et\'es
des groupes de Weyl affines.  Le point est que ces orbites sont en
g\'en\'erales tr\`es grosses : le plus petit groupe d'isom\'etrie d'un
r\'eseau de Niemeier, celui de $({\rm A}_3^8)^+$, a un cardinal de l'ordre
de $3.10^{14}$ !  Ceci dit, ceci demanderait un travail titanesque, et de
toutes fa\c{c}ons serait inaccessible \`a l'ordinateur d\`es que $p$ est
plus grand.  Vraisemblablement le cas $p=17$ me semble
inatteignable quelles que soient les ruses.}

\subsection{Cas du rang $16$} Nous nous concentrons d\'esormais sur le probl\`eme de d\'eterminer $N_p(L,M)$ pour deux r\'eseaux unimodulaires pairs de m\^eme rang $n\leq 24$.
\'Etant donn\'e que $N_p({\rm E}_8,{\rm E}_8)=c_8(p)$, le premier cas
int\'eressant est celui de la dimension $16$. Il sera commode d'introduire
l'op\'erateur $\Z$-lin\'eaire $$T_p : \Z[\mathfrak{X}_n] \longrightarrow
\Z[\mathfrak{X}_n]$$
d\'efini par $T_p([L])=\sum_M [M]$, la somme portant sur les 
$p$-voisins $M$ du r\'eseau $L$. Autrement dit, $T_p([L])=\sum_{[M] \in \mathfrak{X}_n}
N_p(L,M)[M]$. 

\begin{thm}\label{voisins16} Dans la base, $\mathrm{E}_8\oplus
\mathrm{E}_8, \mathrm{E}_{16}$ la matrice de $\mathrm{T}_p$ vaut 
$$\mathrm{c}_{16}(p) \begin{bmatrix}1 & 0 \\ 0 &
1\end{bmatrix} + (1+p+p^2+p^3)\hspace{2pt}\frac{1+p^{11}-\tau(p)}{691}
\begin{bmatrix} -405 & 286 \\ 405 & -286\end{bmatrix},$$
o\`u $\tau$ d\'esigne la fonction de Ramanujan. 
\end{thm}

Bien que ce r\'esultat \'etait sans doute connu des sp\'ecialistes, nous n'en 
avons pas trouv\'e de trace dans la litt\'erature. Je rappelle que la fonction 
$\tau$ de Ramanujan est d\'efinie par $$\Delta = q\prod_{n\geq 1}(1-q^n)^{24}=\sum_{n\geq 1} \tau(n) q^n.$$ 
L'apparition de cette fonction dans les questions de r\'eseaux n'\'etonnera sans doute pas certains lecteurs. 
En effet, il est par exemple bien connu que, si l'on d\'esigne par $$r_L(n)=|\{ x \in L, x \cdot x =
2n\}|,$$ 
alors pour tout nombre premier $p$ on a 
		$$r_{\rm Leech}(p)= \frac{65520}{691} (1+p^{11}-\tau(p)),$$
une formule d'apparence tr\`es similaire \`a celle du th\'eor\`eme ci-dessus ! \ps

Cette identit\'e r\'esulte "simplement" de ce que si $L$ est unimodulaire pair de rang $n$, alors sa s\'erie 
th\'eta $\theta_L= \sum_{n\geq 0} r_L(n) q^n$ est une forme modulaire de poids $n/2$ pour le groupe 
${\rm SL}(2,\Z)$. C'est alors un exercice que de conclure car l'espace des formes modulaires de poids 
$12$ est engendr\'e par $\Delta$ et la s\'erie d'Eisenstein de poids $12$, et car le r\'eseau de Leech 
n'a pas de racine, i.e. $r_{\rm Leech}(1)=0$. Comme nous le verrons, il semblerait que l'apparition de 
$\Delta$ dans le probl\`eme de comptage plus haut de $p$-voisin semble sensiblement plus subtile
que dans cette formule, car elle s'av\`erera \'equivalente \`a un cas non trivial de fonctorialit\'e d'Arthur-Langlands. Mentionnons enfin que la comparaison du th\'eor\`eme et de la formule
ci-dessus conduit \`a la relation "purement quadratique" 
$$N_p({\rm E}_8\oplus {\rm E}_8,{\rm E}_{16}) = \frac{9}{1456} \cdot r_{\rm Leech}(p)
\cdot \frac{p^4-1}{p-1},$$
que nous ne savons pas d\'emontrer directement ! \ps

Discutons bri\`evement la d\'emonstration du Th\'eor\`eme~\ref{voisins16}, ce qui nous permettra notamment d'introduire des objets utiles pour la suite.
Commen\c{c}ons par deux observations bien connues valables en tout rang $n$ \cite{nebevenkov}. Tout d'abord, les
op\'erateurs $T_p$ commutent deux \`a deux. De plus, ce sont des
op\'erateurs auto-adjoints pour le produit scalaire sur $\R[\mathfrak{X}_n]$ d\'efini par $\langle [L], \rangle
[M] \rangle = \delta_{[L],[M]} |{\rm O}(L)|$. Cela signifie concr\`etement
que l'on a l'identit\'e $$N_p(L,M)|O(M)|=N_p(M,L)|O(L)|.$$
Notre probl\`eme g\'en\'eral, \`a savoir d\'eterminer les $N_p(L,M)$, est
donc \'equivalent \`a trouver une base de $\R[\mathfrak{X}_n]$ constitu\'ee
de vecteurs propres commun \`a tous les $T_p$ et de d\'eterminer pour chacun
d'entre eux le syst\`eme $(\lambda_p)$ associ\'e des valeurs
propres des $(T_p)$. Mentionnons la pr\'esence importante de
$(c_n(p))$ comme syst\`eme de valeurs propres "trivial" : l'identit\'e ci-dessus assure en effet que le vecteur
$\sum_{L \in \mathfrak{X}_n} |{\rm O}(L)|^{-1} [L]$
est un vecteur propre de $T_p$ pour la valeur propre $c_n(p)$. C'est en fait le seul
vecteur propre \'evident. \ps

Retournons au cas $n=16$. Notre d\'emonstration repose sur les s\'eries th\'eta de Siegel
$$\vartheta_g : \Z[X_n] \rightarrow \mathrm{M}_{\frac{n}{2}}({\mathrm{Sp}}_{2g}(\Z)),$$
o\`u l'arriv\'ee d\'esigne l'espace des formes modulaires de Siegel de poids
$\frac{n}{2}$ pour le groupe ${\rm Sp}(2g,\Z)$ ("genre $g$"). Les relations
de commutation d'Eichler g\'en\'eralis\'ees (\cite{rallis}, \cite{walling})
assurent que $\vartheta_g$ entrelace $\mathrm{T}_p$ avec un certain
op\'erateur de Hecke explicite du c\^ot\'e des formes de Siegel. Un
r\'esultat classique de Witt, Kneser et Igusa (\cite{kneser}) assure que 
$$\vartheta_g(\mathrm{E}_8 \oplus \mathrm{E}_8)=
\vartheta_g(\mathrm{E}_{16}) \, \, \, \, {\rm
if} \, \, \, g\leq 3.$$
Ceci est bien connu en genre $1$ car ${\rm M}_8({\rm SL}(2,\Z))$ est de dimension $1$. Ces identit\'es absolument remarquables affirment que $\mathrm{E}_8 \oplus
\mathrm{E}_8$ et ${\rm E}_{16}$ sont relativement difficile \`a
diff\'erencier du point de vue de la repr\'esentation des formes\footnote{En revanche, ${\rm E}_{16}$ n'est pas engendr\'e $\Z$-lin\'eairement par ses racines, contrairement \`a ${\rm E}_8^2$.} : ils repr\'esentent toutes les formes quadratiques enti\`eres de rang
$\leq 3$ exactement le m\^eme nombre de fois ! En revanche, il n'en va pas
de m\^eme de la forme ${\rm D}_4$, de sorte que la forme modulaire $$\mathrm{F}=\vartheta_4(\mathrm{E}_8 \oplus
\mathrm{E}_8)-\vartheta_4(\mathrm{E}_{16})$$ qui est parabolique par le
r\'esultat pr\'ec\'edent, ne s'annule pas. C'est pourquoi 
la valeur propre de $T_p$ que nous recherchons est reli\'ee au syst\`eme de
valeurs propres de Hecke de la forme $$F \in \mathrm{S}_8({\rm
Sp}_8(\Z)).$$ Le lien avec la fonction $\Delta$ appara\^it alors de la mani\`ere suivante. D'une part, un r\'esultat de Poor et Yuen \cite{py} affirme que cet espace
est en fait de dimension $1$ (engendr\'e par la fameuse forme de Schottky!).
D'autre part, un autre \'el\'ement non-trivial de cet espace est le
rel\`evement d'Ikeda $I(\Delta)$ de la fonction $\Delta$~ \cite{ikeda1}, dont les
valeurs porpres de Hecke en $p$ sont reli\'ees \`a $\tau(p)$. Il en r\'esulte que $I(\Delta)$ et $F$ sont proportionnels. Le
Th\'eor\`eme~\ref{voisins16} s'en d\'eduit en mettant bout \`a bout les
formules. \ps

L'histoire ne s'arr\^ete pas tout \`a fait ici car nous avons en fait aussi trouv\'e un moyen de se passer de la d\'emonstration difficile d'Ikeda pour construire $I(\Delta)$. En effet, nous avons
observ\'e que dans ce cas tr\`es
pr\'ecis il peut \^etre construit directement \`a
partir de $\Delta$ et de deux constructions de s\'eries th\'eta et de la
trialit\'e pour le groupe ${\rm PGO}_{E_8}^+$ semisimple $\Z$, l'une \`a coefficient dans l'unique invariant harmonique de degr\'e $8$ du groupe de Weyl de ${\rm E}_8$, l'autre \`a coefficient dans la 
repr\'esentation triale de cette derni\`ere. \ps\medskip

\subsection{Cas du rang $24$ et conjecture de Nebe-Venkov}

Nous en venons \`a notre r\'esultat principal, \`a savoir la d\'etermination
des $N_p(L,M)$ pour $L$ et $M$ deux r\'eseaux de Niemeier. Nous aurons
besoin pour le formuler d'introduire certains nouveaux acteurs. Rappelons que si $k \in \{12,16,18,20,22\}$,
l'espace $${\rm S}_k({\rm SL}(2,\Z))$$ des formes modulaires paraboliques de poids $k$ pour
le groupe ${\rm SL}(2,\Z)$ est de dimension $1$. Nous d\'esignerons alors par
$$\Delta_k= \sum_{n\geq 1} \tau_k(n) q^n$$ l'unique \'el\'ement de cet espace tel
que $\tau_k(1)=1$, en particulier $\Delta=\Delta_{12}$. Rappelons aussi que d'apr\`es Tsushima~\cite{tsushima}, si $(j,k) \in \{(6,8), (4,10), (8,8),
(12,6)\}$, l'espace des formes modulaires vectorielles de Siegel de poids
${\rm Sym}^j \otimes \det^k$ pour le groupe ${\rm Sp}(4,\Z)$ est aussi de
dimension $1$. Nous d\'esignerons par $$\Delta_{j,k}$$ un g\'en\'erateur de cet espace et par  $\tau_{j,k}(p)$ la valeur propre de
l'op\'erateur de Hecke de Siegel en $p$ correspondant \`a la "trace enti\`ere du
Frobenius" agissant sur $\Delta_{j,k}$. 

\begin{thmetoile}\label{voisins24} Soient $L$ et $M$ des r\'eseaux de Niemeier. Il existe $11$ polyn\^omes $P, (P_k), (P_{j,k}), P_S \in \Q[X]$ tels que pour tout nombre premier $p$,  
$$N_p(L,M)= P(p)+\sum_k P_k(p) \tau_k(p) + \sum_{j,k} P_{j,k}(p)
\tau_{j,k}(p) + P_S(p) (\tau(p)^2-p^{11}).$$
\end{thmetoile}

Les $11 \cdot 24^2$ polyn\^omes en questions sont explicites : nous renvoyons \`a la table 
\url{http://www.math.polytechnique.fr/~chenevier/niemeier/explicit_formula.pdf} pour une liste. Nous \'epargnerons au lecteur la vue ici d'une de 
ces b\^etes car leurs coefficients, typiquement des nombres comme celui-ci 
$$9241391269/6125466240,$$
sont d'apparence insignifiante. En fait, \`a la vue des formules explicites il est tout \`a fait extraordinaire que $N_p(L,M)$ soit un entier ! 
C'est une analyse d\'etaill\'ee de ces formules, utilisant notamment les estim\'ees de Ramanujan sur les $\tau_\ast(p)$, qui nous permet de d\'emontrer les Th\'eor\`emes~\ref{voisinleech} et~\ref{graphep}. \ps

Notre objectif est d\'esormais d'expliquer d'o\`u sortent ces formules, ce
qui nous permettra aussi de les pr\'esenter sous une forme moins barbare. Tout comme dans 
le cas $n=16$, il s'agit d'\'etudier les syst\`emes de valeurs propres des $T_p$ agissant 
sur $\R[\mathfrak{X}_{24}]$. Pour $p=2$, ainsi que l'ont observ\'e Nebe et Venkov~\cite{nebevenkov}, 
la matrice de $T_2$ se d\'eduit des calculs de Borcherds. Par un heureux hasard, cette matrice a des 
valeurs propres distinctes et enti\`eres, et permet donc de d\'eterminer les vecteurs propres communs des $T_p$. Pour d\'eterminer les valeurs propres, on peut \'etudier les applications $\vartheta_g$ comme dans le cas $n=16$. On est alors confront\'e \`a deux probl\`emes. \ps

Tout d'abord, il faut d\'eterminer le noyau des applications $$\vartheta_g : \Z[\mathfrak{X}_{24}] \longrightarrow {\rm M}_{12}({\rm Sp}(2g,\Z)).$$
C'est une suite d\'ecroissante avec $g$ de sous-espaces, dont Erokhin~\cite{erokhin} avait prouv\'e la nullit\'e quand $g=12$ (c'est tautologique pour $g=24$). Ce probl\`eme est l'analogue pour les r\'eseaux de Niemeier  du probl\`eme mentionn\'e en rang $16$, qui \'etait de d\'eterminer le plus petit genre $g$ diff\'erentiant les s\'eries th\'eta de ${\rm E}_{\rm 16}$ et ${\rm E}_8^2$. Ce probl\`eme a \'et\'e \'etudi\'e en d\'etail par Nebe et Venkov dans \cite{nebevenkov}, o\`u ils determinent un certain nombre de ces relations, utilisant notamment des travaux de 
Borcherds-Freitag-Weissauer~\cite{bfw}. Ils proposent aussi une conjecture pr\'ecise, que nous appellerons {\it la conjecture de Nebe-Venkov}. \ps

Admettant ceci, il s'agirait ensuite de d\'eterminer les valeurs propres des op\'erateurs de Hecke du c\^ot\'e des formes de Siegel sur les images des $\vartheta_g$. Nous nous sommes rendu compte tardivement que ceci avait \'et\'e entrepris par Ikeda~\cite{ikeda2}, qui est arrive arriv\'e \`a d\'eterminer $20$ des $24$ syst\`emes de valeurs propres en combinant de mani\`ere astucieuse deux types de rel\`evements de formes de Siegel qu'il a lui-m\^eme construit (dans le jargon, ce sont les "Ikeda lifts" et autres "Miyawaki lifts"~\cite{ikeda1},~\cite{ikeda2}). Ce travail formidable est malheureusement encore insuffisant pour l'application ci-dessus. Notamment, les $4$ formes myst\'erieuses $\Delta_{j,k}$ que nous avons introduites plus haut n'apparaissent pas dans ses formules. \ps

Nous proc\'edons en fait de mani\`ere compl\`etement diff\'erente, sans \'etudier de s\'eries th\'eta, mais en passant plut\^ot par la description r\'ecente par Arthur du spectre discret automorphe 
des groupes classiques. En retour, nous obtiendrons le :

\begin{thmetoile}\label{conjnv} La conjecture de Nebe-Venkov~\cite{nebevenkov} est vraie. \end{thmetoile}

\noindent Le lien avec les formes automorphes, assez standard, est le suivant. D\'esignons par $${\rm O}(n)$$ le sch\'ema en groupes orthogonaux sur $\Z$ d\'efini par le r\'eseau 
${\rm E}_{n}$ muni de sa forme quadratique $q_{{\rm E}_n}$. L'ensemble $\mathfrak{X}_{n}$ s'identifie alors canoniquement \`a  l'ensemble des classes de ${\rm O}(n)$, i.e. 
$$\mathfrak{X}_n = {\rm O}(n,\Q) \backslash {\rm O}(n,\AAA_f)/ {\rm O}(n,\widehat{\Z}).$$
Cet \'enonc\'e tr\`es classique n'est autre que le fait que les r\'eseaux unimodulaires pairs de rang $n$ sont isom\'etriques sur $\Z_p$ pour tout $p$ (en fait, hyperboliques). Bien entendu, $\AAA_f$ d\'esigne l'anneau des ad\`eles finies de $\Q$.\ps

 L'espace des fonctions $\mathfrak{X}_n \rightarrow \C$ s'identifie alors \`a celui des formes automorphes de ${\rm O}(n)$ qui sont de conducteur $1$ et qui engendrent la repr\'esentation triviale sous le groupe compact ${\rm O}(n,\R)$. L'op\'erateur
 $T_p$ est natuellement un op\'erateur de Hecke de cette th\'eorie. Autrement dit, on se ram\`ene \`a d\'eterminer les $24$ repr\'esentations automorphes $\pi$ de ${\rm O}(24)$ telles que $\pi_\infty$ est trivial et $\pi_p$ non ramifi\'e pour tout $p$, i.e. admet des vecteurs invariants par ${\rm O}(24,\Z_p)$. \ps

Il sera en fait l\'eg\`erement plus commode de raisonner plut\^ot sur le groupe ${\rm SO}(24)$ (qui est connexe). Le groupe ${\rm SO}(n)$ est d\'efini comme \'etant le noyau du d\'eterminant de Dieudonn\'e ${\rm O}(n) \rightarrow \Z/2\Z$,  il est semisimple sur $\Z$ au sens usuel (de SGA 3). Le nombre de classes de ${\rm SO}(24)$ est $25$, et non plus $24$,  car seul le r\'eseau de Leech n'admet pas d'isom\'etrie de d\'eterminant $-1$.

\subsection{Digression : la classification d'Arthur dans le cas des groupes
classiques semisimples sur $\Z$}\label{pararthur} Ne serait-ce que pour
formuler nos r\'esultats, il est n\'ecessaire \`a ce point de faire une
digression sur les \'enonc\'es d'Arthur concernant le spectre automorphe
discret des groupes classiques.  \ps

Ces \'enonc\'es sont discut\'es en toute g\'en\'eralit\'e 
dans~\cite{arthurunipotent},~\cite{arthurlivre} et~\cite{arthur}, du moins dans le cas 
des groupes quasi-d\'eploy\'es. Ils se simplifient notablement dans le cas 
particulier des groupes poss\'edant un mod\`ele entier qui est semisimple sur $\Z$ et 
de leurs repr\'esentations partout non ramifi\'ees, qui est le cadre que
nous allons imposer ici, et que nous retrouverons quand je d\'ecrirai mes travaux avec Renard. 
Je renvoie d'ailleurs \`a mon article~\cite[\S 3]{chrenard2} pour une discussion plus d\'etaill\'ee de tout ce
qui suit. Bien que j'aurais aim\'e agr\'ement\'e ce paragraphe de plus
d'exemples, il m'a sembl\'e n\'ecessaire pour rester concis d'aller
temporairement "droit au but" ! \ps

Dans tout ce qui suit, on supposera que $G$ est un $\Z$-groupe semisimple.
Autrement dit, c'est un sch\'ema en groupes lin\'eaire sur $\Z$ qui est plat
et \`a fibre semisimple sur $\overline{\F}_p$ pour tout nombre premier $p$. Par exemple, 
$G$ peut \^etre un groupe de Chevalley, ou encore le groupe ${\rm SO}(n)$ pour $n \equiv
0 \bmod 8$ introduit plus haut. Ainsi que l'a observ\'e Gross, le dual de Langlands de $G_\Q$ admet alors une action 
triviale du groupe de Galois absolu de $\Q$~\cite{grossinv}. On peut donc simplement le voir comme
un groupe alg\'ebrique semisimple complexe $\widehat{G}$ dont la donn\'ee
radicielle est en dualit\'e avec celle de $G(\C)$. \ps

D\'esignons par $\Pi(G)$ l'ensemble des repr\'esentations complexes $\pi=\bigotimes_v' \pi_v$ admissibles 
irr\'eductibles de $G(\AAA)$ telles que $\pi_p$ est non ramifi\'ee pour tout premier $p$, i.e. $\pi_p^{G(\Z_p)} \neq 0$.  L'objectif ultime est de d\'ecrire ses sous-ensembles $$\Pi_{\rm cusp}(G) \subset \Pi_{\rm disc}(G) \subset
\Pi(G)$$
constitu\'es respectivement des repr\'esentations  qui
sont automorphes cuspidales ou discr\`etes.\footnote{On rappelle, suivant Weil, que 
l'espace homog\`ene $G(\Q)\backslash G(\AAA)$ est muni d'une mesure de Radon
invariante \`a droite sous l'action de $G(\AAA)$, unique \`a un scalaire
$>0$ pr\`es. Elle est de masse totale finie (Borel, Harish-Chandra).
L'espace des formes automorphes de carr\'e int\'egrable de $G$ est l'espace
$\mathcal{L}(G)={\rm L}^2(G(\Q) \backslash G(\AAA))$, qui est muni d'une repr\'esentation
unitaire de $G(\AAA)$ par translations \`a droite. On d\'esigne par
$\mathcal{L}_{\rm disc}(G)$ le sous-espace de $\mathcal{L}(G)$ obtenu en
prenant l'adh\'erence des sous-repr\'esentations irr\'eductibles ferm\'ees.
Suivant Gelfand et Piatetski-Shapiro, le sous-espace $\mathcal{L}_{\rm
cusp}(G)$ constitu\'e des formes cuspidales est $G(\AAA)$-stable et inclus dans
$\mathcal{L}_{\rm disc}(G)$. Par d\'efinition, $\pi \in \Pi(G)$ est dite
automorphe cuspidale (resp. discr\`ete) si elle apparait comme sous-repr\'esentation de
$\mathcal{L}_{\rm cusp}(G)$ (resp. $\mathcal{L}_{\rm disc}(G)$).} 
Selon Langlands et Arthur, pour chaque $\Z$-groupe semisimple $G$, $\Pi_{\rm disc}(G)$ devrait \^etre 
"reconstitu\'e" selon une recette
bien pr\'ecise \`a partir des $\Pi_{\rm cusp}({\rm PGL}(n))$, pour $n\geq
1$ variable. \ps

Soit $\pi \in \Pi(G)$. Pour tout nombre premier $p$, la th\'eorie de Satake associe \`a $\pi_p$ une 
classe de conjugaison semisimple $$c_p(\pi) \subset \widehat{G},$$ 
qui la d\'etermine enti\`erement. De mani\`ere similaire, le caract\`ere infinit\'esimal de $\pi_\infty$ 
d\'efinit selon Harish-Chandra une classe de conjugaison semisimple 
	$$c_\infty(\pi) \subset \widehat{\mathfrak{g}}={\rm Lie}_\C \widehat{G},$$
qui d\'etermine (le module de Harish-Chandra) $\pi_\infty$ parmi un ensemble fini.
\ps

Si $H$ est un groupe alg\'ebrique lin\'eaire complexe, il sera commode
d'introduire l'ensemble $\mathcal{X}(H)$
des collections $(c_v)$, $v$ parcourant les places de $\Q$, telles que $c_p$ (resp. $c_\infty$) 
est une classe de conjugaison semisimple dans $H$ (resp. ${\rm Lie}(H)$). On a d\'ecrit ci-dessus 
une application de param\'etrisation \`a fibres finies $$c : \Pi(G) \rightarrow \mathcal{X}(\widehat{G}), \, \, \pi \mapsto
(c_v(\pi)),$$
qui n'est autre que la param\'etrisation de Langlands (affaiblie \`a l'infini). Quand $G(\R)$ est compact, 
il est n\'ecessairement connexe (Chevalley), et cette application $c$ est m\^eme une injection. \ps

Concentrons nous maintenant sur le cas des groupes {\it classiques}. Dans ce qui
suit un groupe complexe $H$ sera dit classique s'il est isomorphe \`a 
${\rm Sp}(n,\C)$ pour un entier $n\geq 2$ pair, ou \`a ${\rm SO}(n,\C)$ pour un entier $n\geq 1$ et $n\neq 2$. Un tel groupe est semisimple. 
Il est uniquement d\'etermin\'e par son "type", \`a savoir le couple
$$(s(H),n(H))$$
o\`u $s(H)$ vaut $-1$ si $H$ est un groupe symplectique et $1$ sinon, et o\`u $n(H)$ est
l'entier $n$ ci-dessus, i.e. la dimension de la repr\'esentation naturelle,
aussi dite {\it standard}, du groupe $H$. Un $\Z$-groupe
semisimple $G$ sera dit classique si $G(\C)$ l'est, auquel cas on constate
que $\widehat{G}$ l'est aussi, \'eventuellement de type diff\'erent en g\'en\'eral.\ps
Les exemples principaux\footnote{La classification sur $\Q$ de ces groupes est bien connue : 
ils ont pour mod\`eles entiers les ${\rm Sp}(2g)$ ou ${\rm SO}_L$, o\`u $q : L
\rightarrow \Z$ est une forme quadratique dont la forme bilin\'eaire associ\'ee $(x,y) \mapsto q(x+y)-q(x)-q(y)$
est de d\'eterminant $\pm 1$ en rang pair, $\pm 2$ en rang impair. L'unique obstruction \`a
l'existence d'une telle forme est la congruence $p-q \equiv -1,0,1 \bmod 8$
o\`u $(p,q)$ est la signature de $q \otimes \R$. }
 dans nos applications seront le groupe ${\rm Sp}(2g)$
et certains $\Z$-groupes sp\'eciaux orthogonaux, dont les ${\rm SO}(n)$. 
Si $G$ est un $\Z$-groupe semisimple classique, et si $N$ d\'esigne l'entier $n(\widehat{G})$, la repr\'esentation standard de $\widehat{G}$
d\'efinit une application $${\rm St} : \mathcal{X}(\widehat{G}) \rightarrow
\mathcal{X}(SL(N,\C))=\mathcal{X}(\widehat{\PGL(N)}).$$
La th\'eorie d'Arthur et Langlands permet de d\'ecrire l'image de
$c(\Pi_{\rm disc}(G))$ par ${\rm St}$, nous allons rappeler comment. \ps

Soit $\Pi_{\rm cusp}^\bot(\PGL(n))$ le sous-ensemble de
$\Pi_{\rm cusp}(\PGL(n))$ constitu\'e des $\pi$ isomorphes \`a leur
duale $\pi^\vee$ (contragr\'ediente). Un premier r\'esultat fondamental,
sugg\'er\'e par le formalisme de Langlands et d\'emontr\'e par Arthur, est que pour tout $\pi \in \Pi_{\rm
cusp}^\bot(\PGL(n))$, il existe un unique $\Z$-groupe de
Chevalley\footnote{Le groupe non semisimple $\mathbb{G}_m={\rm SO}(1,1)$ est
ici exclu.} classique
$G^\pi$ tel que $n(\widehat{G^\pi})=n$, et tel qu'il existe un $\pi' \in
\Pi_{\rm disc}(G^\pi)$ satisfaisant
	$${\rm St} ( c (\pi') ) = c(\pi).$$
On dit que $\pi$ est symplectique si $\widehat{G^\pi}$ est un groupe
symplectique, et sinon que $\pi$ est orthogonal. L'alternative
symplectique/orthogonal peut aussi \^etre lue sur les fonctions $L$ des
carr\'es symm\'etriques ou altern\'es de $\pi$, mais cela ne nous sera pas
utile ici. \'Evidemment, si $n$ est impair alors toute $\pi \in \Pi_{\rm cusp}(\PGL(n))$ est orthogonale. \ps

Nous arrivons maintenant \`a une d\'efinition importante pour la suite. 
Soit $G$ un $\Z$-groupe semisimple classique, nous appellerons {\it
param\`etre d'Arthur global} de $G$ un quadruplet
$$(k,(n_i),(d_i),(\pi_i))$$ 
tels que :\ps 
 \begin{itemize} 
\item[(i)] $k$ est un entier tel que $1 \leq k \leq n(\widehat{G})$, \ps

\item[(ii)] pour tout entier $1\leq i \leq k$ alors $n_i$ est un entier
$\geq 1$, et de plus $\sum_{i=1}^k n_i = n(\widehat{G})$,\ps

\item[(iii)] pour tout $1\leq i \leq k$, $d_i$ est un diviseur de $n_i$, $\pi_i \in \Pi_{\rm
cusp}^\bot(\PGL(n_i/d_i))$, et $s(\widehat{G^\pi_i})(-1)^{d_i+1}=s(\widehat{G})$,\ps

\item[(iv)] Si $(n_i,d_i)=(n_j,d_j)$ pour $i\neq j$, alors $\pi_i \neq
\pi_j$. \ps

\end{itemize}\bigskip

Nous d\'esignerons par $\Psi_{\rm glob}(G)$ l'ensemble des param\`etres d'Arthur
globaux de $G$. Il ne d\'epend que de $s(\widehat{G})$ et $n(\widehat{G})$. 
Deux param\`etres $(k,(n_i),(d_i),(\pi_i))$ et $(k',(n'_i),(d'_i),(\pi'_i))$ 
sont dits \'equivalents si $k=k'$ et s'il existe une permutation $\sigma \in \got{S}_k$
telle que pour tout $i$, $(n'_i,d'_i,\pi'_i)=(n_{\sigma(i)},d_{\sigma(i)},\pi_{\sigma(i)})$.
La classe d'\'equivalence du param\`etre $(k,(n_i),(d_i),(\pi_i))$ sera en
g\'en\'erale not\'ee symboliquement $$\pi_1[d_1]\oplus \cdots \oplus \pi_r[d_r]$$
et le sous-symbole $\pi_i[d_i]$ sera remplac\'e par $[d_i]$ si $n_i=d_i$ (i.e.
$\pi_i=1$), et par $\pi_i$ si $d_i=1$ et $n_i>1$. \ps

On dispose d'une application de param\'etrisation $$c : \Psi_{\rm glob}(G)
\longrightarrow \mathcal{X}({\rm SL}(n(\widehat{G}),\C))$$
d\'efinie comme suit. Soit $\psi=(k,(n_i),(d_i),(\pi_i)) \in \Psi(G)$. Les
hypoth\`eses sur les $n_i$ et $d_i$ assurent l'existence d'une repr\'esentation
$$\rho_\psi : {\rm SL}(2,\C) \times \prod_{i=1}^k {\rm SL}(n_i/d_i,\C) \longrightarrow {\rm SL}(n(\widehat{G}),\C)$$
obtenue en consid\'erant la somme directe des repr\'esentations $\nu_{d_i}
\otimes \C^{n_i/d_i}$ de $\SL(2,\C) \times {\rm SL}(n_i/d_i,\C)$, o\`u $\nu_d$ 
d\'esigne la $\C$-repr\'esentation irr\'eductible de dimension $d$ de ${\rm SL}(2,\C)$.
Soit $$e \in \mathcal{X}({\rm SL}(2,\C))$$
l'\'el\'ement d'Arthur, i.e. le param\`etre de la repr\'esentation triviale
de $\PGL(2)$. Rappelons que $e=(e_v)$ avec $e_\infty={\rm
diag}(\frac{1}{2},-\frac{1}{2})$, et pour tout nombre premier $p$, $$e_p={\rm
diag}(p^{\frac{1}{2}},p^{-\frac{1}{2}})$$ o\`u $p^{\frac{1}{2}}$ d\'esigne la racine carr\'e positive de
$p$. Un r\'esultat essentiel, aussi d\^u \`a Arthur, est le suivant.

\begin{thm} (Arthur) Pour tout $\pi \in
\Pi_{\rm disc}(G)$, il existe un param\`etre global $\psi(\pi)=(k,(n_i),(d_i),(\pi_i)) \in \Psi_{\rm
glob}(G)$, unique \`a \'equivalence pr\`es, tel que $${\rm
St}(c(\pi))=\rho_{\psi(\pi)} ( e \times \prod_{i=1}^k c(\pi_i)).$$
\end{thm}

La formule de multiplicit\'e Arthur d\'ecrit aussi l'image de $\pi \mapsto
\psi(\pi)$, nous y reviendrons au~\S~\ref{formmultarth}. Nous pouvons d'ores et d\'ej\`a
reformuler le Th\'eor\`eme~\ref{voisins24} en terme de formes automorphes.

\subsection{Les $25$ repr\'esentations $\pi \in \Pi_{\rm disc}({\rm SO}(24))$ telles que $\pi_\infty$ soit la repr\'esentation triviale}\label{disc24} Il est bien connu qu'une forme modulaire $F \in {\rm S}_k({\rm
SL}(2,\Z))$ qui est vecteur propre des op\'erateurs de Hecke ${\rm T}_p$
usuels engendre de mani\`ere naturelle une unique repr\'esentation $\pi_F \in \Pi_{\rm
cusp}(\PGL(2))=\Pi_{\rm cusp}^\bot(\PGL(2))$. Elle a la propri\'et\'e caract\'eristique que 
$p^{\frac{k-1}{2}}{\rm trace}(c_p(\Pi_F))$ co\"incide avec la valeur propre de ${\rm T}_p$
sur $F$, et elle v\'erifie de plus $c_\infty(\Pi_F)={\rm
diag}(\frac{k-1}{2},-\frac{k-1}{2})$. Les $\pi_F$ sont trivialement
symplectiques au sens d'Arthur, car $G^{\pi_F}={\rm SO}(1,1)=\mathbb{G}_m$ est exclus.
Pour ne pas alourdir les notations, on \'ecrira en g\'en\'eral $F$ pour la
repr\'esentation $\pi_F$. Par exemple, les $\Delta_k$ introduits
plus haut peuvent \^etre vus comme des \'el\'ements symplectiques de $\Pi_{\rm
cusp}^\bot(\PGL(2))$. \ps

Plus g\'en\'eralement, consid\'erons l'espace vectoriel ${\rm S}_\rho({\rm Sp}(2g,\Z))$
des formes modulaires paraboliques de Siegel pour le groupe ${\rm Sp}(2g,\Z)$ et
\`a coefficients dans la repr\'esentation $\rho$ de ${\rm GL}_g(\C)$. Si $F
\in {\rm S}_\rho({\rm Sp}(2g,\Z))$ est un vecteur propre pour tous les
op\'erateurs de Hecke, il engendre aussi une unique repr\'esentation $\pi_F
\in \Pi_{\rm cusp}({\rm PGSp}(2g))$. Supposons maintenant $g=2$, alors ${\rm PGSp}_4$ est
$\Z$-isomorphe au groupe de Chevalley ${\rm SO}(3,2)$, de groupe dual ${\rm Sp}(4,\C)$, et $\pi_F$ admet donc
un param\`etre global $\psi(\pi_F) \in \Psi({\rm SO}(3,2))$. L'op\'erateur de Hecke not\'e ${\rm T}_p$ sur ${\rm S}_\rho({\rm Sp}(2g,\Z))$ a la propri\'et\'e que sa valeur propre sur $F$ co\"incide avec la trace 
de la classe de conjugaison $c_p(\pi_F) \subset {\rm Sp}(4,\C)$ multipli\'ee par une puissance de $p$ explicite.\footnote{Si $\rho={\rm Sym}^j \otimes \det^k$ cette puissance est $p^{\frac{j+2k-1}{2}}$.} Les travaux
d'Arthur permettent de d\'emontrer  que si $F$ n'est pas une forme de Saito-Kurokawa, alors $\psi(\pi_F)=(k,(n_i),(d_i),(\pi_i))$ 
est cuspidal, i.e. $k=1=d_1$ et $\psi(\pi_F)=\pi_1 \in \Pi_{\rm cusp}^\bot(\PGL(4))$ (par d\'efinition, symplectique). On notera en g\'en\'eral par abus simplement $F$ la repr\'esentation automorphe $\pi_1$ ainsi d\'efinie.
Observons que comme les formes de Saito-Kurokawa sont \`a valeurs scalaires, les $4$ formes $\Delta_{j,k}$ introduites plus haut ne sont pas de Saito-Kurokawa. \ps

La derni\`ere forme automorphe dont nous aurons besoin est le carr\'e sym\'etrique de la repr\'esentation $\Delta$, qui est une repr\'esentation (orthogonale) not\'ee ${\rm Sym}^2 \Delta$ de $\Pi_{\rm cusp}(\PGL(3))$ d\'efinie par Gelbart et Jacquet. \ps

\begin{thmetoile}\label{24niemeier} Les $25$ repr\'esentations automorphes $\pi$ de ${\rm SO}(24)$
qui sont non ramifi\'ees \`a toutes les places finies, et telles que $\pi_\infty$ est
la repr\'esentation triviale, ont pour param\`etres globaux ceux de la
table~\ref{table24}. \ps
Seul $\Delta[12]$ est le param\`etre de deux telles repr\'esentations, qui
sont par ailleurs conjugu\'ees l'une de l'autre sous l'action de ${\rm O}(24,\Z)$. 
\end{thmetoile}

\begin{table}[h]
\caption{Les $24$ param\`etres des $\pi \in \Pi_{\rm disc}({\rm SO}(24))$ telles que $\pi_\infty$ est triviale. }

{\tiny
\begin{multicols}{3}

$$[23]\oplus [1],$$
\ps
$${\rm Sym}^2 \Delta \oplus [21],$$
\ps
$$\Delta_{22}[2] \oplus [1] \oplus [19],$$
\ps
$${\rm Sym}^2 \Delta \oplus \Delta_{20}[2] \oplus [17],$$
\ps
$$\Delta_{22}[2] \oplus \Delta_{18}[2] \oplus [1] \oplus [15],$$
\ps
$$\Delta_{20}[4] \oplus [1] \oplus [15],$$
\ps
$${\rm Sym}^2 \Delta \oplus \Delta_{20}[2] \oplus \Delta_{16}[2]
\oplus [13],$$
\ps
$${\rm Sym}^2 \Delta \oplus \Delta_{18}[4] \oplus [13],$$
\ps
$$\Delta_{18}[6] \oplus [1] \oplus [11],$$
\ps
$$\Delta_{22}[2] \oplus \Delta_{16}[4] \oplus [1] \oplus[11],$$
\ps
$$\Delta_{12,6}[2] \oplus \Delta_{18}[2] \oplus [1] \oplus [11],$$
\ps
$${\rm Sym}^2 \Delta \oplus \Delta_{20}[2] \oplus \Delta_{16}[2]
\oplus \Delta[2] \oplus [9],$$
\ps
$${\rm Sym}^2 \Delta \oplus \Delta_{18}[4] \oplus \Delta[2] \oplus [9],$$
\ps
$${\rm Sym}^2 \Delta \oplus \Delta_{16}[6] \oplus [9],$$
\ps
$$\Delta_{16}[8] \oplus [1] \oplus [7],$$
\ps
$$\Delta_{22}[2] \oplus \Delta_{18}[2] \oplus \Delta[4] \oplus [1] \oplus [7],$$
\ps
$$\Delta_{20}[4] \oplus \Delta[4] \oplus [1] \oplus [7],$$
\ps
$$\Delta_{8,8}[2] \oplus \Delta_{16}[4] \oplus [1] \oplus [7],$$
\ps
$${\rm Sym}^2 \Delta \oplus \Delta_{20}[2] \oplus \Delta[6] \oplus [5],$$
\ps
$${\rm Sym}^2 \Delta \oplus \Delta_{6,8}[2] \oplus \Delta_{16}[2] \oplus
\Delta[2] \oplus [5],$$
\ps
$$\Delta_{22}[2] \oplus \Delta[8] \oplus [1] \oplus [3],$$
\ps
$$\Delta_{4,10}[2] \oplus \Delta_{18}[2] \oplus \Delta[4] \oplus [1] \oplus
[3],$$
\ps
$${\rm Sym}^2 \Delta \oplus \Delta[10] \oplus [1],$$
\ps
$$\Delta[12].$$
\end{multicols}
}
\label{table24}
\end{table}

Le lecteur v\'erifiera sans peine que les $24$ param\`etres en question sont bien dans $\Psi_{\rm glob}({\rm SO}(24))$. Ils poss\`edent une propri\'et\'e suppl\'ementaire fondamentale que nous allons d\'egager maintenant. 
Si $G$ est un $\Z$-groupe semisimple classique, et si $\psi=(k,(n_i),(d_i),(\pi_i))$ est un param\`etre global de $G$, appelons {\it caract\`ere infinit\'esimal} de $\psi$ la classe de conjugaison semisimple 
$$z_\psi:=\rho_\psi( e \times \prod_{i=1}^k c_\infty(\pi_i)) \subset M_N(\C)$$
o\`u $N=n(\widehat{G})$. Si $G(\R)$ est compact, et si $\psi=\psi(\pi)$ pour un $\pi \in \Pi_{\rm disc}(G)$, il suit que $z_\psi$ est l'image par ${\rm St}$ du caract\`ere infinit\'esimal de $\pi_\infty$. Par exemple, si cette repr\'esentation est triviale, c'est simplement la classe de conjugaison de la demi-somme des co-racines positives de $\widehat{G}$ relativement \`a une paire $(B,T)$ quelconque. Concr\`etement, si $\widehat{G}={\rm SO}(24,\C)$ c'est la classe de conjugaison semisimple ayant pour valeurs propres $$11,10, \cdots,3,2,1,0,0,-1,-2,-3,\cdots,-10,-11.$$
Le lecteur aura besoin de conna\^itre le caract\`ere infinit\'esimal de $(\pi_{\Delta_{j,k}})_\infty$ : c'est la classe de conjugaison semisimple de $\mathfrak{sp}(4,\C)$ de valeurs propres $$\pm (\frac{j+1}{2}+k-2), \pm \frac{j+1}{2}.$$
Cette relation barbare montre d'ailleurs que la param\'etrisation par le couple $(j,k)$ n'est pas tr\`es heureuse pour ces questions, et lorsque nous ferons des choses combinatoirement plus compliqu\'ees il sera primordial de tout param\'etrer par le caract\`ere infinit\'esimal. \ps

En guise d'exemple, le caract\`ere infinit\'esimal de 
$\Delta_{4,10}[2] \oplus \Delta_{18}[2] \oplus \Delta[4] \oplus [1] \oplus [3]$ est la classe de conjugaison {\small
$$\diag(\frac{21}{2},\frac{5}{2},-\frac{5}{2},-\frac{21}{2})\otimes \diag(\frac{1}{2},-\frac{1}{2}) \oplus $$ $$\diag(\frac{17}{2},-\frac{17}{2}) \otimes \diag(\frac{1}{2},-\frac{1}{2}) \oplus 
\diag(\frac{11}{2},-\frac{11}{2}) \otimes \diag(\frac{3}{2},\frac{1}{2},-\frac{1}{2},-\frac{3}{2})\oplus$$ $$\diag(0) \oplus \diag(1,0,-1),$$}
qui a bien les valeurs propres requises. \ps

Le Th\'eor\`eme~\ref{24niemeier} se d\'eduit de la formule de multiplicit\'e d'Arthur, que nous expliquerons en d\'etail dans ce cadre au~\S~\ref{formmultarth}. 

\subsection{Cons\'equences de la liste}  Il n'est pas difficile ici de passer de ${\rm SO}(24)$ \`a ${\rm O}(24)$.
concr\`etement, d\'esignons par $\widetilde{\mathfrak{X}_n}$ l'ensemble des
classes d'isom\'etrie de r\'eseaux unimodulaires pairs orient\'es de $\R^n$,
il est muni d'une projection canonique vers $\mathfrak{X}_n$.  Si $M$ est un
$p$-voisin de $L$, l'orientation de $L$ d\'efinit une orientation canonique
sur $M$, de sorte que l'op\'erateur de Hecke ${\rm T}_p$ sur
$\Z[\mathfrak{X}_n]$ se raffine en un op\'erateur $\widetilde{{\rm T}}_p$ de
$\Z[\widetilde{\mathfrak{X}_n}]$.  On v\'erifie que cet op\'erateur correspond
\`a $p^{\frac{n}{2}-1}$ fois la trace du param\`etre de Satake.  C'est aussi
un op\'erateur de Hecke minuscule au sens de Gross~\cite{grossatake}.  La projection
$\Q[\widetilde{\mathfrak{X}_{n}}] \rightarrow \Q[\mathfrak{X}_{n}]$
entrelace $\widetilde{{\rm T}_p}$ et ${\rm T}_p$.  \ps

Quand $n=24$, seul le r\'eseau de Leech admet deux orientations non \'equivalentes comme on l'a d\'ej\`a dit, de sorte que la projection canonique $$\Q[\widetilde{\mathfrak{X}_{24}}] \rightarrow \Q[\mathfrak{X}_{24}]$$
a son noyau de dimension $1$, engendr\'e par $w=[{\rm Leech}^+]-[{\rm Leech}^-]$. Il n'est pas tr\`es difficile de voir, en \'etudiant l'op\'erateur de changement d'orientation et un calcul pour $p=2$, que le syst\`eme de vecteurs propres des $\widetilde{{\rm T}}_p$ sur $w$ se retrouve comme syst\`eme de valeurs propres des ${\rm T}_p$ sur $\Q[\mathfrak{X}_{24}]$. Ce ph\'enom\`ene est reli\'e \`a l'existence de deux repr\'esentations de param\`etre $\Delta[12]$ ou encore \`a la construction de Borcherds-Freitag-Weissauer~\cite{bfw}. Il en r\'esulte que les syst\`emes de valeurs propres des ${\rm T}_p$ sur $\Q[\mathfrak{X}_{24}]$ sont simplement les traces des param\`etres de Satake 
associ\'es aux $24$ param\`etres de la liste ci-dessus, multipli\'es par $p^{11}$. \ps

Le probl\`eme imm\'ediat qui se pr\'esente alors est que nous ne connaissons pas de formule simple donnant les param\`etres de Satake des $\Delta_{j,k}$, contrairement au cas des $\Delta_k$. C'est en fait un probl\`eme bien connu et assez non trivial. Par chance, pour les $4$ formes en question ils ont \'et\'e r\'ecemment calcul\'es par C. Faber et G. van der Geer pour des petits nombres premiers $p$. Par exemple, leur trace est donn\'ee 
dans~\cite{vdg} pour $p\leq 7$. Leur m\'ethode consiste \`a \'enum\'erer des courbes de genre $2$ sur les corps finis et d'exploiter le lien existant entre l'espace de module des telles courbes et la vari\'et\'e de Siegel en genre $2$. Ils doivent notamment \'etudier la cohomologie des syst\`emes locaux ${\rm Sym}^j \otimes \det^k$ sur cette derni\`ere et leur appliquer la formule de Grothendieck-Lefschetz. En utilisant leurs valeurs num\'eriques pour $p=2$, nous avons pu ainsi v\'erifier que les $24$ valeurs propres de ${\rm T}_2$ obtenues sont bien celles de la matrice de Borcherds-Nebe-Venkov !  \ps

Cette v\'erification constituait pour nous un indice tr\`es fort que la liste~\ref{table24} ci-dessus, obtenue initialement par t\^atonnements, \'etait correcte. Il faut dire qu'\`a ce point de notre travail j'\'etais incapable de comprendre quelle devrait \^etre la forme exacte de la formule de multiplicit\'e annonc\'ee par Arthur dans ce contexte. Je reviendrai sur cette formule au~\S\ref{formmultarth} plus loin. Je voudrais d'abord signaler qu'\`a ce stade nous disposons, pour chaque couple $(L,M)$ de r\'eseaux de Niemeier, d'une expression de $N_p(L,M)$ comme combinaison lin\'eaire explicite \`a coefficients rationnels des quatre inconnues $$\tau_{6,8}(p), \tau_{4,10}(p), \tau_{8,8}(p), \tau_{12,6}(p) \in \Z.$$
Nous avons alors fait l'oservation suivante : admettons que nous sachions pour un nombre premier $p$ donn\'e d\'eterminer $N_p(L,M)$ pour $4$ paires distinctes de r\'eseaux de Niemeier $\{L,M\}$. Alors tr\`es probablement nous pouvons inverser ce syst\`eme et en d\'eduire les $4$ nombres ci-dessus, ainsi donc  en retour tous les $N_p(L,M)$ pour ce nombre premier $p$. Nous sommes ici sauv\'es par le r\'eseau Leech. En effet, il r\'esulte du sens "facile" du Th\'eor\`eme~\ref{voisinleech} que pour $p\leq 23$ et $R \in \{{\rm D}_{24}, {\rm E}_8^3, {\rm E}_8\oplus {\rm D}_{16}, {\rm A}_{24}\}$ 
alors $$N_p({\rm Leech},R^+)=0.$$
On d\'eduit de ceci, essentiellement sans effort, toutes les valeurs des $\tau_{j,k}(p)$ pour $p\leq 23$ : voir la table~\ref{taujk}. Elles co\"incident, dans tous les cas o\`u nous avons pu les comparer, avec les valeurs trouv\'ees par Faber et van der Geer. En peaufinant notre m\'ethode, nous avons en fait r\'eussi, toujours gr\^ace au r\'eseau de Leech, \`a \'etendre cette table jusqu'\`a $p=71$, alors que les travaux de Faber et van der Geer ne donnent des tables que jusqu'\`a $p=37$. On en d\'eduit tous les ${\rm N}_p(L,M)$ pour $p\leq 79$. 

\begin{table}
\caption{Quelques valeurs propres d'op\'erateurs de Hecke en genre $2$}
\label{taujk}

\renewcommand{\arraystretch}{1.5}

\begin{tabular}{|c||c|c|c|c|}
\hline $p$   & $\tau_{6,8}(p)$ & $\tau_{8,8}(p)$ &
$\tau_{12,6}(p)$ &
$\tau_{4,10}(p)$ \\
\hline $2$ & $0$ & $1344$&$-240$&$-1680$\cr
\hline $3$ & $-27000$&$-6408$&$68040$&$55080$\cr
\hline $5$ & $2843100$&$-30774900$&$14765100$&$-7338900$\cr
\hline $7$ & $-107822000$&$451366384$&$-334972400$&$609422800$\cr
\hline $11$ & $3760397784$&$13030789224$&$3580209624$&$25358200824$\cr
\hline $13$ & $9952079500$&$-328006712228$&$91151149180$&$-263384451140$\cr
\hline $17$ &$243132070500$&$5520456217764$&$-11025016477020$&$-2146704955740$\cr
\hline $19$ &
$595569231400$&$-28220918878760$&$-22060913325080$&$43021727413960$\cr
\hline $23$ &
$-6848349930000$&$79689608755152$&$195863810691120$&$-233610984201360$\cr
\hline $29$ &
$53451678149100$&$-1105748270340$&$-1743496339579620$&$-545371828324260$\cr
\hline $31$ & $234734887975744$&$1851264166857664$&$
1979302106496064$&$830680103136064$\cr
\hline $37$ & $448712646713500$ & $22115741387845324$ & $-3685951226317460$ & $11555498201265580$\cr
\hline $41$ & $-1267141915544076$ & $-29442241674311916$ & $106065086529460884$ & $-56208480716702316$\cr
\hline $43$ & $-1828093644641000$ & $308109789751260712$ & $74859021001125400$ & $160336767963955000$\cr
\hline $47$ & $-6797312934516000$ & $43932618784857504$ & $156108802652634720$ & $-116311331328502560$\cr
\hline
\end{tabular}

\end{table}
\ps \bigskip

Une amusette : l'op\'erateur de Hecke ${\rm T}_3$ admet en fait deux valeurs propres \'egales sur $\Z[\mathfrak{X}_{24}]$, contrairement \`a ${\rm T}_2$, \`a savoir la valeur $1827360$. Cela concerne les formes $$\Delta_{4,10}[2] \oplus \Delta_{18}[2] \oplus \Delta[4] \oplus [1] \oplus
[3] \, \, \, \, \, {\rm et}\, \, \, \, \, {\rm Sym}^2 \Delta \oplus \Delta_{6,8}[2] \oplus \Delta_{16}[2] \oplus
\Delta[2] \oplus [5].$$
Je ne parviens pas \`a expliquer cette co\"incidence autrement que par le calcul.

\ps

L'application \`a la conjecture de Nebe-Venkov, i.e. le Th\'eor\`eme~\ref{conjnv}, est une cons\'equence imm\'ediate de la Table~\ref{table24} et de la g\'en\'eralisation par Rallis des relations de commutation d'Eichler~\cite{rallis}, et de l'unicit\'e des param\`etres globaux d'Arthur 
(due en fait \`a Jacquet-Shalika~\cite{jasha}). \ps

Je voudrais terminer ce paragraphe par une derni\`ere application concernant la conjecture de Harder dans~\cite{harder}. Cette conjecture, \'enonc\'ee ci-dessous, est une version en dimension sup\'erieure de la c\'el\`ebre congruence de Ramanujan $\tau(p) \equiv 1+ p^{11} \bmod 691$, ou encore des congruences de Ribet dans sa preuve de la r\'eciproque du th\'eor\`eme de Herbrand~\cite{ribet}. \`A ma connaissance, c'est le premier cas "d\'emontr\'e" des conjectures de Harder. 

\begin{thmetoile}\label{harder} {\rm (Conjecture d'Harder)} Pour tout nombre premier $p$, $$\tau_{4,10}(p) \equiv \tau_{22}(p) + p^8 + p^{13} \bmod 41.$$ \end{thmetoile}

Je me permets de donner la d\'emonstration de ce r\'esultat, qui est tr\`es simple. Consid\'erons les deux repr\'esentations
automorphes de ${\rm SO}(24)$ dont les param\`etres globaux sont
$$\Delta_{4,10}[2] \oplus \Delta_{18}[2] \oplus \Delta[4] \oplus [1] \oplus
[3] \, \, \, \, \, {\rm et} \, \, \, \, \,\Delta_{22}[2] \oplus
\Delta_{18}[2] \oplus \Delta[4] \oplus [1] \oplus [7].$$ Soient $e, f \in
\Z[\mathfrak{X}_{24}]$ des vecteurs propres de ${\rm T}_2$ associ\'es qui
soient de plus primitifs dans le $\Z$-module $\Z[\mathfrak{X}_{24}]$.  Ces
vecteurs sont explicites par le calcul de ${\rm T}_2$ par
Borcherds-Nebe-Venkov.  L'ordinateur nous affirme que le groupe ab\'elien
$$\Z e \oplus \Z f$$ est d'indice exactement $1968=2^4 \cdot 3 \cdot 41$
dans son satur\'e.  Il en r\'esulte que les valeurs propres de ${\rm T}_p$
sur $e$ et $f$ co\"incident modulo $1968$.  Mais on observe sur les
param\`etres que la diff\'erence entre ces deux valeurs propres n'est autre
que l'entier $$ (1+p) \cdot \left[\tau_{4,10}(p) - (\tau_{22}(p)+ p^8 +
p^{13})\right],$$ de sorte que cet entier est $\equiv 0 \bmod 1968$. 
On en d\'eduit la congruence de Harder si $p \not \equiv-1 \bmod 41$.
\ps

Dans le cas g\'en\'eral, on travaille dans l'anneau de Grothendieck $A$ des repr\'esentations
de dimension finie de ${\rm Gal}(\overline{\Q}/\Q)$ \`a coefficients dans
$\overline{\F}_{41}$. Soient les repr\'esentations semisimples $\rho_{22}$ et
$\rho_{4,10}$ \`a coefficients dans $\overline{\F}_{41}$
respectivement associ\'ees \`a $\Delta_{22}$ (Deligne) et
$\Delta_{4,10}$ (Weissauer). Si $\omega$ d\'esigne le caract\`ere
cyclotomique modulo $41$, le th\'eor\`eme de Cebotarev et la congruence
pr\'ec\'edente entra\^inent l'\'egalit\'e dans $A$ 
$$(1+ \omega)\cdot \rho_{4,10} = (1+\omega)\cdot (\rho_{22} + 1 + \omega).$$
C'est un simple exercice alors de voir que l'on peut bien simplifier par
$1+\omega$ cette \'egalit\'e, c'est m\^eme formel \`a partir du fait que l'ordre de
$\omega$ est $>8$ (or il est d'ordre $40$). $\square$ 

\newpage

\section{Motifs de conducteurs $1$ et quelques probl\`emes
associ\'es}\label{cond1chrenard}

\subsection{R\'esum\'e et perspectives} Dans cette partie, j'expose mon
travail~\cite{chrenard2} en collaboration avec David Renard. Je renvoie \`a mon
site~\url{http://www.math.polytechnique.fr/~chenevier/levelone.html} pour de nombreuses tables de r\'esultats. \ps

Les succ\`es des travaux pr\'ec\'edents avec Lannes, notamment la question de la d\'etermination des
$24$ repr\'esentations automorphes de niveau $1$ et coefficient trivial pour
le groupe ${\rm O}(24)$ gr\^ace aux travaux d'Arthur, m'ont permis d'entrevoir de nouveaux horizons sur le 
probl\`eme de l'\'enum\'eration, voire de la classification, des formes
automorphes alg\'ebriques de niveau $1$ pour les groupes classiques sur
$\Z$. J'entendrai par l\`a les $\Z$-formes "semisimples sur $\Z$" des groupes de Lie
r\'eels ${\rm SO}(p,q)$ et ${\rm Sp}(2g,\R)$, par exemples les groupes de
Chevalley correspondants. C'est un probl\`eme \'evidemment tr\`es ancien,
dont l'int\'er\^et n'est plus \`a d\'emontrer, qui a \'et\'e \'etudi\'e par de nombreux auteurs, majoritairement dans le cas du 
groupe ${\rm Sp}(2g,\R)$ et de la th\'eorie des formes de modulaires
Siegel. \ps

La question traditionnelle est de d\'eterminer la dimension de l'espace ${\rm
S}_V({\rm Sp}(2g,\Z))$ des formes
modulaires de Siegel paraboliques pour le groupe ${\rm Sp}(2g,\Z)$ et \`a
coefficients
vectoriels $V$ quelconque (voir par exemple~\cite{vdg}). L'approche classique,
d'apparence inextricable quand le genre $g$ grandit, consiste \`a appliquer le th\'eor\`eme de 
Riemann-Roch au faisceau coh\'erent correspondant sur la vari\'et\'e de
Siegel, ce qui requiert non seulement un calcul de caract\`eristique d'Euler
non trivial mais aussi la d\'etermination de groupes de cohomologie annexes. 
Le cas $g=1$ est bien connu, et appara\^it au moins dans le cours d'arithm\'etique de Serre : la dimension de ${\rm S}_k({\rm
SL}(2,\Z))$ est $[k/12]$, \`a moins que $k>2$ et $k\equiv 2 \bmod 12$ auquel
cas c'est $[k/12]-1$. Le cas du genre $g=2$ est d\^u \`a
Igusa dans le cas des formes "scalaires", et \`a Tsushima~\cite{tsushima} en g\'en\'eral.
La formule de Tsushima, bien qu'explicite, est suffisament obsc\`ene pour \^etre \'epargn\'ee ici au lecteur. Ce n'est que tr\`es r\'ecemment que 
le cas $g=3$ a \'et\'e obtenu par J. Bergstr\"om, C. Faber  et
G. van der Geer~\cite{BFVdG} (partiellement conditionnellement), le cas des formes scalaires remontant \`a
Tsuyumine~\cite{tsuyumine}. \ps

Un probl\`eme tout \`a fait reli\'e \`a ces questions d'apr\`es Langlands est celui
d'\'enum\'erer les motifs purs sur $\Q$ ayant bonne r\'eduction partout, par
exemple en fonction de leurs poids de Hodge.\footnote{Dans le contexte pr\'esent, on se limitera
majoritairement dans notre \'etude aux motifs polaris\'es et dont les nombres de Hodge $h_{p,q}$ sont $0$ ou $1$.} C'est un probl\`eme
classique qui me fascine depuis
longtemps, et qui est par exemple dans l'esprit de mes travaux sur le
probl\`eme de la construction de corps de nombres peu ramifi\'es. Il a
plusieurs facettes qui sont reli\'ees par des conjectures "standards". La forme la plus tractable est de classifier les
repr\'esentations automorphes cuspidales de ${\rm GL}(n)$ sur $\Q$ qui sont non ramifi\'ees \`a toutes les places
finies et qui sont alg\'ebriques \`a l'infini. Ces repr\'esentations sont en bijection
conjecturale naturelle avec les motifs purs sur $\Q$ absolument simples
de rang $n$ qui ont bonne r\'eduction partout (Langlands), avec les
repr\'esentations $\ell$-adiques 
de ${\rm Gal}(\overline{\Q}/\Q)$ irr\'eductibles de rang $n$ qui sont
cristallines en $\ell$ et non ramifi\'ees
hors $\ell$ (Fontaine-Mazur~\cite{fontainemazur}), ou encore avec un sous-ensemble ad\'equat de fonctions $L$ qui sont dans
la classe de Selberg~\cite{selbergclass}. Il convient de mentionner qu'il est connu depuis
Harish-Chandra qu'il n'existe qu'un nombre fini de repr\'esentations
automorphes cuspidales de $\GL(n)$ sur $\Q$ qui soient non ramifi\'ees \`a
toutes les places finies et de caract\`ere infinit\'esimal donn\'e, ce qui
donne un sens au probl\`eme de comptage. En revanche, le comptage dans les
autres mondes semble tr\`es difficile actuellement, surtout dans celui des
motifs o\`u il est ouvert m\^eme en rang $1$ ! \ps

Un exemple embl\'ematique de cette interaction conjecturale est la
conjecture par Shaffarevich, d\'emontr\'ee ensuite ind\'ependamment par
Fontaine et Abrashkin, du fait que la seule courbe projective lisse sur $\Z$
est la droite projective, et plus g\'en\'eralement qu'il n'y a pas de
vari\'et\'e ab\'elienne sur $\Z$ de dimension $>0$.  Admettant que la
fonction $L$ d'une telle vari\'et\'e a les propri\'et\'es attendues, cela
avait \'et\'e v\'erifi\'e par Serre par une m\'ethode inspir\'ee des
formules explicites de Weil~\cite{mestre}.  Un autre exemple classique est
le probl\`eme tr\`es particulier de l'existence de $\pi$ de niveau $1$ pour
$\GL(n)$ telles que $\pi_\infty$ a le caract\`ere infinit\'esimal de la
repr\'esentation triviale.  Il est en effet \'equivalent de d\'eterminer les
entiers $n$ pour lesquels la cohomologie parabolique du groupe ${\rm
SL}(n,\Z)$ \`a coefficients dans $\Q$ est non triviale.  Les formules
explicites de Weil montrent que $n$ doit \^etre $\geq 27$
(\cite{fermigier},\cite{miller}) mais on ne le sait pour aucun autre entier
$n$, un surprenant probl\`eme ouvert.  Comme l'a remarqu\'e Khare
dans~\cite{khareniveau1}, ce probl\`eme est reli\'e aux g\'en\'eralisations
en dimensions sup\'erieures de la conjecture de modularit\'e de Serre.  \ps

Le point de d\'epart de mon travail avec Renard est que les travaux
r\'ecents d'Arthur permettent un point de vue tout nouveau sur ces
questions.  Grosso-modo, l'id\'ee est que la formule de multiplicit\'e
d'Arthur permet th\'eoriquement d'exprimer de mani\`ere tr\`es pr\'ecise les
espaces de formes automorphes alg\'ebriques pour les groupes classiques sur
$\Z$ en fonction de certaines "briques de base".  Cela peut s'utiliser dans
les deux sens : soit pour d\'eterminer ces briques de base, soit pour
reconstruire l'espace des formes en question si l'on conna\^it ces briques
de bases.  Ces briques sont les repr\'esentations automorphes
$\pi$ des $\GL(n)$ pour $n \geq 1$ variable, telles que : \medskip

\begin{itemize}
\item[(i)] (autodualit\'e) $\pi^\vee \simeq \pi$,\ps
\item[(ii)] (niveau $1$) $\pi_p$ est non ramifi\'ee pour tout nombre premier $p$,\ps
\item[(iii)] $\pi_\infty$ est demi-alg\'ebrique r\'eguli\`ere.\medskip
\end{itemize}

Pr\'ecisons cette derni\`ere condition. Si $w \geq 0$ est un entier d\'esignons par ${\rm I}_w$ la repr\'esentation de dimension $2$ du groupe
de Weil ${\rm W}_\R$ de $\R$ qui est 
induite du caract\`ere $$z \mapsto
(z/\overline{z})^{\frac{w}{2}}:=z^w/|z|^w$$ de son
sous-groupe ${\rm W}_{\C}=\C^\ast$ qui est d'indice $2$. La condition (iii)
signifie qu'il existe $[n/2]$ entiers $w_1 > w_2 >
\cdots > w_{[n/2]} \geq 0$ tels que le param\`etre de Langlands ${\rm
L}(\pi_\infty)$ de $\pi_\infty$ ait la propri\'et\'e suivante :  \ps
\begin{itemize}
\item[(iii, $n$ pair)] ${\rm L}(\pi_\infty) \simeq \bigoplus_{i=1}^{n/2}
{\rm I}_{w_i}$, \ps
\item[(iii, $n$ impair)]  ${\rm L}(\pi_\infty) \simeq \chi \oplus
\bigoplus_{i=1}^{[n/2]} {\rm I}_{w_i}$, o\`u $\chi$ est d'ordre $2$  
et $w_{[n/2]} >0$.  \ps
\end{itemize}

Les entiers $w_i$ sont alors uniquement d\'etermin\'es par $\pi$ et seront
appel\'es les poids de Hodge de $\pi$. L'alternative symplectique/orthogonal
d'Arthur assure que tous les $w_i$ sont congrus modulo $2$ : ils sont pairs
si $\pi$ est orthogonale et impairs sinon. L'entier $w(\pi):=w_1$ joue un r\^ole
particulier et sera appel\'e poids motivique\footnote{Au sens du~\S\ref{repgalaut} Ch. I, mentionnons \'egalement que la repr\'esentation $\pi|\cdot|^{-w(\pi)/2}$ est alg\'ebrique.} de $\pi$. Ces d\'enominations
sont naturelles car le $\Q$-motif de rang $n$ (\`a coefficients dans
$\overline{\Q}$ et simple) conjecturalement associ\'e \`a $\pi|\cdot|^{-w(\pi)/2}$ est effectif, pure de poids
$w(\pi)$, et a ses $h_{p,w(\pi)-p}$ non nuls (et \'egaux \`a $1$) quand $p$
parcourt les $\frac{\pm w_i+w(\pi)}{2}$,
$i=1,\dots,[n/2]$, avec en plus $\frac{w(\pi)}{2}$ si $n$ est impair. Observons enfin que le caract\`ere central de $\pi$ est n\'ecessairement
trivial, ce qui entra\^ine notamment que $\chi = \varepsilon_{\C/\R}^{[n/2]}$
quand $n$ est impair. Le r\'esultat principal de mon travail avec Renard est le suivant.

\begin{thmetoile}\label{chrenard} Pour tout entier $n\leq 8$, et toute suite 
$w_1 > w_2 > \dots > w_{[n/2]}$, il existe une formule explicite et
impl\'ement\'ee sur ordinateur\footnote{Il faut au plus une dizaine de minutes \`a
mon ordinateur pour \'evaluer cette formule pour un $r$-uple
$(w_1,\cdots,w_r)$ donn\'e satisfaisant $w_1 <100$.}, donnant le nombre de repr\'esentations automorphes
cuspidales $\pi$ de $\GL(n)$ satisfaisant les conditions (i), (ii) et (iii)
plus haut, de poids de Hodge les $w_i$, et qui sont de plus orthogonales (resp.
symplectiques).
\end{thmetoile}

L'\'etoile qui d\'ecore le terme th\'eor\`eme ci-dessus signifie que ce r\'esultat est conditionnel \`a certains r\'esultats non encore d\'emontr\'es concernant la classification d'Arthur pour 
les formes int\'erieures des groupes classiques~\cite[Chap. 9]{arthur}. Nous discuterons certaines de ces hypoth\`eses au~\S\ref{formmultarth}.\ps

Je renvoie \`a notre article et \`a ma page web~\url{http://www.math.polytechnique.fr/~chenevier/levelone.html}
pour un ensemble de tables. Dans ces tables j'utilise la notation suivante : si
$w_1>\cdots>w_{r}>$ sont des entiers impairs je d\'esigne par
$$S(w_1,\cdots,w_r)$$
le nombre des $\pi$ symplectiques de $\GL(2r)$ satisfaisant (i), (ii) et (iii) de poids
de Hodge les $w_i$. Par exemple, $S(w)=\dim {\rm S}_{w+1}({\rm SL}(2,\Z))$. 
Si $w_1>\cdots>w_r \geq 0$ sont des entiers pairs je d\'esigne par
$$O(w_1,\cdots,w_r) \, \, \, \, {\rm et}\, \, \, \, O^\ast(w_1,\cdots,w_r)$$
le nombre des $\pi$ orthogonaux respectivement de $\GL(2r)$ et $\GL(2r+1)$ satisfaisant (i), (ii) et (iii), de poids
de Hodge les $w_i$. \ps

Pour toutes les valeurs $w_1>\cdots>w_{[n/2]}$ de l'\'enonc\'e nous donnons
plus pr\'ecis\'ement le nombre conjectural des $\pi$ ayant ces poids de
Hodge et qui de plus ont un groupe de Langlands-Sato-Tate donn\'e. 
J'entends par l\`a le groupe compact image du param\`etre global conjectural
du groupe de Langlands de $\Q$ associ\'e \`a $\pi$, qui est un sous-groupe
compact {\it connexe}\footnote{Car ${\rm Spec}(\Z)$ est simplement connexe!}
de ${\rm SL}(n,\C)$ agissant irr\'eductiblement sur $\C^n$. Je reviendrai sur ce groupe au~\S\ref{conjLZ}. Donnons
quelques exemples (on sous-entendra toujours ci-dessous que $\pi$ est une
repr\'esentation automorphe cuspidale satisfaisant
(i), (ii) et (iii)) : \medskip \begin{itemize}
{\it
\item[(1)] Des $\pi$ de $\GL(6)$ de groupe de Langlands-Sato-Tate
le groupe symplectique compact ${\rm Sp}(6)$ apparaissent en poids
motivique $23$ et pas moins. Leurs poids de Hodge en poids motivique $23$ sont $$(23, 13, 5), (23,
15, 3), (23, 15, 7), (23, 17, 5), (23, 17, 9), (23, 19, 3), (23, 19, 11).$$

\item[(2)] Des $\pi$ de $\GL(7)$ de groupe de
Langlands-Sato-Tate le groupe compact ${\rm SO}(7,\R)$ appparaissent en
poids motivique $26$ et pas moins. Leurs poids de Hodge en poids motivique
$26$ sont 
$$(26, 20, 10), (26, 20, 14), (26, 24, 10), (26, 24, 14), (26, 24, 18).$$

\item[(3)] Des $\pi$ de $\GL(8)$ de groupe de
Langlands-Sato-Tate le groupe symplectique compact ${\rm Sp}(8)$ apparaissent en poids motivique
$25$ et pas moins. En poids motivique $25$ il y en a exactement $33$. \ps

\item[(4)] Le premier $\pi$ de $\GL(8)$ de groupe de Langlands-Sato-Tate
le groupe compact ${\rm SO}(8)$ appara\^it en poids
motivique $24$. Il y en a un seul de ce poids, de poids de Hodge 
$(24, 20, 14, 2)$.
}
\end{itemize}

\medskip \medskip Il serait int\'eressant de d\'eterminer des param\`etres de Satake de ces
repr\'esentations. Je voudrais rajouter que je n'ai trouv\'e aucune trace
dans la litt\'erature des repr\'esentations ci-dessus, hormis dans le cas
(2) o\`u de telles repr\'esentations ont \'et\'e mises \'egalement en \'evidence
ind\'ependamment par Bergstr\"om, Faber et van der Geer dans leur analyse de
la vari\'et\'e de Siegel de genre $3$~\cite{BFVdG}. Lorsque $n=7$, nous donnons aussi le nombre conjectural exact 
des $\pi$ de $\GL(7)$ dont le groupe de Langlands-Sato-Tate est le groupe ${\rm G}_2$
en fonction de leurs poids de Hodge. On retrouve ainsi comme cas
particuliers quelques pr\'edictions de~\cite{BFVdG}. La question de
l'existence de motifs sur $\Q$ de groupes de Galois de type ${\rm G}_2$ remonte au
moins \`a Serre~\cite{serremotives}. Elle a \'et\'e \'etudi\'ee par plusieurs auteurs dont Gross
et Savin dans~\cite{grosssavin} par des m\'ethodes automorphes, qui
consid\`erent cependant des repr\'esentations qui sont Steinberg en une
place.\ps \ps 

\medskip \begin{itemize}
{\it
\item[(5)] Le premier $\pi$ de $\GL(7)$ de groupe de Langlands-Sato-Tate
le groupe compact ${\rm G}_2$ appara\^it en poids motivique $24$. Il n'y en n'a
qu'un de ce poids qui est de poids de Hodge $(24,16,8)$. 
}
\end{itemize}

\medskip

Je termine ce paragraphe en donnant quelques cons\'equences qui me semblent int\'eressantes de nos
r\'esultats. L'application qui m'est la plus ch\`ere concerne la question de la d\'etermination, dans
l'esprit de mon travail avec Lannes, des $121$ repr\'esentations automorphes
de conducteur $1$ et \`a coefficients triviaux d'une $\Z$-forme semisimple du groupe ${\rm SO}(25,\R)$.
Ces repr\'esentations sont en bijection avec les r\'eseaux pairs de covolume
$\sqrt{2}$
de l'espace euclidien $\R^{25}$, et Borcherds a d\'emontr\'e
dans~\cite{borcherdsthese} qu'il y a
exactement $121$ classes d'isom\'etrie de tels r\'eseaux (on en obtient des
exemples en consid\'erant la somme orthogonale de n'importe quel r\'eseau
de Niemeier avec la forme quadratique $2 x^2$). Soit $\rho$ le caract\`ere infinit\'esimal de la repr\'esentation triviale de $G(\R)$. 

\begin{thmetoile}\label{so25} Soit ${\rm G}$ un $\Z$-groupe classique tel que ${\rm G}(\R)$ soit le
groupe compact ${\rm SO}(25,\R)$. Il existe exactement $121$ param\`etres d'Arthur pour $G$ de
caract\`ere infinit\'esimal $\rho$ dont les constituants sont les repr\'esentations suivantes satisfaisant
(i), (ii) et (iii) : 
\begin{itemize}
\item[-] la repr\'esentation triviale de $\GL(1)$,\ps
\item[-] les $7$ repr\'esentations de $\GL(2)$ de poids motivique $\leq
23$, \ps
\item[-] les $7$ repr\'esentations symplectiques de $\GL(4)$ de poids
motivique $\leq 23$, \ps
\item[-] les $7$ repr\'esentations symplectiques de $\GL(6)$ de poids
motivique $23$ mentionn\'ees dans le (1) ci-dessus,\ps
\item[-] le carr\'e sym\'etrique de la repr\'esentation de $\GL(2)$ associ\'ee \`a
la forme $\Delta \in {\rm S}_{12}({\rm SL}(2,\Z))$.\ps
\end{itemize}
De plus, ces param\`etres sont exactement les param\`etres des
repr\'esentations automorphes
de conducteur $1$ de $G$ de caract\`ere infinit\'esimal $\rho$.  \end{thmetoile}

Je renvoie \`a notre article pour une liste des param\`etres en question. Le m\^eme ph\'enom\`ene se produit
donc que dans le travail avec Lannes : il y a non seulement 
exactement $121$ param\`etres d'Arthur "acceptables" et une analyse au cas par cas
de ces param\`etres montre que la formule de multiplicit\'e
d'Arthur vaut toujours $1$. Cette co\"incidence me semble assez formidable
d'un point de vue num\'erologique, encore plus que celle dans le travail avec Lannes. Je voudrais rajouter
qu'en dimension plus grande, par exemple dans le cas juste apr\`es o\`u
${\rm G}$ un $\Z$-groupe classique tel que ${\rm
G}(\R)$ soit le groupe compact ${\rm SO}(31,\R)$, il existe des param\`etres
d'Arthur pour ${\rm G}$ de caract\`ere infinit\'esimal $\rho$ dont la multiplicit\'e est
nulle : c'est par exemple le cas pour le param\`etre 
$$\Delta_{29}^{(2)} \oplus \Delta_{27}^{(2)} \oplus \Delta_{17}[9] \oplus [8]$$
o\`u $\Delta_w^{(k)}$ d\'esigne l'une quelconque des $k$ repr\'esentations
de $\GL(2)$ associ\'ee \`a une forme propre de poids $w+1$ pour ${\rm
SL}(2,\Z)$, not\'ee simplement $\Delta_w$ s'il n'y a qu'une telle forme. \ps

Une deuxi\`eme application concerne la construction d'une fonction $L$
dans la classe de Selberg qui me semble remarquable. Si $\underline{w}=(w_1,\cdots,w_r)$
est une suite d\'ecroissante d'entiers impairs positifs, j'utiliserai la notation
$\Delta_{\underline{w}}$ pour d\'esigner l'unique repr\'esentation
automorphe cuspidale de $\GL(2r)$ sur $\Q$ satisfaisant (i), (ii), (iii),
symplectique, et de poids de Hodge les $w_i$, sous r\'eserve qu'il existe
une et une seule telle repr\'esentation. Une mise en garde : ces indexations ne sont pas compatibles avec les notations utilis\'ees au chapitre pr\'ec\'edent, o\`u l'on notait plut\^ot $\Delta_{22}$ pour ce qui est ici $\Delta_{21}$, ou encore $\Delta_{4,10}$ pour ce qui est ici $\Delta_{21,5}$, mais ce sont les seules vraiment raisonnables dans leurs contextes respectifs. \ps

\begin{thmetoile} Il existe une repr\'esentation automorphe temp\'er\'ee non cuspidale de
$\GL(28)$, de caract\`ere infinit\'esimal $\rho$, satisfaisant (i), (ii) et (iii) plus
haut, \`a savoir : $$\Delta_{27,23,9,1} \oplus \Delta_{25,13,3}
\oplus \Delta_{21,5} \oplus \Delta_{19,7} \oplus \Delta_{17} \oplus
\Delta_{15} \oplus \Delta_{11}.$$
\end{thmetoile}

Observer en effet que tous les nombres impairs de $1$ \`a $27$ apparaissent une et
une seule fois dans la "formule" ci-dessus. L'assertion non-triviale, qui
r\'esulte de nos tables, est l'existence de $\Delta_{27,23,9,1}$ et
$\Delta_{25,13,3}$, qui sont des repr\'esentations de $\GL(8)$ et $\GL(6)$
respectivement. Un autre miracle se produit : en utilisant tout ce qu'il y a
dans mes tables, je ne peux construire aucune autre repr\'esentation symplectique $\pi$ de $\GL(2m)$
ayant les propri\'et\'es ci-dessus quand $2m \leq 28$ (donc aucune quand $2m
\leq 26$). Il est int\'eressant de comparer ce r\'esultat avec ce que
donnent les formules explicites de Weil. En effet, les travaux de
Fermigier~\cite{fermigier} (contrairement \`a ceux de Miller~\cite{miller}) s'appliquent aussi aux formes
temp\'er\'ees non cuspidales, et
montrent que l'existence d'un $\pi$ de $\GL(2m)$ comme dans le th\'eor\`eme
ci-dessus n\'ecessite $2m\geq 26$, qui est tr\`es proche de $28$ ! \ps

Une derni\`ere application de notre travail concerne la dimension des
espaces ${\rm S}_V({\rm Sp}(2g,\Z))$ de formes de Siegel. Nous donnons en effet une recette
explicite pour la dimension de cet espace en fonction des divers ${\rm
X}(w_1,\cdots,w_r)$. Je ne voudrais pas faire croire ici que nous retrouvons par nos calculs les dimensions des espaces de formes de Siegel vectorielles de genre pour $g=1,2$ et $3$. Nous utilisons ces dimensions  pour d\'emontrer le th\'eor\`eme principal, pr\'ecis\'ement pour d\'eterminer les quantit\'es $S(w)$, $S(w_1,w_2)$ et $O^\ast(w_1,w_2,w_3)$. En fait, par un "coup de chance" notre m\'ethode permet quand m\^eme de d\'eterminer $O^\ast(w_1,w_2,w_3)$ pour toute une collection de poids de Hodge $(w_1,w_2,w_3)$, par exemple tout ceux $<28$.    
Nos r\'esultats confirment les calculs d\'elicats et ind\'ependants de Bergstr\"om, Faber et
Van der Geer d\'ej\`a mentionn\'es. \ps

\subsection{Id\'ees de la d\'emonstration} La d\'emonstration du
Th\'eor\`eme~\ref{chrenard} se d\'ecoupe en trois \'etapes que je vais
maintenant d\'etailler : \begin{itemize} \medskip 

\item[(a)] Exprimer, gr\^ace \`a la th\'eorie d'Arthur, la dimension de
l'espace des formes automorphes de niveau $1$ pour un $\Z$-groupe classique
$G$ donn\'e, et de composante archim\'edienne une s\'erie discr\`ete
$\pi_\infty$ donn\'ee, en fonction des $S(-)$, $O(-)$ et $O^\ast(-)$. \ps

\item[(b)] Calculer ind\'ependamment, pour des groupes $G$ \'eventuellement bien choisis, la
dimension de l'espace des formes automorphes en questions. \ps

\item[(c)] En extraire, par r\'ecurrence, les valeurs de $S$, $O$ et
$O^\ast$. \ps 
\end{itemize}

\ps

Pour le (a), le probl\`eme est de comprendre, si $\psi \in \Psi_{\rm
glob}(G)$ est un param\`etre d'Arthur de $G$ donn\'e, \`a quelle condition
il existe une repr\'esentation automorphe discr\`ete $\pi$ de $G$ ayant ce
param\`etre et dont la composante archim\'edienne est le $\pi_\infty$
prescrit. Il s'agit donc d'examiner la formule de multiplicit\'e d'Arthur.  Nous en aurons notamment besoin notamment quand $G(\R)$ est
compact, et aussi quand $G={\rm Sp}(2g)$ et $\pi_\infty$ est une s\'erie
discr\`ete holomorphe.  Je soup\c{c}onne que c'est la difficult\'e \`a comprendre ce point qui a fait que notre approche n'a pas \'et\'e consid\'er\'ee plus t\^ot historiquement. 
En effet, bien que la litt\'erature sur les formes de Siegel semble profond\'ement impreign\'ee des conjectures d'Arthur (qui datent des ann\'ees $80$), le lien entre les divers "rel\`evements" consid\'er\'es (Saito-Kurokawa, Ikeda, Miyawaki, et d'autres qui ont suivi) et la forme exacte de la formule de multiplicit\'e d'Arthur exacte n'est jamais consid\'er\'ee \`a ma connaissance (hors du cadre de ${\rm PGSp}(4)$ ou des cas stables), et ce m\^eme conjecturalement.  La difficult\'e \`a extraire ces formules des r\'esultats tr\`es g\'en\'eraux d'Arthur est sans doute responsable de cette "paresse". \ps

Pour le (b), l'id\'ee est de n'effectuer ce calcul
de dimensions, dont on a d\'ej\`a expliqu\'e la difficult\'e quand $G={\rm
Sp}(2g)$, que dans le cas o\`u $G(\R)$ est un groupe compact, ce qui se
ram\`ene \`a un probl\`eme de th\'eorie des invariants. Dans les dimensions
qui nous int\'eressent, le groupe $G$ est n\'ecessairement le groupe sp\'ecial
orthogonal des r\'eseaux de racines ${\rm E}_7$, ${\rm E}_8$ et ${\rm E}_7
\times {\rm A}_1$. Ces r\'eseaux sont uniques dans leur genre (pairs de
d\'eterminants respectifs $2$, $1$ et $2$), de sorte que le groupe $G$ est
de nombre de classes $1$, et on est donc ramen\'e \`a d\'eterminer la
dimension des $G(\Z)$-invariants dans chaque repr\'esentation
irr\'eductible du groupe compact $G(\R)$, i.e. ${\rm SO}(n,\R)$ avec
respectivement $n=7,8$ et $9$. J'expliquerai ci-apr\`es comment nous avons
proc\'ed\'e. \ps

Comme nous le verrons, les calculs de dimension effectu\'es au (b) donnent des
relations lin\'eaires entre certains $S(-)$, $O(-)$ et $O^\ast(-)$. Ces
relations ne permettent pas \`a elles seules de d\'emontrer le th\'eor\`eme,
car elles renferment des inconnues. J'expliquerai enfin comment nous
d\'eterminons ces inconnues, en se ramenant d'une part \`a d'autres calculs de
dimensions connues pour les ${\rm Sp}(2g)$ avec $g=1, 2$ et $3$, et aussi en
d\'emontrant des cas "nouveaux" de fonctorialit\'e \`a la Langlands. Ces
r\'esultats se d\'eduisent en fait simplement des r\'esultats d'Arthur et
sont bas\'es sur l'id\'ee que le groupe de Langlands de $\Z$, sur lequel je reviendrai au~\S\ref{conjLZ}, est simplement connexe (ce que nous "d\'emontrons" !). 
Par exemple, nous prouvons que toute repr\'esentation de
$\GL(5)$ satisfaisant (i), (ii) et (iii) est un $\Lambda^2$ "r\'eduit" d'une
et une seule repr\'esentation symplectique de $\GL(4)$ satisfaisant (i),
(ii) et (iii). \ps

\subsubsection{\'Etape (b) : calculs de dimensions d'espaces d'invariants} Je
commence par cette \'etape car c'est la plus facile.  Le cadre est le
suivant.  Soit $H$ un groupe de Lie compact connexe et soit $\Gamma$ un
sous-groupe fini de $H$.  Il est bien connu que les repr\'esentations
irr\'eductibles de $H$ sont param\'etr\'ees par leur poids dominant. Cela
suppose d'avoir fix\'e un tore maximal $T$ et une chambre de Weyl dans ${\rm
X}^\ast(T) \otimes \R$, et je noterai $V_\lambda$ la repr\'esentation
irr\'eductible de $H$ de plus haut poids $\lambda$. On se pose la question de d\'eterminer la
dimension $$d(\lambda)=\dim(V_\lambda^\Gamma)$$ de l'espace des invariants de $\Gamma$
agissant sur $V_\lambda$. Comme d\'ej\`a dit plus haut, les cas particuliers de paires $(H,\Gamma)$
qui nous int\'eressent sont les paires :
$$({\rm SO}(7,\R), {\rm W}^+({\rm E}_7)), \, \, \, ({\rm SO}(8,\R), {\rm
W}^+({\rm E}_8)), \, \, \, {\rm et} \, \, \, ({\rm SO}(9,\R), {\rm    
W}({\rm E}_8)),$$
o\`u $W(R) \subset {\rm O}(V)$ d\'esigne le groupe de Weyl d'un syst\`eme de racines $R$ de
l'espace euclidien $V$ et o\`u $W^+(R)=W(R) \cap {\rm SO}(V)$. Pour les
applications \`a ${\rm G}_2$ on consid\`ere aussi la paire $({\rm G}_2(\R),{\rm G}_2(\Z))$ o\`u 
${\rm G}_2$ est l'unique $\Z$-groupe semisimple tel que ${\rm G}_2(\R)$ est
compact (Gross~\cite{grossinv}). \ps

Dans tous ces exemples, $\Gamma$ est
intuitivement un tr\`es gros sous-groupe de $H$, on s'attend donc \`a ce que $d(\lambda)$ soit petit pour
des petites valeurs de $|\lambda|$. En revanche, m\^eme pour $\lambda$ petit
$\dim(V_\lambda)$ peut \^etre assez gros, auquel cas $V_\lambda$ est vraisemblablement
ing\'erable pour un ordinateur. Par exemple, dans le cas $H={\rm SO}(8,\R)$ et $\lambda=(15,6,5,4)$ dans les notations standards, $V_\lambda$ est de
dimension $143503360$ alors que $d(\lambda)=0$ (est-ce le record?). Pour calculer $d(\lambda)$ nous partons
plut\^ot de la relation triviale $$d(\lambda)=\frac{1}{|\Gamma|} \sum_{\gamma \in
\Gamma} {\rm Trace}(\gamma,V_\lambda),$$ et utilisons une version
d\'eg\'en\'er\'ee de la formule du caract\`ere de Weyl permettant
d'\'evaluer ${\rm Trace}(t,V_\lambda)$ pour tout $t \in T$, y compris quand $t$ n'est pas r\'egulier. Nous avions \'etabli et utilis\'e une telle formule dans mon
travail~\cite{chclo} avec Clozel dont j'ai d'ailleurs d\'ej\`a parl\'e. Le second ingr\'edient
pour appliquer cette formule est de d\'eterminer
des repr\'esentants dans $T$ de toutes les classes de conjugaison de
$\Gamma$. Heureusement, ce travail d\'elicat avait d\'ej\`a \'et\'e fait par
Carter dans~\cite{carter}, et ce pour tous les groupes de Weyl dans leur repr\'esentation de
reflexion. Dans le cas suppl\'ementaire du groupe $\Gamma={\rm G}_2(\Z)$, nous utilisons
plut\^ot des
g\'en\'erateurs de ce dernier comme sous-groupe de ${\rm SO}(7,\R)$ donn\'es par A. Cohen, G. Nebe et W. Plesken
dans~\cite{cn}. \ps
Les formules finales sont monstrueuses (voir par exemple la discussion page
$9$ de l'introduction de notre article pour quelques pr\'ecisions). Cependant, 
nous les avons impl\'ement\'ees sur ordinateur et au final il faut par exemple environ $5$ minutes \`a l'ordinateur pour
calculer une dimension pour $H={\rm SO}(7)$ pour un poids dominant standard $a \geq b\geq
c\geq 0$ avec $a \leq 100$, ce qui va d\'ej\`a bien au del\`a des
 valeurs qui m'int\'eressent pour l'instant! Je renvoie \`a la
Table~\ref{tableSO7nue} qui suit pour les premi\`eres valeurs. \ps

\begin{table}[htp]
\begin{center} $G={\rm SO}(7,\R)$, $\Gamma={\rm W}^+({\rm E}_7)$.
\end{center}
\begin{tabular}{|c|c||c|c||c|c||c|c||c|c|}
\hline  $\lambda$ & $d(\lambda)$  & $\lambda$ & $d(\lambda)$  & $\lambda$ & $d(\lambda)$  & $\lambda$ & $d(\lambda)$  & $\lambda$ & $d(\lambda)$ \\ 
\hline
(0, 0, 0) & 1& (9, 6, 3) & 2& (10, 7, 2) & 1& (10, 10, 10) & 2& (11, 9, 0) & 2 \\
\hline
(4, 4, 4) & 1& (9, 6, 4) & 1& (10, 7, 3) & 3& (11, 3, 0) & 1& (11, 9, 1) & 1 \\
\hline
(6, 0, 0) & 1& (9, 6, 6) & 1& (10, 7, 4) & 2& (11, 3, 2) & 1& (11, 9, 2) & 4 \\
\hline
(6, 4, 0) & 1& (9, 7, 2) & 1& (10, 7, 5) & 2& (11, 4, 1) & 1& (11, 9, 3) & 4 \\
\hline
(6, 6, 0) & 1& (9, 7, 3) & 1& (10, 7, 6) & 2& (11, 4, 3) & 2& (11, 9, 4) & 5 \\
\hline
(6, 6, 6) & 1& (9, 7, 4) & 2& (10, 7, 7) & 1& (11, 4, 4) & 1& (11, 9, 5) & 4 \\
\hline
(7, 4, 3) & 1& (9, 7, 6) & 1& (10, 8, 0) & 3& (11, 5, 0) & 2& (11, 9, 6) & 5 \\
\hline
(7, 6, 3) & 1& (9, 8, 1) & 1& (10, 8, 2) & 3& (11, 5, 2) & 2& (11, 9, 7) & 3 \\
\hline
(7, 7, 3) & 1& (9, 8, 3) & 1& (10, 8, 3) & 1& (11, 5, 3) & 1& (11, 9, 8) & 2 \\
\hline
(7, 7, 7) & 1& (9, 8, 4) & 1& (10, 8, 4) & 4& (11, 5, 4) & 1& (11, 9, 9) & 1 \\
\hline
(8, 0, 0) & 1& (9, 8, 5) & 1& (10, 8, 5) & 1& (11, 6, 1) & 2& (11, 10, 1) & 3 \\
\hline
(8, 4, 0) & 1& (9, 8, 6) & 1& (10, 8, 6) & 3& (11, 6, 2) & 1& (11, 10, 2) & 3 \\
\hline
(8, 4, 2) & 1& (9, 9, 0) & 1& (10, 8, 7) & 1& (11, 6, 3) & 4& (11, 10, 3) & 5 \\
\hline
(8, 4, 4) & 1& (9, 9, 3) & 1& (10, 8, 8) & 1& (11, 6, 4) & 2& (11, 10, 4) & 4 \\
\hline
(8, 6, 0) & 1& (9, 9, 4) & 1& (10, 9, 1) & 2& (11, 6, 5) & 2& (11, 10, 5) & 6 \\
\hline
(8, 6, 2) & 1& (9, 9, 6) & 1& (10, 9, 2) & 1& (11, 6, 6) & 2& (11, 10, 6) & 5 \\
\hline
(8, 6, 4) & 1& (9, 9, 9) & 1& (10, 9, 3) & 3& (11, 7, 0) & 1& (11, 10, 7) & 5 \\
\hline
(8, 6, 6) & 1& (10, 0, 0) & 1& (10, 9, 4) & 2& (11, 7, 1) & 1& (11, 10, 8) & 3 \\
\hline
(8, 7, 2) & 1& (10, 2, 0) & 1& (10, 9, 5) & 3& (11, 7, 2) & 4& (11, 10, 9) & 2 \\
\hline
(8, 7, 4) & 1& (10, 4, 0) & 2& (10, 9, 6) & 2& (11, 7, 3) & 3& (11, 10, 10) & 2 \\
\hline
(8, 7, 6) & 1& (10, 4, 2) & 1& (10, 9, 7) & 2& (11, 7, 4) & 4& (11, 11, 1) & 1 \\
\hline
(8, 8, 0) & 1& (10, 4, 3) & 1& (10, 9, 8) & 1& (11, 7, 5) & 3& (11, 11, 2) & 2 \\
\hline
(8, 8, 2) & 1& (10, 4, 4) & 2& (10, 9, 9) & 1& (11, 7, 6) & 3& (11, 11, 3) & 3 \\
\hline
(8, 8, 4) & 1& (10, 5, 1) & 1& (10, 10, 0) & 2& (11, 7, 7) & 2& (11, 11, 4) & 2 \\
\hline
(8, 8, 6) & 1& (10, 5, 3) & 1& (10, 10, 2) & 2& (11, 8, 1) & 3& (11, 11, 5) & 3 \\
\hline
(8, 8, 8) & 1& (10, 6, 0) & 2& (10, 10, 3) & 2& (11, 8, 2) & 2& (11, 11, 6) & 3 \\
\hline
(9, 3, 0) & 1& (10, 6, 2) & 2& (10, 10, 4) & 4& (11, 8, 3) & 5& (11, 11, 7) & 3 \\
\hline
(9, 4, 3) & 1& (10, 6, 3) & 1& (10, 10, 5) & 2& (11, 8, 4) & 4& (11, 11, 8) & 2 \\
\hline
(9, 4, 4) & 1& (10, 6, 4) & 3& (10, 10, 6) & 4& (11, 8, 5) & 5& (11, 11, 9) & 1 \\
\hline
(9, 5, 0) & 1& (10, 6, 5) & 1& (10, 10, 7) & 2& (11, 8, 6) & 4& (11, 11, 10) & 1 \\
\hline
(9, 5, 2) & 1& (10, 6, 6) & 2& (10, 10, 8) & 2& (11, 8, 7) & 3& (11, 11, 11) & 1 \\
\hline
(9, 6, 1) & 1& (10, 7, 1) & 2& (10, 10, 9) & 2& (11, 8, 8) & 1& (12, 0, 0) & 2 \\
\hline
\end{tabular}

\caption{{\small Valeurs non nulles de $d(\lambda)=\dim V_\lambda^\Gamma$
pour
$\lambda=(n_1,n_2,n_3)$ avec $n_1 \leq 11$.}}
\label{tableSO7nue}
\end{table}

Mentionnons pour finir que les erreurs possibles lors de
l'impl\'ementation sont multiples, mais il y a de nombreux moyens de
v\'erifier que nos calculs sont exacts. Par exemple, les dimensions doivent \^etre des
entiers (alors qu'elles sont dans un corps cyclotomique gigantesque par
notre proc\'ed\'e de calcul), elles doivent \^etre compatibles avec la s\'erie de Poincar\'e des groupes de
Coxeter~\cite[Chap. V \S 6]{bourbaki} dans le cas des repr\'esentations sur les polyn\^omes harmoniques, ou
encore \^etre compatibles aux formules de multiplicit\'e d'Arthur... \ps \medskip

\subsubsection{\'Etape (a) : explicitation de la formule de multiplicit\'e
d'Arthur}\label{formmultarth} Fixons $G$ un $\Z$-groupe classique semisimple. Nous avons d\'ej\`a
rappel\'e comment Arthur associe \`a chaque $\pi \in \Pi_{\rm disc}(G)$ une
unique classe d'\'equivalence de param\`etres globaux $\psi(\pi)$ dans $\Psi_{\rm
glob}(G)$. Fixons donc r\'eciproquement un param\`etre d'Arthur global
$\psi=(k,(n_i),(d_i),(\pi_i))$. La formule de multiplicit\'e d'Arthur donne
une condition pr\'ecise sous laquelle ce param\`etre est celui d'une
repr\'esentation. \ps

La discussion de cette formule dans le cas g\'en\'eral
nous conduirait ici bien trop loin. C'est pourquoi je propose de se
restreindre dans ce texte au cas des groupes $G$ tels que $G(\R)$ est
compact, et de renvoyer \`a mon article~\cite{chrenard2} pour plus de
g\'en\'eralit\'e. Une condition alors n\'ecessaire pour que $\psi$ soit de la forme $\psi(\pi)$
est que le caract\`ere infinit\'esimal $z_\psi$ de $\psi$ (voir~\S~\ref{disc24}) soit dans l'image, par la
repr\'esentation standard $${\rm St}: \widehat{G} \rightarrow
\SL(n(\widehat{G}),\C),$$ 
du caract\`ere infinit\'esimal $\mu$ d'une repr\'esentation irr\'eductible
de dimension finie du groupe compact $G(\R)$. Le cas $\widehat{G}={\rm
Sp}(2g,\C)$ est un peu plus simple et je vais d'abord me concentrer sur ce
cas. Par exemple, $G$ peut \^etre le groupe sp\'ecial orthogonal de la forme
quadratique ${\rm E}_7$ ou ${\rm E}_8 \times {\rm A}_1$, des cas tr\`es
importants pour la discussion ici. Sous l'hypoth\`ese $\widehat{G}={\rm Sp}(2g,\C)$, la repr\'esentation ${\rm St}$ induit une injection de l'ensemble des classes
de conjugaison semisimples de $\widehat{G}$ dans celles de ${\rm
SL}(n(\widehat{G}),\C)$. Ainsi, Satake et Harish-Chandra assurent qu'il
existent {\it un et un seul} $\pi \in \Pi(G)$ tel que~(voir~\S~\ref{pararthur})
$${\rm St}(c(\pi))= \rho_\psi(e \times \prod_i c(\pi_i)).$$ 
 \ps

La question est alors de d\'eterminer la multiplicit\'e $m(\pi)$ de
ce $\pi$ dans le spectre automorphe du $G$. La formule d'Arthur, faite pour
cela, comporte trois ingr\'edients. Le premier est le groupe fini ${\rm C}_\psi \subset
\widehat{G}$. On observe pour cela que le morphisme $\rho_\psi$ d\'efini
au~\S~\ref{pararthur} se factorise par ${\rm St}$, de sorte que l'on peut m\^eme
supposer qu'il est \`a valeurs dans $\widehat{G}$. Le groupe ${\rm C}_\psi$
est alors le centralisateur de ${\rm Im}(\rho_\psi)$ dans $\widehat{G}$. Il
s'identifie ici naturellement au groupe $\{ \pm 1\} ^{\{1,\cdots,k\}}$. \ps

Arthur d\'efinit ensuite deux caract\`eres ${\rm C}_\psi \rightarrow \{ \pm 1\}$. Le premier,
not\'e $\varepsilon_\psi$, d\'epend des facteurs epsilon des paires
$\{\pi_i,\pi_j\}$ pour $1\leq i < j \leq k$. Ces facteurs ici sont des
signes explicitement calculables en terme des poids de Hodge des $\pi_i$ et
$\pi_j$ (car on est en conducteur $1$), je renvoie \`a~\cite[\S 3.20]{chrenard2} pour un formulaire. Un examen de la recette d'Arthur
montre que l'on a la formule 
$$\varepsilon_\psi(s_i)= \prod_{j \neq i}\varepsilon(\pi_i \times \pi_j)^{{\rm Min}(d_i,d_j)}$$
o\`u $s_i \in {\rm C}_\psi=\{\pm1 \}^{\{1,\cdots,k\}}$ d\'esigne l'\'el\'ement tel que $s_i(x)=-1$ si
$x=i$ et $1$ sinon. \ps

Le dernier ingr\'edient, d'une certaine mani\`ere le plus subtile, est un caract\`ere de ${\rm
C}_\psi$ associ\'e \`a $(G(\R),\pi_\infty)$. Dans le cas qui nous int\'eresse, il
s'est av\'er\'e tr\`es simple. Partons du caract\`ere infinit\'esimal
$z_\psi$ de $\psi$. Par hypoth\`ese, il est l'image d'une classe de
conjugaison semisimple bien d\'efinie dans ${\rm Lie}({\rm Im}(\rho_\psi))$. Fixons
$\mu$ dans cette derni\`ere classe et d\'esignons par $\widehat{T}$ son
centralisateur dans $\widehat{G}$, qui est donc un tore maximal de ce dernier. On a une inclusion \'evidente 
$${\rm C}_\psi \subset \widehat{T}.$$ On consid\`ere
l'unique syst\`eme de racines positives pour $(\widehat{G},\widehat{T})$ pour lequel $\mu \in {\rm
X}_\ast(\widehat{T})$ est un co-caract\`ere dominant, et $\rho^\vee$ la
demi-somme des racines associ\'ee. Ici, $\rho^\vee \in {\rm X}^\ast(\widehat{T})$, il y a donc un sens \`a restreindre $\rho^\vee$ \`a ${\rm C}_\psi$. J'affirme que la formule de multiplicit\'e d'Arthur s'\'ecrit simplement 
\begin{equation}\label{formulemagique} m(\pi) \neq 0 \, \, \, \Leftrightarrow \, \, \, \varepsilon_\psi = \rho^\vee_{|{\rm C}_\psi}\end{equation}
et de plus que la multiplicit\'e est $1$ si elle est non nulle. Avant de "justifier" cet \'enonc\'e, je voudrais signaler d'ores et d\'ej\`a que l'\'etablissement de cette formule est encore conditionnel \`a certains r\'esultats non d\'emontr\'es de la th\'eorie d'Arthur pour les formes int\'erieures des groupes classiques~\cite[Chap. 9]{arthur}, et que tous nos th\'eor\`emes "\'etoil\'es" en d\'ependent notamment.\ps
 Avec un peu d'entrainement il est quasiment instantan\'e de v\'erifier la relation $\varepsilon_\psi = \rho^\vee_{|{\rm C}_\psi}$ pour un param\`etre $\psi \in \Psi_{\rm glob}(G)$ donn\'e. De plus, les conditions obtenues (d'ailleurs, quelque soit la formule raisonnable pr\'ecise comme on pourrait le voir!) portent toujours uniquement sur les positions relatives des poids de Hodge des diff\'erents $\pi_i$ (\`a cause de $\rho^\vee$) et de congruences modulo $4$ sur ces nombres (\`a cause de $\varepsilon_\psi$).  \ps
 
 Le cas particulier le plus important est celui o\`u $k=1$, i.e. $C_{\rm \psi}=Z(\widehat{G})=\{\pm 1\}$ (cas "stable"). 
Dans ce cas la multiplicit\'e est trivialement $1$ (quelque soit la formule de multiplicit\'e exacte) et en particulier toutes les repr\'esentations cuspidales symplectiques de $\GL(n(\widehat{G}))$ satisfaisant (i), (ii) et (iii) vont contribuer au spectre discret de $G$. Cependant, la compr\'ehension de tous les cas est n\'ecessaire pour pouvoir estimer le nombre de ces repr\'esentations en fonction de la dimension des espaces de formes automorphes ! Je renvoie \`a notre article pour des formules explicites, dont le nombre augmente rapidement avec $n(\widehat{G})$ (il y en a au moins autant que de partitions de l'entier $\frac{n(\widehat{G})}{2}$). Dans le cas de ${\rm SO}(9)$, il y a par exemple $16$ formules de multiplicit\'es diff\'erentes. Je renvoie par exemple \`a la figure~\ref{tableSO7} pour une description des param\`etres d'Arthur obtenus en petits poids dans le cas de ${\rm SO}(7)$. \ps
\ps

\begin{table}[htp]
\ps\ps
\renewcommand{\arraystretch}{1.5}
\caption{{\small Param\`etres d'Arthur des repr\'esentations automorphes de niveau $1$ de ${\rm SO}(7)$ dont le caract\`ere infinit\'esimal, dans $\mathfrak{sp}(6,\C)$, a pour valeurs propres 
$\pm w_1, \pm w_2, \pm w_3$ avec $23 \geq w_1 > w_2 > w_3 >0$.}}
\begin{tabular}{|c|c||c|c|}
\hline  $(w_1,w_2,w_3)$ & $\Pi_{w_1,w_2,w_3}({\rm SO}(7))$ & $(w_1,w_2,w_3)$ & $\Pi_{w_1,w_2,w_3}({\rm SO}(7))$ \\
\hline (5,3,1) & $[6]$ & (21,19,17) & $\Delta_{19}[3]$\\
\hline (13,11,9) & $\Delta_{11}[3]$ &  (23,9,1) & $\Delta_{23,9} \oplus [2]$\\

\hline (17,3,1) & $\Delta_{17}\oplus [4]$ &(23,11,7) & $\Delta_{23,7} \oplus \Delta_{11}$\\

\hline (17,11,1) & $\Delta_{17}\oplus \Delta_{11} \oplus [2]$ & (23,11,9) & $\Delta_{23,9} \oplus \Delta_{11}$\\

\hline (17,15,1) & $\Delta_{17}\oplus \Delta_{15} \oplus [2]$ & (23,13,1) & $\Delta_{23,13} \oplus [2]$\\

\hline (17,15,13) & $\Delta_{15}[3]$& (23,13,5) & $\Delta_{23,13,5}$\\

\hline (19,11,7) & $\Delta_{19,7}\oplus \Delta_{11}$& (23,15,3) & $\Delta_{23,15,3}$\\
\hline (19,15,7) & $\Delta_{19,7} \oplus \Delta_{15}$&  (23,15,7) & $\Delta_{23,7} \oplus \Delta_{15}$, $\Delta_{23,15,7}$\\
\hline (19,17,7) & $\Delta_{19,7} \oplus \Delta_{17}$ &(23,15,9) & $\Delta_{23,9} \oplus \Delta_{15}$\\

\hline (19,17,15) & $\Delta_{17}[3]$& (23,15,13) & $\Delta_{23,13} \oplus \Delta_{15}$\\
\hline (21,3,1) & $\Delta_{21}\oplus [4]$&  (23,17,5) & $\Delta_{23,17,5} $\\
\hline (21,11,1) & $\Delta_{21}\oplus \Delta_{11}\oplus [2]$&(23,17,7) & $\Delta_{23,7} \oplus \Delta_{17}$\\

\hline (21,11,5) & $\Delta_{21,5}\oplus \Delta_{11}$& (23,17,9) & $\Delta_{23,9} \oplus \Delta_{17}$, $\Delta_{23,17,9}$\\

\hline (21,11,9) & $\Delta_{21,9}\oplus \Delta_{11}$&  (23,17,13) & $\Delta_{23,13} \oplus \Delta_{17}$\\
\hline (21,15,1) & $\Delta_{21} \oplus \Delta_{15} \oplus [2]$&(23,19,3) & $\Delta_{23,19,3}$\\
\hline (21,15,5) & $\Delta_{21,5} \oplus \Delta_{15}$ & (23,19,7) & $\Delta_{23,7} \oplus \Delta_{19}$\\
\hline (21,15,9) & $\Delta_{21,9} \oplus \Delta_{15}$ & (23,19,9) & $\Delta_{23,9} \oplus \Delta_{19}$\\
\hline (21,15,13) & $\Delta_{21,13} \oplus \Delta_{15}$&  (23,19,11) & $\Delta_{23,19,11}$\\

\hline (21,17,5)&  $\Delta_{21,5} \oplus \Delta_{17}$ &(23,19,13) & $\Delta_{23,13} \oplus \Delta_{19}$\\

\hline (21,17,9) & $\Delta_{21,9} \oplus \Delta_{17}$&  (23,21,1) & ${\rm Sym}^2 \Delta_{11}[2]$\\
\hline (21,17,13) & $\Delta_{21,13}\oplus \Delta_{17}$&  (23,21,7) & $\Delta_{23,7} \oplus \Delta_{21}$\\

\hline (21,19,1) & $\Delta_{21}\oplus \Delta_{19} \oplus [2]$&(23,21,9) & $\Delta_{23,9} \oplus \Delta_{21}$\\
\hline (21,19,5) & $\Delta_{21,5} \oplus \Delta_{19}$& (23,21,13) & $\Delta_{23,13}\oplus \Delta_{21}$\\

\hline (21,19,9) & $\Delta_{21,9} \oplus \Delta_{19}$& (23,21,19) & $\Delta_{21}[3]$\\

\hline  (21,19,13) & $\Delta_{21,13} \oplus \Delta_{19}$ & &\\ \hline

\end{tabular}
\label{tableSO7}
\end{table}

	La d\'ecouverte de la formule~\eqref{formulemagique} a \'et\'e l'un des aspects les plus excitants de ce travail. Pour \^etre honn\^ete, nous l'avons tout d'abord devin\'ee "\`a l'envers" par extrapolation du cas temp\'er\'e en examinant les valeurs que nous avions obtenues en calculant (\'etape pr\'ec\'edente) les dimensions des espaces de formes automorphes pour ${\rm SO}(7)$ et ${\rm SO}(9)$, le point \'etant que quand la dimension est nulle alors tous les param\`etres d'Arthur candidats doivent \^etre de multiplicit\'e nulle. \ps
	
	En effet, je rappelle que si $\psi$ est temp\'er\'e, ce qui correspond au cas tr\`es particulier o\`u $d_i=1$ pour tout $i$, alors les s\'eries discr\`etes de caract\`ere infinit\'esimal $\mu$ de toutes les formes int\'erieures de $G(\R)$ sont en bijection canonique avec les caract\`eres de la $2$-torsion $\widehat{T}[2]$ de $\widehat{T}$ (nous renvoyons par exemple aux articles d'Adams~\cite{adams} et de Shelstad~\cite{shelstadbanff} \`a ce sujet). Observons qu'un tel caract\`ere induit bien un caract\`ere de ${\rm C}_\psi$ car ${\rm C}_\psi \subset \widehat{T}[2]$. 
La canonicit\'e dans la bijection pr\'ec\'edente est subtile et faite de sorte que l'unique s\'erie discr\`ete g\'en\'erique de caract\`ere infinit\'esimal $\mu$ de la forme int\'erieure d\'eploy\'ee de $G(\R)$ corresponde au caract\`ere trivial de $\widehat{T}[2]$ (on utilise ici que ${\rm SO}(2m+1,\C)$ est adjoint 
sinon il faut modifier un peu les \'enonc\'es en utilisant la notion de forme int\'erieure forte~\cite{abv},~\cite{adams}). Ceci \'etant fait, on v\'erifie que la repr\'esentation de dimension finie du groupe compact correspond bien au caract\`ere de $\widehat{T}[2]$ induit par $\rho^\vee$. \ps

Dans le cas d'un $\psi$ g\'en\'eral, il fallait comprendre la param\'etrisation analogue \`a celle ci-dessus pour le param\`etre d'Arthur archim\'edien $\psi_\infty$ de caract\`ere infinit\'esimal $\mu$ associ\'e \`a $\psi$, qui est aussi un param\`etre d'Adams-Johnson~\cite{AJ}. Nous avons travaill\'e sous l'hypoth\`ese standard que les paquets associ\'es par Arthur \`a ces param\`etres dans~\cite{arthur}, pour toutes les formes int\'erieures de $G(\R)$, co\"incident avec ceux d\'efinis par Adams et Johnson dans~\cite{AJ}, ce qui rend nos r\'esultats conditionnels \`a ce fait. \`A chaque tel param\`etre $\psi_\infty$ est associ\'e un unique param\`etre temp\'er\'e $\psi'_\infty$ de m\^eme caract\`ere infinit\'esimal, et on dispose d'une inclusion naturelle sur les centralisateurs ${\rm C}_{\psi_\infty} \subset {\rm C}_{\psi'_\infty}=\widehat{T}[2]$. Le fait fondamental est alors que si une s\'erie discr\`ete $\pi_0$ d'une forme int\'erieure $H$ de $G(\R)$ est dans le paquet d'Arthur associ\'e \`a $\psi_\infty$, le caract\`ere de ${\rm C}_{\psi_\infty}$ qui lui correspond est la restriction du caract\`ere de $\widehat{T}[2]$ associ\'e \`a $(H,\pi_0)$ par le paragraphe pr\'ec\'edent. Cela d\'emontre~\eqref{formulemagique} dans le cas tr\`es particulier de la forme compacte et de sa repr\'esentation de dimension finie de caract\`ere infinit\'esimal $\mu$. Nous avons tard\'e \`a r\'ealiser que ce fait, qui me semble tout \`a fait fondamental aux applications arithm\'etiques, a \'et\'e compris depuis longtemps par Kottwitz dans~\cite[\S 9]{kottwitz}. Cela nous a permis aussi, toujours dans~\cite{chrenard2}, d'\'etendre la discussion de~\cite{adams} au cas de tous les param\`etres d'Adams-Johnson. \ps
	Lorsque $\widehat{G}={\rm SO}(2m,\C)$ (auquel cas en fait $4 | m$) la situation est presque la m\^eme que dans le cas pr\'ec\'edent, si ce n'est que la repr\'esentation standard de $\widehat{G}$ ne permet que de s\'eparer les classes de ${\rm O}(2m,\C)$-conjugaison de ${\rm SO}(2m,\C)$, ce qui introduit quelques modifications mineures. Au final, la formule de multiplicit\'e est la m\^eme que la formule~\eqref{formulemagique}, \`a ceci pr\`es que les multiplicit\'es non nulles peuvent \^etre $1$ ou $2$ d'apr\`es Arthur. Plus g\'en\'eralement, cette m\^eme formule est la formule de multiplicit\'e d'Arthur explicite pour tout $\Z$-groupe semisimple $G$ tel que $G(\R)$ est compact. Nous l'avons par exemple aussi utilis\'ee comme point de d\'epart dans notre discussion conjecturale sur $G_2$. Je renvoie enfin \`a mon article pour le cas de ${\rm Sp}(2g)$, pour lequel ce n'est pas le caract\`ere $\rho^\vee$ qui intervient. \ps
	
	On d\'eduit de la formule~\eqref{formulemagique} la v\'erification manquante du Th\'eor\`eme~\ref{24niemeier} avec Lannes, ainsi que le Th\'eor\`eme~\ref{so25}. Des crit\`eres assez simples permettent d'\'eviter de faire cette v\'erification au cas par cas : voir par exemple la d\'emonstration du Thm. 1.13 dans~\cite{chrenard2}.\ps

\subsubsection{\'Etape (c) : d\'etermination des termes manquants et simple
connexit\'e du groupe de Langlands de $\Z$}\label{scLZ} Expliquons maintenant bri\`evement comment d\'eduire le Th\'eor\`eme~\ref{chrenard} des \'etapes pr\'ec\'edentes. Un r\^ole important est jou\'e par les $\Z$-isog\'enies 
centrales exceptionnelles ${\rm SL}(2) \rightarrow {\rm SO}(2,1)={\rm PGL}(2)$, ${\rm SO}(2,2) \rightarrow {\rm PGL}(2) \times {\rm PGL}(2)$, et ${\rm Sp}(4) \rightarrow {\rm SO}(3,2)={\rm PGSp}(4)$.
Gr\^ace aux travaux d'Arthur et \`a un lemme \'el\'ementaire de rel\`evement de formes automorphes pour les isog\'enies entre groupes de Chevalley, ces isog\'enies nous permettent de d\'emontrer les \'egalit\'es suivantes, qui sont des manifestations de cas particuliers de fonctorialit\'e de Langlands :
\begin{itemize}\medskip
\item[(i)] $O^\ast(w)=S(\frac{w}{2})$, \ps \ps
\item[(ii)] $O(w,v)=S(\frac{w+v}{2})\cdot S(\frac{w-v}{2})$ si $v \neq 0$, et $O(w,0)= \frac{S(w/2)\cdot (S(w/2)-1)}{2}$, \ps\ps
\item[(iii)] $O^\ast(w,v)=S(\frac{w+v}{2},\frac{w-v}{2})$.\ps\ps
\end{itemize}
Le (i) n'est autre que la fonctorialit\'e carr\'e sym\'etrique de Gelbart-Jacquet~\cite{GJ}, le (ii) est la fonctorialit\'e "carr\'e tensoriel" et le (iii) est le "carr\'e altern\'e r\'eduit" d\'ej\`a mentionn\'e. \ps

Je dois rajouter aussi que la trivialit\'e des facteurs epsilon des repr\'esentations orthogonales (Arthur) entra\^ine que les termes $O(w_1,\cdots,w_r)$ et $O^\ast(w_1,\cdots,w_r)$ sont nuls hors d'une certaine congruences sur $\frac{1}{2}(\sum_i w_i)$, qui fait notamment que les \'enonc\'es ci-dessus ont un sens. Enfin, un argument de caract\`ere central montre que $O(w_1,\cdots,w_r)$ est nul si $r \equiv 1 \bmod 2$. \ps

J'ai d\'ej\`a dit que $S(w)$ est la dimension de l'espace des formes modulaires paraboliques pour ${\rm SL}(2,\Z)$. La formule de multiplicit\'e d'Arthur et le travail de Tsushima d\'ej\`a cit\'e permet d'en d\'eduire $S(w,v)$. 
Les \'etapes (a) et (b) dans le cas du groupe ${\rm SO}(7)$ permettent alors de tirer la valeur des 
$S(w_1,w_2,w_3)$. Ces m\^eme \'etapes dans le cas de ${\rm SO}(8)$ permettent de d\'eterminer $O^\ast(w_1,w_2,w_3,w_4)$ si $w_4 \neq 0$, ainsi que les 
$$O^\ast(w_1,w_2,w_3)+ 2 O(w_1,w_2,w_3,0).$$
Quand ce dernier est $\leq 1$, on en d\'eduit $O^\ast(w_1,w_2,w_3)$ et $O(w_1,w_2,w_3,0)=0$. Cela se produit tout le temps quand $w_1 \leq 26$ et c'est ce qui nous a permit de confirmer les calculs 
de~\cite{BFVdG}, qui eux permettent de d\'eterminer $O^\ast(w_1,w_2,w_3)$ en g\'en\'eral. Enfin, les \'etapes (a) et (b) pour ${\rm SO}(9)$ conduisent \`a la d\'etermination de ${\rm S}(w_1,w_2,w_3,w_4)$ !

\subsection{Une ouverture : structure du groupe de Langlands de $\Z$ motivique} \label{conjLZ}
Je voudrais terminer cette partie en revenant un peu sur ce groupe conjectural, qui est au coeur des consid\'erations
pr\'ec\'edentes.  \ps

Soit $\mathcal{L}_\Q$ le groupe conjectural de Langlands, que l'on verra \`a la
Kottwitz~\cite{kottwitzctt} comme un groupe topologique localement compact
muni de classes de conjugaison d'homomorphismes de groupes topologiques $\iota_v : {\rm W}'_{\Q_v} \rightarrow
\mathcal{L}_\Q$, o\`u $v$ parcourt les places de $\Q$ (on rappelle que ${\rm
W}'_{\Q_v}$ est le groupe de Weil-Deligne de $\Q_v$ sous sa forme ${\rm
SU}(2)$, cf.~\S\ref{notations}). Une construction conjecturale axiomatique de $\mathcal{L}_\Q$ a
notamment \'et\'e propos\'ee par Arthur dans~\cite{arthurconjlan}. D\'esignons par $\mathcal{L}_\Z$ le plus grand quotient de $\mathcal{L}_\Q$ dans lequel :
\begin{itemize}\ps
\item[(i)] (conducteur 1) l'image de $\iota_p(I_p \times {\rm SU}(2))$ est triviale pour
tout nombre premier $p$, $I_p$ d\'esignant le groupe d'inertie de ${\rm
W}_{\Q_p}$, 
\item[(ii)] (essentielle alg\'ebricit\'e) l'image de $\iota_\infty(Z(W_\R))$ est centrale, o\`u
$Z(W_\R) \simeq \R^\ast$ d\'esigne le
centre de ${\rm W}_\R$. \ps
\end{itemize}

Pour tout entier $n\geq 1$ les classes d'isomorphie de repr\'esentations irr\'eductibles continues $\rho :
\mathcal{L}_\Z \rightarrow {\rm GL}(n,\C)$ sont en bijection avec les
repr\'esentations automorphes cuspidales $\pi$ de $\GL(n)$ sur $\Q$ qui sont
non ramifi\'ees \`a toutes les places finies et telles que $\pi \otimes
|.|^s$ est alg\'ebrique pour un certain $s \in \C$ (\S\ref{conjfolk}). Le
lecteur voulant v\'erifier cette affirmation observera que d'apr\`es le
lemme de puret\'e de Clozel, si $\pi$ est une repr\'esentation automorphe
alg\'ebrique cuspidale de $\GL(n)$ sur $\Q$, alors $Z(W_\R)$ agit par des
scalaires dans ${\rm L}(\pi_\infty)$. \ps 

L'\'epimorphisme norme $|.| : \mathcal{L}_\Q \rightarrow \R_{>0}$ se factorise par
$\mathcal{L}_\Z$ et induit un isomorphisme $\mathcal{L}_\Z^{\rm ab} \isomo \R_{>0}$. 
De plus, $|.|$ admet une section centrale canonique au dessus de $\mathcal{L}_\Z$ par la condition (ii). On
dispose donc d'un isomorphisme canonique ("d\'ecomposition polaire")
	$$\mathcal{L}_\Z \isomo  \mathcal{L}_\Z^1 \times \R_{>0}$$
o\`u $\mathcal{L}_\Z^1 = \{ g \in \mathcal{L}_\Z, |g|=1\} \simeq
\mathcal{L}_\Z/\iota_\infty(\R_{>0})$, de sorte que le groupe
parfait $\mathcal{L}_\Z^1$ param\`etre aussi les $\pi$ comme plus haut de
caract\`ere central trivial.\ps

Une structure suppl\'ementaire importante que $\mathcal{L}^1_\Z$,
identifi\'e \`a $\mathcal{L}_\Z/\iota_\infty(\R_{>0})$, h\'erite de
$\mathcal{L}_\Q$, est la collection des $\iota_v$. D'une part, si ${\rm W}_\R^1$ d\'esigne le sous-groupe des \'el\'ements de ${\rm W}_\R$
de norme $1$, i.e. une extension de $\Z/2\Z$ par $\mathbb{S}^1$, on dispose
enfin d'une classe de conjugaison canonique $\iota_\infty: {\rm W}_\R^1 \rightarrow \mathcal{L}^1_\Z$. 
La classe de conjugaison de cocaract\`eres induite $$\mathbb{S}^1 \rightarrow
{\mathcal L}^1_\Z$$
est particuli\`erement importante ("cocaract\`ere de Hodge"). D'autre part, on dispose pour chaque
premier $p$, d'une classe de conjugaison de Frobenius dans $\mathcal{L}^1_\Z$.
On conjecture depuis Sato, Tate et Langlands qu'elles sont \'equi-r\'eparties
dans l'espace des classes de conjugaison de $\mathcal{L}_\Z^1$ muni de sa
mesure de probabilit\'e invariante (voir aussi l'expos\'e de Serre~\cite[Ch. 1, appendice]{serreabelian}). Ce dernier appara\^it donc comme le
groupe de Sato-Tate universel pour les "motifs de conducteur $1$". \ps

Concr\`etement, si $\rho$ et $\pi$ se correspondent
comme plus haut, alors le groupe compact $\rho(\mathcal{L}^1_\Z)$ est aussi   
le groupe de Sato-Tate conjectural de $\pi$, nous le noterons
$\mathcal{L}_\pi$. Par d\'efinition, c'est donc un sous-groupe compact
connexe et irr\'eductible de $\GL_n(\C)$ bien d\'efini modulo ${\rm
GL}(n,\C)$-conjugaison. 
\ps 

On peut d\'emontrer que le groupe $\mathcal{L}_\Z^1$ est compact, connexe, et
m\^eme simplement connexe. La compacit\'e est en effet essentiellement dans
les axiomes, la connexit\'e est un analogue automorphe du th\'eor\`eme de
Minkowski se d\'eduisant des formules explicites de Weil (Serre,
Mestre~\cite{mestre}), et la simple connexit\'e a \'et\'e discut\'ee au~\S\ref{scLZ} (voir aussi~\cite{arthurconjlan}). C'est donc un produit direct 
de groupes de Lie compacts connexes, simplement connexes, quasi-simples (en
quantit\'e d\'enombrable). \ps

\`A ma connaissance, on ne conna\^it pas actuellement la structure exacte 
de~$\mathcal{L}^1_\Z$. Une conjecture na\"ive, quoiqu'assez naturelle, est d'esp\'erer que
tous les groupes de Lie compacts, connexes, parfaits, et simplement connexes (mais pas
forc\'ement quasi-simples) sont facteurs direct de $\mathcal{L}^1_\Z$. Par
manque d'arguments totalement probants je la pose ici comme une question.

\begin{question} Est-ce que tout groupe de Lie compact, connexe, \'egal \`a
son groupe d\'eriv\'e, et simplement connexe, est facteur direct de $\mathcal{L}^1_\Z$ ?
\end{question}

Bien que je n'ai pas trouv\'e de trace de cette question dans la
litt\'erature, elle est probablement du folklore, et au moins assez proche de questions que posent Serre
dans~\cite{serremotives}. L'existence d'une infinit\'e de formes modulaires
paraboliques propres et normalis\'ees pour ${\rm SL}(2,\Z)$ montre par exemple que
${\rm SU}(2)^m$ est facteur direct de $\mathcal{L}_\Z^1$ pour tout $m\geq
1$. \ps

Plus g\'en\'eralement, les r\'esultats de ce chapitre montrent que la
question ci-dessus admet une r\'eponse affirmative pour les groupes dont les facteurs simples sont de type $A_1$,
${\rm B}_2$, ${\rm G}_2$, ${\rm B}_3$, ${\rm C}_3$, ${\rm C}_4$ et ${\rm D}_4$. Mieux, pour chaque groupe quasi-simple
$H$ de ce type, nous donnons le nombre des $\pi$ tels que $\mathcal{L}_\pi
\simeq H$ et dont le cocaract\`ere de Hodge est un cocaract\`ere
r\'egulier donn\'e quelconque de $H$. \ps

Il me semble qu'il ne devrait pas
\^etre trop difficile de d\'emontrer que la question admet une r\'eponse
affirmative plus g\'en\'eralement pour les groupes $H$ dont les facteurs
simples ne sont ni de type ${\rm A}_n$ ou ${\rm D}_{2n+1}$ avec $n>1$, ni de type ${\rm E}_6$, autrement dit pour les $H$ de centre un $2$-groupe
ab\'elien \'el\'ementaire. C'est quelque chose que j'aimerais v\'erifier
dans un future proche. Le
point est que sous cette condition, $-1$ est un \'el\'ement dans le groupe
de Weyl de $H$ et donc le $\Z$-groupe semisimple $G$ de Chevalley de dual de
Langlands le complexifi\'e de $H$ a la propri\'et\'e que $G(\R)$ admet des s\'eries 
discr\`etes. \ps

La vraie question int\'eressante me semble plut\^ot de construire des
facteurs de type ${\rm A}_n$, ${\rm D}_{2n+1}$ avec $n>1$, ou de type ${\rm
E}_6$.  Je ne connais aucun $\pi$ tel que $\mathcal{L}_\pi$ est de ce type. 
Seul le cas du type ${\rm A}_2$ a \'et\'e \'etudi\'e \`a ma connaissance,
par Ash et Pollack dans~\cite{ashpollack}.  Dans le language adopt\'e ici,
ces auteurs ont cherch\'e des $\pi$ alg\'ebriques de conducteur $1$ telles
que $\mathcal{L}_\pi \simeq {\rm SU}(3)$. Leurs cocaract\`eres de Hodge sont
de la forme $u \mapsto (u^{2k},1,u^{-{2k}})$   
avec $k$ entier. Il est remarquable que ces auteurs n'aient trouv\'e aucun tel
$\pi$ pour $|k| \leq 121$, ce qui \'etait contraire \`a leurs attentes. Ils ont alors
avanc\'e l'hypoth\`ese qu'il n'existerait aucun tel $\pi$, ce qui est contraire \`a la croyance
na\"ive avanc\'ee plus haut. Comme le font remarquer ces auteurs, il est
remarquable que si l'on exclut l'hypoth\`ese de conducteur $1$, le
probl\`eme analogue admet des solutions pour des petits $k$, comme l'ont montr\'e van Geemen et
Top dans~\cite{geemen} (voir aussi~\cite{dettweiler}), qui donnent des
exemples o\`u le conducteur est une puissance de $2$ et o\`u $k=1$.\ps

\newpage

\section{Quelques cas des conjectures de Bloch-Kato}\label{thmlivre}

Dans cette derni\`ere partie, je vais exposer de mani\`ere tr\`es succinte les th\'eor\`emes principaux de mon livre avec Bella\"iche~\cite[Ch. 8,9]{bchlivre}. Ces r\'esultats 
sont dans la lign\'ee des travaux de Ribet~\cite{ribet} sur la r\'eciproque du fameux th\'eor\`eme de Herbrand, ou encore de la d\'emonstration par Mazur-Wiles~\cite{mazurwiles} de la conjecture principale d'Iwasawa, et plus particuli\`erement de l'approche de Wiles~\cite{wilesmc}. Je renvoie \`a l'article d'exposition de Mazur~\cite{mazurribet} pour un tour d'horizon de ce sujet. \ps

Soit $F$ un corps quadratique imaginaire et soit $\pi$ une repr\'esentation automorphe cuspidale alg\'ebrique r\'eguli\`ere de $\GL_n$ sur $F$ telle que\footnote{On se place donc en "poids motivique $-1$".}  $\pi^\vee \simeq \pi^c |\cdot|^{-1}$. Soit $L(s,\pi)$ La fonction $L$ de $\pi$, qui est une fonction enti\`ere de $s$, et soit  $L^\ast(s,\pi)$ sa fonction $L$ compl\'et\'ee \`a la place archim\'edienne (voir~\cite{serrezeta},\cite{tate},\cite{cogdell}). 
Elle satisfait $$L^\ast(s,\pi)=\varepsilon(s,\pi)L^\ast(-s,\pi).$$
En particulier, le facteur $\varepsilon(\pi):=\varepsilon(0,\pi)$ est un signe $\pm 1$. Sous nos hypoth\`eses, $L(s,\pi)$ et $L^\ast(s,\pi)$ ont m\^eme ordre d'annulation en $0$, ce signe dicte donc la parit\'e de l'ordre d'annulation de $L(s,\pi)$ en $s=0$. \ps

Bien que la repr\'esentation $\pi$ ne soit pas tout-\`a-fait polaris\'ee au sens du~\S\ref{repgalaut}, elle le devient apr\`es torsion par un caract\`ere de Hecke alg\'ebrique de $F$ bien choisi, de sorte que l'on peut \'egalement lui associer des repr\'esentations galoisiennes $\rho_{\pi,\iota}$. Fixons un nombre premier $p$ ainsi qu'une r\'ealisation $p$-adique $\rho_{\pi,\iota}$ de $\pi$. Bloch et Kato d\'efinissent dans~\cite{bkato} un sous-espace de dimension finie 
$$H^1_f(F,\rho_{\pi,\iota}) \subset H^1({\rm Gal}(\overline{F}/F),\rho_{\pi,\iota}).$$
Ce sous-espace est par d\'efinition celui param\'etrant les classes d'extensions $U$ de la repr\'esentation triviale $1$ par $\rho_{\pi,\iota}$, dont la restriction $U_v$ \`a un groupe de d\'ecomposition en $v$ a les propri\'et\'es suivantes : $U_v^{I_v} \rightarrow 1$ est surjectif si $v$ ne divise pas $p$, $I_v$ d\'esignant le groupe d'inertie en $v$, et ${\rm D}_{\rm cris}(U_v) \rightarrow {\rm D}_{\rm cris}(1)$ est surjectif si $v|p$. \ps

Les conjectures de Bloch et Kato~\cite{bkato}, reprises par Fontaine et Perrin-Riou dans~\cite{FP}, sont des g\'en\'eralisations d'une variante ... de la conjecture de Birch et Swinnerton-Dyer.\footnote{Elles entra\^inent la conjecture de Birch-Swinnerton-Dyer "modulo" la finitude du groupe de Tate-Shaffarevich.}  Je renvoie \`a \cite[Ch. 5]{bchlivre} pour une discussion d\'etaill\'ee de ces conjectures dans le contexte adopt\'e ici. 

\begin{conjecture} (Bloch-Kato) ${\rm ord}_{s=0}\, L(s,\pi) = \dim H^1_f(F,\rho_{\pi,\iota})$.\end{conjecture}

\noindent En particulier, la version tr\`es affaiblie ci-dessous devrait \^etre \'egalement vraie. 

\begin{conjecture} (Conjecture du signe) Si $\varepsilon(\pi)=-1$ alors $H^1_f(F,\rho_{\pi,\iota}) \neq 0$. \end{conjecture}

Pour r\'esumer mes travaux avec Bella\"iche, nous pouvons dire que nous avons fait un progr\`es certain vers la conjecture du signe, et un tout petit progr\`es vers la conjecture de Bloch-Kato. \ps

Commen\c{c}ons par d\'ecrire bri\`evement notre r\'esultat concernant la conjecture du signe. Nous faisons d'abord les hypoth\`eses suppl\'ementaires suivantes : \ps
\begin{itemize}
\item[(i)] $n \not \equiv 0 \bmod 4$,\ps
\item[(ii)] $p$ est d\'ecompos\'e dans $F$ et $\pi$ est non ramifi\'ee aux places divisant $p$, \ps
\item[(iii)]  les entiers $p_{i,\sigma}$ associ\'es \`a $\pi_\infty$ (voir~\S\ref{repgalaut} Ch. 1) sont diff\'erents de $0$ et $1$,\ps
\item[(iv)]  si $v$ est une place finie non d\'ecompos\'ee de $F$ telle que $\pi_v$ est ramifi\'ee, alors $\pi_v$ est une s\'erie principale suffisamment r\'eguli\`ere au sens 
de~\cite[\S 9.2.1]{bchlivre}.\ps
\end{itemize}
	
	Nous avons \'egalement deux autres hypoth\`eses importantes, que l'on note ${\rm AC}(\pi)$ et ${\rm Rep}(n+2)$ {\it loc.cit.}, et que je d\'ecris maintenant bri\`evement. Sous l'hypoth\`ese (i), il existe une forme hermitienne sur $F^{n+2}$ relative \`a $F/\Q$ qui est d\'efinie positive sur $\R$ et d'indice $[n/2+1]$ \`a toutes les places finies de $\Q$, notons $U(n+2)$ son $\Q$-groupe unitaire. L'hypoth\`ese ${\rm Rep}(n+2)$ est la g\'en\'eralisation au groupe $U(n+2)$ des r\'esultats de Labesse~\cite{labesse} 
(qui suppose "pour simplifier" que son groupe unitaire est d\'efini sur un corps totalement r\'eel de degr\'e $\geq 2$) concernant le changement de base des repr\'esentations automorphes de $U(n+2)$ \`a $\GL(n+2)$ sur $F$. Cette hypoth\`ese assure notamment l'existence de repr\'esentations galoisiennes associ\'ees aux repr\'esentations automorphes de $U(n+2)$. \ps

 L'hypoth\`ese ${\rm AC}(\pi)$ est un cas particulier des conjectures d'Arthur~\cite{arthurunipotent}.
Elle affirme l'existence d'une certaine repr\'esentation automorphe $\pi'$ de $U(n+2)$, choisie "aussi peu ramifi\'ee que possible" (je renvoie \`a~\cite[\S 6.9]{bchlivre} pour la recette pr\'ecise), 
dont la repr\'esentation galoisienne $\rho_{\pi',\iota}$ vaut $$\rho_{\pi,\iota}\oplus 1 \oplus \omega$$ apr\`es torsion convenable par un caract\`ere, $\omega$ d\'esignant le caract\`ere cyclotomique $p$-adique. Les conjectures g\'en\'erales d'Arthur, tout-\`a-fait dans l'esprit de celles d\'ej\`a discut\'ees dans mes travaux avec Lannes et avec Renard aux~\S\ref{niemeierchlannes} et~\ref{cond1chrenard} de ce chapitre, montrent qu'une telle repr\'esentation $\pi'$ doit exister pr\'ecis\'ement sous l'hypoth\`ese $\varepsilon(\pi)=-1$. \ps

Le r\'esultat suivant est d\'emontr\'e dans le chapitre 8 de mon livre avec Bella\"iche. 

\begin{thm}\label{thmsignebk} Supposons que $\pi$ satisfasse les conditions (i) \`a (iv), ainsi que les hypoth\`eses ${\rm AC}(\pi)$ et ${\rm Rep}(n+2)$. Alors la conjecture du signe est vraie pour  $\pi$.
\end{thm}

Quelques remarques s'imposent. Tout d'abord, lorsque $n=1$ les hypoth\`eses (i), (iv) sont trivialement v\'erifi\'ees, ainsi que ${\rm AC}(\pi)$ et ${\rm Rep}(3)$ d'apr\`es Rogawski~\cite{roglivre},~\cite{picard}.
Le th\'eor\`eme dans ce cas \'etait le r\'esultat principal de mon article avec Bella\"iche~\cite{bchens}. Il \'etait d'ailleurs d\'ej\`a connu ! car il se d\'eduit des r\'esultats bien ant\'erieurs de Rubin sur la "conjecture principale d'Iwasawa" pour les corps quadratiques imaginaires~\cite{rubinmc}. L'int\'er\^et de~\cite{bchens} consistait donc plut\^ot dans l'approche que dans le r\'esultat.  \ps

L'id\'ee de la d\'emonstration, pour $n$ g\'en\'eral, est de construire des congruences modulo une puissance de $p$ entre la repr\'esentation $\pi'$ de $U(n+2)$ et des repr\'esentations automorphes $\pi''$ de repr\'esentation galoisienne irr\'eductible, soit encore de d\'eformer irr\'eductiblement la repr\'esentation galoisienne 
\begin{equation}\label{repgalpi'} \rho_{\pi,\iota}\oplus 1 \oplus \omega. \end{equation} Il s'agit ensuite de comprendre comment \'etendre "l'argument de r\'eseaux de Ribet"~\cite{ribet} dans ce contexte pour produire une extension globalement non triviale de $1$ par $\rho_{\pi,\iota}$ ayant les propri\'et\'es locales requises. \ps

Cette strat\'egie avait \'et\'e men\'ee \`a bien par Bella\"iche dans sa th\`ese~\cite{Bellaichethese}, toujours quand $n=1$. Il construit ses congruences par la m\'ethode "d'augmentation du niveau" de Ribet et observe qu'un contr\^ole soigneux de la ramification en toutes les places permet de construire une extension avec toutes les propri\'et\'es d\'esir\'ees. Le fait que $F$ soit quadratique imaginaire plut\^ot que CM quelconque joue un r\^ole tr\`es s\'erieux dans son argument, et explique aussi cette hypoth\`ese ici, car il doit \`a un moment \'eliminer la possibilit\'e de construire \`a la place une unit\'e d'ordre infinie de $\OO_F$ (par la th\'eorie de Kummer!), ce qui n'existe pas si $F$ est imaginaire. Le probl\`eme est l'apparition du caract\`ere cyclotomique dans~\eqref{repgalpi'}, qui n'a pas d'analogue chez Ribet, qui travaille avec des s\'eries d'Eisenstein plut\^ot que des repr\'esentations automorphes discr\`etes non-temp\'er\'ees. Une faiblesse de l'approche de Bella\"iche est qu'elle ne construit l'extension cherch\'ee que modulo certains nombres premiers, et donc n'obtient pas tout-\`a-fait l'\'enonc\'e ci-dessus m\^eme quand $n=1$. \ps

Un second progr\`es substantiel a \'et\'e alors r\'ealis\'e par Skinner et Urban~\cite{sku1},~\cite{sku2}, dans un contexte proche o\`u le $\Q$-groupe $U(3)$ est remplac\'e par ${\rm Sp}(4)$. Ces auteurs, inspir\'es sans doute par la m\'ethode de Wiles~\cite{wilesmc}, montrent comment produire des extensions "en caract\'eristique $0$" en utilisant des familles de formes modulaires \`a la Coleman comme argument de d\'eformation. Les travaux de Kisin~\cite{kisinoc} que nous avons rappel\'es au~\S\ref{appraffinee} Ch. 1 jouent alors un r\^ole important dans l'\'etablissement des propri\'et\'es en $p$ des extensions construites. \ps

Dans l'article~\cite{bchens}, nous mettions ensemble les deux approches ci-dessus pour d\'emontrer le cas $n=1$ du Th\'eor\`eme~\ref{thmsignebk}, \`a l'aide notamment de mon travail ant\'erieur~\cite{chcrelle} sur les familles de Coleman pour $U(3)$~(\S\ref{fougereunitaire} Ch. 1). Le th\'eor\`eme ci-dessus est une g\'en\'eralisation de cette m\'ethode. Les difficult\'es principales sont d'ordre technique. Nous ne supposons pas, notamment, que $\rho_{\pi,\iota}$ est irr\'eductible. En revanche, comme le lecteur l'aura remarqu\'e, notre travail est soumis \`a un certain nombre d'autres d'hypoth\`eses ! \ps

Il y a lieu cependant de ne pas \^etre trop pessimiste sur ces hypoth\`eses. Tout d'abord, les travaux r\'ecents d'Arthur~\cite{arthur}, et notamment leurs  g\'en\'eralisations par Mok~\cite{mok} au cas des groupes unitaires, autorisent un certain optimisme sur ${\rm AC}(\pi)$. De plus, il semblerait que ${\rm Rep}(n+2)$ soit connue des sp\'ecialistes (voir par exemple les travaux de Morel~\cite{morel}). L'hypoth\`ese (i) d'apparence importante pourrait sans doute \^etre contourn\'ee via d'autres cas des conjectures d'Arthur (cela vient de ce que l'on n'a pas besoin de supposer que $\rho_{\pi,\iota}$ est irr\'eductible dans notre argument). L'hypoth\`ese (ii) est faite par commodit\'e et l'hypoth\`ese (iv) nous a \'et\'e n\'ecessaire principalement pour formuler de mani\`ere non ambig\"ue ${\rm AC}(\pi)$ aux places ramifi\'ees dans $F$ :  on pourrait vraisemblablement proc\'eder maintenant autrement \'etant donn\'e les r\'esultats r\'ecents de Mok suscit\'es sur la conjecture de Langlands locale pour les groupes unitaires. Enfin, Bella\"iche a compris dans~\cite{joeltrans} comment s'affranchir de (iii) \'egalement. Je voudrais rajouter pour finir cette interminable discussion qu'il ne sera peut-\^etre pas n\'ecessaire de faire tout cela, car Skinner et Urban ont annonc\'e \`a l'I.C.M.~\cite{skuicm} une g\'en\'eralisation de notre th\'eor\`eme. \ps

Je voudrais terminer ce m\'emoire en expliquant le second r\'esultat principal de notre livre, qui est un raffinement du r\'esultat pr\'ec\'edent. Je suppose d\'esormais que $\pi=\chi$ est un caract\`ere de Hecke alg\'ebrique du corps quadratique imaginaire $F$ tel que  $\chi(z\overline{z})=|z|$ pour tout $z \in \AAA_F^\ast$, et  tel que $\chi_\infty(z)=z^q \overline{z}^{1-q}$ avec 
$q \in \Z \backslash \{0,1\}$. Je suppose aussi que le nombre premier $p$ est d\'ecompos\'e dans $F$ et que $\chi$ est non-ramifi\'e en $p$. \ps

On peut alors consid\'erer la vari\'et\'e de Hecke unitaire 
$$\mathcal{E}_3 \subset Y_3^\bot \times \mathcal{T}^3$$ 
que nous avons introduite au \S\ref{fougereunitaire} Ch. 1. Elle est associ\'ee au groupe $U(3)$ relativement \`a $F$ d\'efini ci-dessus et \`a l'ensemble $S$ des nombres premiers en lesquels soit $F$, soit $\chi$, est ramifi\'e. Comme nous l'avons d\'ej\`a dit plus haut, sous l'hypoth\`ese $\varepsilon(\chi)=-1$ l'assertion ${\rm AC}(\chi)$ d\'emontr\'ee par Rogawski assure l'existence d'une repr\'esentation automorphe $\pi'$ de $U(3)$ "minimalement ramifi\'ee" et telle que $\rho_{\pi',\iota}$ soit \`a une torsion pr\`es la repr\'esentation 
	$$\chi_{\iota} \oplus 1 \oplus \omega,$$
o\`u $\chi_\iota$ est la r\'ealisation $\iota$-adique de $\chi$. Fixons $v$ l'une des deux places de $F$ divisant $p$ et choisissons en cette place le raffinement de $\pi'$ tout particulier suivant : \`a la m\^eme torsion pr\`es que ci-dessus c'est le raffinement 
$\Phi_v=(1,\ast,p^{-1})$. La repr\'esentation raffin\'ee $(\pi',\{\Phi_v\})$ nous fournit donc un point
$$x \in \mathcal{E}_3$$
ainsi que nous l'avons expliqu\'e au~\S\ref{fougereunitaire} du Chapitre 1. Nous l'appelons le {\it point d'Arthur anti-ordinaire} associ\'e \`a $\chi$. 
L'hypoth\`ese de ramification minimale dans le choix de $\pi'$ assure que ce point appartient \'egalement au sous-espace
$$\mathcal{E}_3^{\rm min} \subset \mathcal{E}_3$$
obtenu en prenant l'adh\'erence Zariski des points automorphes raffin\'es qui param\'etrent les repr\'esentations $\Pi$ de $\GL_3$ sur $F$ de param\`etre de Langlands ${\rm L}(\Pi_u)$ trivial sur le facteur ${\rm SU}(2)$ pour toute place finie $u$ de $F$. Des arguments de $K$-type en th\'eorie des vari\'et\'es de Hecke montrent que ce sous-espace est encore d'\'equidimension $3$, autrement dit que c'est une r\'eunion de composantes irr\'eductibles de $\mathcal{E}_3$. Nous d\'emontrons alors le r\'esultat suivant~\cite[Ch. 9]{bchlivre}.

\begin{thm} On suppose $\varepsilon(\chi)=-1$. Soit $x \in \mathcal{E}_3^{\rm min}$ le point d'Arthur anti-ordinaire associ\'e \`a $\chi$. Soit $t$ la dimension de l'espace tangent de $\mathcal{E}_3^{\rm min}$ en $x$ et soit $h=\dim H^1_f(F,\chi_\lambda)$. On a alors l'in\'egalit\'e 
$$ t \leq \frac{(h+1)(h+2)}{2}.$$
En particulier $h\geq 1$, et si $\mathcal{E}_3^{\rm min}$ n'est pas lisse en $x$ alors $h\geq 2$.
\end{thm}

\`A notre connaissance, ce type d'\'enonc\'e est le premier de la sorte. Il n'a notamment pas d'analogue en th\'eorie d'Iwasawa. Notre d\'emonstration consiste en une \'etude des diff\'erents lieux de r\'eductibilit\'e, au voisinage du point $x$, de la famille de pseudo-caract\`eres galoisiens (locaux et globaux) port\'ee par $\mathcal{E}_3^{\rm min}$. Nos r\'esultats mentionn\'es au Chapitre 1~\S\ref{prelimvcar} et~\S\ref{appraffinee} y jouent un r\^ole important. Il n'est pas clair actuellement comment relier la g\'eom\'etrie de $\mathcal{E}_3^{\rm min}$ en $x$ aux fonctions $L$ (vraisemblablement $p$-adiques), encore moins \`a la conjecture exacte de Bloch-Kato, qui rappelons-le n'est m\^eme pas connue dans ce cas particulier. Il serait notamment int\'eressant de savoir, par exemple par des calculs dans des cas concrets, s'il faut s'attendre \`a ce que la r\'eciproque de la derni\`ere assertion du th\'eor\`eme soit \'egalement vraie !

\chapter*{Remerciements}

	Je remercie chaleureusement Michael Harris de m'avoir initi\'e \`a la th\'eorie des formes automorphes, et de m'avoir g\'en\'ereusement fait partag\'e, tout au long de ces dix derni\`eres ann\'ees, ses id\'ees stimulantes et sa vision unique sur ce magnifique sujet.  	Je pense \'egalement \`a ma longue et enrichissante collaboration avec Jo\"el Bella\"iche, notamment \`a nos innombrables discussions  sur les conjectures de Langlands et Arthur  : merci Jo\"el !  \ps

	Je voudrais exprimer ma sinc\`ere gratitude \`a Guy Henniart, Richard Taylor et Jacques Tilouine qui ont accept\'e d'\^etre rapporteurs
de ce travail, et ce malgr\`e les d\'elais contraignants qui leurs \'etaient impos\'es. Je remercie tout particuli\`erement Guy Henniart d'avoir accept\'e avec enthousiasme d'\^etre le coordinateur de cette habilitation, et aussi pour ses encouragements.
Je suis \'egalement tr\`es honor\'e que Laurent Clozel, Michael Harris, Jean Lannes 
et Jean-Pierre Serre aient accept\'e de faire partie de mon jury. \ps

	C'est avec grand plaisir que je remercie mes collaborateurs Teodor Banica, Jo\"el Bella\"iche, Laurent Berger, Julien Bichon, Laurent Clozel, Michael Harris, Chandrashekhar Khare, Jean Lannes, Michael Larsen, et David Renard, pour les joies que j'ai eues \`a travailler avec eux et pour ce qu'ils ont rendu possible. Je suis particuli\`erement redevable \`a Laurent Clozel et Jean Lannes de ce qu'ils m'ont appris. Merci aussi aux autres, qui ont r\'eguli\`erement r\'epondu \`a mes questions, ou avec qui j'ai eu des discussions math\'ematiques inspirantes; je pense notamment \`a Ahmed Abbes, Nicolas Bergeron, Christophe Breuil, Emmanuel Breuillard, Kevin Buzzard, Pierre Colmez, Christophe Cornut, Jean-Fran\c{c}ois Dat, Matthew Emerton, Laurent Fargues, Jean-Pierre Labesse, Vincent Lafforgue, Julien March\'e, Barry Mazur, Colette Moeglin, Fabrice Orgogozo, et Claudio Procesi. \ps
	
	Un survol de ce m\'emoire suffit pour se rendre compte de toute l'importance dans ma recherche,  parfois m\^eme d\'ependance, des travaux de nombreux math\'ematiciens. Il serait maladroit ici de tenter de les \'enum\'erer tous, mais leurs apparitions au fil du texte rendra, je l'esp\`ere, claire au lecteur, toute l'admiration que j'ai pour leur travail. \ps	
	
	Depuis ma th\`ese, j'ai eu la chance de b\'en\'eficier des ambiances chaleureuses et motivantes de l'\'equipe de g\'eom\'etrie arithm\'etique du L.A.G.A., de l'I.H.\'E.S. et du C.M.L.S. J'adresse un grand merci \`a tous mes coll\`egues ou ex-coll\`egues de ces institutions. Merci aussi \`a tous ceux ayant rendu techniquement possible cette soutenance ! \ps

Merci enfin \`a tous mes amis que j'ai plaisir \`a saluer ici. Ma pens\'ee la plus ch\`ere va \`a Valeria et Julia : merci pour tout le bonheur que vous m'apportez.\ps

\backmatter

\address{Ga\"etan Chenevier \\ Centre de Math\'ematiques Laurent Schwartz \\ \'Ecole Polytechnique \\ 91128 Palaiseau Cedex \\ FRANCE }

\end{document}